\documentclass[a4paper,12pt,amsart,frenchb]{article}

\usepackage{amsmath,amsbsy,amsfonts,amssymb}
\usepackage[french]{babel}
\newcommand{\ass}[2]{\vskip0.3cm\noindent
{\bf {#1}}. { \sl {#2}}\vskip0.3cm\noindent
}

\begin{document}

\title{  Une formule int\'egrale reli\'ee \`a la conjecture locale de Gross-Prasad }
\author{J.-L. Waldspurger}
\date{10 f\'evrier 2009}
\maketitle

\bigskip

{\bf Introduction}

Soit $F$ un corps local non archim\'edien de caract\'eristique nulle. Soit $V$ un espace vectoriel sur $F$, de dimension finie $d$, muni d'une forme quadratique non d\'eg\'en\'er\'ee $q$. On suppose donn\'ee une d\'ecomposition en somme directe de sous-espaces orthogonaux deux \`a deux $V=W\oplus D\oplus Z$. On suppose que $D$ est une droite et que $Z$ est muni d'une base $\{ v_{i}; i=\pm 1,...,\pm r\}$ telle que $q(v_{i},v_{j})=\delta_{i,-j}$ pour tous $i,j$, o\`u $\delta_{i,-j}$ est le symbole de Kronecker. On note $G$, resp. $H$, le groupe sp\'ecial orthogonal de $V$, resp. $W$, et $U$ le radical unipotent du sous-groupe parabolique de $G$ qui conserve le drapeau de sous-espaces isotropes
$$Fv_{r}\subset Fv_{r}\oplus Fv_{r-1}\subset ...\subset Fv_{r} \oplus ...\oplus Fv_{1}.$$
Fixons un \'el\'ement non nul $v_{0}\in D$ et un caract\`ere continu non trivial $\psi$ de $F$. D\'efinissons un caract\`ere $\xi$ de $U(F)$ par la formule
$$\xi(u)=\psi(\sum_{i=0,...,r-1}q(uv_{i},v_{-i-1})).$$
Le groupe $H$ est le sous-groupe des \'el\'ements de $G$ qui agissent par l'identit\'e sur $D\oplus Z$. Il normalise $U$ et la conjugaison par $H(F)$ conserve $\xi$ ($\xi$ est essentiellement le caract\`ere de $U(F)$ le plus r\'egulier possible qui soit conserv\'e par cette conjugaison). Soient $\pi$, resp. $\sigma$, une repr\'esentation admissible irr\'eductible de $G(F)$, resp. $H(F)$, dans un espace (complexe) $E_{\pi}$, resp. $E_{\sigma}$. Notons $Hom_{H,\xi}(\pi,\sigma)$ l'espace des applications lin\'eaires $\varphi:E_{\pi}\to E_{\sigma}$ telles que
$$\varphi(\pi(hu)e)=\xi(u)\sigma(h)\varphi(e)$$
pour tous $u\in U(F)$, $h\in H(F)$, $e\in E_{\pi}$. D'apr\`es [AGRS] th\'eor\`eme 1' et [GGP] corollaire 20.4, cet espace est de dimension $0$ ou $1$. On note $m(\sigma,\pi)$ cette dimension.

Supposons maintenant que $G$ et $H$ sont quasi-d\'eploy\'es sur $F$ et affectons les donn\'ees $V$, $W$, $q$, $G$ et $H$ d'un indice $i$ (pour "isotrope"). Dans cette introduction, supposons pour simplifier $dim(W_{i})\geq3$. On sait qu'\`a isomorphisme pr\`es, il existe un unique espace $V_{a}$ ($a$ pour "anisotrope") de m\^eme dimension $d$ que $V_{i}$, muni d'une forme quadratique $q_{a}$ de m\^eme discriminant que $q_{i}$ mais qui n'est pas isomorphe \`a $q_{a}$ (c'est-\`a-dire d'indice de Witt oppos\'e). Il existe de m\^eme un unique espace $W_{a}$ de m\^eme dimension que $W_{i}$, muni d'une forme quadratique de m\^eme discriminant que la restriction de $q_{i}$ \`a $W_{i}$, mais qui n'est pas isomorphe \`a cette restriction. On v\'erifie que $V_{a}$ est encore isomorphe \`a la somme directe orthogonale $W_{a}\oplus D\oplus Z$, o\`u les formes quadratiques sur $D$ et $Z$ sont les m\^ emes que pr\'ec\'edemment. On note $G_{a}$, resp. $H_{a}$, le groupe sp\'ecial orthogonal de $V_{a}$, resp. $W_{a}$. C'est une forme int\'erieure de $G_{i}$, resp. $H_{i}$.

 La conjecture locale de Gross-Prasad suppose l'existence des $L$-paquets et certaines de leurs propri\'et\'es. On y reviendra ci-dessous.  Gross et Prasad \'enoncent leur conjecture pour les $L$-paquets g\'en\'eriques. On se limite ici aux $L$-paquets temp\'er\'es. Soit $\Pi_{i}$, resp.
$\Sigma_{i}$, un $L$-paquet de repr\'esentations temp\'er\'ees de $G_{i}(F)$, resp. $H_{i}(F)$. Il peut lui correspondre un $L$-paquet $\Pi_{a}$, resp. $\Sigma_{a}$, de repr\'esentations temp\'er\'ees de $G_{a}(F)$, resp. $H_{a}(F)$. Ce $L$-paquet est alors unique. Ou bien, il n'y a pas de tel $L$-paquet $\Pi_{a}$, resp. $\Sigma_{a}$. Dans ce cas, on pose $\Pi_{a}=\emptyset$, resp. $\Sigma_{a}=\emptyset$. En tout cas, pour $(\sigma,\pi)\in (\Sigma_{i}\times \Pi_{i})\cup (\Sigma_{a}\times \Pi_{a})$, la dimension
$m(\sigma,\pi)$ est d\'efinie.

\ass{Conjecture (Gross-Prasad)}{Il existe un unique couple $(\sigma,\pi)\in(\Sigma_{i}\times \Pi_{i})\cup (\Sigma_{a}\times \Pi_{a})$ tel que $m(\sigma,\pi)=1$.}

C'est une partie de la conjecture 6.9 de [GP]. D\'ecrivons les propri\'et\'es des $L$-paquets temp\'er\'es que nous admettrons (on les \'enonce pour le couple $(\Pi_{i},\Pi_{a})$, mais on admet les propri\'et\'es similaires pour le couple $(\Sigma_{i},\Sigma_{a})$). Notons $\sharp$ l'un des indices $i$ ou $a$. Rappelons qu'\`a toute repr\'esentation admissible irr\'eductible $\pi$ de $G_{\sharp}(F)$  est associ\'e un caract\`ere $\theta_{\pi}$ que l'on peut consid\'erer comme une distribution ou comme une fonction localement int\'egrable sur $G_{\sharp}(F)$.  Dans le cas du groupe $G_{i}$, on sait d\'efinir la notion de mod\`ele de Whittaker de $\pi$. Plus exactement, il y a une  notion de mod\`ele de Whittaker relatif \`a ${\cal O}$ pour chaque orbite nilpotente r\'eguli\`ere ${\cal O}\subset \mathfrak{g}_{i}(F)$ (pour tout groupe r\'eductif $L$, on note $\mathfrak{l}$ son alg\`ebre de Lie). On suppose

(1) pour $\sharp=i$ ou $a$, $\Pi_{\sharp}$ est un ensemble fini, non vide si $\sharp=i$, et la distribution $\theta_{\Pi_{\sharp}}=\sum_{\pi\in \Pi_{\sharp}}\theta_{\pi}$ sur $G_{\sharp}(F)$ est stable;

(2) le transfert \`a $G_{a}(F)$ de la distribution $\theta_{\Pi_{i}}$ est $(-1)^d\theta_{\Pi_{a}}$ (en particulier est nul si $\Pi_{a}=\emptyset$);

(3) pour toute orbite nilpotente r\'eguli\`ere ${\cal O}\subset \mathfrak{g}_{i}(F)$, il existe un et un seul \'el\'ement de $\Pi_{i}$ qui admet un mod\`ele de Whittaker relatif \`a ${\cal O}$.

On reviendra sur ces propri\'et\'es en 13.2. Notre r\'esultat est le suivant.

\ass{Th\'eor\`eme}{ Supposons v\'erifi\'ees les propri\'et\'es ci-dessus. Supposons  de plus que $\Pi_{i}$ et $\Pi_{a}$ soient form\'es uniquement de repr\'esentations supercuspidales. Alors la conjecture ci-dessus est v\'erifi\'ee.}

Ce th\'eor\`eme r\'esulte d'une formule int\'egrale qui calcule la dimension $m(\sigma,\pi)$ \`a l'aide des caract\`eres de $\sigma$ et $\pi$, dans le cas o\`u $\pi$ est supercuspidale. Revenons aux notations du d\'ebut en abandonnant les indices $i$ et $a$. Consid\'erons l'ensemble des sous-tores $T\subset H$  pour lesquels il existe une d\'ecomposition en somme directe orthogonale $W=W'\oplus W''$ de sorte que

- la dimension de $W'$ est paire et les groupes sp\'eciaux orthogonaux $H''$ de $W''$ et $G''$ de $V''=W''\oplus D\oplus Z$ sont quasi-d\'eploy\'es sur $F$;

- le tore $T$ est un sous-tore maximal du groupe sp\'ecial orthogonal de $W'$ et il ne contient aucun sous-tore d\'eploy\'e non trivial.

On fixe un ensemble de repr\'esentants ${\cal T}$ des classes de conjugaison par $H(F)$ dans cet ensemble de tores. Soit $T\in {\cal T}$. On lui associe des groupes $H''$ et $G''$ comme ci-dessus. Soit $\pi$  une repr\'esentation admissible irr\'eductible de $G(F) $.   Harish-Chandra a d\'ecrit le comportement local du caract\`ere $\theta_{\pi}$. Soit $x$ un \'el\'ement semi-simple de $G(F)$. Notons $G_{x}$ la composante neutre du commutant de $x$ dans $G$. Alors, pour toute orbite nilpotente ${\cal O}$ dans $\mathfrak{g}_{x}(F)$, il existe un coefficient $c_{\pi,{\cal O}}(x)\in {\mathbb C}$ de sorte que, pour toute fonction $f\in C_{c}^{\infty}(\mathfrak{g}_{x}(F))$ dont le support soit contenu dans un voisinage assez petit de $0$, on ait l'\'egalit\'e
$$(4) \qquad \int_{\mathfrak{g}_{x}(F)}\theta_{\pi}(xexp(X))f(X)dX=\sum_{{\cal O}}c_{\pi,{\cal O}}(x)\int_{{\cal O}}\hat{f}(X)dX.$$
La somme porte sur les orbites nilpotentes dans $\mathfrak{g}_{x}(F)$ et le dernier terme est la transform\'ee de Fourier de l'int\'egrale orbitale sur ${\cal O}$. Bien s\^ur, les mesures et la transformation de Fourier doivent \^etre d\'efinies pr\'ecis\'ement. Supposons que $x$ est un \'el\'ement de $T(F)$ en position g\'en\'erale. Alors $G_{x}=T\times G''$, en particulier les orbites nilpotentes de $\mathfrak{g}_{x}(F)$ sont celles de $\mathfrak{g}''(F)$. Supposons d'abord $d$ impair. Par hypoth\`ese, $G''$ est quasi-d\'eploy\'e. En dimension impaire, cela implique qu'il est d\'eploy\'e. Son alg\`ebre de Lie $\mathfrak{g}''(F)$ poss\`ede une unique orbite nilpotente r\'eguli\`ere, on la note ${\cal O}_{reg}$ et on pose $c_{\pi}(x)=c_{\pi,{\cal O}_{reg}}(x)$. Supposons maintenant $d$ pair. Alors $\mathfrak{g}''(F)$ poss\`ede (en g\'en\'eral) plusieurs orbites nilpotentes r\'eguli\`eres. On peut les param\'etrer par un sous-ensemble de $F^{\times}/F^{\times 2}$. Posons $\nu_{0}=q(v_{0})$. On montre que $\nu_{0}$ appartient \`a l'ensemble de param\`etres, on lui associe une orbite ${\cal O}_{\nu_{0}}$ et on pose $c_{\pi}(x)=c_{\pi,{\cal O}_{\nu_{0}}}(x)$. On a ainsi d\'efini une fonction $c_{\pi}$ sur un ouvert de Zariski de $T(F)$. Soit $\sigma$ une repr\'esentation admissible irr\'eductible de $H(F)$. On d\'efinit de fa\c{c}on similaire une fonction $c_{\sigma}$ sur un ouvert de Zariski de $T(F)$. Posons
$$(5) \qquad m_{geom}(\sigma,\pi)=\sum_{T\in {\cal T}}w(T)^{-1}\int_{T(F)}c_{\check{\sigma}}(x)c_{\pi}(x)D^H(x)\Delta(x)^rdx.$$
Les fonctions $D^H$ et $\Delta$ sont des d\'eterminants \'el\'ementaires et $w(T)$ est le nombre d'\'el\'ements d'un certain normalisateur. La mesure sur $T(F)$ est de masse totale $1$. La repr\'esentation $\check{\sigma}$ est la contragr\'ediente de $\sigma$.

\ass{Th\'eor\`eme}{(i) Pour des repr\'esentations admissibles irr\'eductibles $\sigma$ de $H(F)$ et $\pi$ de $G(F)$, l'expression ci-dessus est absolument convergente.

(ii) Si $\pi$ est supercuspidale, on a l'\'egalit\'e $m(\sigma,\pi)=m_{geom}(\sigma,\pi)$.}

Ce th\'eor\`eme est, lui, ind\'ependant de toute hypoth\`ese sur les $L$-paquets. Indiquons comment on d\'eduit le premier th\'eor\`eme du second. R\'etablissons les indices $i$ et $a$, posons 
$$m(\Sigma_{i},\Pi_{i})=\sum_{(\sigma,\pi)\in \Sigma_{i}\times \Pi_{i}}m(\sigma,\pi)$$
et d\'efinissons de m\^eme $m(\Sigma_{a},\Pi_{a})$. A l'aide du second th\'eor\`eme, ces termes se calculent comme des sommes index\'ees par des ensembles de tores ${\cal T}_{i}$ et ${\cal T}_{a}$. On peut regrouper ces tores selon leur classe de conjugaison stable. Il y a une correspondance entre classes de conjugaison stable dans ${\cal T}_{a}$ et classes de conjugaison stable dans ${\cal T}_{i}$. Cette correspondance est en fait une injection du premier ensemble de classes dans le second et c'est presqu'une surjection: l'unique classe dans ${\cal T}_{i}$ qui n'est pas dans l'image est la classe r\'eduite au tore $T=\{1\}\in {\cal T}_{i}$. Les formules de transfert de caract\`eres de $L$-paquets contiennent des signes, dont le produit est $-1$.  On en d\'eduit que, pour toute  classe de conjugaison stable  $\{T\}_{a}\subset{\cal T}_{a}$, la contribution de cette classe \`a $m(\Sigma_{a},\Pi_{a})$ est l'oppos\'e de la contribution \`a $m(\Sigma_{i},\Pi_{i})$ de la classe de conjugaison stable dans ${\cal T}_{i}$ image de $\{T\}_{a}$. Alors seul le tore $\{1\}\in {\cal T}_{i}$ contribue de fa\c{c}on non nulle \`a la somme $m(\Sigma_{a},\Pi_{a})+m(\Sigma_{i},\Pi_{i})$. A l'aide d'un r\'esultat de Rodier, cette contribution du tore $\{1\}$ s'interpr\`ete comme le produit des nombres d'\'el\'ements de $\Sigma_{i}$, resp. $\Pi_{i}$ qui admettent un mod\`ele de Whittaker relatif \`a une certaine orbite nilpotente r\'eguli\`ere. D'apr\`es  (3), ces nombres sont \'egaux \`a $1$. On obtient
$$m(\Sigma_{a},\Pi_{a})+m(\Sigma_{i},\Pi_{i})=1$$
d'o\`u le premier th\'eor\`eme. Remarquons que l'apparition d'un signe n\'egatif dans les formules de transfert, qui est cruciale pour le calcul ci-dessus, est probablement r\'eminiscente de fait que le produit des $L$-groupes ${^LH}\times{^LG}$  a une repr\'esentation naturelle qui est symplectique.

La preuve du second th\'eor\`eme est plus compliqu\'ee. Appelons quasi-caract\`ere sur $G(F)$ une fonction $\theta$ d\'efinie presque partout sur $G(F)$, invariante par conjugaison et poss\'edant un d\'eveloppement de la forme (4) au voisinage de tout point semi-simple. Pour une fonction $f\in C_{c}^{\infty}(G(F))$, disons que $f$ est tr\`es cuspidale si, pour tout sous-groupe parabolique propre $P=MU$ de $G$ (avec des notations standard), et pour tout $m\in M(F)$, on a l'\'egalit\'e
$$\int_{U(F)}f(mu)du=0.$$
Pour tout entier $N\in {\mathbb N}$, on d\'efinit une fonction $\kappa_{N}$ sur $G(F)$. C'est l'image r\'eciproque de la fonction caract\'eristique d'un sous-ensemble compact de $H(F)U(F)\backslash G(F)$, qui devient de plus en plus grand quand $N$ tend vers l'infini. Soient $\theta$ un quasi-caract\`ere sur $H(F)$ et $f\in C_{c}^{\infty}(G(F))$ une fonction tr\`es cuspidale. On pose
$$I_{N}(\theta,f)=\int_{H(F)U(F)\backslash G(F)}\int_{H(F)}\int_{U(F)}\theta(h)f(g^{-1}hug)\xi(u)du\,dh\,\kappa_{N}(g)dg.$$
La plus grande partie de l'article consiste \`a prouver que cette expression a une limite quand $N$ tend vers l'infini et \`a calculer cette limite. Celle-ci est, comme l'expression (5) ci-dessus, une somme sur les tores $T\in {\cal T}$ d'int\'egrales sur $T(F)$ de fonctions d\'eduites de $\theta$ et $f$. Cf. 7.8 pour un \'enonc\'e pr\'ecis. L'expression $I_{N}(\theta,f)$ ressemble beaucoup \`a celles qui interviennent dans la partie g\'eom\'etrique de la formule des traces locale d'Arthur ([A3]). D'ailleurs, pour l'\'etudier, on s'inspire largement des m\'ethodes d'Arthur. Il y a toutefois une diff\'erence importante entre les deux situations. Dans la formule des traces locale, il n'y a pas de probl\`eme de singularit\'es. La formule finale ne fait intervenir que des points r\'eguliers du groupe. En particulier, si on se limite \`a des fonctions dont le support est form\'e d'\'el\'ements elliptiques r\'eguliers, la partie g\'eom\'etrique de la formule des traces locale est essentiellement triviale. Ici, il y a des singularit\'es. Pour un \'el\'ement semi-simple $x\in H(F)$, le groupe $G_{x}$ est en g\'en\'eral plus gros que $H_{x}$ et on peut dire que la singularit\'e du probl\`eme cro\^{\i}t en m\^eme temps que $dim(G_{x})-dim(H_{x})$. L'\'etude  de $I_{N}(\theta,f)$ passe donc par une \'etude locale. On commence par se ramener au cas o\`u $\theta$ et $f$ ont des supports concentr\'es dans des voisinages invariants par conjugaison d'un point semi-simple  $x\in H(F)$. Une m\'ethode de descente imit\'ee d'Harish-Chandra ram\`ene alors le probl\`eme a un probl\`eme similaire, o\`u les fonctions $\theta$ et $f$ vivent cette fois sur les alg\`ebres de Lie $\mathfrak{h}_{x}(F)$ et $\mathfrak{g}_{x}(F)$. Parce que $\theta$ est un quasi-caract\`ere, on peut ensuite exprimer  l'avatar de $I_{N}(\theta,f)$ en fonction de la transform\'ee de Fourier de $f$. Il s'av\`ere qu'apr\`es cette transformation, l'expression converge beaucoup mieux. On peut maintenant prouver l'existence d'une limite et calculer celle-ci par des m\'ethodes similaires \`a celles d'Arthur. Le second th\'eor\`eme ci-dessus s'en d\'eduit en rempla\c{c}ant $\theta$ par $\theta_{\check{\sigma}}$ et $f$ par un coefficient de $\pi$. On montre en effet facilement que $m(\sigma,\pi)$ est essentiellement la limite de $I_{N}(\theta,f)$ quand $N$ tend vers l'infini.

Les trois premi\`eres sections sont consacr\'ees aux notations et \`a divers rappels d'analyse harmonique. Les sections 4 \`a 6 \'etablissent les propri\'et\'es qui nous seront utiles des quasi-caract\`eres et des fonctions tr\`es cuspidales. Les sections 7 \`a 12 sont consacr\'ees \`a l'\'etude de l'expression $I_{N}(\theta,f)$ d\'efinie ci-dessus et au calcul de sa limite. La preuve des deux th\'eor\`emes est donn\'ee dans la section 13.

 \bigskip

\section{Notations et premi\`eres d\'efinitions}

\bigskip

\subsection{Groupes}

Soit $F$ un corps local non archim\'edien de caract\'eristique nulle. On en  fixe une cl\^oture alg\'ebrique $\bar{F}$. On note $val_F$ et $\vert .\vert _F$ les valuation et  valeur absolue usuelles de $F$ et on note de la m\^eme fa\c{c}on leurs prolongements \`a $\bar{F}$. On note$\mathfrak{o}_{F}$ l'anneau des entiers de $F$, ${\mathbb F}_{q}$ son corps r\'esiduel et on fixe une uniformisante $\varpi_{F}$.

Tous les groupes alg\'ebriques sont suppos\'es d\'efinis sur $F$. Soit $G$ un groupe alg\'ebrique r\'eductif connexe. On note aussi $G$ son groupe de points sur $\bar{F}$, c'est-\`a-dire $G=G(\bar{F})$. On note $A_G$ le plus grand tore d\'eploy\'e central dans $G$, $X(G)$ le groupe des caract\`eres d\'efinis sur $F$ de $G$, ${\cal A}_G=Hom(X(G),{\mathbb R})$ et ${\cal A}^*=X(G)\otimes_{\mathbb Z}{\mathbb R}$ le dual de ${\cal A}$. On d\'efinit l'homomorphisme $H_G:G(F)\to {\cal A}_G$ par $H_G(g)(\chi)=log(\vert \chi(g)\vert _F)$ pour tous $g\in G(F)$ et $\chi\in X(G)$. On note $\mathfrak{g}$ l'alg\`ebre de Lie de $G$ et
$$\begin{array}{ccc}G\times\mathfrak{g}&\to&\mathfrak{g}\\(g,X)&\mapsto&gXg^{-1}\\ \end{array}$$
l'action adjointe. On appelle L\'evi de $G$ un sous-groupe $M$ tel qu'il existe un sous-groupe parabolique $P$ de $G$ (d\'efini sur $F$) de sorte que $M$ soit une composante de L\'evi de $P$. Pour un tel L\'evi, on note ${\cal P}(M)$ l'ensemble des sous-groupes paraboliques de $G$ de composante de L\'evi $M$, ${\cal L}(M)$ celui des L\'evi de $G$ contenant $M$ et ${\cal F}(M)$ celui des sous-groupes paraboliques de $G$ contenant $M$. Pour $Q\in {\cal F}(M)$, on notera sans plus de commentaire $Q=LU$ la d\'ecomposition de $Q$ en sa composante de L\'evi $L$ contenant $M$ et son radical unipotent $U$. Il y a une d\'ecomposition naturelle ${\cal A}_M={\cal A}_M^G\oplus{\cal A}_G$.  On note $proj_M^G$ et $proj_G$ les projections sur chacun des facteurs. Le sous-espace ${\cal A}_M^G$ est engendr\'e par l'ensemble $\check{\Sigma}_M$ des coracines indivisibles. A un \'el\'ement $P$ de ${\cal P}(M)$ est associ\'e une chambre positive ${\cal A}_P^+\subset {\cal A}_M$ et un sous-ensemble de coracines simples $\check{\Delta}_P\subset \check{\Sigma}_M$.  Bruhat et Tits ont d\'efini la notion de sous-groupe compact sp\'ecial de $G(F)$. Si $K$ est un tel sous-groupe et $M$ est un L\'evi de $G$, on dit que $K$ est en bonne position relativement \`a $M$ s'il existe un sous-tore d\'eploy\'e maximal $A\subset M$ de sorte que $K$ fixe un point sp\'ecial de l'appartement associ\'e \`a $A$ dans l'immeuble de $G$. Supposons qu'il en soit ainsi et soit $P=MU\in {\cal P}(M)$. On d\'efinit la fonction $H_P:G(F)\to {\cal A}_M$ par $H_P(g)=H_M(m)$ pour $g=muk\in G(F)$, avec $m\in M(F)$, $u\in U(F)$, $k\in K$. Suivant Harish-Chandra, on d\'efinit une fonction hauteur $\vert \vert .\vert \vert $ sur $G(F)$, \`a valeurs dans ${\mathbb R}_{\geq 1}=\{x\in {\mathbb R}; x\geq 1\}$ et (en modifiant l\'eg\`erement la d\'efinition d'Harish-Chandra) une fonction $\sigma$ par $\sigma(g)=sup(1,log(\vert \vert g\vert \vert ))$. On d\'efinit une fonction sur $\mathfrak{g}(F)$, \'egalement not\'ee $\sigma$, de la fa\c{c}on suivante. On fixe une base de $\mathfrak{g}(F)$ sur $F$. Pour $X\in \mathfrak{g}(F)$, on pose $\sigma(X)=sup(1,sup\{-val_{F}(X_{i})\})$, o\`u les $X_{i}$ sont les coordonn\'ees de $X$.  Le cas \'ech\'eant, on ajoutera des exposants $^{G}$ aux notations que l'on vient d'introduire pour pr\'eciser le groupe ambiant.

Soit $G$ un groupe. On note $Z_G$ son centre. Soit $A$ un ensemble muni d'une action de $G$. Pour un sous-ensemble $B\subset A$, on note $Z_G(B)$ le centralisateur de $B$ dans $G$ et $Norm_G(B)$ le normalisateur. Si $B=\{x\}$, on note simplement $Z_G(x)=Z_G(\{x\})$. Quand $A=G$, on suppose implicitement que l'action de $G$ est l'action par conjugaison. De m\^eme si $G$ est un groupe alg\'ebrique lin\'eaire et $A=\mathfrak{g}$ est son alg\`ebre de Lie. Pour une fonction $f$ sur $A$ et pour $g\in G$, on note $^gf$ la fonction $a\mapsto f(g^{-1}(a))$. 

Quand $G$ est un groupe alg\'ebrique lin\'eaire, on note $G^0$ sa composante neutre. Pour $x\in G$, resp. $X\in \mathfrak{g}$, on note $G_x=Z_G(x)^0$, resp. $G_{X}=Z_{G}(X)^0$, la composante neutre du centralisateur de $x$, resp. $X$.

Soit $G$ un groupe r\'eductif connexe. On note $G_{ss}$ l'ensemble de ses \'el\'ements semi-simples et $G_{reg}$ le sous-ensemble des \'el\'ements semi-simples r\'eguliers. On d\'efinit de m\^eme $\mathfrak{g}_{ss}$ et $\mathfrak{g}_{reg}$. Pour $x\in G_{ss}(F)$, l'op\'erateur $ad(x)-1$ est d\'efini et inversible sur $\mathfrak{g}(F)/\mathfrak{g}_{x}(F)$, on pose:
$$D^G(x)=\vert det(ad(x)-1)_{\vert  \mathfrak{g}(F)/\mathfrak{g}_{x}(F)})\vert _{F}.$$
De m\^eme, pour $X\in \mathfrak{g}_{ss}(F)$, on pose:
$$D^G(X)=\vert det(ad(X)_{\vert \mathfrak{g}(F)/\mathfrak{g}_{X}(F)})\vert _{F}.$$
Pour tout sous-ensemble $\Gamma\subset G(F)$, on pose $\Gamma^G=\{g^{-1}\gamma g; g\in G(F), \gamma\in \Gamma\}$.
On dit qu'un sous-ensemble $\Omega\subset G(F)$ est compact modulo conjugaison s'il existe un sous-ensemble compact $\Gamma\subset G(F)$ tel que $\Omega\subset \Gamma^G$

\bigskip

\subsection{ Mesures}

On fixe pour tout l'article un caract\`ere continu et non trivial $\psi:F\to {\mathbb C}^{\times}$.  Soit $G$ un groupe r\'eductif connexe. On munit $\mathfrak{g}(F)$ d'une forme bilin\'eaire sym\'etrique non d\'eg\'en\'er\'ee $<.,.>$ invariante par conjugaison par $G(F)$.   Pour tout ensemble topologique $X$ totalement discontinu, on note $C_{c}^{\infty}(X)$ l'espace des fonctions sur $X$, \`a valeurs dans ${\mathbb C}$, localement constantes et \`a support compact. On d\'efinit la transformation de Fourier $f\mapsto \hat{f}$ de $C_{c}^{\infty}(\mathfrak{g}(F))$ dans lui-m\^eme par
$$\hat{f}(X)=\int_{\mathfrak{g}(F)} f(Y)\psi(<X,Y>)dY,$$
o\`u $dY$ est la mesure de Haar autoduale, c'est-\`a-dire telle que $\hat{\hat{f}}(X)=f(-X)$. L'espace $\mathfrak{g}(F)$ sera toujours muni de cette mesure. Si $H$ est un sous-groupe r\'eductif de $G$, le m\^eme proc\'ed\'e munit $\mathfrak{h}(F)$ d'une mesure.

On note $Nil(\mathfrak{g})$ l'ensemble des orbites nilpotentes. Soit ${\cal O}$ une telle orbite. Pour $X\in {\cal O}$, la forme bilin\'eaire $(Y,Z)\mapsto <X, [Y,Z] >$ sur $\mathfrak{g}(F)$ se descend en une forme symplectique sur $\mathfrak{g}(F)/\mathfrak{g}_{X}(F)$, c'est-\`a-dire sur l'espace tangent \`a ${\cal O}$ au point $X$. Ainsi, ${\cal O}$ est muni d'une structure de vari\'et\'e $F$-analytique symplectique et on en d\'eduit une mesure "autoduale" sur ${\cal O}$. Cette mesure est invariante par conjugaison par $G(F)$.

Appelons $G$-domaine dans $G(F)$, resp. $\mathfrak{g}(F)$, un sous-ensemble de $G(F)$, resp. $\mathfrak{g}(F)$, qui est ouvert, ferm\'e et invariant par conjugaison. On sait que l'on peut d\'efinir une application exponentielle $exp:\omega\to \Omega$, o\`u $\omega$ est un certain $G$-domaine dans $\mathfrak{g}(F)$ contenant $0$ et $\Omega$ un certain $G$-domaine dans $G(F)$ contenant $1$. Cette application est un hom\'eomorphisme \'equivariant pour les actions de $G(F)$. On a d\'ej\`a muni $\mathfrak{g}(F)$ d'une mesure et on munit $G(F)$ de la mesure de Haar telle que le Jacobien de l'exponentielle soit \'egal \`a $1$ au point $0\in \mathfrak{g}(F)$. On d\'efinit de m\^eme une mesure de Haar sur $H(F)$ pour tout sous-groupe r\'eductif $H$ contenu dans $G$. Si $K$ est un sous-groupe compact sp\'ecial de $G(F)$, on munit $K$ de la mesure de Haar  de masse totale $1$. Soient $M$ un L\'evi de $G$ et $P=MU\in {\cal P}(M)$. On doit munir $U(F)$ d'une mesure de Haar. On sera toujours dans l'une des situations suivantes. Ou bien le choix de la mesure sera sans importance et on ne la pr\'ecisera pas. Ou bien sera fix\'e un sous-groupe compact sp\'ecial $K$ de $G(F)$ en bonne position relativement \`a $M$. Dans ce cas on choisira la mesure telle que, pour toute $f\in C_{c}^{\infty}(G(F))$, on ait l'\'egalit\'e:
$$\int_{G(F)}f(g)dg=\int_{K}\int_{U(F)}\int_{M(F)}f(muk)dm\,du\,dk.$$
Autrement dit, de sorte que l'on ait l'\'egalit\'e:
$$mes(K,dg)=mes(K\cap M(F),dm)mes(K\cap U(F),du),$$
avec une notation \'evidente. En inversant le proc\'ed\'e ci-dessus, on munit aussi $\mathfrak{u}(F)$ d'une mesure. Dans la situation ci-dessus, pour $f\in C_{c}^{\infty}(G(F))$, on d\'efinit $f_{P}\in C_{c}^{\infty}(M(F))$ par:
$$f_{P}(m)=\delta_{P}(m)^{1/2} \int_{K}\int_{U(F)}f(muk)du\,dk,$$
o\`u $\delta_{P}$ est le module usuel.

Soit $T$ un sous-tore de $G$. Le groupe $T(F)$ est muni d'une mesure par la d\'efinition ci-dessus, notons-la $dt$. Il y a une autre mesure de Haar qui intervient naturellement dans la th\'eorie, que l'on  note $d_{c}t$, et qui est d\'efinie de la fa\c{c}on suivante. Si $T$ est d\'eploy\'e, le sous-groupe compact maximal de $T(F)$ est de volume $1$ pour $d_{c}t$. En g\'en\'eral, $d_{c}t$ est compatible avec la mesure que l'on vient de d\'efinir sur $A_{T}(F)$ et avec la mesure sur $T(F)/A_{T}(F)$ de masse totale $1$. Pour \'eviter les confusions, nous n'utiliserons que la mesure $dt$, mais il sera n\'ecessaire d'introduire dans nos formules la constante $\nu(T)$ d\'efinie par $d_{c}t=\nu(T)dt$.

Soit $M$ un L\'evi de $G$. On munit ${\cal A}_{M}^G$ de la mesure pour laquelle le quotient 
$${\cal A}_{M}^G/proj_{M}^G(H_{M}(A_{M}(F)))$$
 est de volume $1$.

\bigskip

\section{ Int\'egrales orbitales pond\'er\'ees}

\bigskip

\subsection{ $(G,M)$-familles}

Un groupe r\'eductif connexe $G$ est fix\'e pour toutes les sections 2 \`a 6. On fixe aussi une forme bilin\'eaire sur $\mathfrak{g}(F)$ comme en 1.2. Soit $M$ un L\'evi de $G$. Arthur a introduit la notion de $(G,M)$-famille: c'est une famille $(c_{P})_{P\in {\cal P}(M)}$ de fonctions $C^{\infty}$ sur $i{\cal A}_{M}^*$ (o\`u $i=\sqrt{-1}$) v\'erifiant une certaine condition de compatibilit\'e ([A1] p.36). Consid\'erons une telle $(G,M)$-famille. On sait lui associer un nombre complexe $c_{M}$ ([A1] p.37). On a besoin pour cela d'une mesure sur ${\cal A}_{M}$: on l'a fix\'ee dans la section pr\'ec\'edente. Soit $L\in {\cal L}(M)$. On d\'eduit de notre $(G,M)$-famille une $(G,L)$-famille, on note $c_{L}$ le nombre qui lui est associ\'e. Soit $Q\in {\cal P}(L)$. On d\'eduit aussi de la famille de d\'epart une $(L,M)$-famille dont on note $c_{L}^Q$ le nombre associ\'e.

Soit $(Y_{P})_{P\in {\cal P}(M)}$ une famille d'\'el\'ements de ${\cal A}_{M}$. On dit qu'elle est $(G,M)$-orthogonale, resp. et positive, si elle v\'erifie la condition suivante. Soient $P$ et $P'$ deux \'el\'ements adjacents de ${\cal P}(M)$. Il y a une unique coracine $\check{\alpha}$ telle que $\check{\alpha}\in \check{\Delta}_{P}$ et $-\check{\alpha}\in \check{\Delta}_{P'}$. On demande que $Y_{P}-Y_{P'}\in {\mathbb R}\check{\alpha}$, resp. $Y_{P}-Y_{P'}\in {\mathbb R}_{\geq 0}\check{\alpha}$. Pour $P\in {\cal P}(M)$, d\'efinissons une fonction $c_{P}$ sur $i{\cal A}_{M}^*$ par $c_{P}(\lambda)=e^{-\lambda(Y_{P})}$. Supposons que la famille $(Y_{P})_{P\in {\cal P}(M)}$ soit $(G,M)$-orthogonale. Alors la famille $(c_{P})_{P\in {\cal P}(M)}$ est une $(G,M)$-famille. Soit $L\in {\cal L}(M)$. La $(G,L)$-famille d\'eduite de cette $(G,M)$-famille est associ\'ee \`a la famille de points $(Y_{Q})_{Q\in {\cal P}(L)}$ ainsi d\'efinie: $Y_{Q}=proj_{L}(Y_{P})$ pour n'importe quel $P\in {\cal P}(M)$ tel que $P\subset Q$. De m\^eme, soit $Q\in {\cal P}(L)$. Alors la $(L,M)$-famille d\'eduite de notre $(G,M)$-famille est associ\'ee \`a la famille de points $(Y_{P'})_{P'\in {\cal P}^L(M)}$ ainsi d\'efinie: $Y_{P'}=Y_{P}$, o\`u $P$ est l'unique \'el\'ement de ${\cal P}(M)$ tel que $P\subset Q$ et $P\cap L=P'$.

\bigskip

\subsection{Formules de descente}

Soient  $M$ un L\'evi de $G$, $(c_{P})_{P\in {\cal P}(M)}$ et  $(d_{P})_{P\in {\cal P}(M)}$ deux $(G,M)$-familles. Pour $P\in {\cal P}(M)$, posons $(cd)_{P}=c_{P}d_{P}$. Alors  $((cd)_{P})_{P\in {\cal P}(M)}$
est encore une $(G,M)$-famille. On a une \'egalit\'e ([A2] corollaire 7.4):
$$(1)\qquad  (cd)_{M}=\sum_{L,L'\in {\cal L}(M)}d_{M}^G(L,L')c_{M}^Qd_{M}^{Q'}.$$
Le terme $d_{M}^G(L,L')$ est un r\'eel positif ou nul, qui est non nul si et seulement si
$${\cal A}_{M}^G={\cal A}_{L}^G\oplus {\cal A}_{L'}^G.$$
On a $d_{M}^G(M,G)=d_{M}^G(G,M)=1$.
On doit fixer un param\`etre auxiliaire $\xi\in {\cal A}_{M}^G$, en position g\'en\'erale. Pour $L,L'$ v\'erifiant la condition pr\'ec\'edente, notons $\xi_{L}$ et $\xi_{L'}$ les projections de $\xi$ sur chacun des facteurs. Alors $Q$ est l'unique \'el\'ement de ${\cal P}(L)$ tel que $\xi_{L}\in {\cal A}_{Q}^+$ et $Q'$ est d\'efini de fa\c{c}on similaire.

Supposons que  $(d_{P})_{P\in {\cal P}(M)}$ est associ\'ee \`a une famille $(G,M)$-orthogonale de points $(Y_{P})_{P\in {\cal P}(M)}$. Alors:
$$(2) \qquad (cd)_{M}=\sum_{Q\in {\cal F}(M)}c_{M}^Qu_{Q}(Y_{Q}),$$
o\`u, pour $Q=LU$, $u_{Q}$ est une fonction sur ${\cal A}_{L}$, bien s\^ur ind\'ependante de nos $(G,M)$-familles. La fonction $u_{G}$ est constante de valeur $1$. La formule r\'esulte de [A1] (6.3) et lemme 6.3.

Pour une seule $(G,M)$-famille  $(c_{P})_{P\in {\cal P}(M)}$ et pour $L\in {\cal L}(M)$, on a aussi:
$$(3) \qquad c_{L}=\sum_{L'\in {\cal L}(M)}d_{M}^G(L,L')c_{M}^{Q'},$$
avec les m\^emes d\'efinitions qu'en (1).

\bigskip

\subsection{ Int\'egrales orbitales pond\'er\'ees}

Soient $M$ un L\'evi de $G$ et $K$ et sous-groupe compact sp\'ecial de $G(F)$ en bonne position relativement \`a $M$. Pour $g\in G(F)$, la famille de points $(H_{P}(g))_{P\in {\cal P}(M)}$ est $(G,M)$-orthogonale et positive. On note $(v_{P}(g))_{P\in {\cal P}(M)}$ la $(G,M)$-famille associ\'ee et $v_{M}(g)$ le nombre associ\'e \`a cette $(G,M)$-famille. La fonction $g\mapsto v_{M}(g)$ est invariante \`a gauche par $M(F)$ et \`a droite par $K$.

Soient $f\in C_{c}^{\infty}(G(F))$ et $x\in M(F)\cap G_{reg}(F)$. On d\'efinit l'int\'egrale orbitale pond\'er\'ee
$$J_{M}(x,f)=D^G(x)^{1/2}\int_{G_{x}(F)\backslash G(F)}f(g^{-1}xg)v_{M}(g)dg.$$
L'int\'egrale a un sens puisque $G_{x}=M_{x}\subset M$.

\ass{Lemme}{(i) Soit $f\in C_{c}^{\infty}(G(F))$. La fonction $x\mapsto J_{M}(x,f)$ d\'efinie sur $M(F)\cap G_{reg}(F)$ est localement constante et invariante par conjugaison par $M(F)$. L'adh\'erence dans $M(F)$ de son support est compacte modulo conjugaison.

(ii) Il existe un entier $k\geq 0$ et, pour toute $f\in C_{c}^{\infty}(G(F))$, il existe $c>0$ de sorte que l'on ait l'in\'egalit\'e:
$$\vert J_{M}(x,f)\vert  \leq c(1+\vert log\,D^G(x)\vert )^k$$
pour tout $x\in M(F)\cap G_{reg}(F)$.}

Preuve. Le (i) est \'evident. Le (ii) est d\^u \`a Arthur mais nous allons rappeler la d\'emonstration car nous l'utiliserons plus loin. D'apr\`es le (i), on peut fixer un sous-tore maximal $T$ de $M$, un sous-ensemble compact $\omega\subset T(F)$ et se contenter de majorer $\vert J_{M}(x,f)\vert $ pour $x\in \omega$.  Fixons une norme sur ${\cal A}_{M}$. Il existe $c>0$ tel que, pour tout $P\in {\cal P}(M)$ et tout $g\in G(F)$, on ait l'in\'egalit\'e $\vert H_{P}(g)\vert \leq c\sigma(g)$. Par construction, $v_{M}(g)$ est polynomial en les $H_{P}(g)$, il y a donc un entier $k\geq 0$ et $c>0$ tel que $v_{M}(g)\leq c\sigma(g)^k$. Posons $\sigma_{T}(g)=inf\{\sigma(tg); t\in T(F)\}$. Puisque $v_{M}(g)$ est invariante \`a gauche par $T(F)\subset M(F)$, on a m\^eme  $v_{M}(g)\leq c\sigma_{T}(g)^k$. Rappelons le lemme 4.2 de [A3], qui pr\'ecise un r\'esultat de Harish-Chandra. Pour tous sous-ensembles compacts $\Omega\subset T(F)$ et $\Gamma\subset G(F)$, il existe $c>0$ de sorte que, pour tout $x\in \Omega$ et tout $g\in G(F)$ tels que $g^{-1}xg\in \Gamma$, on ait l'in\'egalit\'e:
$$(1) \qquad \sigma_{T}(g)\leq c(1+\vert log\,D^G(x)\vert ).$$
On applique cela \`a $\Omega=\omega$ et au support $\Gamma$ de $f$. Alors pour $x\in \omega\cap G_{reg}(F)$, on peut majorer le terme $v_{M}(g)$ intervenant dans la d\'efinition de $J_{M}(x,f)$ par $c(1+\vert log\,D^G(x)\vert )^k$, o\`u $c$ d\'epend de $f$ mais pas $k$. On obtient:
$$\vert  J_{M}(x,f)\vert\leq c (1+\vert log\,D^G(x)\vert )^k D^G(x)^{1/2}\int_{G_{x}(F)\backslash G(F)}\vert f(g^{-1}xg)\vert dg$$
$$\qquad  \leq c(1+\vert log\,D^G(x)\vert )^k J_{G}(x,\vert f\vert ).$$
D'apr\`es [HCvD] th\'eor\`eme 13, $J_{G}(x,\vert f\vert )$ est born\'e sur $\omega\cap G_{reg}(F)$ et cela conclut. $\square$

\bigskip

\subsection{ Formule des traces locale}

Soient $M_{min}$ un L\'evi minimal de $G$ et $K$ un sous-groupe compact sp\'ecial de $G(F)$ en bonne position relativement \`a $M_{min}$. Les d\'efinitions du paragraphe pr\'ec\'edent se descendent \`a l'alg\`ebre de Lie: pour tous $M\in {\cal L}(M_{min})$, $X\in \mathfrak{m}(F)\cap \mathfrak{g}_{reg}(F)$, $f\in C_{c}^{\infty}(\mathfrak{g}(F))$, on d\'efinit l'int\'egrale orbitale pond\'er\'ee $J_{M}(X,f)$.

Soient $f,f'\in C_{c}^{\infty}(\mathfrak{g}(F))$. Pour $M\in {\cal L}(M_{min})$ et $X\in\mathfrak{m}(F)\cap \mathfrak{g}_{reg}(F)$, posons:
$$J_{M}(X,f,f')=\sum_{L,L'\in {\cal L}(M)}d_{M}^G(L,L')J_{M}^L(X,f_{\bar{Q}})J_{M}^{L'}(X,f'_{Q'})$$
Les d\'efinitions sont les m\^emes qu'en 2.2(1); $\bar{Q}$ est le sous-groupe parabolique oppos\'e \`a $Q$.  On note $W^M=Norm_{M}(M_{min})/M_{min}$, $a_{M}=dim(A_{M})$.  Pour un sous-tore maximal $T$ de $M$, on pose $W(M,T)=Norm_{M(F)}(T)/T(F)$. On dit que $T$ est elliptique dans $M$ si $A_{T}=A_{M}$. On fixe un ensemble ${\cal T}_{ell}(M)$ de repr\'esentants des classes de conjugaison de sous-tores maximaux de $M$, elliptiques dans $M$. Posons:
$$J(f,f')=\sum_{M\in {\cal L}(M_{min})}\vert W^M\vert \vert W^G\vert ^{-1}(-1)^{a_{G}-a_{M}}\sum_{T\in {\cal T}_{ell}(M)}\vert W(M,T)\vert ^{-1}\nu(T)^{-1}\int_{\mathfrak{t}(F)}J_{M}(X,f,f') dX.$$
Cette expression est absolument convergente en vertu du lemme 2.3(ii) et du lemme suivant.

\ass{Lemme}{Soient $V$ un espace vectoriel de dimension finie sur $F$ et $(R_{i})_{i=1,...,n}$ une famille finie de polyn\^omes non nuls sur $V$. Alors la fonction $v\mapsto \prod_{i=1,...,n}log(\vert R_{i}(v)\vert _{F})$ est localement int\'egrable sur $V$. $\square$}

\ass{Th\'eor\`eme}{Pour toutes $f,f'\in C_{c}^{\infty}(\mathfrak{g}(F))$, on a l'\'egalit\'e $J(\hat{f},f')=J(f,\hat{f}')$.}

Cf. [W1] th\'eor\`eme 5.2, qui reprenait [A3]. Il n'y a pas de $\nu(T)$ dans [W1], ce qui est d\^u au fait que les mesures sur les tores n'y sont pas les m\^emes que les n\^otres (il y a d'ailleurs aussi dans cette r\'ef\'erence une erreur dans la d\'efinition des mesures sur les espaces ${\cal A}_{M}$).

\bigskip

\subsection{La condition (H)}

On conserve les m\^emes hypoth\`eses. Pour $\varphi\in C_{c}^{\infty}(\mathfrak{g}(F))$, consid\'erons la condition:

\ass{(H)}{pour tout $M\in {\cal L}(M_{min})$, il existe $\varphi_{M}\in C_{c}^{\infty}(\mathfrak{m}(F))$ telle que $\varphi_{P}=\varphi_{M}$ pour tout $P\in {\cal P}(M)$.}

En vertu de l'\'egalit\'e $(\hat{\varphi})_{P}=(\varphi_{P}){\hat{}}$, $\varphi$ v\'erifie (H) si et seulement $\hat{\varphi}$ v\'erifie (H). En g\'en\'eral, pour tout sous-ensemble $B$ d'un ensemble $A$, notons ${\bf 1}_{B}$ la fonction caract\'eristique de $B$ dans $A$. Si $\Omega$ est un $G$-domaine dans $\mathfrak{g}(F)$ et si $\varphi$ v\'erifie (H), alors $\varphi{\bf 1}_{\Omega}$ v\'erifie aussi (H) et on a $(\varphi{\bf 1}_{\Omega})_{M}=\varphi_{M}{\bf 1}_{\Omega\cap \mathfrak{m}(F)}$ pour tout $M$.

Pour tout sous-ensemble $B\subset \mathfrak{g}(F)$, posons $B^K=\{k^{-1}Xk; k\in K, X\in B\}$. Posons:
$$\Omega=\bigcup_{M\in {\cal L}(M_{min})}\bigcup_{T\in{\cal T}_{ell}(M)}(\mathfrak{t}(F)\cap\mathfrak{g}_{reg}(F))^K.$$
C'est un ouvert de $\mathfrak{g}(F)$.

\ass{Lemme}{Soit $\varphi\in C_{c}^{\infty}(\mathfrak{g}(F))$.

(i) Supposons $Supp(\varphi)\subset \Omega$. Alors $\varphi$ v\'erifie (H).

(ii) Supposons $Supp(\varphi)\subset \mathfrak{g}_{reg}(F)$. Alors il existe une familles finie $(\varphi_{i})_{i=1,...,n}$ d'\'el\'ements de $C_{c}^{\infty}(\mathfrak{g}(F))$ et  une famille finie $(g_{i})_{i=1,...,n}$ d'\'el\'ements de $G(F)$ telles que $\varphi_{i} $ v\'erifie (H) pour tout $i$ et $\varphi=\sum_{i=1,...,n}{^{g_{i}}\varphi_{i}}$.}

Preuve. Supposons $Supp(\varphi)\subset \Omega$. Soient $P=MU\in {\cal F}(M_{min})$ et $X\in \mathfrak{m}(F)\cap \mathfrak{g}_{reg}(F)$. Par un calcul familier, on a l'\'egalit\'e:
$$\varphi_{P}(X)=D^G(X)^{1/2}D^M(X)^{-1/2}\int_{K}\int_{U(F)}\varphi(k^{-1}u^{-1}Xuk)du\,dk.$$
Soit $u\in U(F)$ pour lequel il existe $k\in K$ tel que $\varphi(k^{-1}u^{-1}Xuk)\not=0$. Alors $u^{-1}Xu\in \Omega$ et on peut fixer $L\in {\cal L}(M_{min})$, $T\in {\cal T}_{ell}(L)$, $Y\in \mathfrak{t}(F)\cap \mathfrak{g}_{reg}(F)$ et $k\in K$ de sorte que $u^{-1}Xu=kYk^{-1}$. Posons $g=uk$. On a $gYg^{-1}=X$, donc $gTg^{-1}$ est le commutant $G_{X}$ de $X$. Le plus grand tore d\'eploy\'e de $gTg^{-1}$ est $gA_{L}g^{-1}$. Celui de $G_{X}$ contient $A_{M}$. Donc $A_{M}\subset gA_{L}g^{-1}\subset gA_{M_{min}}g^{-1}$. Puisque $M$ est le commutant de $A_{M}$, on a $gA_{M_{min}}g^{-1}\subset M$. Alors $A_{M_{min}}$ et $gA_{M_{min}}g^{-1}$ sont deux tores d\'eploy\'es maximaux de $M$, ils sont donc conjugu\'es par un \'el\'ement de $M(F)$. Fixons $m\in M(F)$ tel que $mgA_{M_{min}}g^{-1}m^{-1}=A_{M_{min}}$. D'apr\`es Bruhat et Tits, le normalisateur $Norm_{G(F)}(A_{M_{min}})$ est contenu dans $M_{min}(F)K$. Donc $mg\in M_{min}(F)K$, puis $g\in M(F)K$ et enfin $u\in M(F)K$. Parce que $K$ est en bonne position relativement \`a $M$, cela entra\^{\i}ne $u\in U(F)\cap K$. Donc:
 $$\varphi_{P}(X)=D^G(X)^{1/2}D^M(X)^{-1/2}\int_{K}\int_{U(F)\cap K}\varphi(k^{-1}u^{-1}Xuk)du\,dk.$$
L'int\'egrale sur $U(F)\cap K$ est absorb\'ee par celle sur $K$, on obtient:
$$\varphi_{P}(X)=D^G(X)^{1/2}D^M(X)^{-1/2}(mes(U(F)\cap K))\int_{K}\varphi(k^{-1}Xk)dk.$$
Comme on l'a remarqu\'e en 1.2, $mes(U(F)\cap K)$ ne d\'epend pas du sous-groupe parabolique $P\in {\cal P}(M)$. L'expression ci-dessus n'en d\'epend donc pas non plus et c'est la condition pour que $\varphi$ v\'erifie (H).

Supposons maintenant que $Supp(\varphi)\subset \mathfrak{g}_{reg}(F)$. Pour tout $X\in Supp(\varphi)$, fixons $g_{X}\in G(F)$ tel que
$$g_{X}^{-1}Xg_{X}\in \bigcup_{M\in {\cal L}(M_{min})}\bigcup_{T\in{\cal T}_{ell}(M)}\mathfrak{t}(F)\cap\mathfrak{g}_{reg}(F),$$
puis un voisinage $\omega_{X}$ de $X$ tel que $g_{X}^{-1}\omega_{X}g_{X}\subset \Omega$. Une partition de l'unit\'e nous ram\`ene au cas o\`u $Supp(\varphi)$ est contenu dans un tel voisinage $\omega_{X}$. Dans ce cas, $\varphi={^{g_{X}}\varphi'}$, o\`u $\varphi'={^{g_{X}^{-1}}\varphi}$. Mais $Supp(\varphi')\subset \Omega$, donc $\varphi'$ v\'erifie (H). $\square$

\bigskip

\subsection{Transform\'ees de Fourier d'int\'egrales orbitales, germes de Shalika}

Pour tout ${\cal O}\in Nil(\mathfrak{g})$, on d\'efinit l'int\'egrale orbitale nilpotente 
$$J_{{\cal O}}(f)=\int_{{\cal O}}f(X) dX$$
pour $f\in C_{c}^{\infty}(\mathfrak{g}(F))$, et sa transform\'ee de Fourier
$$\hat{J}_{{\cal O}}(f)=J_{{\cal O}}(\hat{f}).$$
 Pour $\lambda\in F^{\times}$, d\'efinissons $f^{\lambda}$ par $f^{\lambda}(X)=f(\lambda X)$. Notons $F^{\times 2} $ le groupe des carr\'es dans $F^{\times}$. On a l'\'egalit\'e
$$J_{{\cal O}}(f^{\lambda})=\vert \lambda\vert _{F}^{-dim({\cal O})/2}J_{{\cal O}}(f)$$
pour tout $\lambda\in F^{\times 2} $.

On pose $\delta(G)=dim(G)-dim(T)$, o\`u $T$ est n'importe quel sous-tore maximal de $G$. On sait qu'il existe une unique fonction $\Gamma_{{\cal O}}$ sur $\mathfrak{g}_{reg}(F)$, le germe de Shalika associ\'e \`a ${\cal O}$, v\'erifiant les deux conditions suivantes:

$$\Gamma_{{\cal O}}(\lambda X)=\vert \lambda\vert _{F}^{(\delta(G)-dim({\cal O}))/2}\Gamma_{{\cal O}}(X)$$
pour tous $X\in \mathfrak{g}_{reg}(F)$ et $\lambda\in F^{\times 2} $;

pour toute $f\in C_{c}^{\infty}(\mathfrak{g}(F))$, il existe un voisinage $\omega$ de $0$ dans $\mathfrak{g}(F)$ tel que:
$$J_{G}(X,f)=\sum_{{\cal O}\in Nil(\mathfrak{g})}\Gamma_{{\cal O}}(X)J_{{\cal O}}(f)$$
pour tout $X\in \omega\cap \mathfrak{g}_{reg}(F)$.

Remarquons que $\Gamma_{{\cal O}}$ co\"{\i}ncide sur tout compact avec une int\'egrale orbitale, en particulier y est born\'e.

Il existe une unique fonction $\hat{j}$ sur $\mathfrak{g}(F)\times \mathfrak{g}(F)$, localement int\'egrable, localement constante sur $\mathfrak{g}_{reg}(F)\times \mathfrak{g}_{reg}(F)$, telle que, pour toute $f\in C_{c}^{\infty}(\mathfrak{g}(F))$ et tout $X\in \mathfrak{g}_{reg}(F)$, on ait l'\'egalit\'e:
$$J_{G}(X,\hat{f})=\int_{\mathfrak{g}(F)}f(Y)\hat{j}(X,Y)dY.$$
De m\^eme, pour ${\cal O}\in Nil(\mathfrak{g})$, il existe une unique fonction $Y\mapsto\hat{j}({\cal O},Y)$ sur $\mathfrak{g}(F)$, localement int\'egrable, localement constante sur $\mathfrak{g}_{reg}(F)$, telle que, pour toute fonction $f\in C_{c}^{\infty}(\mathfrak{g}(F))$ , on ait l'\'egalit\'e:
$$\hat{J}_{{\cal O}}(f)=\int_{\mathfrak{g}(F)}f(Y)\hat{j}({\cal O},Y) dY.$$

On a les \'egalit\'es:
$$(1) \qquad \hat{j}(\lambda X,Y)=\vert \lambda\vert_{F} ^{\delta(G)/2}\hat{j}(X,\lambda Y),\,\, \hat{j}({\cal O},\lambda Y)=\vert \lambda\vert _{F}^{-dim({\cal O})/2}\hat{j}({\cal O},Y)$$
pour tous $X,Y\in \mathfrak{g}_{reg}(F)$, ${\cal O}\in Nil(\mathfrak{g})$ et $\lambda\in F^{\times 2} $.

Soient $\omega$ et $\Omega$ deux $G$-domaines dans $\mathfrak{g}(F)$ compacts modulo conjugaison. La conjecture de Howe entra\^{\i}ne l'existence d'une famille finie $(X_{i})_{i=1,...,n}$ d'\'elements de $\Omega\cap \mathfrak{g}_{reg}(F)$ et d'une famille finie $(f_{i})_{i=1,...,n}$ d'\'el\'ements de $C_{c}^{\infty}(\omega)$  v\'erifiant la condition suivante. Pour toute distribution invariante $D$ dont la transform\'ee de Fourier est \`a support dans $\Omega$ et toute $f\in C_{c}^{\infty}(\omega)$, on a l'\'egalit\'e
$$(2) \qquad D(f)=\sum_{i=1,...,n}J_{G}(X_{i},\hat{f})D(f_{i}).$$
Il en r\'esulte que, pour $X\in \Omega\cap \mathfrak{g}_{reg}(F)$ et $Y\in \omega\cap \mathfrak{g}_{reg}(F)$, on a l'\'egalit\'e
$$(3) \qquad \hat{j}(X,Y)=\sum_{i=1,...,n}\hat{j}(X_{i},Y)J_{G}(X,\hat{f}_{i}).$$
Au voisinage de $0$, on a un r\'esultat plus pr\'ecis. Soit $\omega$ un $G$-domaine de $\mathfrak{g}(F)$ compact modulo conjugaison et contenant $0$. Alors il existe un $G$-domaine $\Omega$ de $\mathfrak{g}(F)$ compact modulo conjugaison et contenant $0$ tel que, pour $X\in \Omega\cap \mathfrak{g}_{reg}(F)$ et $Y\in \omega\cap \mathfrak{g}_{reg}(F)$, on ait l'\'egalit\'e
$$(4) \qquad\hat{j}(X,Y)=\sum_{{\cal O}\in Nil(\mathfrak{g})}\Gamma_{{\cal O}}(X)\hat{j}({\cal O},Y).$$

Soient $M$ un L\'evi de $G$ et $X\in \mathfrak{m}(F)\cap \mathfrak{g}_{reg}(F)$. Pour $Y\in \mathfrak{g}_{reg}(F)$, fixons un ensemble de repr\'esentants $(Y_{i})_{i=1,...,r}$ des classes de conjugaison par $M(F)$ dans l'ensemble des \'el\'ements de $\mathfrak{m}(F)$ qui sont conjugu\'es \`a $Y$ par un \'el\'ement de $G(F)$. On v\'erifie l'\'egalit\'e
$$(5) \qquad \hat{j}^G(X,Y)D^G(Y)^{1/2}=  \sum_{i=1,...,r}\hat{j}^M(X,Y_{i})D^M(Y_{i})^{1/2} .$$

\bigskip

\section{Voisinages d'\'el\'ements semi-simples}

\bigskip

\subsection{Bons voisinages}

On fixe pour toute la section un \'el\'ement $x\in G_{ss}(F)$. On dira qu'un sous-ensemble $\omega\subset \mathfrak{g}_{x}(F)$ est un bon voisinage de $0$ s'il v\'erifie les conditions (1) \`a (7) ci-dessous.

(1) L'ensemble $\omega$ est un $G_{x}$-domaine compact modulo conjugaison, invariant par $Z_{G}(x)(F)$ et contenant $0$.

(2) L'exponentielle est d\'efinie sur $\omega$; c'est un hom\'eomorphisme \'equivariant pour la conjugaison par $Z_{G}(x)(F)$ de $\omega$ sur un $G_{x}$-domaine $exp(\omega)$ de $G_{x}(F)$.

(3) Pour tout $\lambda\in F^{\times}$ tel que $\vert \lambda\vert _{F}\leq 1$, on a $\lambda\omega\subset \omega$.

(4) On a l'\'egalit\'e 
$$\{g\in G(F); g^{-1}xexp(\omega)g\cap xexp(\omega)\neq\emptyset\}=Z_{G}(x)(F).$$

(5) Pour tout sous-ensemble compact $\Gamma\subset G(F)$, il existe un sous-ensemble compact $\Gamma'\subset G(F)$ tel que l'on ait l'inclusion:
$$\{g\in G(F); g^{-1}xexp(\omega)g\cap \Gamma\neq\emptyset\} \subset G_{x}(F)\Gamma'.$$

Fixons un r\'eel $c_{F}>0$ tel que $c_{F}^k<\vert (k+1)!\vert _{F}$ pour tout entier $k\geq 1$.

(6) Pour tout sous-tore maximal $T\subset G_{x}$, tout caract\`ere alg\'ebrique $\chi$ de $T$ et tout \'el\'ement $X\in \mathfrak{t}(F)\cap \omega$, on a l'in\'egalit\'e $\vert \chi(X)\vert_{F}<c_{F}$. 

Consid\'erons un sous-espace propre $W\subset \mathfrak{g}(F)$ pour l'op\'erateur $ad(x)$. Notons $\lambda$ la valeur propre. Soit $X\in \omega$. Alors $ad(X)$ conserve $W$. Soit $W_{X}$ un sous-espace propre de $W$ pour l'op\'erateur $ad(X)$, de valeur propre $\mu$. Alors $W_{X}$ est aussi un espace propre pour l'op\'erateur $ad(xexp(X))$, de valeur propre $\lambda exp(\mu)$.

(7) Supposons $\lambda\neq 1$. Alors $\vert \lambda exp(\mu)-1\vert _{F}=\vert \lambda-1\vert _{F}$.

De bons voisinages de $0$ existent, aussi petits que l'on veut en ce sens que, si $\omega_{0}$ est un voisinage de $0$ dans $\mathfrak{g}_{x}(F)$, il existe un bon voisinage $\omega$ de $0$ tel que $\omega\subset \omega_{0}^{G_{x}}$. Consid\'erons un bon voisinage $\omega$ de $0$. Les conditions (1) et (2) entra\^{\i}nent que l'ensemble $\Omega=(xexp(\omega))^G$ est un $G$-domaine dans $G(F)$, compact modulo conjugaison. La condition (4) entra\^{\i}ne que, pour $X\in \omega$, $Z_{G}(xexp(X))(F)\subset Z_{G}(x)(F)$ et $G_{xexp(X)}=(G_{x})_{X}\subset G_{x}$. On note simplement $G_{x,X}=(G_{x})_{X}$. La condition (6) entra\^{\i}ne que l'exponentielle de $\omega$ sur $exp(\omega)$ pr\'eserve les mesures. Plus g\'en\'eralement, elle entra\^{\i}ne qu'un certain nombre de jacobiens qui interviendront plus tard sont \'egaux \`a $1$. Les conditions (6) et (7) entra\^{\i}nent que, pour tout $X\in \omega$, on a l'\'egalit\'e
$$D^G(xexp(X))=D^G(x)D^{G_{x}}(X).$$

Une cons\'equence de la propri\'et\'e (4) est que, pour toute fonction $\varphi$ sur $\omega$, invariante par conjugaison par $Z_{G}(x)(F)$, il existe une unique fonction $f$ sur $G(F)$, invariante par conjugaison par $G(F)$, \`a support dans $\Omega$ et telle que $f(xexp(X))=\varphi(X)$ pour tout $\varphi\in \omega$. Si $\varphi$ est localement constante sur $\omega\cap\mathfrak{g}_{reg}(F)$, $f$ est localement constante sur $G_{reg}(F)$.

Le cas \'ech\'eant, on peut renforcer les conditions impos\'ees aux bons voisinages. Supposons par exemple que $G_{x}$ se d\'ecompose en le produit de deux groupes r\'eductifs connexes $G_{x}=G'\times G''$, conserv\'es chacun par $Z_{G}(x)(F)$. On peut imposer

(8) $\omega=\omega'\times \omega''$, o\`u $\omega'\subset \mathfrak{g}'(F)$ et $\omega''\subset \mathfrak{g}''(F)$ v\'erifient des conditions analogues \`a (1) et (2).

Ou bien, supposons fix\'ee une repr\'esentation alg\'ebrique $\rho$ de $G$ dans un espace vectoriel $V$ de dimension finie sur $F$. Les conditions pr\'ec\'edant (7) s'appliquent en rempla\c{c}ant $\mathfrak{g}(F)$ par $V$ et les op\'erateurs $ad(x)$, $ad(X)$ et $ad(xexp(X))$ par $\rho(x)$, $\rho(X)$, $\rho(xexp(X))$. On peut imposer une condition $(7)_{\rho}$ analogue \`a (7) pour les ensembles de valeurs propres ainsi d\'efinis.

\bigskip

\subsection{Correspondance des L\'evi}

Soit $M$ un L\'evi de $G$ contenant $x$. On a l'\'egalit\'e $M_{x}=M\cap G_{x}$ et ce groupe est un L\'evi de $G_{x}$: c'est le commutant de $A_{M}\subset G_{x}$. On a $A_{M}\subset A_{M_{x}}$. Pour $P=MU\in {\cal P}(M)$, on a les \'egalit\'es $P_{x}=P\cap G_{x}$, $U_{x}=U\cap G_{x}$, $P_{x}=M_{x}U_{x}$ et $P_{x}$ appartient \`a ${\cal P}(M_{x})={\cal P}^{G_{x}}(M_{x})$. Inversement, soit $R$ un L\'evi de $G_{x}$. Notons ${\bf R}$ le commutant de $A_{R}$ dans $G$. C'est un L\'evi de $G$. On a:

(1) $x\in {\bf R}(F)$, ${\bf R}_{x}=M$ et $A_{{\bf R}}=A_{R}$.

Preuve. Puisque $x$ commute \`a $A_{R}$, $x$ appartient \`a ${\bf R}(F)$. On sait d\'ej\`a que $A_{{\bf R}}\subset A_{{\bf R}_{x}}$. Le groupe $R$ commute \`a $A_{R}$, donc est inclus dans ${\bf R}$, puis dans ${\bf R}\cap G_{x}={\bf R}_{x}$. Cela entra\^{\i}ne $A_{{\bf R}_{x}}\subset A_{R}$. Enfin, par construction de ${\bf R}$, $A_{R}$ est inclus dans $A_{{\bf R}}$. Alors $A_{{\bf R}}=A_{{\bf R}_{x}}=A_{R}$, ce qui entra\^{\i}ne les deux derni\`eres assertions de (1) $\square$

L'application $R\mapsto {\bf R}$ est une bijection de l'ensemble des L\'evi de $G_{x}$ sur celui des L\'evi $M$ de $G$ contenant $x$ et tels que $A_{M}=A_{M_{x}}$.

\bigskip

\subsection{Descente des poids}

Fixons un L\'evi minimal $R_{min}$ de $G_{x}$, un sous-groupe compact sp\'ecial $K_{x}$ de $G_{x}(F)$ en bonne position relativement \`a $R_{min}$ et un sous-groupe compact sp\'ecial $K$ de $G(F)$ en bonne position relativement \`a ${\bf R}_{min}$. Remarquons que la bijection du paragraphe pr\'ec\'edent se restreint en une bijection de ${\cal L}^{G_{x}}(R_{min})$ sur le sous-ensemble des $M\in {\cal L}^G({\bf R}_{min})$ tels que $A_{M}=A_{M_{x}}$. Le L\'evi ${\bf R}_{min}$ n'a bien s\^ur aucune raison d'\^etre minimal dans $G$.

\ass{Lemme}{Soient $R\in {\cal L}(R_{min})$, $g\in G_{x}(F)$ et $y\in G(F)$. On a l'\'egalit\'e:
$$v_{{\bf R}}(gy)=\sum_{S\in{\cal L}(R)}\sum_{Q\in {\cal P}({\bf S})}v_{R}^{Q_{x}}(g)u_{Q}(H_{Q}(gy)-H_{Q_{x}}(g)).$$}

Preuve. Pour $P\in {\cal P}({\bf R})$ et $\lambda\in i{\cal A}_{{\bf R}}=i{\cal A}_{R}$, posons $c_{P}(g)(\lambda)=v_{P_{x}}(g)(\lambda)$. D'apr\`es [A4] p.233, la famille $(c_{P}(g))_{P\in{\cal P}({\bf R})}$ est une $(G,{\bf R})$-famille. D\'efinissons la $(G,{\bf R})$-famille $(d_{P}(g,y))_{P\in {\cal P}({\bf R})}$ par $d_{P}(g,y)=c_{P}(g)^{-1}v_{P}(gy)$. Elle est associ\'ee \`a la famille de points $(H_{P}(gy)-H_{P_{x}}(g))_{P\in {\cal P}({\bf R})}$. On a $v_{P}(gy)=c_{P}(g)d_{P}(g,y)$. D'apr\`es 2.2(2),
$$v_{{\bf R}}(gy)=\sum_{Q=LU\in {\cal F}({\bf R})}c_{{\bf R}}^Q(g)u_{Q}(H_{Q}(gy)-proj_{L}(H_{Q_{x}}(g))).$$
Soit $Q=LU\in {\cal F}({\bf R})$. D'apr\`es [A4], lemme 4.1,
$$c_{{\bf R}}^Q(g)=\sum_{S\in {\cal L}(R)}d_{R}^L({\bf R},S)v_{R}^{Q_{S}}(g).$$
  Le groupe $Q_{S}$ appartient \`a ${\cal P}(S)$ et est contenu dans $Q_{x}$. La constante $d_{R}^L({\bf R},S)$ est similaire \`a celle de 2.2. La condition $d_{R}^L({\bf R},S)\neq 0$ impose ${\cal A}_{R}^L={\cal A}_{{\bf R}}^L\oplus {\cal A}_{S}^L$. Or ${\cal A}_{R}={\cal A}_{{\bf R}}$. Donc ${\cal A}_{S}={\cal A}_{L}$. Mais alors $L={\bf S}$ et $Q_{S}=Q_{x}$. On obtient:
$$c_{{\bf R}}^Q(g)=\left\lbrace\begin{array}{cc}v_{R}^{Q_{x}}(g),&\,\,{\rm si}\,\,L={\bf S}\,\,{\rm pour\,\,un}\,\,S\in {\cal P}(R),\\ 0,&\,\,{\rm sinon.}\\ \end{array}\right.$$
Dans le cas o\`u $L={\bf S}$ pour un $S\in {\cal P}(R)$, on a $H_{Q_{x}}(g)\in {\cal A}_{S}={\cal A}_{L}$, donc $proj_{L}(H_{Q_{x}}(g))=H_{Q_{x}}(g)$. Toutes ces formules conduisent \`a celle de l'\'enonc\'e. $\square$

\bigskip

\section{Quasi-caract\`eres}

\bigskip

\subsection{Quasi-caract\`eres de $G(F)$}

Soit $\theta$ une fonction d\'efinie presque partout sur $G(F)$ et invariante par conjugaison. On dit que c'est un quasi-caract\`ere si et seulement si, pour tout $x\in G_{ss}(F)$, il existe un bon voisinage $\omega$ de $0$ dans $\mathfrak{g}_{x}(F)$ et, pour tout ${\cal O}\in Nil(\mathfrak{g}_{x})$, il existe $c_{\theta,{\cal O}}(x)\in{\mathbb C}$ de sorte que l'on ait l'\'egalit\'e
$$(1) \qquad \theta(xexp(X))=\sum_{{\cal O}\in Nil(\mathfrak{g}_{x})}c_{\theta,{\cal O}}(x)\hat{j}({\cal O},X)$$
presque partout pour $X\in \omega$. Autrement dit, pour toute $f\in C_{c}^{\infty}(\mathfrak{g}_{x}(F))$ \`a support dans $\omega$, on a l'\'egalit\'e
$$\int_{\mathfrak{g}_{x}(F)}\theta(xexp(X))f(X)dX=\sum_{{\cal O}\in Nil(\mathfrak{g}_{x})}c_{\theta,{\cal O}}(x) \hat{J}_{{\cal O}}(f).$$
Les coefficients $c_{\theta,{\cal O}}(x)$ sont uniquement d\'etermin\'es.

Soit $\theta$ un quasi-caract\`ere. Alors $\theta$ est localement int\'egrable sur $G(F)$ et localement constant sur $G_{reg}(F)$. Pour tout $G$-domaine $\Omega$ dans $G(F)$, $\theta{\bf 1}_{\Omega}$ est un quasi-caract\`ere. Soient $x\in G_{ss}(F)$ et $\omega$ un bon voisinage de $0$ dans $\mathfrak{g}_{x}(F)$. On dit que $\theta$ est d\'eveloppable dans $xexp(\omega)$ si l'\'egalit\'e (1) est v\'erifi\'ee pour $X\in \omega\cap \mathfrak{g}_{reg}(F)$.

\subsection{Quasi-caract\`eres de $\mathfrak{g}(F)$}

La d\'efinition pr\'ec\'edente s'adapte aux alg\`ebres de Lie. Soit $\theta$ une fonction d\'efinie presque partout sur $\mathfrak{g}(F)$. On dit que c'est un quasi-caract\`ere si et seulement si, pour tout $X\in \mathfrak{g}_{ss}(F)$, il existe un $G_{X}$-domaine $\omega$ dans $\mathfrak{g}_{X}(F)$, contenant $0$, et, pour tout ${\cal O}\in Nil(\mathfrak{g}_{X})$, il existe $c_{\theta,{\cal O}}(X)\in{\mathbb C}$ de sorte que l'on ait l'\'egalit\'e
$$(1) \qquad \theta(X+Y)=\sum_{{\cal O}\in Nil(\mathfrak{g}_{X})}c_{\theta,{\cal O}}(X)\hat{j}({\cal O},Y)$$
presque partout pour $Y\in \omega$. Les quasi-caract\`eres de $\mathfrak{g}(F)$ ont des propri\'et\'es similaires \`a celles des quasi-caract\`eres de $G(F)$.

Soit $\theta$ un quasi-caract\`ere. On note simplement $c_{\theta,{\cal O}}=c_{\theta,{\cal O}}(0)$ les coefficients du d\'eveloppement de $\theta$ au point $0$. Soit $\lambda\in F^{\times 2} $. Alors $\theta^{\lambda}$ est un quasi-caract\`ere. Soient $X\in \mathfrak{g}_{ss}(F)$ et $\omega$ comme ci-dessus, tel que $\theta$ soit d\'eveloppable dans $X+\omega$. Alors $\theta^{\lambda}$ est d\'eveloppable dans $\lambda^{-1}X+\lambda^{-1}\omega$. Remarquons que $\mathfrak{g}_{\lambda^{-1}X}=\mathfrak{g}_{X}$.  Pour tout ${\cal O}\in Nil(\mathfrak{g}_{X})$, on a l'\'egalit\'e
$$(2) \qquad c_{\theta^{\lambda},{\cal O}}(\lambda^{-1}X)=\vert \lambda\vert _{F}^{-dim({\cal O})/2}c_{\theta,{\cal O}}(X).$$
Cela r\'esulte imm\'ediatement des formules de 2.6.

\ass{Th\'eor\`eme}{Soit $D$ une distribution sur $\mathfrak{g}(F)$, invariante par conjugaison et \`a support compact modulo conjugaison. Alors sa transform\'ee de Fourier est la distribution associ\'ee \`a une fonction localement int\'egrable $\theta$ qui est un quasi-caract\`ere.}

C'est un r\'esultat d'Harish-Chandra. La premi\`ere assertion r\'esulte du th\'eor\`eme 4.4 de [HCDS]. La derni\`ere n'est pas tr\`es clairement \'enonc\'ee dans cette r\'ef\'erence mais r\'esulte de la preuve du th\'eor\`eme 4.4, en particulier du th\'eor\`eme 5.11 et du corollaire 6.10.

\subsection{Localisation}

On fixe un \'el\'ement $x\in G_{ss}(F)$ et un bon voisinage $\omega$ de $0$ dans $\mathfrak{g}_{x}(F)$. Soit $\theta$ un quasi-caract\`ere de $G(F)$. On d\'efinit une fonction $\theta_{x,\omega}$ sur $\mathfrak{g}_{x}(F)$ par
$$\theta_{x,\omega}(X)=\left\lbrace\begin{array}{cc}\theta(xexp(X)),&\,\,{\rm si}\,\,X\in \omega,\\ 0,&\,\,{\rm sinon}.\\ \end{array}\right.$$
Alors $\theta_{x,\omega}$ est un quasi-caract\`ere de $\mathfrak{g}_{x}(F)$. On a les \'egalit\'es $c_{\theta,{\cal O}}(xexp(X))=c_{\theta_{x,\omega},{\cal O}}(X)$ pour tout $X\in \omega\cap \mathfrak{g}_{x,ss}(F)$ et tout ${\cal O}\in Nil(\mathfrak{g}_{x,X})$. En particulier $c_{\theta,{\cal O}}(x)=c_{\theta_{x,\omega},{\cal O}}$ pour tout ${\cal O}\in Nil(\mathfrak{g}_{x})$.

Inversement, soit $\theta$ un quasi-caract\`ere de $\mathfrak{g}_{x}(F)$. Supposons que $\theta$ soit invariant par conjugaison par $Z_{G}(x)(F)$ et \`a support dans $\omega$. Soit $\boldsymbol{\theta}$ la fonction sur $G(F)$, invariante par conjugaison par $G(F)$, \`a support dans $\Omega=(xexp(\omega))^G$ et telle que $\boldsymbol{\theta}(xexp(X))=\theta(X)$ pour $X\in \omega$, cf. 3.1. Alors $\boldsymbol{\theta}$ est un quasi-caract\`ere sur $G(F)$. On a l'\'egalit\'e $\boldsymbol{\theta}_{x,\omega}=\theta$.

\bigskip

\section{Fonctions tr\`es cuspidales}

\bigskip

\subsection{D\'efinition}

Soit $f\in C_{c}^{\infty}(G(F))$. On dit que $f$ est tr\`es cuspidale si et seulement si, pour tout sous-groupe parabolique propre $P=MU$ de $G$ et pour tout $x\in M(F)$, on a l'\'egalit\'e
$$\int_{U(F)}f(xu)du=0.$$
Remarquons que l'int\'egrale ci-dessus est localement constante en $x$. Sa nullit\'e sur $M(F)$ \'equivaut \`a sa nullit\'e sur $M(F)\cap G_{reg}(F)$. Mais, pour $x\in M(F)\cap G_{reg}(F)$, on a l'\'egalit\'e
$$\delta_{P}(x)^{1/2}\int_{U(F)}f(xu) du=D^G(x)^{1/2}D^M(x)^{-1/2}\int_{U(F)}f(u^{-1}xu)du.$$
Alors $f$ est tr\`es cuspidale si et seulement si, pour tout sous-groupe parabolique propre $P=MU$ de $G$ et pour tout $x\in M(F)\cap G_{reg}(F)$, on a l'\'egalit\'e
$$\int_{U(F)}f(u^{-1}xu)du=0.$$

 Disons qu'un \'el\'ement $x\in G_{reg}(F)$ est elliptique si $A_{G_{x}}=A_{G}$. Si le support de $f$ est contenu dans l'ensemble des \'el\'ements r\'eguliers elliptiques de $G(F)$, alors $f$ est tr\`es cuspidale.
Si $\Omega$ est un $G$-domaine dans $G(F)$ et si $f$ est tr\`es cuspidale, alors $f{\bf 1}_{\Omega}$ est tr\`es cuspidale. Si $f$ est tr\`es cuspidale, $^gf$ l'est pour tout $g\in G(F)$.

\bigskip

\subsection{Int\'egrales orbitales pond\'er\'ees de fonctions tr\`es cuspidales}

Soient $M$ un L\'evi de $G$ et $K$ un sous-groupe compact sp\'ecial de $G(F)$ en bonne position relativement \`a $M$.

\ass{Lemme}{Soient $f\in C_{c}^{\infty}(G(F))$ et $x\in M(F)\cap G_{reg}(F)$. On suppose $f$ tr\`es cuspidale.

(i) L'int\'egrale orbitale pond\'er\'ee $J_{M}(x,f)$ ne d\'epend pas de $K$.

(ii) Pour tout $y\in G(F)$, on a l'\'egalit\'e $J_{M}(x,{^yf})=J_{M}(x,f)$.

(iii) Si $A_{G_{x}}\neq A_{M}$, alors $J_{M}(x,f)=0$.}

Preuve. Soit $\tilde{K}$ un autre sous-groupe compact sp\'ecial en bonne position relativement \`a $M$. Soit $g\in G(F)$. En utilisant $K$, on a d\'efini la famille de points $(H_{P}(g))_{P\in {\cal P}(M)}$ et la $(G,M)$-famille $(v_{P}(g))_{P\in {\cal P}(M)}$. En utilisant $\tilde{K}$, on d\'efinit de m\^eme une famille de points $(\tilde{H}_{P}(g))_{P\in {\cal P}(M)}$ et une $(G,M)$-famille $(\tilde{v}_{P}(g))_{P\in {\cal P}(M)}$. On d\'efinit la $(G,M)$-famille $(d_{P})_{P\in {\cal P}(M)}$ par $d_{P}(g)=v_{P}(g)\tilde{v}_{P}(g)^{-1}$. Elle est associ\'ee \`a la famille de points $(H_{P}(g)-\tilde{H}_{P}(g))_{P\in {\cal P}(M)}$. On a $v_{P}(g)=\tilde{v}_{P}(g)d_{P}(g)$, d'o\`u, d'apr\`es 2.2(2)
$$v_{M}(g)=\sum_{Q\in {\cal F}(M)}\tilde{v}_{M}^Q(g)u_{Q}(H_{Q}(g)-\tilde{H}_{Q}(g)).$$
Alors
$$J_{M}(x,f)= D^G(x)^{1/2}\int_{G_{x}(F)\backslash G(F)}f(g^{-1}xg)v_{M}(g) dg$$
$$\qquad =D^G(x)^{1/2}\sum_{Q\in {\cal F}(M)}\int_{G_{x}(F)\backslash G(F)}f(g^{-1}xg)\tilde{v}_{M}^Q(g)u_{Q}(H_{Q}(g)-\tilde{H}_{Q}(g))dg$$
$$\qquad =D^G(x)^{1/2}\sum_{Q=LU\in {\cal F}(M)}\int_{L_{x}(F)\backslash L(F)}\int_{\tilde{K}}\int_{U(F)}f(k^{-1}u^{-1}l^{-1}xluk)$$
$$\qquad \qquad \tilde{v}_{M}^L(l)u_{Q}(H_{Q}(k))du\,dk\,dl.$$
Si $Q\neq G$, l'int\'egrale int\'erieure sur $U(F)$ est nulle puisque $f$ est tr\`es cuspidale. Le terme pour $Q=G$ est l'int\'egrale orbitale pond\'er\'ee calcul\'ee \`a l'aide de $\tilde{K}$. Cela prouve (i).

Par changement de variable
$$J_{M}(x,{^yf})=D^G(x)^{1/2}\int_{G_{x}(F)\backslash G(F)}f(g^{-1}xg)v_{M}(gy^{-1})dg.$$
La relation 2.2(2) permet d'exprimer $v_{M}(gy^{-1})$ \`a l'aide des $v_{M}^Q(g)$ pour $Q\in {\cal F}(M)$. Comme ci-dessus, les termes index\'es par $Q\neq G$ ont une contribution nulle. Le terme pour $Q=G$ donne l'int\'egrale $J_{M}(x,f)$. Cela prouve (ii).

Soit $M(x)$ le commutant de $A_{G_{x}}$ dans $G$. C'est un L\'evi de $M$. Quitte \`a changer de groupe $K$, ce qui est loisible d'apr\`es (i), on peut supposer que $K$ est en bonne position relativement \`a $M(x)$. La formule 2.2(3) appliqu\'ee aux poids conduit ais\'ement \`a l'\'egalit\'e
$$J_{M}(x,f)=\sum_{L\in {\cal L}(M(x))}d_{M(x)}^G(M,L)J_{M(x)}^L(x,f_{Q}).$$
Supposons $A_{G_{x}}\neq A_{M}$. Alors $M(x)\neq M$ et tous les L\'evi $L$ intervenant dans cette somme sont diff\'erents de $G$. Les fonctions correspondantes $f_{Q}$ sont nulles et on obtient l'assertion (iii). $\square$

\bigskip

\subsection{La distribution invariante associ\'ee \`a une fonction tr\`es cuspidale}

Soit $f\in C_{c}^{\infty}(G(F))$. On suppose $f$ tr\`es cuspidale. Soit $x\in G_{reg}(F)$. Notons $M(x)$ le commutant de $A_{G_{x}}$ dans $G$. C'est un L\'evi de $G$. On pose
$$\theta_{f}(x)=(-1)^{a_{M(x)}-a_{G}}\nu(G_{x})^{-1}D^G(x)^{-1/2}J_{M(x)}(x,f),$$
l'int\'egrale orbitale pond\'er\'ee \'etant calcul\'ee \`a l'aide d'un sous-groupe compact sp\'ecial de $G(F)$ en bonne position relativement \`a $M(x)$. Cette d\'efinition est loisible d'apr\`es le (i) du lemme pr\'ec\'edent.

\ass{Lemme}{La fonction $\theta_{f}$ est invariante par conjugaison, \`a support compact modulo conjugaison, localement int\'egrable sur $G(F)$ et localement constante sur $G_{reg}(F)$.}

Preuve. Soient $y\in G(F)$ et $x\in G_{reg}(F)$. On a $M(yxy^{-1})=yM(x)y^{-1}$. Par transport de structure,
$$J_{M(x)}(x,f)=J_{yM(x)y^{-1}}(yxy^{-1},{^yf}).$$
Le second terme est \'egal \`a $J_{yM(x)y^{-1}}(yxy^{-1},f)$ d'apr\`es le (ii) du lemme pr\'ec\'edent. Alors $\theta_{f}(x)=\theta_{f}(yxy^{-1})$, d'o\`u la premi\`ere assertion. Le support de $\theta_{f}$ est contenu dans $(Supp(f))^G$, d'o\`u la deuxi\`eme assertion. Puisque $\theta_{f}$ est invariante par conjugaison, la locale int\'egrabilit\'e et la locale constance se testent sur les tores maximaux de $G$. Soit $T$ un tel tore, notons $M(T)$ le commutant de $A_{T}$ dans $G$. Pour $x\in T(F)\cap G_{reg}(F)$, on a $M(x)=M(T)$. Les assertions \`a prouver r\'esultent des assertions similaires pour la fonction $x\mapsto J_{M(T)}(x,f)$ sur $T(F)\cap G_{reg}(F)$. Celles-ci r\'esultent du lemme 2.3. $\square$

\bigskip

\subsection{Localisation: premi\`eres propri\'et\'es}

Soit $x\in G_{ss}(F)$. On suppose que $G_{x}$ est le produit de deux groupes r\'eductifs connexes $G_{x}=G'\times G''$, conserv\'es chacun par $Z_{G}(x)(F)$. Tout \'el\'ement $X\in\mathfrak{g}_{x}(F)$ se d\'ecompose en somme d'un \'el\'ement de $\mathfrak{g}'(F)$ et d'un \'el\'ement de $\mathfrak{g}''(F)$. On notera sans plus de commentaire $X=X'+X''$ cette d\'ecomposition. De m\^eme, tout L\'evi $R$ de $G_{x}$ se d\'ecompose en $R=R'\times R''$, o\`u $R'$ est un L\'evi de $G'$ et $R''$ un L\'evi de $G''$. On note $f\mapsto f^{\sharp}$ la transformation de Fourier partielle dans $C_{c}^{\infty}(\mathfrak{g}_{x}(F))$ relative \`a la deuxi\`eme variable. C'est-\`a-dire que, pour $X=X'+X''\in \mathfrak{g}_{x}(F)$,
$$f^{\sharp}(X)=\int_{\mathfrak{g}''(F)}f(X'+Y'')\psi(<Y'',X''>)dY''.$$
Si $R$ est un L\'evi de $G_{x}$, on d\'efinit de m\^eme une transformation de Fourier partielle $f\mapsto f^{\sharp}$ dans $C_{c}^{\infty}(\mathfrak{r}(F))$. Soit $\omega$ un bon voisinage de $0$ dans $\mathfrak{g}_{x}(F)$, auquel on impose la condition (8) de 3.1. Cette situation sera conserv\'ee jusqu'en 5.8 inclus.

Soit $f\in C_{c}^{\infty}(G(F))$. Pour $g\in G(F)$, on d\'efinit $^gf_{x,\omega}\in C_{c}^{\infty}(\mathfrak{g}_{x}(F))$ par
$$^gf_{x,\omega}(X)=\left\lbrace\begin{array}{cc}0,&\,\,{\rm si}\,\,X\not\in \omega,\\ f(g^{-1}xexp(X)g),&\,\,{\rm si}\,\,X\in \omega.\\ \end{array}\right.$$
On pose $^gf_{x,\omega}^{\sharp}=(^gf_{x,\omega})^{\sharp}$. Pour $y\in Z_{G}(x)(F)$ et $X\in \mathfrak{g}_{x}(F)$, on a les \'egalit\'es
$$(1) \qquad {^{yg}f}_{x,\omega}(X)={^gf}_{x,\omega}(y^{-1}Xy),\,\,^{yg}f^{\sharp}_{x,\omega}(X)={^gf}_{x,\omega}^{\sharp}(y^{-1}Xy).$$

Soient $M$ un L\'evi de $G$ tel que $x\in M(F)$. On fixe un sous-groupe compact sp\'ecial $K$ de $G(F)$ en bonne position relativement \`a $M$. Soit $P=MU\in {\cal P}(M)$. Pour $f\in C_{c}^{\infty}(G(F))$, d\'efinissons des fonctions $\varphi[P,f]$, $\varphi^{\sharp}[P,f]$ et $J^{\sharp}_{M,x,\omega}(.,f)$ sur $\mathfrak{m}_{x}(F)\cap \mathfrak{g}_{x,reg}(F)$ par
$$\varphi[P,f](X)=D^{G_{x}}(X)^{1/2}D^{M_{x}}(X)^{-1/2}\int_{U(F)}{^{u}f}_{x,\omega}(X)du,$$
$$\varphi^{\sharp}[P,f](X)=D^{G_{x}}(X)^{1/2}D^{M_{x}}(X)^{-1/2}\int_{U(F)}{^{u}f}^{\sharp}_{x,\omega}(X)du,$$
 $$J^{\sharp}_{M,x,\omega}(X,f)=D^{G_{x}}(X)^{1/2}\int_{G_{x,X}(F)\backslash G(F)}{^gf}_{x,\omega}^{\sharp}(X)v_{M}(g) dg.$$
 
 \ass{Lemme}{(i) Ces trois int\'egrales sont absolument convergentes.
 
 (ii) Les fonctions $\varphi[P,f]$ et $\varphi^{\sharp}[P,f]$ se prolongent en des \'el\'ements de $C_{c}^{\infty}(\mathfrak{m}_{x}(F))$ et on a l'\'egalit\'e $\varphi^{\sharp}[P,f]=(\varphi[P,f])^{\sharp}$.
 
 (iii) La fonction $X\mapsto J^{\sharp}_{M,x,\omega}(X,f)$ est invariante par conjugaison par $M_{x}(F)$. Son support est compact modulo conjugaison. Elle est localement constante sur $\mathfrak{m}_{x}(F)\cap \mathfrak{g}_{x,reg}(F)$. Il existe $c>0$ et un entier $k\geq 0$ (d'ailleurs ind\'ependant de $f$) tels que l'on ait l'in\'egalit\'e
 $$\vert J^{\sharp}_{M,x,\omega}(X,f)\vert  \leq c(1+\vert log(D^{G_{x}}(X))\vert )^k$$
 pour tout $X\in \mathfrak{m}_{x}(F)\cap \mathfrak{g}_{x,reg}(F)$.}
 
 Preuve. Ecrivons
 $$(2) \qquad \varphi[P,f](X)=\int_{U_{x}(F)\backslash U(F)}D^{G_{x}}(X)^{1/2}D^{M_{x}}(X)^{-1/2}\int_{U_{x}(F)}{^{uv}f}_{x,\omega}(X)du\,dv,$$
  $$(3) \qquad \varphi^{\sharp}[P,f](X)=\int_{U_{x}(F)\backslash U(F)}D^{G_{x}}(X)^{1/2}D^{M_{x}}(X)^{-1/2}\int_{U_{x}(F)}{^{uv}f}^{\sharp}_{x,\omega}(X)du\,dv.$$
  D'apr\`es 2.1(5), on peut fixer un sous-ensemble compact $\Gamma$ de $G(F)$ tel que ${^gf}_{x,\omega}=0$ pour $g\in G(F)$, $g\not\in G_{x}(F)\Gamma$. On a donc aussi ${^gf}^{\sharp}_{x,\omega}=0$ pour un tel $g$. On v\'erifie que l'application 
 $$U_{x}(F)\backslash U(F) \to G_{x}(F)\backslash G(F)$$
 est d'image ferm\'ee et est un hom\'eomorphisme de sa source sur son image.  Il en r\'esulte que les int\'egrales en $v\in U_{x}(F)\backslash U(F)$ dans les \'egalit\'es (2) et (3) sont \`a support compact. Les expressions que l'on int\`egre \'etant localement constantes, il suffit de fixer $v\in U(F)$ et de prouver les assertions pour les expressions en question. Gr\^ace \`a (1), celles-ci s'\'ecrivent
 $$D^{G_{x}}(X)^{1/2}D^{M_{x}}(X)^{-1/2}\int_{U_{x}(F)}{^vf}_{x,\omega}(u^{-1}Xu)du,$$
 $$D^{G_{x}}(X)^{1/2}D^{M_{x}}(X)^{-1/2}\int_{U_{x}(F)}{^vf}^{\sharp}_{x,\omega}(u^{-1}Xu)du.$$
 Par un calcul familier, elles sont \'egales \`a
 $$\int_{\mathfrak{u}_{x}(F)}{^vf}_{x,\omega}(X+N)dN\,\,{\rm et}\,\,\int_{\mathfrak{u}_{x}(F)}{^vf}^{\sharp}_{x,\omega}(X+N)dN.$$
 Les assertions sont maintenant faciles \`a prouver.
 
 De m\^eme, pour d\'emontrer les assertions relatives \`a la fonction $J^{\sharp}_{M,x,\omega}(.,f)$, on peut fixer $\gamma\in \Gamma$ et prouver les m\^emes assertions pour la fonction
 $$X\mapsto D^{G_{x}}(X)^{1/2}\int_{G_{x,X}(F)\backslash G_{x}(F)}{^{g\gamma}f}^{\sharp}_{x,\omega}(X)v_{M}(g\gamma) dg,$$
 ou encore
 $$X\mapsto D^{G_{x}}(X)^{1/2}\int_{G_{x,X}(F)\backslash G_{x}(F)}{^{\gamma}f}^{\sharp}_{x,\omega}(g^{-1}Xg)v_{M}(g\gamma)dg.$$
 Il suffit alors de reprendre la preuve du lemme 2.3. $\square$
 
 \bigskip
 
 \subsection{Localisation pour une fonction tr\`es cuspidale}
 
 Les groupes $M$, $P$ et $K$ sont comme dans le paragraphe pr\'ec\'edent.
 
 \ass{Lemme}{Soit $f\in C_{c}^{\infty}(G(F))$ une fonction tr\`es cuspidale.
 
 (i) Si $P\not=G$, les fonctions $\varphi[P,f]$ et $\varphi^{\sharp}[P,f]$ sont nulles.
 
 (ii) La fonction $J^{\sharp}_{M,x,\omega}(.,f)$ ne d\'epend pas du choix de $K$. Elle s'annule aux points $X\in \mathfrak{m}_{x}(F)\cap \mathfrak{g}_{x,reg}(F)$ tels que $A_{G_{x,X}}\not=A_{M}$. En particulier, elle est partout nulle si $A_{M_{x}}\neq A_{M}$.
 Pour tous $y\in G(F)$ et $X\in \mathfrak{m}_{x}(F)\cap \mathfrak{g}_{x,reg}(F)$, on a l'\'egalit\'e 
 $$J^{\sharp}_{M,x,\omega}(X,f)=J^{\sharp}_{M,x,\omega}(X,{^yf}).$$}
 
 Preuve. Soit $X\in \mathfrak{m}_{x}(F)\cap \mathfrak{g}_{x,reg}(F)$. Par d\'efinition, $\varphi[P,f](X)=0$ si $X\not\in \omega$. Supposons $X\in \omega$. Alors
 $$\varphi[P,f](X)=D^{G_{x}}(X)^{1/2}D^{M_{x}}(X)^{-1/2}\int_{U(F)}f(u^{-1}xexp(X)u)du.$$
 On a $xexp(X)\in M(F)\cap G_{reg}(F)$. Si $P\neq G$, cette int\'egrale est nulle puisque $f$ est tr\`es cuspidale. Donc $\varphi[P,f]=0$. Puisque $\varphi^{\sharp}[P,f]=(\varphi[P,f])^{\sharp}$, on a aussi $\varphi^{\sharp}[P,f]=0$.
 
 On d\'emontre (ii) en reprenant la preuve du lemme 5.2. On y avait utilis\'e l'hypoth\`ese de forte cuspidalit\'e de $f$ pour annuler certaines int\'egrales. Maintenant, les int\'egrales similaires s'annulent d'apr\`es la nullit\'e des fonctions $\varphi^{\sharp}[Q,f]$ pour tout $Q\in {\cal F}(M)$, $Q\neq G$. Cela conduit aux m\^emes r\'esultats. $\square$
 
 \bigskip
 
 \subsection{Les distributions locales associ\'ees \`a une fonction tr\`es cuspidale}
   
   Soit $f\in C_{c}^{\infty}(G(F))$ une fonction tr\`es cuspidale. On d\'efinit une fonction $\theta_{f,x,\omega}$ sur $\mathfrak{g}_{x,reg}(F)$ par
   $$\theta_{f,x,\omega}(X)=\left\lbrace\begin{array}{cc}0,&\,\,{\rm si}\,\,X\not\in \omega,\\ \theta_{f}(xexp(X)),&\,\,{\rm si}\,\,X\in \omega.\\ \end{array}\right.$$
   
   Soit $X\in \mathfrak{g}_{x,reg}(F)$. Notons    ${\bf M}(X)$ le commutant de  $A_{G_{x,X}}$  dans $G$. On pose
   $$\theta^{\sharp}_{f,x,\omega}(X)=(-1)^{a_{{\bf M}(X)}-a_{G}}\nu(G_{x,X})^{-1}D^{G_{x}}(X)^{-1/2}J^{\sharp}_{{\bf M}(X),x,\omega}(X,f),$$
   le dernier terme \'etant calcul\'e \`a l'aide d'un sous-groupe compact sp\'ecial de $G(F)$ en bonne position relativement \`a ${\bf M}(X)$. Cette d\'efinition est loisible d'apr\`es le (ii) du lemme pr\'ec\'edent.
   
   \ass{Lemme}{Les fonctions $\theta_{f,x,\omega}$ et $\theta^{\sharp}_{f,x,\omega}$ sont invariantes par conjugaison par $G_{x}(F)$, \`a support compact modulo conjugaison, localement int\'egrables sur $\mathfrak{g}_{x}(F)$ et localement constantes sur $\mathfrak{g}_{x,reg}(F)$.}
   
   Preuve. Pour la fonction $\theta_{f,x,\omega}$, les assertions r\'esultent du lemme 5.3. Pour la fonction $\theta^{\sharp}_{f,x,\omega}$, elles se prouvent comme dans ce lemme, en utilisant les lemmes 5.4(iii) et 5.5(ii). $\square$
   
   \bigskip
   
   \subsection{Descente des int\'egrales orbitales pond\'er\'ees}

On fixe un L\'evi minimal $R_{min}$ de $G_{x}$ et un sous-groupe compact sp\'ecial $K_{x}$ de $G_{x}(F)$ en bonne position relativement \`a $R_{min}$. Rappelons l'application $R\mapsto {\bf R}$ de 3.2. On fixe un sous-groupe compact sp\'ecial $K$ de $G(F)$ en bonne position relativement \`a ${\bf R}_{min}$.

\ass{Lemme}{Soit $f\in C_{c}^{\infty}(G(F))$ une fonction tr\`es cuspidale. Pour tout $S\in {\cal L}(R_{min})$, il existe une fonction $f^S\in C_{c}^{\infty}(\mathfrak{s}(F))$, \`a support dans $\omega\cap \mathfrak{s}(F)$, telle que, pour tout $R\in {\cal L}(R_{min})$, on ait les \'egalit\'es

(i) $J_{{\bf R}}(xexp(X),f)=D^G(x)^{1/2}\sum_{S\in {\cal L}(R)}J_{R}^S(X,f^S)$ pour tout $X\in \mathfrak{r}(F)\cap \mathfrak{g}_{x,reg}(F)\cap \omega$;

(ii) $J^{\sharp}_{{\bf R},x,\omega}(X,f)=\sum_{S\in{\cal L}(R)}J_{R}^S(X,f^{S,\sharp})$ pour tout $X\in \mathfrak{r}(F)\cap \mathfrak{g}_{x,reg}(F)$, o\`u $f^{S,\sharp}=(f^S)^{\sharp}$.}

Preuve. Soient $R\in {\cal L}(R_{min})$ et $X\in \mathfrak{r}(F)\cap \mathfrak{g}_{x,reg}(F)\cap \omega$. On a les \'egalit\'es
$$J_{{\bf R}}(xexp(X),f)=D^G(xexp(X))^{1/2}\int_{G_{x,X}(F)\backslash G(F)}{^gf}_{x,\omega}(X)v_{{\bf R}}(g)dg,$$
$$(1) \qquad J_{{\bf R}}(xexp(X))=D^G(x)^{1/2}D^{G_{x}}(X)^{1/2}\int_{G_{x,X}(F)\backslash G(F)}{^gf}_{x,\omega}(X)v_{{\bf R}}(g)dg.$$
Fixons un sous-groupe ouvert compact $K'$ de $K$ tel que $f$ soit invariante par conjugaison par $K'$. Soit $\Delta$ un ensemble de repr\'esentants de $G_{x}(F)\backslash G(F)/K'$. On a l'\'egalit\'e
$$\int_{G_{x,X}(F)\backslash G(F)}{^gf}_{x,\omega}(X)v_{{\bf R}}(g) dg=\sum_{\delta\in \Delta}m(\delta)\int_{G_{x,X}(F)\backslash G_{x}(F)}{^{g\delta}f}_{x,\omega}(X)v_{{\bf R}}(g\delta)dg,$$
o\`u $m(\delta)=mes(K')mes(G_{x}(F)\cap \delta K'\delta^{-1})^{-1}$. D'apr\`es 3.1(5), on peut fixer un sous-ensemble fini $\Delta_{0}\subset \Delta$ tel que ${^{g\delta}f}_{x,\omega}=0$ pour tout $g\in G_{x}(F)$ et tout $\delta\in \Delta$ tel que $\delta\not\in \Delta_{0}$. Le lemme 3.3 conduit \`a l'\'egalit\'e suivante
$$(2) \qquad J_{{\bf R}}(xexp(X),f)=D^G(x)^{1/2}\sum_{S\in{\cal L}(R)}\sum_{\delta\in \Delta_{0}}\sum_{Q\in {\cal P}({\bf S})}m(\delta)c^{Q,\delta}(X),$$
o\`u
$$c^{Q,\delta}(X)=D^{G_{x}}(X)^{1/2}\int_{G_{x,X}(F)\backslash G_{x}(F)}{^{\delta}f}_{x,\omega}(g^{-1}Xg)v_{R}^{Q_{x}}(g)u_{Q}(H_{Q}(g\delta)-H_{Q_{x}}(g))dg.$$
Un calcul familier remplace cette expression par
$$c^{Q,\delta}(X)=D^S(X)^{1/2}\int_{G(x,X)(F)\backslash S(F)}\int_{K_{x}}\int_{\mathfrak{u}_{x}(F)}{^{k\delta}f}_{x,\omega}(l^{-1}Xl+N)v_{R}^S(l)u_{Q}(H_{Q}(k\delta))dN\,dk\,dl,$$
o\`u $U_{x}$ est le radical unipotent de $Q_{x}$. D\'efinissons une fonction $f^{Q,\delta}$ sur $\mathfrak{s}(F)$ par
$$f^{Q,\delta}(Y)=\int_{K_{x}}\int_{\mathfrak{u}_{x}(F)}{^{k\delta}f}_{x,\omega}(Y+N)u_{Q}(H_{Q}(k\delta))dN\,dk.$$
Cette fonction appartient \`a $C_{c}^{\infty}(\mathfrak{s}(F))$ et on a l'\'egalit\'e
$$c^{Q,\delta}(X)=J_{R}^S(X,f^{Q,\delta}).$$
En posant 
$$f^S=\sum_{\delta\in \Delta_{0}}\sum_{Q\in {\cal P}({\bf S})}m(\delta)f^{Q,\delta},$$
l'\'egalit\'e (2) devient celle du (i) de l'\'enonc\'e.

Pour tout $X\in \mathfrak{r}(F)\cap \mathfrak{g}_{x,reg}(F)$, le terme $J^{\sharp}_{{\bf R},x,\omega}(X,f)$ s'exprime par une formule analogue \`a (1): il suffit de supprimer le facteur $D^G(x)^{1/2}$ et de remplacer les fonctions $^gf_{x,\omega}$ par $^gf^{\sharp}_{x,\omega}$. On peut refaire le calcul ci-dessus. Les fonctions $f^{Q,\delta}$ sont remplac\'ees par les fonctions
$$Y\mapsto \int_{K_{x}}\int_{\mathfrak{u}_{x}(F)}{^{k\delta}f}_{x,\omega}^{\sharp}(Y+N)u_{Q}(H_{Q}(k\delta))dn\,dk.$$
On v\'erifie ais\'ement que ce sont les images $f^{Q,\delta,\sharp}$ de $f^{Q,\delta}$ par transformation de Fourier partielle. Le (ii) de l'\'enonc\'e s'ensuit. $\square$
 
 \bigskip

 \subsection{Transform\'ees de Fourier des distributions locales}
 
 \ass{Proposition}{Soit $f\in C_{c}^{\infty}(G(F))$ une fonction tr\`es cuspidale. Alors la fonction $\theta^{\sharp}_{f,x,\omega}$ est la transform\'ee de Fourier partielle de $\theta_{f,x,\omega}$, c'est-\`a-dire que, pour toute $\varphi\in C_{c}^{\infty}(\mathfrak{g}_{x}(F))$, on a l'\'egalit\'e
 $$\int_{\mathfrak{g}_{x}(F)}\theta^{\sharp}_{f,x,\omega}(X)\varphi(X)dX=\int_{\mathfrak{g}_{x}(F)}\theta_{f,x,\omega}(X)\varphi^{\sharp}(X)dX.$$}
 
 Preuve. Soit $\varphi\in C_{c}^{\infty}(\mathfrak{g}_{x}(F))$. Notons $\theta^{\sharp}(\varphi)$ le membre de gauche de l'\'egalit\'e de l'\'enonc\'e et $\theta(\varphi^{\sharp})$ celui de droite. D'apr\`es la formule de Weyl, on a
 $$\theta(\varphi^{\sharp})=\sum_{R\in {\cal L}(R_{min})}\vert W^R\vert \vert W^{G_{x}}\vert ^{-1}\sum_{T\in {\cal T}_{ell}(R)}\vert W(R,T)\vert ^{-1}\int_{\mathfrak{t}(F)}J_{G_{x}}(X,\varphi^{\sharp})\theta_{f,x,\omega}(X)D^{G_{x}}(X)^{1/2}dX.$$
 Soient $R\in {\cal L}(R_{min})$, $T\in {\cal T}_{ell}(R)$ et $X\in \mathfrak{t}(F)\cap \mathfrak{g}_{x,reg}(F)$.  Supposons d'abord $X\in \omega$. Alors
 $$\theta_{f,x,\omega}(X)=(-1)^{a_{M(xexp(X))}-a_{G}}\nu(T)^{-1}D^G(xexp(X))^{-1/2} J_{M(xexp(X))}(xexp(X),f).$$
 On a $D^G(xexp(X))=D^G(x)D^{G_{x}}(X)$. On a aussi $M(xexp(X))={\bf R}$. Enfin $J_{{\bf R}}(xexp(X),f)$ est calcul\'e par le lemme 5.7 et on obtient:
  $$\theta_{f,x,\omega}(X)=(-1)^{a_{R}-a_{G}}\nu(T)^{-1}D^{G_{x}}(X)^{-1/2}\sum_{S\in {\cal L}(R)}J^S_{R}(X,f^S).$$
 Cette formule reste vraie pour $X\not\in \omega$ car ses deux membres sont nuls: par d\'efinition de $\theta_{f,x,\omega}$ pour celui de gauche; parce que les fonctions $f^S$ sont \`a support dans $\mathfrak{s}(F)\cap \omega$ pour celui de droite. D'o\`u l'\'egalit\'e
 $$\theta(\varphi^{\sharp})= \sum_{S\in {\cal L}(R)}\vert W^S\vert \vert W^{G_{x}}\vert ^{-1}(-1)^{a_{S}-a_{G}}\theta^S(\varphi^{\sharp},f^S),$$
 o\`u
 $$(1) \qquad \theta^S(\varphi^{\sharp},f^S)=\sum_{R\in {\cal L}^S(R_{min})}\vert W^R\vert \vert W^{G_{x}}\vert ^{-1}(-1)^{a_{R}-a_{S}}$$
 $$\qquad \sum_{T\in {\cal T}_{ell}(R)}\vert W(R,T)\vert ^{-1}\nu(T)^{-1}\int_{\mathfrak{t}(F)}J_{G_{x}}(X,\varphi^{\sharp})J_{R}^S(X,f^S)dX.$$
 Soient $S\in {\cal L}(R_{min})$, $\alpha$ et $\beta$ deux \'el\'ements de $C_{c}^{\infty}(\mathfrak{s}(F))$. Pour $R\in {\cal L}^S(R_{min})$, $T\in {\cal T}_{ell}(R)$ et $X\in \mathfrak{t}(F)\cap \mathfrak{g}_{x,reg}(F)$, posons
 $$(2) \qquad j_{R}^S(X,\alpha,\beta)=\sum_{S''_{1},S''_{2}\in {\cal L}^{S''}(R'')}d_{R''}^{S''}(S''_{1},S''_{2}) J_{S'\times R''}^{S'\times S''_{1}}(X,\alpha_{S'\times \bar{Q}''_{1}})J_{R}^{S'\times S''_{2}}(X,\beta_{S'\times Q''_{2}}).$$
Les sous-groupes paraboliques $Q''_{1}$ et $Q''_{2}$ sont d\'etermin\'es par un param\`etre auxiliaire $\xi''\in {\cal A}_{R''}^{S''}$. On pose
$$j^S(\alpha,\beta)=\sum_{R\in {\cal L}^{S}(R_{min})}\vert W^R\vert \vert W^S\vert ^{-1}(-1)^{a_{R}-a_{S}}\sum_{T\in {\cal T}_{ell}(R)}\vert W(R,T)\vert ^{-1}\nu(T)^{-1}\int_{\mathfrak{t}(F)}j_{R}^S(X,\alpha,\beta)dX.$$
Revenons \`a la formule (2) et supposons que $\alpha$ v\'erifie l'hypoth\`ese (H) de 2.5. On peut alors remplacer $\alpha_{S'\times \bar{Q}''_{1}}$ par $\alpha_{S'\times S''_{1}}$. Soit $S''_{1}\in {\cal L}^{S''}(M'')$. Pour tout $\gamma\in C_{c}^{\infty}(\mathfrak{s}'(F)\times \mathfrak{s}''_{1}(F))$, on a la formule de descente
$$J_{S'\times R''}^{S'\times S''_{1}}(X, \gamma)=J_{R}^{R'\times S''_{1}}(X,\gamma_{Q'\times S''_{1}})$$
o\`u $Q'$ est un \'el\'ement quelconque de ${\cal P}^{S'}(R')$. Appliqu\'ee \`a $\gamma=\alpha_{S'\times S''_{1}}$, cette formule devient
$$J_{S'\times R''}^{S'\times S''_{1}}(X,\alpha_{S'\times S''_{1}})=J_{R}^{R'\times S''_{1}}(X,\alpha_{R'\times S''_{1}}).$$
Alors
$$j_{R}^S(X,\alpha,\beta)=\sum_{S''_{1}\in {\cal L}^{S''}(R'')}J_{R}^{R'\times S''_{1}}(X,\alpha_{R'\times S''_{1}})\sum_{S''_{2}\in {\cal L}^{S''}(R'')}d_{R''}^{S''}(S''_{1},S''_{2})J_{R}^{S'\times S''_{2}}(X,\beta_{S'\times Q''_{2}}).$$
Fixons $S''_{1}\in {\cal L}^{S''}(R'')$. L'application $S_{2}\mapsto S''_{2}$ est une bijection de l'ensemble des $S_{2}\in {\cal L}^S(R)$ tels que $d_{R}^S(R'\times S''_{1},S_{2})\neq 0$ sur l'ensemble des $S''_{2}\in {\cal L}^{S''}(R'')$ tels que $d_{R''}^{S''}(S''_{1},S''_{2})\neq 0$. La bijection r\'eciproque est $S''_{2}\mapsto S'\times S''_{2}$. Fixons un param\`etre auxiliaire $\xi\in {\cal A}_{R}^S$ dont la seconde composante soit $\xi''$. Pour $S_{2}$ et $S''_{2}$ se correspondant par la bijection ci-dessus, on a $d_{R}^S(R'\times S''_{1},S_{2})=d_{R''}^{S''}(S''_{1},S''_{2})$ et le parabolique $Q_{2}$ associ\'e \`a $R'\times S''_{1}$, $S_{2}$ et $\xi$ n'est autre que $S'\times Q''_{2}$. Cela conduit \`a l'\'egalit\'e
$$j_{R}^S(X,\alpha,\beta)=\sum_{S''_{1}\in {\cal L}^{S''}(R'')}J_{R}^{R'\times S''_{1}}(X,\alpha_{R'\times S''_{1}})\sum_{S_{2}\in {\cal L}^S(R)}d_{R}^S(R'\times S''_{1},S_{2})J_{R}^{S_{2}}(X,\beta_{Q_{2}}),$$
ou encore, d'apr\`es 2.2(3),
$$(3) \qquad j_{R}^S(X,\alpha,\beta)=\sum_{S''_{1}\in {\cal L}^{S''}(R'')}J_{R}^{R'\times S''_{1}}(X,\alpha_{R'\times S''_{1}})J_{R'\times S''_{1}}^S(X,\beta).$$
Supposons que la fonction $\varphi$ v\'erifie l'hypoth\`ese (H). Une g\'en\'eralisation imm\'ediate de ce que l'on a dit en 2.5 montre que $\varphi^{\sharp}$ v\'erifie aussi cette hypoth\`ese et que $(\varphi_{S})^{\sharp}=(\varphi^{\sharp})_{S}$. On peut noter sans ambig\"uit\'e $\varphi^{\sharp}_{S}$ cette fonction. Elle v\'erifie aussi (H) et $j_{R}^S(X,\varphi^{\sharp}_{S},f^S)$ se calcule par la formule (3). Remarquons que le terme de cette formule index\'e par $S''_{1}=R''$ est $J_{R}^R(X,\varphi^{\sharp}_{R})J_{R}^S(X,f^S)$, \'egal \`a $J_{G_{x}}(X,\varphi^{\sharp})J_{R}^S(X,f^S)$, qui est pr\'ecis\'ement le terme intervenant dans (1). On en d\'eduit
$$j^S(\varphi^{\sharp}_{S},f^S)-\theta^S(\varphi^{\sharp},f^S)=\sum_{R\in {\cal L}^S(R_{min})}\vert W^R\vert \vert W^S\vert ^{-1}(-1)^{a_{R}-a_{S}}\sum_{T\in {\cal T}_{ell}(R)}\vert W(R,T)\vert ^{-1}$$
$$\qquad \nu(T)^{-1}\int_{\mathfrak{t}(F)}\sum_{S''_{1}\in {\cal L}^{S''}(R''),S''_{1}\not=R''}J_{R}^{R'\times S''_{1}}(X,\varphi^{\sharp}_{R'\times S''_{1}})J_{R'\times S''_{1}}^S(X,f^S)dX.$$
Posons
$$j(\varphi^{\sharp})= \sum_{S\in {\cal L}(R_{min})}\vert W^S\vert \vert W^{G_{x}}\vert ^{-1}(-1)^{a_{S}-a_{G}}j^S(\varphi^{\sharp}_{S},f^S).$$
Alors
$$j(\varphi^{\sharp})-\theta(\varphi^{\sharp})=\sum_{S\in {\cal L}(R_{min})}\vert W^S\vert \vert W^{G_{x}}\vert ^{-1}(-1)^{a_{S}-a_{G}}(j^S(\varphi^{\sharp}_S,f^S)-\theta^S(\varphi^{\sharp},f^S)),$$
$$j(\varphi^{\sharp})-\theta(\varphi^{\sharp})=\sum_{R\in {\cal L}(R_{min})}\vert W^R\vert \vert W^{G_{x}}\vert ^{-1}(-1)^{a_{R}-a_{G}}\sum_{T\in {\cal T}_{ell}(R)}\vert W(R,T)\vert ^{-1}\nu(T)^{-1}$$
$$\qquad \int_{\mathfrak{t}(F)}\sum_{S''_{1}\in {\cal L}^{G''}(R''), S''_{1}\neq R''}J_{R}^{R'\times S''_{1}}(X,\varphi^{\sharp}_{R'\times S''_{1}})\sum_{S\in {\cal L}^{G_{x}}(R'\times S''_{1})}J_{R'\times S''_{1}}^S(X,f^S)dX.$$
D'apr\`es le lemme 5.7, la derni\`ere somme dans l'expression ci-dessus est
$$D^G(x)^{-1/2}J_{{\bf S}_{1}}(xexp(X),f),$$
o\`u on a pos\'e $S_{1}=R'\times S''_{1}$. Or $xexp(X)\in {\bf R}\subsetneq {\bf S}_{1}$ puisque $R''\subsetneq S''_{1}$. D'apr\`es le lemme 5.2(iii), l'expression ci-dessus est nulle. D'o\`u l'\'egalit\'e
$$\theta(\varphi^{\sharp})=j(\varphi^{\sharp}).$$
D\'efinissons $j^{\sharp}(\varphi)$ en rempla\c{c}ant $j^S(\varphi^{\sharp}_{S},f^S)$ par $j^S(\varphi_{S},f^{S,\sharp})$ dans la formule (4). Le m\^eme calcul conduit \`a l'\'egalit\'e
$$\theta^{\sharp}(\varphi)=j^{\sharp}(\varphi).$$
Il suffit en effet de remplacer l'usage du lemme 5.2(iii) par celui du lemme 5.5(ii).

L'\'egalit\'e $\theta(\varphi^{\sharp})=\theta^{\sharp}(\varphi)$ que l'on veut prouver \'equivaut donc \`a $j(\varphi^{\sharp})=j^{\sharp}(\varphi)$. Il suffit de fixer $S\in {\cal L}(R_{min})$ et de prouver l'\'egalit\'e $j^S(\varphi^{\sharp}_{S},f^S)=j^S(\varphi_{S},f^{S,\sharp})$. On va plus g\'en\'eralement prouver l'\'egalit\'e $j^S(\alpha^{\sharp},\beta)=j^L(\alpha,\beta^{\sharp})$ pour toutes $\alpha,\beta\in C_{c}^{\infty}(\mathfrak{s}(F))$. Par lin\'earit\'e, on peut supposer $\alpha=\alpha'\otimes \alpha''$, $\beta=\beta'\otimes \beta''$, o\`u $\alpha',\beta'\in C_{c}^{\infty}(\mathfrak{s}'(F))$, $\alpha'',\beta''\in C_{c}^{\infty}(\mathfrak{s}''(F))$. Consid\'erons l'expression (2). On a
$$j_{R}^S(X,\alpha^{\sharp},\beta)=J_{S'}(X',\alpha')J_{R'}^{S'}(X',\beta')\sum_{S''_{1},S''_{2}\in {\cal L}^{S''}(R'')}d_{R''}^{S''}(S''_{1},S''_{2})J_{R''}^{S''_{1}}(X'',\hat{\alpha}''_{\bar{Q}''_{1}})J_{R''}^{S''_{2}}(X'',\beta''_{Q''_{2}}),$$
puis
$$j_{R}^S(X,\alpha^{\sharp},\beta)=J_{S'}(X',\alpha')J_{R'}^{S'}(X',\beta')J_{R''}^{S''}(X'',\hat{\alpha}'',\beta''),$$
avec la notation de 2.4. On en d\'eduit ais\'ement
$$j^S(\alpha^{\sharp},\beta)=k^{S'}(\alpha',\beta')J^{S''}(\hat{\alpha}'',\beta''),$$
o\`u
$$k^{S'}(\alpha',\beta')=\sum_{R'\in {\cal L}^{S'}(R'_{min})}\vert W^{R'}\vert \vert W^{S'}\vert ^{-1}(-1)^{a_{R'}-a_{S'}}\sum_{T'\in {\cal T}_{ell}(R')}\vert W(R',T')\vert ^{-1}$$
$$\qquad \nu(T')^{-1}\int_{\mathfrak{t}'(F)}J_{S'}(X',\alpha')J_{R'}^{S'}(X',\beta')dX'.$$
De m\^eme
$$j^S(\alpha,\beta^{\sharp})=k^{S'}(\alpha',\beta')J^{S''}(\alpha'',\hat{\beta}'').$$
On d\'eduit alors l'\'egalit\'e cherch\'ee $j^S(\alpha^{\sharp},\beta)=j^S(\alpha,\beta^{\sharp})$ du th\'eor\`eme 2.4.

Cela prouve l'\'egalit\'e de l'\'enonc\'e pour les fonctions $\varphi$ v\'erifiant l'hypoth\`ese (H). Puisque les deux distributions $\varphi\mapsto \theta^{\sharp}(\varphi)$ et $\varphi\mapsto \theta(\varphi^{\sharp})$ sont invariantes par conjugaison par $G_{x}(F)$, le lemme 2.5(ii) g\'en\'eralise leur \'egalit\'e \`a toute fonction $\varphi$ \`a support dans $\mathfrak{g}_{x,reg}(F)$. Pour obtenir l'\'egalit\'e pour toute $\varphi$, il suffit de prouver que les deux distributions sont localement int\'egrables. C'est vrai pour la premi\`ere d'apr\`es le lemme 5.6. Consid\'erons la seconde. Fixons deux $G$-domaines $\Gamma'\subset \mathfrak{g}'(F)$ et $\Gamma''\subset \mathfrak{g}''(F)$, compacts modulo conjugaison. D'apr\`es la conjecture de Howe (cf. 2.6(2) ), on peut fixer une famille finie $(X''_{i})_{i=1,...,n}$ d'\'el\'ements de $\omega''\cap \mathfrak{g}''_{reg}(F)$ et une famille finie $(\beta_{i})_{i=1,...,n}$ d'\'el\'ements de $C_{c}^{\infty}(\Gamma'')$ de sorte que, pour tout $X''\in \omega''\cap \mathfrak{g}''_{reg}(F)$ et toute $\beta\in C_{c}^{\infty}(\Gamma'')$, on ait l'\'egalit\'e
$$J_{G}(X'',\hat{\beta})=\sum_{i=1,...,n}J_{G''}(X''_{i},\hat{\beta})J_{G''}(X'',\hat{\beta}_{i}).$$
 Soient $\alpha\in C_{c}^{\infty}(\Gamma')$ et $\beta\in C_{c}^{\infty}(\Gamma'')$. Posons
$$\gamma=\alpha\otimes(\beta-\sum_{i=1,...,n}J_{G''}(X''_{i},\hat{\beta})\beta_{i}).$$
Alors $J_{G_{x}}(X,\gamma^{\sharp})=0$ pour tout $X\in \omega\cap\mathfrak{g}_{x,reg}(F)$. Puisque $\theta_{f,x,\omega}$ est \`a support dans $\omega$, il en r\'esulte que $\theta(\gamma^{\sharp})=0$. Donc
$$\theta(\alpha\otimes\hat{\beta})=\sum_{i=1,...,n}J_{G''}(X''_{i},\hat{\beta})\theta(\alpha\otimes \beta_{i}).$$
La fonction $\theta_{f,x,\omega}$ est localement int\'egrable. Il en r\'esulte que la distribution $\alpha\mapsto \theta(\alpha\otimes \beta_{i})$ l'est aussi. D'apr\`es [HCvD] th\'eor\`emes 13 et 15, la distribution $\beta\mapsto J_{G''}(X''_{i},\hat{\beta})$ est aussi localement int\'egrable. Donc la distribution $\alpha\otimes \beta\mapsto \theta(\alpha\otimes \hat{\beta})$ est int\'egrable sur $\Gamma'\times \Gamma''$. Cela ach\`eve la d\'emonstration. $\square$

\bigskip

\subsection{Le quasi-caract\`ere associ\'e \`a une fonction tr\`es cuspidale}

\ass{Corollaire}{Soit $f\in C_{c}^{\infty}(G(F))$. Supposons $f$ tr\`es cuspidale. Alors la fonction $\theta_{f}$ est un quasi-caract\`ere.}

Preuve. Soit $x\in G_{ss}(F)$. On applique la proposition pr\'ec\'edente au cas $G'=\{1\}$, $G''=G_{x}$. On obtient que la transform\'ee de Fourier de $\theta_{f,x,\omega}$ est la fonction $\theta^{\sharp}_{f,x,\omega}$, que l'on peut d'ailleurs plut\^ot noter $\hat{\theta}_{f,x,\omega}$. Le support de cette fonction est compact modulo conjugaison. D'apr\`es le th\'eor\`eme 4.2, $\theta_{f,x,\omega}$ est donc un quasi-caract\`ere sur $\mathfrak{g}_{x}(F)$. Alors $\theta_{f}$ co\"{\i}ncide au voisinage de $x$ avec un quasi-caract\`ere. Cela \'etant vrai pour tout $x\in G_{ss}(F)$, la conclusion s'ensuit. $\square$

\bigskip

\section{Fonctions tr\`es cuspidales sur les alg\`ebres de Lie} 

\bigskip

\subsection{Premi\`eres propri\'et\'es}

La d\'efinition des fonctions tr\`es cuspidales s'adapte aux fonctions sur l'alg\`ebre de Lie $\mathfrak{g}(F)$. Soit $f\in C_{c}^{\infty}(\mathfrak{g}(F))$. On dit qu'elle est tr\`es cuspidale si et seulement si, pour tout sous-groupe parabolique propre $P=MU$ de $G$ et tout $X\in \mathfrak{m}(F)$, on a l'\'egalit\'e 
$$(1) \qquad \int_{\mathfrak{u}(F)}f(X+N)dN=0.$$
Ou encore si et seulement si, pour tout sous-groupe parabolique propre $P=MU$ de $G$ et tout $X\in \mathfrak{m}(F)\cap \mathfrak{g}_{reg}(F)$, on a l'\'egalit\'e 
$$\int_{U(F)}f(u^{-1}Xu)du=0.$$
Les propri\'et\'es \'enonc\'ees en 5.1 et 5.2 restent vraies. Il y a une propri\'et\'e suppl\'ementaire: si $f$ est tr\`es cuspidale, $\hat{f}$ l'est aussi. En effet, notons $f_{U}(X)$ l'int\'egrale (1). On v\'erifie que $(\hat{f})_{U}=(f_{U})\hat{}$ et l'assertion s'ensuit.

Soit $f\in C_{c}^{\infty}(\mathfrak{g}(F))$, supposons $f$ tr\`es cuspidale. On d\'efinit une fonction $\theta_{f}$ sur $\mathfrak{g}_{reg}(F)$ par
$$\theta_{f}(X)=(-1)^{a_{M(X)}-a_{G}}\nu(G_{X})^{-1}D^G(X)^{-1/2}J_{M(X)}(X,f),$$
o\`u $M(X)$ est le commutant de $A_{G_{X}}$ dans $G$. Elle a des propri\'et\'es similaires \`a celles \'enonc\'ees au lemme 5.3.

\ass{Lemme}{Soit $f\in C_{c}^{\infty}(\mathfrak{g}(F))$ une fonction tr\`es cuspidale.

(i) La fonction $\theta_{\hat{f}}$ est la transform\'ee de Fourier  de $\theta_{f}$.

(ii) La fonction $\theta_{f}$ est un quasi-caract\`ere.}

Preuve. La d\'emonstration de la proposition 5.8 se simplifie grandement. En effet, fixons  un L\'evi minimal $M_{0}$ de $G$ et un sous-groupe compact sp\'ecial $K$ de $G(F)$ en bonne position relativement \`a $M_{0}$. Pour toute $\varphi\in C_{c}^{\infty}(\mathfrak{g}(F))$, on a simplement
$$(2)\qquad J(\hat{\varphi},f)=\int_{\mathfrak{g}(F)}\hat{\varphi}(X)\theta_{f}(X)dX.$$
L'assertion (i) s'ensuit en appliquant le th\'eor\`eme 2.4. Pour prouver (2), soient $M\in {\cal L}(M_{0})$, $T\in {\cal T}_{ell}(M)$ et $X\in \mathfrak{t}(F)\cap \mathfrak{g}_{reg}(F)$. Parce que $f_{Q}=0$ pour tout sous-groupe parabolique propre $Q$ de $G$, on a
$$J_{M}(X,\hat{\varphi},f)=J_{G}(X,\hat{\varphi})J_{M}(X,f).$$
Mais alors $J(\hat{\varphi},f)$ co\"{\i}ncide avec le membre de droite de (2) exprim\'e \`a l'aide de la formule de Weyl.

L'assertion (ii) r\'esulte de (i) et du th\'eor\`eme 4.2. $\square$

\bigskip

\subsection{Rel\`evement au groupe}

\ass{Lemme}{Soient $x\in G_{ss}(F)$, $\omega$ un bon voisinage de $0$ dans $\mathfrak{g}_{x}(F)$ et $\varphi\in C_{c}^{\infty}(\mathfrak{g}_{x}(F))$. On suppose $\varphi$ tr\`es cuspidale et \`a support dans $\omega$. On suppose que $\theta_{\varphi}$ est invariant par $Z_{G}(x)(F)$ et que $x$ est elliptique, c'est-\`a-dire que $A_{G_{x}}=A_{G}$. Alors il existe $f\in C_{c}^{\infty}(G(F))$ telle que $f$ soit tr\`es cuspidale et $\theta_{f,x,\omega}=\theta_{\varphi}$.}

Preuve. Fixons un L\'evi minimal $R_{min}$ de $G_{x}(F)$ et un sous-groupe compact sp\'ecial $K$ de $G(F)$ en bonne position relativement \`a ${\bf R}_{min}$. Fixons un  sous-groupe ouvert compact $K'$ de $K$ tel que $\varphi$ soit invariante par conjugaison par $K'\cap Z_{G}(x)(F)$. Posons $\Sigma=\{k^{-1}xexp(X)k; k\in K', X\in \omega\}$. C'est un sous-ensemble ouvert et ferm\'e de $G(F)$. Comme en 3.1, on peut d\'efinir une fonction $f$ sur $G(F)$ par
$$f(g)=\left\lbrace\begin{array}{cc}0,&\,\,{\rm si}\,\,g\not\in \Sigma,\\ \varphi(X),&\,\,{\rm si}\,\,g=k^{-1}xexp(X)k,\,\,{\rm avec}\,\,k\in K',\,X\in \omega.\\ \end{array}\right.$$
Montrons que $f$ est tr\`es cuspidale. Soient $P=MU$ un sous-groupe parabolique propre de $G$ et $m\in M(F)\cap G_{reg}(F)$.  Posons 
$$f(U,m)=\int_{U(F)}f(u^{-1}mu)du.$$
On veut prouver que $f(U,m)=0$. C'est \'evident si $\{u^{-1}mu; u\in U(F)\}\cap \Sigma=\emptyset$. Supposons cette intersection non vide. On peut fixer $v\in U(F)$, $k\in K'$ et $X\in \omega$ tels que $v^{-1}mv=k^{-1}xexp(X)k$. Posons $g=vk^{-1}$, $P'=g^{-1}Pg$, $M'=g^{-1}Mg$, $U'=g^{-1}Ug$ et $m'=g^{-1}mg$. Gr\^ace au changement de variable $u\mapsto uv $ et \`a l'invariance de $f$ par $K'$, on a l'\'egalit\'e
$$f(U,m)=\int_{U(F)}f(kv^{-1}u^{-1}muvk^{-1})du,$$
puis
$$f(U,m)=\int_{U(F)}f(g^{-1}u^{-1}mug)du=\int_{U'(F)}f(u^{_{'}-1}m'u')du'=f(U',m').$$
Cela nous ram\`ene \`a la m\^eme question pour le sous-groupe parabolique $P'$ et l'\'el\'ement $m'$ de $M'(F)$. Mais $m'=xexp(X)$. En oubliant les donn\'ees $P'$ et $m'$, on peut donc supposer que $m=xexp(X)$, avec $X\in \omega$. Dans ce cas, on a $A_{M}\subset G_{m}\subset G_{x}$, donc $x\in M(F)$. De plus, puisque $P$ est propre, on a $A_{G}\subsetneq A_{M}$. Puisque $x$ est elliptique, on a aussi $A_{G_{x}}\subsetneq A_{M}\subset A_{M_{x}}$, donc $P_{x}=M_{x}U_{x}$ est un sous-groupe parabolique propre de $G_{x}$. Ecrivons
$$f(U,m)=\int_{U_{x}(F)\backslash U(F)}\int_{U_{x}(F)}f(v^{-1}u^{-1}muv)du\,dv.$$
On peut fixer $v\in U(F)$ et prouver que
$$\int_{U_{x}(F)}f(v^{-1}u^{-1}muv)du=0.$$
De nouveau, c'est \'evident si $\{v^{-1}u^{-1}muv; u\in U_{x}(F)\}\cap \Sigma=\emptyset$. Supposons cette intersection non vide. Alors il existe $w\in U_{x}(F)$ et $k\in K'$ tels que $v^{-1}w^{-1}mwv\in k^{-1}xexp(\omega)k$. Puisque $m\in xexp(\omega)$, la condition 3.1(4) entra\^{\i}ne que $wvk^{-1}\in Z_{G}(x)(F)$. Donc $v\in Z_{G}(x)(F)k$. Ecrivons $v=gk$, avec $g\in Z_{G}(x)(F)$. Pour $u\in U_{x}(F)$, on a $v^{-1}u^{-1}muv=k^{-1}xexp(g^{-1}u^{-1}Xug)k$. Donc 
$$f(v^{-1}u^{-1}muv)=\varphi(g^{-1}u^{-1}Xug)={^g\varphi}(u^{-1}Xu),$$
avec une d\'efinition \'evidente de $^g\varphi$. L'int\'egrale \`a calculer est \'egale \`a
$$\int_{U_{x}(F)}{^g\varphi}(u^{-1}Xu)du.$$
Elle est nulle parce que $^g\varphi$ est tr\`es cuspidale et $P_{x}$ est propre.

Posons
$$c=[G_{x}(F)\backslash Z_{G}(x)(F)K'/K']mes(K')mes(G_{x}(F)\cap K')^{-1}.$$
On va prouver  que $\theta_{f,x,\omega}=c\theta_{\varphi}$. En explicitant les d\'efinitions, on voit qu'il s'agit de prouver l'assertion suivante. Soit $X\in \omega\cap \mathfrak{g}_{x,reg}(F)$. Notons $R$ le commutant de $A_{G_{x,X}}$ dans $G_{x}$. Alors on a l'\'egalit\'e
$$(1) \qquad J_{{\bf R}}(xexp(X),f)=cD^G(x)^{1/2}J_{R}(X,\varphi).$$
On peut supposer que $R$ contient $R_{min}$. Reprenons la preuve du lemme 5.7, o\`u on prend pour $K'$ le groupe introduit ci-dessus. Par un raisonnement fait plusieurs fois, ${^gf}_{x,\omega}=0$ si $g\in G(F)$ et $g\not\in Z_{G}(x)(F)K'$. On peut prendre pour ensemble $\Delta_{0}$ un ensemble de repr\'esentants de $G_{x}(F)\backslash Z_{G}(x)(F)K'/K'$, inclus dans $Z_{G}(x)(F)$. Pour $\delta\in \Delta_{0}$, on a $^{\delta}f_{x,\omega}={^{\delta}\varphi}$. L'hypoth\`ese que $\varphi$ est tr\`es cuspidale entra\^{\i}ne que pour $S\in {\cal L}(R_{min})$, $S\neq G_{x}$, la fonction $f^S$ est nulle. Pour $S=G_{x}$, on a ${\bf S}=G$ parce que $x$ est elliptique. Alors
$$f^{G_{x}}=mes(K')mes(G_{x}(F)\cap K')^{-1}\sum_{\delta\in \Delta_{0}} {^{\delta}\varphi}.$$
D'apr\`es l'hypoth\`ese d'invariance de $\theta_{\varphi}$ par $Z_{G}(x)(F)$, on a
$$J_{R}(X,{^{\delta}\varphi})=J_{R}(X,\varphi),$$
pour tout $\delta\in \Delta_{0}$, et l'\'egalit\'e (1) r\'esulte du (i) du lemme 5.7. $\square$

\subsection{Support des distributions associ\'ees aux fonctions tr\`es cuspidales}

\ass{Lemme}{(i) Soient $ \theta$ un  quasi-caract\`ere de $\mathfrak{g}(F))$ et $\omega$ un $G$-domaine dans $\mathfrak{g}(F)$ compact modulo conjugaison. On suppose  la transform\'ee de Fourier de $\theta$  \`a support compact modulo conjugaison. Alors il existe une famille finie $(X_{i})_{i=1,...,n}$ d'\'el\'ements de $\mathfrak{g}_{reg}(F)$ et une famille finie $(c_{i})_{i=1,...,n}$ de nombres complexes telles que, pour tout $Y\in \omega\cap \mathfrak{g}_{reg}(F)$, on ait l'\'egalit\'e
$$\theta(Y)=\sum_{i=1,...,n}c_{i}\hat{j}(X_{i},Y).$$

(ii) Soient $X\in \mathfrak{g}_{reg}(F)$ et $\omega$ un $G$-domaine dans $\mathfrak{g}(F)$ compact modulo conjugaison. Alors il existe une fonction tr\`es cuspidale $f\in C_{c}^{\infty}(\mathfrak{g}(F))$ telle que, pour tout $Y\in \omega\cap \mathfrak{g}_{reg}(F)$, on ait l'\'egalit\'e
$$\theta_{f}(Y)=\hat{j}(X,Y).$$

(iii) Soient $\omega$ un $G$-domaine dans $\mathfrak{g}(F)$ compact modulo conjugaison et ${\cal O}\in Nil(\mathfrak{g})$. Il existe une  fonction $f\in C_{c}^{\infty}(\mathfrak{g}(F)$, tr\`es cuspidale, telle que, pour tout $Y\in \omega\cap \mathfrak{g}_{reg}(F)$, on ait l'\'egalit\'e
$$\theta_{f}(Y)=\hat{j}({\cal O},Y).$$

(iv) Pour tout $Y\in \mathfrak{g}_{reg}(F)$, il existe une fonction tr\`es cuspidale $f\in C_{c}^{\infty}(\mathfrak{g}(F))$ telle que $\theta_{f}(Y)\not=0$.}

Preuve.  Soit $\theta$ comme en (i).  Soit $\Omega$ un $G$-domaine de $\mathfrak{g}(F)$, compact modulo conjugaison et contenant le support de $\hat{\theta}$. Pour $\varphi\in C_{c}^{\infty}(\omega)$, on a, avec les notations de 2.6(2) , 
$$\theta(\varphi)=\sum_{i=1,...,n}J_{G}(X_{i},\hat{\varphi})\theta(f_{i}).$$
 Autrement dit
 $$\int_{\mathfrak{g}(F)}\theta(Y)\varphi(Y)dY=\int_{\mathfrak{g}(F)}\sum_{i=1,...,n}\theta(f_{i})\hat{j}(X_{i},Y)\varphi(Y)dY.$$
 L'assertion (i) s'ensuit.
 
 Fixons un L\'evi minimal $M_{min}$ de $G$. Soit $f'\in C_{c}^{\infty}(\mathfrak{g}(F))$, supposons $f'$ tr\`es cuspidale. Pour toute $\varphi\in \mathfrak{g}(F)$, on a l'\'egalit\'e:
 $$\int_{\mathfrak{g}(F)}\theta_{\hat{f}'}(Y)\varphi(Y)dY=\int_{\mathfrak{g}(F)}\theta_{f'}(X)\hat{\varphi}(X) dX$$
 $$(1) \qquad =\sum_{M\in {\cal L}(M_{min})}\vert W^M\vert \vert W^G\vert ^{-1}\sum_{T\in {\cal T}_{ell}(M)}\vert W(M,T)\vert ^{-1}\int_{\mathfrak{t}(F)}\theta_{f'}(X)J_{G}(X,\hat{\varphi})D^G(X)^{1/2}dX.$$
 Soient $X$ et $\omega$ comme en (ii). On ne perd rien \`a supposer qu'il existe $M\in {\cal L}(M_{min})$ et $T\in {\cal T}_{ell}(M)$ de sorte que $X\in \mathfrak{t}(F)$.  La formule 2.6(3)  montre qu'il existe un voisinage $\omega'_{X}$ de $X$ dans $\mathfrak{g}(F)$ tel que, pour tout $Y\in \omega\cap \mathfrak{g}_{reg}(F)$, la fonction $X'\mapsto \hat{j}(X',Y)$ soit constante dans $\omega'_{X}$. Supposons
 
 (2) il existe $f'\in C_{c}^{\infty}(\mathfrak{g}(F))$, tr\`es cuspidale, telle que $\theta_{f'}(X)\neq 0$.
 
 Fixons une telle fonction. On peut trouver un voisinage $\omega_{X}$ de $X$ dans $\mathfrak{t}(F)$ tel que
 
 - $\omega_{X}$ est ouvert et compact;

 - la fonction $X'\mapsto \theta_{f'}(X')D^G(X')^{1/2}$ soit constante dans $\omega_{X}$;
 
 - pour $w\in Norm_{M(F)}(T)$ tel que $w\not\in T(F)$, $w^{-1}\omega_{X}w\cap \omega_{X}=\emptyset$.
 
 Rempla\c{c}ons $f'$ par son produit avec la fonction caract\'eristique du $G$-domaine $\omega_{X}^G$. La formule (1) devient
 $$\int_{\mathfrak{g}(F)}\theta_{\hat{f}'}(Y)\varphi(Y)dY=\theta_{f'}(X)D^G(X)^{1/2}\int_{\omega_{X}}J_{G}(X',\hat{\varphi})dX'$$
 $$\qquad =\theta_{f'}(X)D^G(X)^{1/2}\int_{\omega_{X}}\int_{\mathfrak{g}(F)}\hat{j}(X',Y)\varphi(Y)dY\,dX'.$$
 La double int\'egrale est bien s\^ur absolument convergente. On en d\'eduit l'\'egalit\'e
 $$(3) \qquad \theta_{\hat{f}'}(Y)=\theta_{f'}(X)D^G(X)^{1/2}\int_{\omega_{X}}\hat{j}(X',Y)dX'$$
 pour tout $Y\in \mathfrak{g}_{reg}(F)$. Imposons de plus la condition $\omega_{X}\subset \omega'_{X}$.  Alors l'\'egalit\'e (3) devient
     $$\theta_{\hat{f'}}(Y)=\theta_{f'}(X)D^G(X)^{1/2}mes(\omega_{X})\hat{j}(X,Y)$$
  pour tout $Y\in \omega\cap\mathfrak{g}_{reg}(F)$. En posant
  $$f=\theta_{f'}(X)^{-1}D^G(X)^{-1/2}mes(\omega_{X})^{-1}\hat{f'},$$
  on obtient (ii), sous l'hypoth\`ese (2).
  
  Soit $Y\in \mathfrak{g}_{reg}(F)$. Appliquons 2.6(4) : on peut fixer un $G$-domaine $\Omega$ de $\mathfrak{g}(F)$ tel que
  $$(4)\qquad    \hat{j}(X,Y)=\sum_{{\cal O}\in Nil(\mathfrak{g})}\Gamma_{{\cal O}}(X)\hat{j}({\cal O},Y)$$
  pour tout $X\in \Omega\cap \mathfrak{g}_{reg}(F)$. On peut \'evidemment supposer $\lambda\Omega\subset \Omega$ pour tout $\lambda\in F^{\times}$ tel que $\vert \lambda\vert _{F}\leq 1$. Soit $X\in \Omega\cap \mathfrak{g}_{reg}(F)$ tel que $X$ soit elliptique. Alors l'hypoth\`ese (2) est v\'erifi\'ee car toute fonction \`a support r\'egulier elliptique est tr\`es cuspidale. On peut donc trouver une fonction $f_{X}$ tr\`es cuspidale telle que $\theta_{f_{X}}(Y)=\hat{j}(X,Y)$. Soit $\lambda\in F^{\times 2} $ tel que $\vert \lambda\vert _{F}\leq 1$. Rempla\c{c}ons $X$ par $\lambda X$. En utilisant (4) et les formules de 2.6, on a
  $$\theta_{f_{\lambda X}}(Y)=\sum_{{\cal O}\in Nil(\mathfrak{g})}\vert \lambda\vert _{F}^{(\delta(G)-dim({\cal O}))/2}\Gamma_{{\cal O}}(X)\hat{j}({\cal O},Y). $$
  En prenant une combinaison lin\'eaire convenable des fonctions $f_{\lambda X}$, on peut s\'eparer les orbites selon leurs dimensions. On peut en particulier isoler l'orbite $\{0\}$ et trouver une fonction $f$ tr\`es cuspidale telle que
  $$\theta_{f}(Y)=\Gamma_{\{0\}}(X)\hat{j}(\{0\},Y).$$
  La fonction $\hat{j}(\{0\},.)$ est constante de valeur $1$. D'apr\`es Harish-Chandra ([HCDS] lemme 9.6), le germe de Shalika $\Gamma_{\{0\}}$ ne s'annule pas sur l'ensemble des points elliptiques r\'eguliers. Donc $\theta_{f}(Y)\neq 0$, ce qui d\'emontre (iv).
  
  On peut maintenant achever la d\'emonstration de (ii): d'apr\`es (iv), l'hypoth\`ese (2) est v\'erifi\'ee pour tout point $X\in \mathfrak{g}_{reg}(F)$.
  
  Soit $\omega$ un $G$-domaine de $\mathfrak{g}(F)$ compact modulo conjugaison. Choisissons $\Omega$ tel que (4) soit v\'erifi\'ee pour tous $X\in \Omega\cap \mathfrak{g}_{reg}(F)$, $Y\in \omega\cap \mathfrak{g}_{reg}(F)$. Gr\^ace \`a (ii), pour tout $X\in \Omega\cap \mathfrak{g}_{reg}(F)$,on peut trouver une fonction $f_{X}$ tr\`es cuspidale telle que $\theta_{f_{X}}(Y)=\hat{j}(X,Y)$ pour tout $Y\in \omega\cap \mathfrak{g}_{reg}(F)$. On sait que les germes de Shalika sont lin\'eairement ind\'ependants ([HCDS] lemme 9.5). Etant homog\`enes, leurs restrictions \`a $\Omega$ le sont aussi.  En utilisant (4), on voit que, pour ${\cal O}\in Nil(\mathfrak{g})$, une combinaison lin\'eaire convenable $f$ de fonctions $f_{X}$ va v\'erifier $\theta_{f}(Y)=\hat{j}({\cal O},Y)$ pour tout $Y\in \omega\cap \mathfrak{g}_{reg}(F)$. Cela prouve (iii). $\square$
  
  \bigskip
  
  \subsection{ Quasi-caract\`eres \`a support compact modulo conjugaison}
  
  \ass{Proposition}{Soit $\theta$ un quasi-caract\`ere de $\mathfrak{g}(F)$ \`a support compact modulo conjugaison. Alors il existe une fonction $f\in C_{c}^{\infty}(\mathfrak{g}(F))$ tr\`es cuspidale telle que $\theta=\theta_{f}$.}
  
  Preuve. On d\'emontre la proposition par r\'ecurrence sur $dim(G)$. On suppose qu'elle est vraie pour tout groupe de dimension strictement inf\'erieure \`a $dim(G)$. Soit $X\in \mathfrak{g}_{ss}(F)$. On va prouver
  
  (1) il existe un $G$-domaine $\omega$ dans $\mathfrak{g}(F)$ et une fonction tr\`es cuspidale $f\in C_{c}^{\infty}(\mathfrak{g}(F))$ tels que $X\in \omega$ et $\theta(Y)=\theta_{f}(Y)$ pour tout $Y\in \omega$.
  
  Autrement dit, il existe une fonction tr\`es cuspidale $f\in C_{c}^{\infty}(\mathfrak{g}(F))$ telle que $\theta_{f}$ ait le m\^eme d\'eveloppement que $\theta$ au voisinage de $X$. Si $X=0$, cela r\'esulte du lemme 6.3(iii). Si $X$ est central, cela r\'esulte du cas $X=0$ par translation. Supposons $X$ non central, donc $dim(G_{X})<dim(G)$. Supposons d'abord $X$ elliptique. La notions de bon voisinage s'adapte au cas de l'alg\`ebre de Lie. Soit $\omega_{X}$ un bon voisinage de $0$ dans  $\mathfrak{g}_{X}(F)$. Soit $\theta_{X}$ la fonction sur $\mathfrak{g}_{X}(F)$ qui est nulle hors de $X+\omega_{X}$ et qui co\"{\i}ncide avec $\theta$ sur $X+\omega_{X}$. Alors $\theta_{X}$ est un quasi-caract\`ere de $\mathfrak{g}_{X}(F)$ \`a support compact modulo conjugaison. Appliquant l'hypoth\`ese de r\'ecurrence, on peut fixer une fonction tr\`es cuspidale $f_{X}\in C_{c}^{\infty}(\mathfrak{g}_{X}(F))$ telle que $\theta_{f_{X}}=\theta_{X}$. On peut reprendre la d\'emonstration du lemme 6.2 et montrer qu'il existe une fonction tr\`es cuspidale $f\in C_{c}^{\infty}(\mathfrak{g}(F))$ telle que $\theta_{f}$ co\"{\i}ncide avec $\theta_{f_{X}}$ dans $X+\omega_{X}$ (remarquons que la condition d'invariance par $Z_{G}(x)(F)$ impos\'ee en 6.2 dispara\^{\i}t  car le commutant d'un \'el\'ement semi-simple d'une alg\`ebre de Lie est toujours connexe). En posant $\omega=(X+\omega_{X})^G$, $f$ et $\omega$ satisfont (1). Supposons maintenant $X$ non elliptique. Notons $M$ le commutant de $A_{G_{X}}$ dans $G$. On a $X\in \mathfrak{m}(F)$ et $G_{X}\subset M$. Soit $\omega_{X}$ un bon voisinage de $0$ dans $\mathfrak{g}_{X}(F)$. Posons
  $$\omega_{M}=(X+\omega_{X})^M,\,\,\omega_{G}=(X+\omega_{X})^G.$$
  Les ensembles $\omega_{M}$ et $\omega_{G}\cap \mathfrak{m}(F)$ sont des $M$-domaines dans $\mathfrak{m}(F)$, compacts modulo conjugaison. Notons $\theta_{M}$ la fonction sur $\mathfrak{m}(F)$ qui est nulle hors de $\omega_{M}$ et co\"{\i}ncide avec $\theta$ sur $\omega_{M}$.   C'est un quasi-caract\`ere de $\mathfrak{m}(F)$, \`a support compact modulo conjugaison. En appliquant l'hypoth\`ese de r\'ecurrence et le (i) du lemme 6.3, on peut fixer des familles finies $(X_{i})_{i=1,...,n}$ d'\'el\'ements de $\mathfrak{m}_{reg}(F)$ et $(c_{i})_{i=1,...,n}$ de nombres complexes de sorte que
 $$\theta_{M}(Y)=\sum_{i=1,...,n}c_{i}\hat{j}^M(X_{i},Y)$$
 pour tout $Y\in \omega_{G}\cap \mathfrak{m}_{reg}(F)$. La preuve du (i) du lemme 6.3 montre que l'on peut aussi bien remplacer chaque $X_{i}$ par tout \'el\'ement suffisamment proche. On peut donc supposer $X_{i}\in \mathfrak{g}_{reg}(F)$. La fonction
 $$Y\mapsto D^G(Y)^{1/2}D^M(Y)^{-1/2}$$
 est constante sur $X+\omega_{X}$. Notons $c$ sa valeur. Montrons que
 $$(2) \qquad \theta(Y)=\sum_{i=1,...,n}cc_{i}\hat{j}^G(X_{i},Y)$$
 pour tout $Y\in (X+\omega_{X})\cap \mathfrak{g}_{reg}(F)$. D'apr\`es 2.6(5), le membre de droite ci-dessus est \'egal \`a
 $$\sum_{Y'\in {\cal Y}}\sum_{i=1,...,n}cc_{i}\hat{j}^M(X_{i},Y')D^G(Y)^{-1/2}D^M(Y')^{1/2},$$
 o\`u ${\cal Y}$ est un ensemble de  repr\'esentants des classes de conjugaison par $M(F)$ dans l'ensemble des \'el\'ements de $\mathfrak{m}(F)$ conjugu\'es \`a $Y$ par un \'el\'ement de $G(F)$. Cet ensemble ${\cal Y}$ est contenu dans $\omega_{G}\cap \mathfrak{m}(F)$. La somme ci-dessus vaut donc
 $$\sum_{Y'\in {\cal Y}}c\theta_{M}(Y')D^G(Y)^{-1/2}D^M(Y')^{1/2}.$$
 On peut supposer $Y\in {\cal Y}$ et le terme index\'e par $Y$ est \'egal \`a $\theta(Y)$. Soit $Y'\in {\cal Y}$, $Y'\not=Y$. Soit $g\in G(F)$ tel que $gYg^{-1}=Y'$. Cet \'el\'ement n'appartient pas \`a $M(F)$. Si $Y'\in \omega_{M}$, il existe $m\in M(F)$ tel que $mgY(mg)^{-1}\in X+\omega_{X}$. Alors $mg\in G_{X}(F)$ puisque $\omega_{X}$ est un bon voisinage de $0$ dans $\mathfrak{g}_{X}(F)$. Puisque $G_{X}\subset M$, on en d\'eduit $g\in M(F)$, contradiction. Donc $Y'\not\in \omega_{M}$ et $\theta_{M}(Y')=0$.
 Cela d\'emontre (2).

  D'apr\`es le lemme 6.3(ii), il existe une fonction tr\`es cuspidale $f\in C_{c}^{\infty}(\mathfrak{g}(F))$ telle que $\theta_{f}(Y)$ co\"{\i}ncide avec le membre de droite de (2) dans $\omega_{G}$. Alors $\theta_{f}(Y)=\theta(Y)$ pour $Y\in \omega_{G}$, ce qui ach\`eve la preuve de (1).
 
 La preuve de la proposition est maintenant \'el\'ementaire. Pour tout $X\in \mathfrak{g}_{ss}(F)$, on fixe $\omega$ et $f$ v\'erifiant (1) et on les note plut\^ot $\omega_{X}$ et $f_{X}$. On fixe un sous-ensemble compact $\Gamma\subset \mathfrak{g}(F)$ tel que $Supp(\theta)\subset \Gamma^G$. On a
 $$\Gamma\subset \mathfrak{g}(F)\subset \bigcup_{X\in \mathfrak{g}_{ss}(F)}\omega_{X}.$$
 On peut donc choisir une famille finie $(X_{i})_{i=1,...,n}$ d'\'el\'ements de $\mathfrak{g}_{ss}(F)$ telle que 
 $$\Gamma\subset \bigcup_{i=1,...n}\omega_{X_{i}}.$$
 Pour tout $i$, notons $\Delta_{i}$ le compl\'ementaire dans $\mathfrak{g}(F)$ de $\bigcup_{j=1,...,i-1}\omega_{X_{i}}$, $\varphi_{i}$ la fonction caract\'eristique de $\omega_{X_{i}}\cap \Delta_{i}$ et posons$f_{i}=f_{X_{i}} \varphi_{i}$. Alors $f_{i}$ est tr\`es cuspidale et on a l'\'egalit\'e $\theta=\sum_{i=1,...,n}\theta_{f_{i}}$. $\square$
 
\bigskip

\section{Enonc\'e du th\'eor\`eme principal}

\bigskip

\subsection{Groupes orthogonaux}

Soit $V$ un espace vectoriel sur $F$ de dimension finie $d$, muni d'une forme bilin\'eaire sym\'etrique et non d\'eg\'en\'er\'ee $q$ (on dira parfois que $V$ est un espace quadratique,la forme $q$ \'etant sous-entendue). Pour $v\in V$, on pose 
$$q(v)=\frac{1}{2}q(v,v).$$ 
 On appelle syst\`eme hyperbolique dans $V$ une famille $(v_{i})_{i=\pm 1,...,\pm n}$ d'\'el\'ements de $V$ telle que $q(v_{i},v_{j})=\delta_{i,-j}$ pour tous $i,j$, o\`u $\delta_{i,-j}$ est le symbole de Kronecker. On dit que $q$ est hyperbolique si et seulement s'il existe un syst\`eme hyperbolique dans $V$ qui soit une base de $V$. On peut d\'ecomposer $V$ (de fa\c{c}on non unique) en somme orthogonale de deux sous-espaces $V_{hyp}$ et $V_{an}$ de sorte que la restriction de $q$ \`a $V_{hyp}$ soit hyperbolique et la restriction $q_{an}$ de $q$ \`a $V_{an}$ soit anisotrope. La classe d'\'equivalence de $q_{an}$ est uniquement d\'efinie, on l'appelle le noyau anisotrope de $q$ et on note $d_{an}(V)$  son rang, c'est-\`a-dire la dimension de $V_{an}$. Evidemment, $d_{an}(V)\equiv d\,\,mod\,\,2$.

On introduit le groupe orthogonal $O(V)$ de $(V,q)$ et son sous-groupe sp\'ecial orthogonal $SO(V)$. Notons-les ici $G^+$ et $G$, ainsi qu'on le fera souvent dans la suite. Le groupe $G^+(F)$ agit sur $V$, on note cette action $(g,v)\mapsto gv$. Le groupe $G$ est d\'eploy\'e sur $F$ si $d_{an}(V)=0$ ou $1$, quasi-d\'eploy\'e et non d\'eploy\'e si $d_{an}(V)=2$ et non quasi-d\'eploy\'e si $d_{an}(V)=3$ ou $4$.  On d\'efinit la forme bilin\'eaire sur $\mathfrak{g}(F)$
$$<X,Y>=\frac{1}{2}trace(XY),$$
la trace \'etant la forme lin\'eaire usuelle sur $End(V)$. C'est cette forme que l'on utilise pour normaliser les constructions de 1.2.

Si $W$ est un sous-espace non d\'eg\'en\'er\'e de $V$, le groupe  $H=SO(W)$ se plonge naturellement dans   $G$: un \'el\'ement de $H$ agit par l'identit\'e sur l'orthogonal de $W$. De m\^eme, $\mathfrak{h}$ se plonge dans $\mathfrak{g}$.

Soient $v',v''\in V$. D\'efinissons $c_{v',v''}\in End(V)$ par
$$c_{v',v''}(v)=q(v,v')v''-q(v,v'')v'.$$
On v\'erifie que $c_{v',v''}$ appartient \`a $\mathfrak{g}(F)$ et que cette alg\`ebre est engendr\'ee en tant qu'espace vectoriel par de tels \'el\'ements.

On peut classifier les orbites unipotentes r\'eguli\`eres de $\mathfrak{g}(F)$. Il n'y en a pas si $d_{an}(V)\geq 3$. Il y en a une seule si $d_{an}(V)=1$ ou si $d\leq 2$ et on la note ${\cal O}_{reg}$. Supposons $d\geq 4$ et $d_{an}(V)=0$ ou $2$. Introduisons le sous-ensemble ${\cal N}^V$ suivant:

- si $d_{an}(V)=0$, ${\cal N}^V=F^{\times}/F^{\times 2}$;

- si $d_{an}(V)=2$, ${\cal N}^V$ est le sous-ensemble des \'el\'ements de $F^{\times}/F^{\times 2}$ qui sont repr\'esent\'es par $q_{an}$. 

Soit $\nu\in {\cal N}^V$, dont on fixe un rel\`evement dans $F^{\times}$. On peut d\'ecomposer $V$ en somme orthogonale $V=D\oplus W$, o\`u $D$ est une droite et la restriction de $q$ \`a $W$ a pour noyau anisotrope la forme $x\mapsto \nu x^2$ de dimension $1$. Notons $H=SO(W)$. Il y a une unique orbite nilpotente r\'eguli\`ere dans $\mathfrak{h}(F)$. Soit $N$ un \'el\'ement de cette orbite, que l'on identifie \`a un \'el\'ement de $\mathfrak{g}(F)$. On note ${\cal O}_{\nu}$ la $G(F)$-orbite de $N$. Elle ne d\'epend pas des choix et l'application $\nu\mapsto {\cal O}_{\nu}$ est une bijection de ${\cal N}^V$ sur l'ensemble des orbites nilpotentes r\'eguli\`eres de $\mathfrak{g}(F)$.

Consid\'erons une d\'ecomposition orthogonale $V=Z\oplus V_{an}$. Supposons la restriction de $q$ \`a $V_{an}$ anisotrope et la restriction de $q$ \`a $Z$ hyperbolique. Fixons une base hyperbolique $(v_{i})_{i=\pm 1,...\pm n}$ de $Z$. Si $V_{an}\not=\{0\}$, soit $c\in {\mathbb Z}$ tel qu'il existe $v\in V_{an}$ pour lequel $val_{F}(q(v))=c$. Si $V_{an}=\{0\}$, soit $c$ un \'el\'ement quelconque de $ {\mathbb Z}$. Notons $R_{an}$ l'ensemble des $v\in V_{an}$ tels que $val_{F}(q(v))\geq c$. C'est un $\mathfrak{o}_{F}$-r\'eseau de $V_{an}$. Notons $R_{Z}$ le $\mathfrak{o}_{F}$-r\'eseau de $Z$ engendr\'e par les \'el\'ements $v_{i}$ pour $i=1,...,n$ et par les $\varpi_{F}^cv_{-i}$ pour $i=1,...,n$. Posons $R=R_{Z}\oplus R_{an}$. Quand on fait varier la d\'ecomposition orthogonale, la base hyperbolique et l'entier $c$, le r\'eseau $R$ parcourt un certain ensemble de $\mathfrak{o}_{F}$-r\'eseaux de $V$. Par d\'efinition, c'est l'ensemble des r\'eseaux sp\'eciaux. Pour un tel r\'eseau $R$ sp\'ecial, notons $K$ le stabilisateur de $R$ dans $G(F)$. C'est un sous-groupe compact sp\'ecial de $G(F)$. Inversement, si $K$ est un sous-groupe compact sp\'ecial de $G(F)$, il est le stabilisateur d'un r\'eseau sp\'ecial $R$ (qui n'est pas unique).

 Pour tout  $\mathfrak{o}_{F}$-r\'eseau $R$ de $V$, on d\'efinit une fonction $val_{R}$ sur $V$, \`a valeurs dans ${\mathbb Z}\cup\{\infty\}$ par $val_{R}(v)=sup\{i\in {\mathbb Z}; v\in \varpi_{F}^{i}R\}$. 

\subsection{La situation}

 On conserve les donn\'ees et notations du paragraphe pr\'ec\'edent. Soit $r\in {\mathbb N}$ tel que $2r+1\leq d$. On suppose donn\'ee une d\'ecomposition orthogonale $V=W\oplus D\oplus Z$, o\`u $D$ est une droite et $Z$ un espace hyperbolique de dimension $2r$.  On  note $q_{W}$ la restriction de $q$ \`a $W$ et $d_{W}=d-2r-1$ la dimension de $W$.   On pose $V_{0}=W\oplus D$. On note $H$, resp. $G_{0}$, le groupe sp\'ecial orthogonal de $W$, resp. $V_{0}$, et $H^+$ le groupe orthogonal de $W$. On identifie $H^+$ \`a un sous-groupe de $G$: un \'el\'ement $h\in H^+$ s'identifie \`a l'\'el\'ement de $G$ qui agit par $h$ sur $W$ et par $det(h)$ sur $D\oplus Z$.  On fixe une base $v_{0}$ de $D$ et un syst\`eme hyperbolique maximal $(v_{i})_{i=\pm 1,...,\pm r}$ de $Z$. On note $Z_{+}$, resp. $Z_{-}$, le sous-espace de $Z$ engendr\'e par les $v_{i}$, $i=1,...,r$, resp. par les $v_{-i}$. On note $A$ le sous-tore maximal de $SO(Z)$ qui conserve chaque droite $Fv_{i}$.  Pour $a\in A(F)$ et $i=1,...,r$, on note $a_{i}\in F^{\times}$ la valeur propre de $a$ sur $v_{i}$, c'est-\`a-dire que $av_{i}=a_{i}v_{i}$. On note $P$ le sous-groupe parabolique de $G$ form\'e des \'el\'ements qui conservent le drapeau
$$Fv_{r}\subset Fv_{r}\oplus Fv_{r-1}\subset...\subset Fv_{r}\oplus...\oplus Fv_{1}$$
de $V$. On note $U$ le radical unipotent de $P$ et $M$ sa composante de L\'evi qui contient $A$. On a $M=AG_{0}$. Remarquons que $A_{M}=A$, sauf dans le cas o\`u $V_{0}$ est hyperbolique de dimension $2$, auquel cas $A_{M}=M$.  Fixons une famille $(\xi_{i})_{i=0,...,r-1}$ d'\'el\'ements de $F^{\times}$. On d\'efinit une fonction $\xi$ sur $U(F)$ par
$$\xi(u)=\psi(\sum_{i=0,...,r-1}\xi_{i}q(uv_{i},v_{-i-1})).$$
On v\'erifie que c'est un caract\`ere de $U(F)$ invariant par conjugaison par le sous-groupe $H^+(F)$ de $M(F)$. 

On fixe un r\'eseau sp\'ecial $R_{0}$ de $V_{0}$. On peut choisir un r\'eseau $R_{Z}$ de $Z$ ayant une base form\'ee de vecteurs proportionnels aux $v_{i}$, de sorte que le r\'eseau $R=R_{0}\oplus R_{Z}$ de $V$ soit sp\'ecial. On note $K_{0}$, resp. $K$,  le stabilisateur de $R_{0}$ dans $G_{0}(F)$, resp. de $R$ dans $G(F)$. Ce sont des  sous-groupes compacts sp\'eciaux de $G_{0}(F)$, resp. $G(F)$.  Le groupe $K$ est en bonne position relativement \`a $M$. On a $K\cap M(F)=(K\cap A(F))K_{0}$ et $K\cap A(F)$ est le plus grand sous-groupe compact de $A(F)$. Pour tout entier $N\in {\mathbb N}$, on d\'efinit une fonction $\kappa_{N}$ sur $G(F)$ de la fa\c{c}on suivante. Elle est invariante \`a droite par $K$, \`a gauche par $U(F)$. Sa restriction \`a $M(F)$ est la fonction caract\'eristique de l'ensemble des $ag_{0}$, avec $a\in A(F)$, $g_{0}\in G_{0}(F)$, tels que $\vert val_{F}(a_{i})\vert \leq N$ pour tout $i=1,...,r$ et $val_{R_{0}}(g_{0}^{-1}v_{0})\geq -N$. La fonction $\kappa_{N}$ est invariante \`a gauche par le sous-groupe $(K\cap A(F))H^+(F)$ de $M(F)$. L'image de son support dans $U(F)H(F)\backslash G(F)$ est compacte. Plus pr\'ecis\'ement

(1) il existe $c>0$ tel que, pour tout entier $N\geq 1$ et tout $g\in G(F)$ pour lequel $\kappa_{N}(g)=1$, il existe $g'\in G(F)$ tel que $g\in U(F)H(F)g'$ et $\sigma(g')\leq cN$.

On peut \'ecrire $g=uag_{0}k$, avec $u\in U(F)$, $a\in A(F)$, $g_{0}\in G_{0}(F)$ et $k\in K$. Les bornes sur les valuations des coordonn\'ees de $a$ entra\^{\i}nent $\sigma(a)\leq cN$, pour $c$ convenable. Un raisonnement analogue \`a celui de la preuve du lemme III.5 de [W2] montre qu'il existe $c'>0$ tel que, pour tout $N\geq 1$ et tout $v\in \varpi_{F}^{-N}R_{0}$ tel que $q(v)=\nu_{0}=q(v_{0})$, il existe $y\in G_{0}(F)$ tel que $y^{-1}v_{0}=v$ et $\sigma(y)\leq c'N$. Appliquons cela \`a $v=g_{0}^{-1}v_{0}$. Alors $g_{0}y^{-1}$ appartient \`a $H(F)$. On a $g=ug_{0}y^{-1}g'$, avec $g'=ayk$ et ce dernier \'el\'ement satisfait la majoration requise.

\bigskip

\subsection{Les ingr\'edients de la formule int\'egrale}

On d\'efinit une fonction $\Delta$ sur $H_{ss}(F)$ par
$$\Delta(t)=\vert det((1-t)_{\vert W/W''(t)})\vert _{F},$$
o\`u $W''(t)$ est le noyau de $1-t$ agissant dans $W$.

Notons $\underline{{\cal T}}$ l'ensemble des sous-tores $T$ de $H$, en g\'en\'eral non maximaux, pour lesquels il existe une d\'ecomposition orthogonale $W=W'\oplus W''$ de sorte que les conditions (1) \`a (4) ci-dessous soient v\'erifi\'ees. On note $H'$ le groupe sp\'ecial orthogonal de $W'$ et on pose $V''=W''\oplus D\oplus Z$.

(1) $A_{T}=\{1\}$.

(2) $dim(W')$ est pair.

(3) $T$ est inclus dans $H'$ et c'en est un sous-tore maximal.

(4) Si $d$ est pair, $d_{an}(W'')=1$; si $d$ est impair, $d_{an}(V'')=1$.  

Evidemment, $W'$ et $W''$ sont d\'etermin\'es par $T$: $W''$ est l'intersection des noyaux de $t-1$ agissant sur $W$, pour $t\in T$.

Pour $T\in \underline{\cal T}$, on pose
$$W(H,T)=Norm_{H(F)}(T)/Z_{H(F)}(T).$$
On note $T_{\natural}$ le sous-ensemble des $t\in T$ tels que les valeurs propres de l'action de $t$ dans $W'$ soient toutes distinctes. C'est un ouvert de Zariski, non vide. Notons $H'=G'$, resp. $H''$, $G''$, les groupes sp\'eciaux orthogonaux de $W'$, resp. $W''$, $V''$ et $H^{_{'}+}$ le groupe orthogonal de $W'$. Pour $t\in T_{\natural}$, $t$ est un \'el\'ement semi-simple r\'egulier dans $H'$ et m\^eme dans $H^{_{'}+}$ en ce sens que son commutant dans $H^{_{'}+}$ est r\'eduit \`a $T$. Alors $Z_{H}(t)=TH''$, $Z_{G}(t)=TG''$. En particulier, les orbites nilpotentes de $\mathfrak{h}_{t}(F)$, resp. $\mathfrak{g}_{t}(F)$ sont les m\^emes que celles de $\mathfrak{h}''(F)$, resp. $\mathfrak{g}''(F)$.  

Soient $\theta$, resp. $\tau$, un quasi-caract\`ere de $H(F)$, resp. $G(F)$. Soit $T\in \underline{\cal T}$, pour lequel on adopte les notations ci-dessus. Soit $t\in T_{\natural}(F)$. Supposons $d$ pair.  Alors $dim(W'')$ est impair et l'hypoth\`ese (4) dit que $H''$ est d\'eploy\'e. Donc $\mathfrak{h}''(F)$ poss\`ede une unique orbite nilpotente r\'eguli\`ere ${\cal O}_{reg}$. On pose $c_{\theta}(t)=c_{\theta,{\cal O}_{reg}}(t)$. La dimension $dim(V'')$ est paire, n\'ecessairement non nulle. Si elle est \'egale \`a $2$, $\mathfrak{g}''(F)$ poss\`ede une unique orbite nilpotente r\'eguli\`ere ${\cal O}_{reg}$ et on pose $c_{\tau}(t)=c_{\tau,{\cal O}_{reg}}(t)$. Supposons $dim(V'')\geq 4$. Notons $q_{W'',an}$ le noyau anisotrope de la restriction de $q$ \`a $W''$ et $q_{D}$ la restriction de $q$ \`a $D$. Notons aussi $\nu_{0}=q(v_{0})$. La forme $q_{W'',an}$ est de rang $1$. Par construction, la restriction de $q$ \`a $V''$ a m\^eme noyau anisotrope que $q_{W'',an}\oplus q_{D}$. Elle est donc de rang $0$ ou $2$ et $G''$ est quasi-d\'eploy\'e. Puisque $q_{W'',an}\oplus q_{D}$ repr\'esente $\nu_{0}$, on a $\nu_{0}\in {\cal N}^{V''}$ et l'orbite ${\cal O}_{\nu_{0}}$ de $\mathfrak{g}''(F)$ est d\'efinie. On pose $c_{\tau}(t)=c_{\tau,{\cal O}_{\nu_{0}}}(t)$. Supposons maintenant $d$ impair. Alors $dim(V'')$ est impair et l'hypoth\`ese (4) dit que $G''$ est d\'eploy\'e. De fa\c{c}on analogue \`a ci-dessus, on pose $c_{\tau}(t)=c_{\tau,{\cal O}_{reg}}(t)$. La dimension de $W''$ est paire. Si elle est inf\'erieure ou \'egale \`a $2$, on pose encore $c_{\theta}(t)=c_{\theta,{\cal O}_{reg}}(t)$. Supposons $dim(W'')\geq 4$. Avec les m\^emes notations que ci-dessus, la restriction de $q$ \`a $V''$ a encore m\^eme noyau anisotrope que $q_{W'',an}\oplus q_{D}$. Si  $d_{an}(W'')$ \'etait \'egal \`a $4$, le noyau anisotrope de $q_{W'',an}\oplus q_{D}$ serait de rang $3$, contrairement \`a (4). Donc $d_{an}(W'')=0$ ou $2$. On a $-\nu_{0}\in {\cal N}^{W''}$: c'est \'evident si  $d_{an}(W'')=0$; si  $d_{an}(W'')=2$, cela r\'esulte du fait que $q_{W'',an}\oplus q_{D}$ n'est pas anisotrope puisque $d_{an}(V'')=1$. Donc l'orbite ${\cal O}_{-\nu_{0}}$ de $\mathfrak{h}''(F)$ est d\'efinie. On pose $c_{\theta}(t)=c_{\theta,{\cal O}_{-\nu_{0}}}(t).$

\ass{Proposition}{(i) Les fonctions $c_{\theta}$ et $c_{\tau}$ sont localement constantes sur $T_{\natural}(F)$.

(ii) La fonction $t\mapsto c_{\theta}(t)c_{\tau}(t)D^{H}(t)\Delta(t)^r$ est localement int\'egrable sur $T(F)$.}

La preuve est donn\'ee dans les quatre paragraphes suivants. Le tore $T\in \underline{\cal T}$  est fix\'e pour ces  paragraphes.

\bigskip

\subsection{Fonctions localement int\'egrables sur un \'el\'ement de $\underline{{\cal T}}$}

   Pour $t\in T(F)$, notons $E''(t)$ le noyau de $t-1$ dans $W'$ et $E'(t)$ son orthogonal dans $W'$. On note $J'(t)$, resp. $J''(t)$ le groupe sp\'ecial orthogonal de $E'(t)$, resp. $E''(t)$,  $J'(t)_{t}$ la composante neutre du commutant de $t$ dans $J'(t)$ et $\mathfrak{z}_{t}$ le centre de l'alg\`ebre de Lie $\mathfrak{j}'(t)_{t}$. Le groupe $H'_{t}$ pr\'eserve n\'ecessairement les espaces propres de $t$, donc est inclus dans $J'(t)J''(t)$. Puisque $J''(t)$ est \'evidemment inclus dans ce commutant $H'_{t}$, on a l'\'egalit\'e
$$H'_{t}=J'(t)_{t}J''(t).$$
On aura besoin des r\'esultats suivants.

(1) $\mathfrak{z}_{t}$ est inclus dans le centre de $\mathfrak{h}_{t}$ et dans le centre de $\mathfrak{g}_{t}$.

(2) Il existe un voisinage $\omega$ de $0$ dans $\mathfrak{t}(F)$ sur lequel l'exponentielle est d\'efinie et tel que, pour $X\in \omega$, on a $\mathfrak{z}_{t}\subset \mathfrak{z}_{texp(X)}$, avec \'egalit\'e si et seulement si $X\in \mathfrak{z}_{t}(F)$.

\bigskip

Preuve. Le noyau de $t-1$ dans $W$ est $E''(t)\oplus W''$ et son orthogonal dans $W$ est $E'(t)$. On en d\'eduit comme ci-dessus que $H_{t}$ est le produit de $J'(t)_{t}$ et du groupe sp\'ecial orthogonal de $E''(t)\oplus W''$.Donc le centre de $\mathfrak{h}_{t}$ est le produit de $\mathfrak{z}_{t}$ et du centre de l'alg\`ebre de Lie du second groupe. Cela prouve que $\mathfrak{z}_{t}$ est inclus dans le centre de $\mathfrak{h}_{t}$. Un raisonnement analogue vaut en rempla\c{c}ant $\mathfrak{h}_{t}$ par $\mathfrak{g}_{t}$.

 Consid\'erons un voisinage $\omega$ de $0$ dans $\mathfrak{t}(F)$ sur lequel l'exponentielle est d\'efinie et tel que $H'_{texp(X)}\subset H'_{t}$ pour tout $X\in \omega$. Soit $X\in \omega$. Puisque $X$ commute \`a $t$, il pr\'eserve les espaces propres de $t$, donc pr\'eserve $E'(t)$ et $E''(t)$. Notons $X'$ et $X''$ les restrictions de $X$ \`a chacun de ces deux espaces et posons $\tilde{t}=texp(X)$. Le groupe $H'_{\tilde{t}}$ est le commutant de $X$ dans $H'_{t}$.  Donc  $H'_{\tilde{t}}$ est le produit du commutant $J'(t)_{t,X'}$ de $X'$ dans $J'(t)_{t}$ et du commutant $J''(t)_{X''}$ de $X''$ dans $J''(t)$.
 En choisissant $\omega$ assez petit, on peut imposer que toutes les valeurs propres de $\tilde{t}$ dans $E'(t)$ soient diff\'erentes de $1$. Alors $E''(\tilde{t})\subset E''(t)$. Le groupe $J'(\tilde{t})_{\tilde{t}}$ est le sous-groupe des \'el\'ements de $H'_{\tilde{t}}$ qui agissent trivialement sur $E''(\tilde{t})$. Ce sous-groupe contient certainement $J'(t)_{t,X'}$ et est donc le produit de ce groupe et d'un certain sous-groupe de $J''(t)_{X''}$, que l'on note $\tilde{J}$. Donc $\mathfrak{z}_{\tilde{t}}$ est le produit du centre de $\mathfrak{j}'(t)_{t,X'}$ et du centre de $\tilde{\mathfrak{j}}$. L'alg\`ebre $\mathfrak{j}'(t)_{t,X'}$ est le commutant dans $\mathfrak{j}'(t)_{t}$ de l'\'el\'ement semi-simple $X'$ de cette alg\`ebre. Sur une extension de $F$, c'est donc une sous-alg\`ebre de L\'evi de $\mathfrak{j}'(t)_{t}$ et son centre contient le centre de cette alg\`ebre, c'est-\`a-dire contient $\mathfrak{z}_{t}$. Cela d\'emontre l'inclusion $\mathfrak{z}_{t}\subset \mathfrak{z}_{\tilde{t}}$. Supposons qu'il y a \'egalit\'e. Une sous-alg\`ebre de L\'evi \'etant le commutant de son centre, cela entra\^{\i}ne que $\mathfrak{j}'(t)_{t,X'}=\mathfrak{j}'(t)_{t}$. Donc $X'\in \mathfrak{z}_{t}$. De plus, il est clair que $X''$ appartient  au centre de $\tilde{\mathfrak{j}}$, donc \`a $\mathfrak{z}_{\tilde{t}}=\mathfrak{z}_{t}$. Mais tout \'el\'ement de $\mathfrak{z}_{t}$ agit trivialement sur $E''(t)$, donc $X''=0$. Alors $X=X'$ appartient \`a $\mathfrak{z}_{t}$. La r\'eciproque est ais\'ee. Cela prouve (2). $\square$
 
 Pour un espace vectoriel $E$ sur $F$, de dimension finie, et pour $i\in {\mathbb Z}$, on note $C_{i}(E)$ l'espace des fonctions $\varphi:E\to {\mathbb C}$ telles que
$$\varphi(\lambda e)=\vert \lambda\vert_{F} ^{i}\varphi(e)$$
pour tout $e\in E$ et tout $\lambda\in F^{\times 2}$.  On note $C_{\geq i}(E)$ l'espace des combinaisons lin\'eaires d'\'el\'ements $C_{j}(E)$ pour $j\geq i$. Remarquons que, si $E=\{0\}$, on a $C_{\geq i}(E)={\mathbb C}$ si  $i\leq 0$, $C_{\geq i}(E)=\{0\}$ si $i>0$.

 Soit $\delta:T(F)\to {\mathbb Z}$ une fonction. On note $C_{\geq \delta}(T)$ l'espace des fonctions $f$ d\'efinies presque partout sur $T(F)$ v\'erifiant la condition suivante. Soit $t\in T(F)$. Alors il existe un voisinage $\omega$ de $0$ dans $\mathfrak{t}(F)$, sur lequel l'exponentielle est d\'efinie, et il existe une fonction $\varphi\in C_{\geq \delta(t)}(\mathfrak{t}(F)/\mathfrak{z}_{t}(F))$ tels que l'on ait l'\'egalit\'e
$$f(texp(X))=\varphi(\bar{X})$$
presque partout pour $X\in \omega$, o\`u $\bar{X}$ d\'esigne la projection de $X$ dans $\mathfrak{t}(F)/\mathfrak{z}_{t}(F)$.Il revient au m\^eme de demander qu'il existe un suppl\'ementaire $\mathfrak{s}$ de $\mathfrak{z}$ dans $\mathfrak{t}$, une fonction $\varphi\in C_{\geq \delta(t)}(\mathfrak{s}(F))$ et des voisinages $\omega_{z}$ de $0$ dans $\mathfrak{z}(F)$ et $\omega_{s}$ de $0$ dans $\mathfrak{s}(F)$ de sorte que l'on ait l'\'egalit\'e
$$(3) \qquad f(texp(X_{z}+X_{s}))=\varphi(X_{s})$$
presque partout pour $X_{z}\in \omega_{z}$ et $X_{s}\in \omega_{s}$.

\ass{Lemme}{Supposons $\delta(t)=inf(dim(\mathfrak{z}_{t})-dim(\mathfrak{t})+1,0)$ pour tout $t\in T(F)$. Alors tout \'el\'ement de $C_{\geq \delta}(T)$ est localement int\'egrable sur $T(F)$.}

Preuve. Pour tout $n\in {\mathbb N}$, posons $T_{n}=\{t\in T; dim(\mathfrak{z}_{t})\geq dim(\mathfrak{t})-n\}$. Cet ensemble est un ouvert de Zariski. On va prouver par r\'ecurrence

$(4)_{n}$ tout \'el\'ement de $C_{\geq \delta}(T)$ est localement int\'egrable sur $T_{n}(F)$.

Soit $n\in {\mathbb N}$. Si $n>0$, on suppose $(4)_{n'}$ vraie pour tout $n'<n$. Soit $f\in C_{\geq\delta}(T)$. Pour prouver $(4)_{n}$, il suffit de fixer $t\in T(F)$  tel que $dim(\mathfrak{z}_{t})= dim(\mathfrak{t})-n$ et de prouver que $f$ est int\'egrable dans un voisinage de $t$. Si $n=0$, on a $\mathfrak{z}_{t}=\mathfrak{t}$ et $f$ est localement constante au voisinage de $t$. L'assertion s'ensuit. Supposons $n>0$. Fixons comme avant l'\'enonc\'e un espace $\mathfrak{s}$, une fonction $\varphi\in C_{\geq \delta(t)}(\mathfrak{s}(F))$ et des voisinages $\omega_{z}$ et $\omega_{s}$ de sorte que l'on ait l'\'egalit\'e (3). On suppose aussi  que $\omega_{z}$ et $\omega_{s}$ sont ouverts et compacts et que le voisinage $\omega=\omega_{z}\times\omega_{s}$ v\'erifie (2).  On suppose enfin que l'exponentielle de $\omega$ sur son image pr\'eserve les mesures. Ecrivons $\varphi=\sum_{i\geq \delta(t)}\varphi_{i}$, o\`u $\varphi_{i}\in C_{i}(\mathfrak{s}(F))$ et $\varphi_{i}=0$ sauf pour un nombre fini d'indices. On peut choisir une base $(e_{j})_{j=1,...,m}$ de $\mathfrak{s}(F)$ de sorte que le r\'eseau engendr\'e par cette base soit inclus dans $\omega_{s}$. On ne perd rien \`a supposer que $\omega_{s}$ est \'egal \`a ce r\'eseau.  

 Pour $i\geq \delta(t)$, d\'efinissons une fonction $f_{i}$ sur $T(F)$ de la fa\c{c}on suivante. Elle est nulle hors de $texp(\omega)$. Pour $X_{z}\in \omega_{z}$ et $X_{s}\in \omega_{s}$, on pose $f_{i}(texp(X_{z}+X_{s}))=\varphi_{i}(X_{s})$. Montrons que
 $$(5) \qquad f_{i}\in C_{\geq \delta}(T).$$
 Pour $\lambda\in F^{\times 2}$ tel que $\vert \lambda\vert _{F}\leq 1$, d\'efinissons une fonction $f[\lambda]$ sur $T(F)$ de la fa\c{c}on suivante. Elle est nulle hors de $texp(\omega)$. Pour  $X\in \omega$, on pose $f[\lambda](texp(X))=f(texp(\lambda X))$. On a $f[\lambda]=\sum_{i\geq \delta(t)}\vert \lambda\vert _{F}^{i}f_{i}$. Par interpolation, chaque $f_{i}$ est combinaison lin\'eaire de fonctions $f[\lambda]$. Il suffit donc de fixer $\lambda$ et de prouver que $f[\lambda]$ appartient \`a $C_{\geq \delta}(T)$. On fixe $X\in \omega$, on pose $t'=texp(X)$ et on doit \'etudier le comportement de la fonction $Y\mapsto f[\lambda](t'exp(Y))$ au voisinage de $0$. Posons $t''=texp(\lambda X)$ et introduisons la fonction $\varphi''\in C_{\geq\delta(t'')}(\mathfrak{t}(F)/\mathfrak{z}_{t''})$ telle que $f(t''exp(Y))=\varphi''(\bar{Y})$ pour $Y$ assez proche de $0$. On a alors
  $$(6) \qquad f[\lambda](t'exp(Y))=\varphi''(\lambda\bar{Y})=\varphi^{_{''}\lambda}(\bar{Y})$$
   pour $Y$ assez proche de $0$. La fonction $ \varphi^{_{''}\lambda}$ appartient \'evidemment \`a $C_{\geq\delta(t'')}(\mathfrak{t}(F)/\mathfrak{z}_{t''})$. De plus, la preuve de (2) montre que $\mathfrak{z}_{t'}=\mathfrak{z}_{t''}$, d'o\`u aussi $\delta(t')=\delta(t'')$. Alors l'\'egalit\'e (6) est le d\'eveloppement requis pour que $f[\lambda]$ appartienne \`a $C_{\geq \delta}(T)$. Cela prouve (5).

Notons $\Omega_{s}$ l'ensemble des \'el\'ements de $\omega_{s}$ dont les coordonn\'ees $(\lambda_{j})_{j=1,...,m}$ dans la base  $(e_{j})_{j=1,...,m}$ v\'erifient la condition $inf\{val_{F}(\lambda_{j}); j=1,...,m\}=0$ ou $1$. C'est un sous-ensemble ouvert  et compact de $\mathfrak{s}(F)$. L'ensemble $\omega$ est r\'eunion disjointe de $\omega_{z}\times\{0\}$, qui est de mesure nulle,  et des ensembles $\omega_{z}\times\varpi_{F}^{2k}\Omega_{s}$, pour $k\in {\mathbb N}$.  Soit $k\in {\mathbb N}$. Puisque $\Omega_{s}$ ne contient pas $0$, tout \'el\'ement  $X\in \omega_{z}\times\varpi_{F}^{2k}\Omega_{s}$ appartient \`a $\omega$ mais pas \`a $\mathfrak{z}_{t}(F)$. D'apr\`es (2), on a donc $\mathfrak{z}_{t}\subsetneq \mathfrak{z}_{texp(X)}$. Il en r\'esulte que $texp(X)\in \bigcup_{n'<n}T_{n'}(F)$. D'apr\`es l'hypoth\`ese de r\'ecurrence et (5), toute fonction $f_{i}$ est int\'egrable sur  $texp(\omega_{z}+\varpi_{F}^{2k}\Omega_{s})$. La fonction $f$ co\"{\i}ncide sur cet ensemble avec $\sum_{i\geq \delta(t)}f_{i}$ et est donc aussi int\'egrable. Pour prouver que $f$ est int\'egrable sur $texp(\omega)$, il reste \`a prouver que la s\'erie
$$(7) \qquad \sum_{k\in {\mathbb N}}\int_{texp(\omega_{z}+\varpi_{F}^{2k}\Omega_{s})}\vert f(t')\vert  dt'$$
est convergente. Elle est major\'ee par
$$\sum_{i\geq\delta(t)}mes(\omega_{z})\sum_{k\in {\mathbb N}}\int_{\varpi_{F}^{2k}\Omega_{s}}\vert \varphi_{i}(X)\vert dX.$$
Cette derni\`ere int\'egrale  est \'egale \`a
$$\int_{\Omega_{s}}\vert \varphi_{i}(\varpi_{F}^{2k}X)\vert q^{-2kdim(\mathfrak{s})}dX$$
ou encore
$$q^{-2k(i+dim(\mathfrak{s}))}\int_{\Omega_{s}}\vert \varphi_{i}(X)\vert dX.$$
On a 
$$i+dim(\mathfrak{s})\geq \delta(t)+dim(\mathfrak{t})-dim(\mathfrak{z}_{t})\geq 1.$$
Donc la s\'erie
$$\sum_{k\in {\mathbb N}}q^{-2k(i+dim(\mathfrak{s}))}$$
est convergente et aussi la s\'erie (7). Cela ach\`eve la d\'emonstration. $\square$

\bigskip

\subsection{Les fonctions d\'eterminants}

On d\'efinit une fonction 
$$\begin{array}{cccc}\delta_{0}:&T(F)&\to &{\mathbb Z}\\ &t&\mapsto& \delta(H_{t})-\delta(H'') +rdim(E''(t))\\ \end{array}$$
o\`u les notations sont celles introduites dans le paragraphe pr\'ec\'edent.

\ass{Lemme}{La fonction $D^{H}\Delta^r$ appartient \`a $C_{\geq\delta_{0}}(T)$.}

Preuve. Soient $t\in T(F)$ et $X\in \mathfrak{t}(F)$. On note $X'$, resp. $X''$ la restriction de $X$ \`a $E'(t)$, resp. $E''(t)$. On suppose $texp(X)\in T_{\natural}(F)$. Si $X$ est assez proche de $0$, on a les \'egalit\'es
$$D^H(texp(X))=D^H(t)D^{H_{t}}(X),\,\,\Delta(texp(X))=\vert det((1-t)_{\vert E'(t)})\vert_{F} \vert det(X''_{\vert E''(t)})\vert_{F} .$$
La fonction $D^{H_{t}}$ est invariante par translations par le centre de $\mathfrak{h}_{t}(F)$, donc aussi par $\mathfrak{z}_{t}(F)$ d'apr\`es le (1) du paragraphe pr\'ec\'edent. La seconde fonction est aussi invariante par translations par $\mathfrak{z}_{t}(F)$, ainsi que l'est toute  fonction de $X$ ne d\'ependant que de $X''$. Cette seconde fonction est homog\`ene de degr\'e $dim(E''(t))$. On a l'\'egalit\'e
$$D^{H_{t}}(X)=lim_{Y\in \mathfrak{h}_{t,X}(F),Y\to 0}D^{H_{t}}(X+Y)D^{H_{t,X}}(Y)^{-1}.$$
On a $H_{t,X}=TH''$. Sur les \'el\'ements r\'eguliers de $\mathfrak{h}_{t}(F)$, $D^{H_{t}}$ est homog\`ene de degr\'e $\delta(H_{t})$. De m\^eme, sur les \'el\'ements r\'eguliers de $\mathfrak{h}_{t,X}(F)$, $D^{H_{t,X}}$ est homog\`ene de degr\'e $\delta(H_{t,X})=\delta(H'')$. Il en r\'esulte que, sur un ouvert dense de $\mathfrak{t}(F)$, $D^{H_{t}}$ est homog\`ene de degr\'e $\delta(H_{t})-\delta(H'')$. Le r\'esultat s'ensuit. $\square$

\subsection{La fonction $c_{\tau}$}

On d\'efinit une fonction $\delta_{G,T}:T(F)\to {\mathbb Z}$ par les formules suivantes, pour $t\in T(F)$:

- si $d$ est impair et $E''(t)\not=\{0\}$, $\delta_{G,T}(t)=\frac{1}{2}(\delta(G'')-\delta(G_{t})+2)$;

- si $d$ est pair, ou si $d$ est impair et $E''(t)=\{0\}$, $\delta_{G,T}(t)=\frac{1}{2}(\delta(G'')-\delta(G_{t}))$.

\ass{Lemme}{La fonction $c_{\tau}$ appartient \`a $C_{\geq \delta_{G,T}}(T)$.}

Preuve. Fixons $t\in T(F)$ et un   bon voisinage $\omega$ de $0$ dans $\mathfrak{g}_{t}(F)$. Posons $\underline{\tau}=\tau_{t,\omega}$, cf. 4.3. C'est un quasi-caract\`ere sur $\mathfrak{g}_{t}(F)$. Pour $X\in \omega\cap \mathfrak{t}(F)$ tel que $texp(X)\in T_{\natural}(F)$, on a l'\'egalit\'e
$$c_{\tau}(texp(X))=c_{\underline{\tau},{\cal O}}(X)$$
o\`u ${\cal O}$ est une certaine orbite nilpotente r\'eguli\`ere appartenant \`a $Nil(\mathfrak{g}_{t,X})=Nil(\mathfrak{g}'')$. Il s'agit donc d'\'etudier la fonction $c_{\underline{\tau},{\cal O}}$ sur un ouvert dense de $\mathfrak{t}(F)$, au voisinage de $0$. On peut g\'en\'eraliser la question \`a un quasi-caract\`ere quelconque $\underline{\tau}$ de $\mathfrak{g}_{t}(F)$. On peut restreindre $\omega$ et supposer que $\underline{\tau}$ est d\'eveloppable dans $\omega$. En utilisant le lemme 6.3(iii), on peut fixer une orbite ${\cal O}_{t}\in Nil(\mathfrak{g}_{t})$ et supposer que $\underline{\tau}(Y)=\hat{j}({\cal O}_{t},Y)$ presque partout pour $Y\in \omega$. Cette fonction est localement invariante par translations par le centre de $\mathfrak{g}_{t}(F)$, donc aussi par $\mathfrak{z}_{t}(F)$ d'apr\`es 7.4(1). La fonction $c_{\underline{\tau},{\cal O}}$ sur $\mathfrak{t}(F)$ l'est donc aussi. Soit $\lambda\in F^{\times 2}$. D'apr\`es 4.2(2), le quasi-caract\`ere $\underline{\tau}^{\lambda}$ co\"{\i}ncide avec $\vert \lambda\vert _{F}^{-dim({\cal O}_{t})/2}\underline{\tau}$ au voisinage de $0$. La m\^eme formule nous dit alors que la fonction $c_{\underline{\tau},{\cal O}}$ sur $\mathfrak{t}(F)$ est homog\`ene de degr\'e $(dim({\cal O})-dim({\cal O}_{t}))/2$. On a $dim({\cal O}_{t})\leq \delta(G_{t})$. Puisque ${\cal O}$ est r\'eguli\`ere, on a $dim({\cal O})=\delta(G'')$. Si $d$ est pair, ou si $d$ est impair et $E''(t)=\{0\}$, le degr\'e pr\'ec\'edent est sup\'erieur ou \'egal \`a $\delta_{G,T}(t)$ et cela ach\`eve la d\'emonstration. Supposons $d$ impair et $E''(t)\not=\{0\}$. Si ${\cal O}_{t}$ n'est pas r\'eguli\`ere, on a $dim({\cal O}_{t})\leq \delta(G_{t})-2$ et on conclut encore. Supposons ${\cal O}_{t}$ r\'eguli\`ere. Le tore $T$ est un sous-tore maximal de $H'_{t}=J'(t)_{t}J''(t)$ et se d\'ecompose donc en $T=T'T''$, o\`u $T'$ est un sous-tore maximal de $J'(t)_{t}$ et $T''$ un sous-tore maximal de $J''(t)$. L'hypoth\`ese $A_{T}=\{1\}$ entra\^{\i}ne $A_{T''}=\{1\}$ et l'hypoth\`ese $E''(t)\not=\{0\}$ entra\^{\i}ne $T''\not=\{1\}$. Notons $\tilde{G}$ le groupe sp\'ecial orthogonal de $E''(t)\oplus V''$. On a $T''\subset J''(t)\subset \tilde{G}$. Comme on l'a vu dans la preuve de 7.4(1), on a $G_{t}=J'(t)_{t}\tilde{G}$.  L'orbite ${\cal O}_{t}$ se d\'ecompose en la somme d'une orbite nilpotente dans $\mathfrak{j}'(t)_{t}(F)$ et d'une orbite nilpotente $\tilde{{\cal O}}$ dans $\tilde{\mathfrak{g}}(F)$. Ces deux orbites sont r\'eguli\`eres. Cela entra\^{\i}ne que $\tilde{G}$ est quasi-d\'eploy\'e. Or $dim(E''(t)\oplus V'')$ est  impair, donc $\tilde{G}$ est d\'eploy\'e. Donc $\tilde{\mathfrak{g}}(F)$ poss\`ede une unique orbite nilpotente r\'eguli\`ere, \`a savoir $\tilde{{\cal O}}$, qui est induite \`a partir de l'orbite $\{0\}$ d'une sous-alg\`ebre de L\'evi minimale, c'est-\`a-dire de l'alg\`ebre de Lie d'un tore d\'eploy\'e maximal. Cela entra\^{\i}ne que la fonction $\hat{j}(\tilde{{\cal O}},.)$ est \`a support dans l'ensemble des \'el\'ements qui appartiennent \`a une sous-alg\`ebre de Borel. Les propri\'et\'es de $T''$ montrent qu'un \'el\'ement de $\mathfrak{t}''(F)$ en position g\'en\'erale poss\`ede un voisinage dans $\tilde{\mathfrak{g}}(F)$ dont aucun \'el\'ement n'appartient \`a une telle alg\`ebre. Donc $\hat{j}(\tilde{{\cal O}},.) $ s'annule au voisinage de presque tout \'el\'ement de $\mathfrak{t}''(F)$. Il en r\'esulte que $\underline{\tau}$ s'annule au voisinage de presque tout \'el\'ement de $\mathfrak{t}(F)$. A fortiori, la fonction $c_{\underline{\tau},{\cal O}}$ est nulle sur $\mathfrak{t}(F)$. Cela ach\`eve la d\'emonstration. $\square$

\subsection{Preuve de la proposition 7.3}

Evidemment, un lemme analogue au lemme 7.6 vaut si l'on remplace $G$ par $H$ et $\delta_{G,T}$ par une fonction $\delta_{H,T}$ d\'efinie de fa\c{c}on similaire (l'entier $d_{W}$ remplace $d$). Il r\'esulte tout d'abord de ces lemmes que les fonctions $c_{\theta}$ et $c_{\tau}$ sont localement constantes sur $T_{\natural}(F)$. Evidemment, si $\delta_{i}$, $i=1,2,3$ sont trois fonctions sur $T(F)$ telles que $\delta_{1}+\delta_{2}\geq \delta_{3}$, on a $f_{1}f_{2}\in C_{\geq\delta_{3}}(T)$ pour toutes $f_{1}\in C_{\geq\delta_{1}}(T)$ et $f_{2}\in C_{\geq\delta_{2}}(T)$. En vertu des lemmes 7.4, 7.5 et 7.6, pour d\'emontrer le (ii) de la proposition, il suffit de prouver que
$$(1) \qquad \delta_{0}(t)+\delta_{G,T}(t)+\delta_{H,T}(t)\geq inf(dim(\mathfrak{z}_{t})-dim(\mathfrak{t})+1,0)$$
pour tout $t\in T(F)$. Posons $e=dim(E''(t))$. Puisque $d$ ou $d_{W}$ est impair, les d\'efinitions entra\^{\i}nent que le membre de gauche est \'egal \`a

 $ \frac{1}{2}(\delta(G'')-\delta(G_{t})-\delta(H'')+\delta(H_{t}))+re$, si $e=0$;

$ \frac{1}{2}(\delta(G'')-\delta(G_{t})-\delta(H'')+\delta(H_{t}))+1+re$, si $e>0$.

On a d\'ej\`a dit que $G_{t}$ \'etait le produit de $J'(t)_{t}$ et du groupe sp\'ecial orthogonal $\tilde{G}$ de $E''(t)\oplus V''$. De m\^eme, $H_{t}$ est le produit de $J'(t)_{t}$ et du groupe sp\'ecial orthogonal $\tilde{H}$ de $E''(t)\oplus W''$. On peut remplacer $-\delta(G_{t})+\delta(H_{t})$ par $-\delta(\tilde{G})+\delta(\tilde{H})$ dans les formules pr\'ec\'edentes. Il est facile de calculer
$$(2) \qquad \delta(G)=\left\lbrace\begin{array}{cc}d(d-2)/2,&\,\,{\rm si}\,\,d\,\,{\rm est\,\,pair},\\ (d-1)^2/2,&\,\,{\rm si}\,\,d\,\,{\rm est\,\,impair}.\\ \end{array}\right.$$
On calcule de m\^eme $\delta(G'')$, $\delta(H'')$, $\delta(\tilde{G})$ et $\delta(\tilde{H})$ en rempla\c{c}ant $d$ par $d''+1+2r$, $d''$, $d''+1+2r+e$, $d''+e$, o\`u $d''=dim(W'')$. On obtient que le membre de gauche de (1) est sup\'erieur ou \'egal \`a 
$$\left\lbrace\begin{array}{c}0,\text{\,\, si\,\,} e=0;\\ -e/2+1,\text{\,\,si\,\,} e>0.\\ \end{array}\right.$$
  Dans le premier cas, il est clairement sup\'erieur au membre de droite de (1). Supposons $e>0$. L'alg\`ebre $\mathfrak{z}_{t}$ est celle d'un sous-tore de $J'(t)$, donc de dimension inf\'erieure ou \'egale \`a $dim(E'(t))/2$, qui est \'egale \`a $dim(\mathfrak{t})-e/2$. Le membre de droite de (1) est donc inf\'erieur ou \'egal \`a $-e/2+1$, donc au membre de gauche. Cela ach\`eve la preuve. $\square$
  
  \bigskip
  
  \subsection{Le th\'eor\`eme}
  
  Soient $\theta$ un quasi-caract\`ere sur $H(F)$ et $f\in C_{c}^{\infty}(G(F))$ une fonction tr\`es cuspidale. Pour $T\in \underline{\cal T}$, on d\'efinit la fonction $c_{\theta_{f}}$ sur $T_{\natural}(F)$ et on la note simplement $c_{f}$. Fixons un ensemble de repr\'esentants ${\cal T}$ des classes de conjugaison par $H(F)$ dans $\underline{\cal T}$. Posons
$$I(\theta,f)=\sum_{T\in {\cal T}}\vert W(H,T)\vert ^{-1}\nu(T)\int_{T(F)}c_{\theta}(t)c_{f}(t)D^H(t)\Delta(t)^rdt.$$
D'apr\`es la proposition 7.3, cette expression est absolument convergente.

Pour $g\in G(F)$, on d\'efinit une fonction $^gf^{\xi}$ sur $H(F)$ par
$$^gf^{\xi}(x)=\int_{U(F)}f(g^{-1}xug)\xi(u)du.$$
Elle appartient \`a $C_{c}^{\infty}(H(F))$. On pose 
$$I(\theta,f,g)=\int_{H(F)}\theta(x){^gf}^{\xi}(x)dx,$$
puis, pour un entier $N\in {\mathbb N}$,
$$I_{N}(\theta,f)=\int_{U(F)H(F)\backslash G(F)}I(\theta,f,g)\kappa_{N}(g) dg.$$
Ces int\'egrales sont \`a supports compacts.

\ass{Th\'eor\`eme}{Pour tout quasi-caract\`ere $\theta$ sur $H(F)$ et toute fonction tr\`es cuspidale $f\in C_{c}^{\infty}(G(F))$, on a l'\'egalit\'e
$$lim_{N\to \infty}I_{N}(\theta,f)=I(\theta,f).$$}

Ce th\'eor\`eme sera d\'emontr\'e en 12.3. Contentons-nous ici de la remarque facile suivante. Supposons $d_{W}\geq 1$. Soit $y\in H^+(F)$ $y\not\in H(F)$, que l'on identifie comme on l'a dit en 7.2 \`a un \'el\'ement de $G(F)$. Posons $\theta^+=(\theta+{^y\theta})/2$. Par de simples changements de variables, on v\'erifie les \'egalit\'es
$$I(\theta^+,f)=I(\theta,f),\,\,I_{N}(\theta^+,f)=I_{N}(\theta,f).$$
On peut donc remplacer $\theta$ par $\theta^+$ pour d\'emontrer le th\'eor\`eme. Autrement dit, on peut supposer  $\theta$ invariant par conjugaison par $H^+(F)$.  

\bigskip

\subsection{Le th\'eor\`eme pour les alg\`ebres de Lie}

Soient $\theta$ un quasi-caract\`ere sur $\mathfrak{h}(F)$ et $f\in C_{c}^{\infty}(\mathfrak{g}(F))$. Les d\'efinitions pos\'ees pour les groupes dans les paragraphes pr\'ec\'edents se descendent aux alg\`ebres de Lie. Ainsi, on a d\'efini en 7.2 un caract\`ere $\xi$ de $U(F)$. Il s'en d\'eduit un caract\`ere de $\mathfrak{u}(F)$, d\'efini par la m\^eme formule qu'en 7.2 et que l'on note encore $\xi$. On d\'efinit une fonction $f^{\xi}$ sur $\mathfrak{h}(F)$ par
$$f^{\xi}(Y)=\int_{\mathfrak{u}(F)}f(Y+N)\xi(N)dN.$$
Pour $g\in G(F)$, on pose
$$I(\theta,f,g)=\int_{\mathfrak{h}(F)}\theta(Y){^gf}^{\xi}(Y)dY,$$
puis, pour un entier $N\in {\mathbb N}$,
$$I_{N}(\theta,f)=\int_{U(F)H(F)\backslash G(F)}I(\theta,f,g)\kappa_{N}(g)dg.$$
On d\'efinit la fonction $\Delta$ sur $\mathfrak{h}(F)$ par
$$\Delta(Y)=\vert det(Y\vert W/W''(Y))\vert _{F},$$
o\`u $W''(Y) $ est le noyau de $Y$ agissant dans $W$. Pour $T\in {\cal T}$, on note $\mathfrak{t}_{\natural}$ le sous-ensemble des $X\in \mathfrak{t}$ tels que les valeurs propres de l'action de $X$ dans $W'$ soient toutes distinctes, o\`u $W'$ est comme en 7.3. Supposons $f$ tr\`es cuspidale. On d\'efinit les fonctions $c_{\theta}$ et $c_{f}=c_{\theta_{f}}$ sur $\mathfrak{t}_{\natural}(F)$. On pose
$$I(\theta,f)=\sum_{T\in {\cal T}}\vert W(H,T)\vert ^{-1}\nu(T)\int_{\mathfrak{t}(F)}c_{\theta}(Y)c_{f}(Y)D^{H}(Y)\Delta(Y)^rdY.$$
Une analogue de la proposition 7.3 entra\^{\i}ne l'absolue convergence de cette expression.

\ass{Th\'eor\`eme}{Pour tout quasi-caract\`ere $\theta$ sur $\mathfrak{h}(F)$ et toute fonction tr\`es cuspidale $f\in C_{c}^{\infty}(\mathfrak{g}(F))$, on a l'\'egalit\'e
$$lim_{N\to \infty}I_{N}(\theta,f)=I(\theta,f).$$}

Ce th\'eor\`eme sera d\'emontr\'e en 12.3.

\bigskip

\section{Localisation}

\bigskip

\subsection{Un cas trivial}

On fixe pour toute la section un quasi-caract\`ere $\theta$ sur $H(F)$, invariant par conjugaison par $H^+(F)$, et une fonction tr\`es cuspidale $f\in C_{c}^{\infty}(G(F))$. Soit $x\in G_{ss}(F)$. Notons $V''$ le noyau de $x-1$ agissant dans $V$. Supposons que $x$ n'est conjugu\'e \`a aucun \'el\'ement de $H(F)$. Par le th\'eor\`eme de Witt, cette hypoth\`ese \'equivaut \`a dire que $V''$ ne contient aucun sous-espace non d\'eg\'en\'er\'e isomorphe (comme espace quadratique) \`a $D\oplus Z$. Soit $\omega$ un bon voisinage de $0$ dans $\mathfrak{g}_{x}(F)$, v\'erifiant la condition $(7)_{\rho}$ de 3.1, o\`u $\rho$ est la repr\'esentation de $G$ dans $V$. Pour $X\in \omega$, le noyau de $xexp(X)-1$ est contenu dans $V''$ et v\'erifie a fortiori la m\^eme condition que $V''$. Donc $xexp(X)$ n'est conjugu\'e \`a aucun \'el\'ement de $H(F)$. Posons $\Omega=(xexp(\omega))^G$. Alors $\Omega\cap H(F)=\emptyset$. Supposons $f$ \`a support dans $\Omega$. Pour tout $t\in H_{ss}(F)$, le compl\'ementaire de $\Omega$ dans $G(F)$ est un voisinage de $t$ invariant par conjugaison par $G(F)$ et sur lequel $f$ est nulle. Donc $\theta_{f}$ y est nul aussi et le d\'eveloppement de $\theta_{f}$ au voisinage de $t$ est nul. Il en r\'esulte que $I(\theta,f)=0$. D'autre part, tout \'el\'ement de $U(F)H(F)$ a pour partie semi-simple un \'el\'ement conjugu\'e \`a un \'el\'ement de $H(F)$. Il en r\'esulte que $^gf^{\xi}=0$ pour tout $g\in G(F)$, donc $I_{N}(\theta,f)=0$. Alors l'\'egalit\'e du th\'eor\`eme est triviale.

\subsection{Localisation de $I_{N}(\theta,f)$}

Soit $x\in H_{ss}(F)$. On note $W''$, resp. $V_{0}''$, $V''$, le noyau de $x-1$ agissant dans $W$, resp. $V_{0}$, resp. $V$. On a $V''_{0}=W''\oplus D$, $V''=W''\oplus D\oplus Z$. On note $W'$ l'orthogonal de $W''$ dans $W$. On note $H'=G'$, resp. $H''$, $G''_{0}$, $G''$, les groupes sp\'eciaux orthogonaux de $W'$, resp. $W''$, $V''_{0}$, $V''$. On a les \'egalit\'es $H_{x}=H'_{x}H''$, $G_{x}=G'_{x}G''$.  On fixe un bon voisinage $\omega$ de $0$ dans $\mathfrak{g}_{x}(F)$, auquel on impose la condition (8) de 3.1, c'est-\`a-dire $\omega=\omega'\times\omega''$, o\`u $\omega'\subset \mathfrak{g}'_{x}(F)$, $\omega''\subset \mathfrak{g}''(F)$. On pose $\Omega=(xexp(\omega))^G$. On suppose

\ass{Hypoth\`ese} {Le support de $f$ est contenu dans $\Omega$.}

La situation ci-dessus, les notations et cette hypoth\`ese seront conserv\'ees jusqu'en 10.9.

On d\'efinit le quasi-caract\`ere $\theta_{x,\omega}$ de $\mathfrak{g}_{x}(F)$, cf. 4.3, et, pour $g\in G(F)$, la fonction $^gf_{x,\omega}$ sur $\mathfrak{g}_{x}(F)$, cf. 5.4.   Pour $g\in G(F)$, on d\'efinit une fonction $^gf_{x,\omega}^{\xi}$ sur $\mathfrak{h}_{x}(F)$ par
$$^gf_{x,\omega}^{\xi}(X)=\int_{\mathfrak{u}_{x}(F)}{^gf}_{x,\omega}(X+N)\xi(N)dN.$$
Remarquons que $x$ appartient \`a $M(F)$ et que l'on a l'inclusion $\mathfrak{h}_{x}\subset \mathfrak{m}_{x}$. Posons
$$I_{x,\omega}(\theta,f,g)=\int_{\mathfrak{h}_{x}(F)}\theta_{x,\omega}(X){^gf}_{x,\omega}^{\xi}(X)dX,$$
puis
$$I_{x,\omega,N}(\theta,f)=\int_{U_{x}(F)H_{x}(F)\backslash G(F)}I_{x,\omega}(\theta,f,g)\kappa_{N}(g)dg.$$
Cette int\'egrale a un sens:  la fonction $g\mapsto I_{x,\omega}(\theta,f,g)$ est invariante \`a gauche par $U_{x}(F)H_{x}(F)$.   Elle est \`a support compact. En effet, d'apr\`es 3.1(5), il existe un sous-ensemble compact $\Gamma\subset G(F)$ tel que $^gf_{x,\omega}$ est nulle pour $g\in G(F)$,$g\not\in G_{x}(F)\Gamma$. D'autre part, on v\'erifie que, pour tout $\gamma\in G(F)$, la fonction $g\mapsto \kappa_{N}(g\gamma)$ sur $G_{x}(F)$ a un support d'image compacte dans $U_{x}(F)H_{x}(F)\backslash G_{x}(F)$. L'assertion en r\'esulte.

Posons
$$C(x)=\vert H^+(F)/H(F)\vert  \vert Z_{H^+}(x)(F)/H_{x}(F)\vert ^{-1}\Delta(x)^r.$$

\ass{Lemme}{On a l'\'egalit\'e
$$I_{N}(\theta,f)=C(x)I_{x,\omega,N}(\theta,f).$$}

Preuve. Pour tout groupe r\'eductif connexe $L$, fixons un ensemble de repr\'esentants ${\cal T}(L)$ des classes de conjugaison par $L(F)$ dans l'ensemble des sous-tores maximaux de $L$. Soit $g\in G(F)$. D'apr\`es la formule de Weyl, on a
$$(1) \qquad I(\theta,f,g)=\sum_{T\in {\cal T}(H)}\vert W(H,T)\vert ^{-1}\int_{T(F)}\theta(t)J_{H}(t,{^gf}^{\xi})D^H(t)^{1/2}dt.$$
Pour deux sous-tores (pas forc\'ement maximaux) $T$ et $T'$ de $H$, notons $W^+(T,T')$ l'ensemble des isomorphismes de $T$ sur $T'$ induits par la conjugaison par un \'el\'ement de $H^+(F)$.  On va prouver les assertions suivantes.

(2) Soient $T\in {\cal T}(H)$ et $t\in T(F)\cap H_{reg}(F)$. Alors $J_{H}(t,{^gf}^{\xi})=0$ si $t$ n'appartient pas \`a 
$$\bigcup_{T_{1}\in {\cal T}(H_{x})}\bigcup_{w\in W^+(T_{1},T)}w(xexp(\mathfrak{t}_{1}(F)\cap \omega)).$$

(3) Soit $T\in {\cal T}(H)$ et, pour $i=1,2$, soient $T_{i}\in {\cal T}(H_{x})$ et $w_{i}\in W^+(T_{i},T)$. Alors les ensembles $w_{1}(xexp(\mathfrak{t}_{1}(F)\cap \omega))$ et $w_{2}(xexp(\mathfrak{t}_{2}(F)\cap \omega))$ sont disjoints ou confondus.

(4) Soient $T\in {\cal T}(H)$, $T_{1}\in {\cal T}(H_{x})$ et $w_{1}\in W^+(T_{1},T)$. Le nombre des couples $(T_{2},w_{2})$ tels que $T_{2}\in {\cal T}(H_{x})$, $w_{2}\in W^+(T_{2},T)$ et $w_{2}(xexp(\mathfrak{t}_{2}(F)\cap \omega))=w_{1}(xexp(\mathfrak{t}_{1}(F)\cap \omega))$ est \'egal \`a
$$\vert W(H_{x},T_{1})\vert \vert Z_{H^+}(x)(F)/H_{x}(F)\vert \vert Z_{H^+}(T_{1})(F)/T_{1}(F)\vert ^{-1}.$$

Soient $T$ et $t$ comme en (2). Supposons $J_{H}(t,{^gf}^{\xi})\not=0$. Alors il existe $u\in U(F)$ tel que la classe de conjugaison par $G(F)$ de $tu$ coupe le support de $f$. Elle coupe donc aussi $xexp(\omega)$. La partie semi-simple de $tu$ est conjugu\'ee \`a $t$ et la partie semi-simple d'un \'el\'ement de $xexp(\omega)$ reste dans cet ensemble. Donc la classe de conjugaison par $G(F)$ de $t$ coupe $xexp(\omega)$. Soient $X\in \omega$ et $y\in G(F)$ tels que $yty^{-1}=xexp(X)$. Le noyau de $t-1$ agissant dans $V$ contient $D\oplus Z$. D'apr\`es l'hypoth\`ese $(7)_{\rho}$ de 3.1, celui de $xexp(X)-1$ est contenu dans $W''$. Donc $W''$ contient $y(D\oplus Z)$. Mais il contient aussi  $D\oplus Z$. Ces deux espaces $D\oplus Z$ et $y(D\oplus Z)$ sont isomorphes et non d\'eg\'en\'er\'es, en tant qu'espaces  quadratiques. D'apr\`es le th\'eor\`eme de Witt, on peut trouver $y''\in G''(F)$ tel que $y''y(D\oplus Z)=D\oplus Z$. On a $G''\subset G_{x}$. Quitte \`a remplacer $y$ par $y''y$ et $X$ par $y''Xy^{_{''}-1}$, on est ramen\'e au cas o\`u $y$ conserve $D\oplus Z$. Dans ce cas, puisque $t$ agit trivialement sur $D\oplus Z$, $xexp(X)$ agit trivialement lui aussi, donc $X\in \mathfrak{h}_{x}(F)$. Quitte \`a multiplier encore $y$ \`a gauche par un \'el\'ement de $H_{x}(F)$, on peut supposer que $X\in \mathfrak{t}_{1}(F)$ pour un \'el\'ement $T_{1}\in {\cal T}(H_{x})$. L'\'el\'ement $y$ conserve $W$. Notons $h$ sa restriction \`a cet espace. Alors $h\in H^+(F)$ et $hth^{-1}=xexp(X)$.  N\'ecessairement, la conjugaison par $h$ envoie le commutant de $t$ dans $H$ sur celui de $xexp(X)$. Mais $t$ est r\'egulier dans $H$, donc ces commutants sont $T$ et $T_{1}$. Si on note $w$ l'\'el\'ement de $W^+(T_{1},T)$ induit par la conjugaison par $y^{-1}$, on a alors $t\in w(xexp(\mathfrak{t}_{1}(F)\cap \omega))$ ce qui prouve (2).

Passons \`a la preuve de (3). Pour $i=1,2$ soit $y_{i}\in H^+(F)$ tel que $w_{i}$ soit induit par la conjugaison par $y_{i}$. On identifie $y_{i}$ \`a un \'el\'ement de $G(F)$. Posons $y=y_{2}^{-1}y_{1}$. Supposons que les ensembles $w_{1}(xexp(\mathfrak{t}_{1}(F)\cap \omega))$ et $w_{2}(xexp(\mathfrak{t}_{2}(F)\cap \omega))$ ne sont pas disjoints. Alors $y(xexp(\omega))y^{-1}\cap (xexp(\omega))\not=\emptyset$. D'apr\`es 3.1(4), $y$ appartient \`a $Z_{G}(x)(F)$. D'apr\`es 3.1(1), la conjugaison par $y$ conserve $\omega$. D'autre part, d'apr\`es la d\'efinition de $y$, cette conjugaison envoie $T_{1}$ sur $T_{2}$, donc aussi $\mathfrak{t}_{1}$ sur $\mathfrak{t}_{2}$. Elle envoie alors $xexp(\mathfrak{t}_{1}(F)\cap \omega)$ sur $xexp(\mathfrak{t}_{2}(F)\cap \omega)$ et les ensembles $w_{1}(xexp(\mathfrak{t}_{1}(F)\cap \omega))$ et $w_{2}(xexp(\mathfrak{t}_{2}(F)\cap \omega))$ sont confondus. Cela prouve (3).

Soient $T$, $T_{1}$ et $w_{1}$ comme en (4). Posons
$${\cal Y}=\{y\in Z_{H^+}(x)(F); yT_{1}y^{-1}\in {\cal T}(H_{x})\}/Z_{H^+}(T_{1})(F).$$
La preuve de (3) montre que l'application $y\mapsto (T_{2}=yT_{1}y^{-1}, w_{2}=w_{1}ad(y^{-1}))$ est une surjection de ${\cal Y}$ sur l'ensemble des couples $(T_{2},w_{2})$ dont on veut calculer le nombre d'\'el\'ements. Cette application est aussi injective, le nombre \`a calculer est donc $\vert {\cal Y}\vert $. On v\'erifie que l'application naturelle
$${\cal Y}\to H_{x}(F)\backslash Z_{H^+}(x)(F)/Z_{H^+}(T_{1})(F)$$
est surjective et que toutes ses fibres ont pour nombre d'\'el\'ements $\vert W(H_{x},T_{1})\vert $. Enfin, parce que $H_{x}(F)$ est un sous-groupe distingu\'e de $Z_{H^+}(x)(F)$ et que $Z_{H^+}(T_{1})\cap H_{x}=T_{1}$,  on a
$$\vert H_{x}(F)\backslash Z_{H^+}(x)(F)/Z_{H^+}(T_{1})(F)\vert =\vert Z_{H^+}(x)(F)/H_{x}(F)\vert \vert Z_{H^+}(T_{1})(F)/T_{1}(F)\vert ^{-1}.$$
Cela prouve (4).

Ces trois propri\'et\'es permettent de transformer l'expression (1) de la fa\c{c}on suivante
$$I(\theta,f,g)=\sum_{T_{1}\in {\cal T}(H_{x})}\sum_{T\in {\cal T}(H)}\sum_{w_{1}\in W^+(T_{1},T)}\vert W(H,T)\vert^{-1}w(T_{1})$$
$$\qquad \int_{\mathfrak{t}_{1}(F)\cap \omega} \theta(w_{1}(xexp(X)))J_{H}(w_{1}(xexp(X)),{^gf}^{\xi})D^H(w_{1}(xexp(X)))^{1/2}dX,$$
o\`u 
$$w(T_{1})=  \vert W(H_{x},T_{1})\vert^{-1} \vert Z_{H^+}(x)(F)/H_{x}(F)\vert^{-1} \vert Z_{H^+}(T_{1})(F)/T_{1}(F)\vert .$$
On a $D^H(w_{1}(xexp(X)))=D^H(xexp(X))$. On a $\theta(w_{1}(xexp(X)))=\theta(xexp(X))$ puisqu'on a suppos\'e $\theta$ invariant par $H^+(F)$.  Si $w_{1}$ \'etait induit par la conjugaison par un \'el\'ement de $H(F)$, on aurait aussi $J_{H}(w_{1}(xexp(X)),{^{g}f}^{\xi})=J_{H}(xexp(X),{^{g}f}^{\xi})$, et $w_{1}$ dispara\^{\i}trait de la formule ci-dessus. En g\'en\'eral, on a seulement $J_{H}(w_{1}(xexp(X)),{^{g}f}^{\xi})=J_{H}(xexp(X),{^{yg}f}^{\xi})$, o\`u $y\in H^+(F)$ d\'epend de $w_{1}$. Mais ce terme $y$ dispara\^{\i}t par changement de variables quand on calcule $I_{N}(\theta,f)$. Ces arguments conduisent \`a l'\'egalit\'e
$$(5) \qquad I_{N}(\theta,f)=\int_{U(F)H(F)\backslash G(F)}\sum_{T_{1}\in {\cal T}(H_{x})}w'(T_{1})$$
$$\qquad \int_{\mathfrak{t}_{1}(F)\cap \omega}\theta(xexp(X))J_{H}(xexp(X),{^gf}^{\xi} )D^H(xexp(X))^{1/2}dX\kappa_{N}(g)dg,$$
o\`u
$$w'(T_{1})=w(T_{1})\sum_{T\in {\cal T}(H)} \vert W^+(T_{1},T)\vert \vert W(H,T)\vert^{-1}.$$
Soit $T_{1}\in {\cal T}(H_{x})$. Remarquons que $W(H,T)$ a m\^eme nombre d'\'el\'ements que $W(H,T_{1})$ pour tout $T$ tel que $W^+(T_{1},T)$ est non vide. On a donc
$$w'(T_{1})=w(T_{1})\vert W(H,T_{1})^{-1}\vert \vert {\cal Y}_{T_{1}}\vert ,$$
o\`u
$${\cal Y}_{T_{1}}=\{(T,w_{1}); T\in {\cal T}(H), w_{1}\in W^+(T_{1},T)\}.$$
Posons 
$${\cal Y}'_{T_{1}}=\{y\in H^+(F)/Z_{H^+}(T_{1})(F); yT_{1}y^{-1}\in {\cal T}(H)\}.$$
L'application $y\mapsto (T=yT_{1}y^{-1}, w_{1}=ad(y))$ est une bijection de ${\cal Y}'_{T_{1}}$ sur ${\cal Y}_{T_{1}}$. L'application naturelle
$${\cal Y}'_{T_{1}}\to H(F)\backslash H^+(F)/Z_{H^+}(T_{1})(F)$$
est surjective et toutes ses fibres ont pour nombre d'\'el\'ements $\vert W(H,T_{1})\vert $. Enfin, parce que $H$ est un sous-groupe distingu\'e de $H^+$ et $Z_{H^+}(T_{1})\cap H=T_{1}$, on a l'\'egalit\'e
$$\vert H(F)\backslash H^+(F)/Z_{H^+}(T_{1})(F)\vert =\vert H^+(F)/H(F)\vert \vert Z_{H^+}(T_{1})(F)/T_{1}(F)\vert ^{-1}.$$
Cela conduit \`a l'\'egalit\'e
$$w'(T_{1})=\vert H^+(F)/H(F)\vert\vert Z_{H^+}(x)(F)/H_{x}(F)\vert^{-1}\vert W(H_{x},T_{1})\vert ^{-1}.$$
Pour $X\in \omega\cap\mathfrak{h}_{x,reg}(F)$ et $g\in G(F)$, on a 
$$J_{H}(xexp(X),{^gf}^{\xi})=D^H(xexp(X))^{1/2}\int_{H_{x}(F)\backslash H(F)}\int_{T_{1}(F)\backslash H_{x}(F)}{^{yg}f}^{\xi}(xexp(h^{-1}Xh))dh\,dy.$$
D'autre part, on a $D^H(xexp(X))=D^H(x)D^{H_{x}}(X)$. Ces \'egalit\'es transforment la formule (5) en
$$(6) \qquad I_{N}(\theta,f)=C'(x)\int_{U(F)H_{x}(F)\backslash G(F)}\Phi(g)\kappa_{N}(g)dg,$$
o\`u 
$$C'(x)=\vert H^+(F)/H(F)\vert\vert Z_{H^+}(x)(F)/H_{x}(F)\vert^{-1} D^H(x)$$
et
$$\Phi(g)=\sum_{T_{1}\in {\cal T}(H_{x})}\vert W(H_{x},T_{1})^{-1}\int_{\mathfrak{t}_{1}(F)\cap \omega}\theta(xexp(X))\int_{T_{1}(F)\backslash H_{x}(F)} {^gf}^{\xi}(xexp(h^{-1}Xh))dh\,D^{H_{x}}(X)dX.$$
D\'efinissons une fonction $\varphi_{g}$ sur $\mathfrak{h}_{x}(F)$ par
$$\varphi_{g}(X)=\left\lbrace\begin{array}{cc}0,&\,\,{\rm si}\,\,X\not\in \omega,\\ \theta(xexp(X)){^gf}^{\xi}(xexp(X)),&\,\,{\rm si}\,\, X\in \omega.\\ \end{array}\right.$$
D'apr\`es la formule de Weyl,
$$\Phi(g)=\int_{\mathfrak{h}_{x}(F)}\varphi_{g}(X) dX.$$
Soient $X\in \omega\cap \mathfrak{h}_{x,reg}(F)$ et $g\in G(F)$. On a
 $${^gf}^{\xi}(xexp(X))= \int_{U(F)}{^gf}(xexp(X)u)\xi(u)du$$
 $$\qquad = \int_{U_{x}(F)\backslash U(F)}\int_{U_{x}(F)}{^gf}(xexp(X)uv)\xi(uv)du \,dv.$$
 Pour $u\in U_{x}(F)$, l'application $v\mapsto (xexp(X)u)^{-1}v^{-1}xexp(X)uv$ est une bijection de $U_{x}(F)\backslash U(F)$ sur lui-m\^eme. Gr\^ace \`a l'hypoth\`ese $(7)_{\rho}$ de 3.1, son jacobien est \'egal \`a la valeur absolue du d\'eterminant de $1-ad(x)^{-1}$ agissant sur $\mathfrak{u}(F)/\mathfrak{u}_{x}(F)$. Remarquons que, avec les notations de 7.1 et 7.2, l'application
 $$\begin{array}{ccc}W'\otimes Z_{+}&\to&\mathfrak{u}(F)\\ (w',z)&\mapsto&c_{w',z}\\ \end{array}$$ 
est une bijection de $W'\otimes Z_{+}$ sur un suppl\'ementaire de $\mathfrak{u}_{x}(F)$ dans $\mathfrak{u}(F)$. Le jacobien ci-dessus est donc \'egal \`a $\Delta(x)^r$. D'autre part, on a
$$\xi( (xexp(X)u)^{-1}v^{-1}xexp(X)uv)=1.$$
Cela conduit \`a l'\'egalit\'e
$${^gf}^{\xi}(xexp(X))=\Delta(x)^r \int_{U_{x}(F)\backslash U(F)}\int_{U_{x}(F)}{^gf}(v^{-1}xexp(X)uv)\xi(u)du\,dv$$
$$\qquad =\Delta(x)^r \int_{U_{x}(F)\backslash U(F)}\int_{U_{x}(F)}{^{vg}f}(xexp(X)u)\xi(u)du\,dv.$$
Gr\^ace \`a la condition (6) de 3.1, l'application 
$$\begin{array}{ccc}\mathfrak{u}_{x}(F)&\to&U_{x}(F)\\N&\mapsto& exp(-X)exp(X+N)\\ \end{array}$$
est bijective et pr\'eserve les mesures. On a $\xi(exp(-X)exp(X+N))=\xi(N)$. On a donc aussi
$$ {^gf}^{\xi}(xexp(X))=\Delta(x)^r \int_{U_{x}(F)\backslash U(F)}\int_{\mathfrak{u}_{x}(F)}{^{vg}f}(xexp(X+N))\xi(N)dN\,dv.$$
Remarquons que la partie semi-simple de $X+N$ est conjugu\'ee \`a $X$ par un \'el\'ement de $G_{x}(F)$, donc $X+N\in \omega$ et ${^{vg}f}(xexp(X+N))={^{vg}f}_{x,\omega}(X+N)$. Alors
$${^gf}^{\xi}(xexp(X))=\Delta(x)^r \int_{U_{x}(F)\backslash U(F)}{^{vg}f}_{x,\omega}^{\xi}(X)dv.$$
Par ailleurs, on a $\theta(xexp(X))=\theta_{x,\omega}(X)$. Donc
$$\varphi_{g}(X)=\Delta(x)^r\theta_{x,\omega}(X)\int_{U_{x}(F)\backslash U(F)}{^{vg}f}_{x,\omega}^{\xi}(X)dv.$$
Cette \'egalit\'e reste vraie si $X\not\in \omega$ puisque les deux membres sont nuls. Alors on reconna\^{\i}t
$$\Phi(g)=\Delta(x)^r\int_{U_{x}(F)\backslash U(F)}I_{x,\omega}(\theta,f,vg)dv.$$
En remarquant que $C'(x)\Delta(x)^r=C(x)$, la formule (6) devient
$$I_{N}(\theta,f)=C(x)\int_{U_{x}(F)H_{x}(F)\backslash G(F)}I_{x,\omega}(\theta,f,g)\kappa_{N}(g)dg$$
$$\qquad =C(x)I_{x,\omega,N}(\theta,f).\,\,\square$$

\bigskip

\subsection{Localisation de $I(\theta,f)$}

Modifions les notations de 7.3: pour $T\in \underline{\cal T}$, on note maintenant $W'_{T}$, $W''_{T}$ et $V''_{T}$ les espaces que l'on avait not\'es $W'$, $W''$ et $V''$ dans ce paragraphe. On note $\underline{\cal T}_{x}$ le sous-ensemble des $T\in \underline{\cal T}$ tels que $T\subset H_{x}$ et $W'\subset W'_{T}$. Remarquons que ces conditions impliquent que $T$ se d\'ecompose en $T'T''$ o\`u $T'$ est un sous-tore maximal de $H'$ et $T''$ est un sous-tore de $H''$. On a $x\in T'$. Pour $T\in \underline{\cal T}_{x}$, on  a $xexp(X)\in T_{\natural}(F)$ pour tout $X\in \mathfrak{t}_{\natural}(F)\cap \omega$. On d\'efinit des fonctions $c_{\theta,x,\omega}$ et $c_{f,x,\omega}$ presque partout sur $\mathfrak{t}(F)$. Elles sont nulles hors de $\mathfrak{t}(F)\cap \omega$. Pour $X\in \mathfrak{t}_{\natural}(F)\cap \omega$,
$$c_{\theta,x,\omega}(X)=c_{\theta}(xexp(X)),\,\,c_{f,x,\omega}(X)=c_{f}(xexp(X)).$$
En fait, les fonctions $\theta_{x,\omega}$ et $\theta_{f,x,\omega}$ sont des quasi-caract\`eres et les fonctions ci-dessus sont associ\'ees \`a ces quasi-caract\`eres comme en 7.9. On fixe un ensemble de repr\'esentants ${\cal T}_{x}$ des classes de conjugaison par $H_{x}(F)$ dans $\underline{\cal T}_{x}$. Enfin, on  d\'efinit une fonction $\Delta''$ sur $\mathfrak{h}_{x}(F)$ par
$$\Delta''(X)=\vert det(X\vert W''/W''(X))\vert _{F},$$
o\`u $W''(X)$ est le noyau de $X$ agissant dans $W''$. Posons
$$I_{x,\omega}(\theta,f)=\sum_{T\in {\cal T}_{x}}\vert W(H_{x},T)\vert ^{-1}\nu(T)\int_{\mathfrak{t}(F)}c_{\theta,x,\omega}(X)c_{f,x,\omega}(X)D^{H_{x}}(X)\Delta''(X)^rdX.$$
On pourrait montrer que cette int\'egrale est absolument convergente de la m\^eme fa\c{c}on qu'en 7.3. Cela va aussi r\'esulter de la preuve suivante.

\ass{Lemme}{On a l'\'egalit\'e $I(\theta,f)=C(x)I_{x,\omega}(\theta,f)$.}

Preuve. On a les propri\'et\'es suivantes.

(1) Soient $T\in {\cal T}$ et $t\in T_{\natural}(F)$. Alors $ c_{f}(t)=0$ si $t$ n'appartient pas \`a 
$$\bigcup_{T_{1}\in {\cal T}_{x}}\bigcup_{w\in W^+(T_{1},T)}w(xexp(\mathfrak{t}_{1}(F)\cap \omega)).$$

(2) Soit $T\in {\cal T}$ et, pour $i=1,2$, soient $T_{i}\in {\cal T}_{x}$ et $w_{i}\in W^+(T_{i},T)$. Alors les ensembles $w_{1}(xexp(\mathfrak{t}_{1}(F)\cap \omega))$ et $w_{2}(xexp(\mathfrak{t}_{2}(F)\cap \omega))$ sont disjoints ou confondus.

(3) Soient $T\in {\cal T}$, $T_{1}\in {\cal T}_{x}$ et $w_{1}\in W^+(T_{1},T)$. Le nombre des couples $(T_{2},w_{2})$ tels que $T_{2}\in {\cal T}_{x}$, $w_{2}\in W^+(T_{2},T)$ et $w_{2}(xexp(\mathfrak{t}_{2}(F)\cap \omega))=w_{1}(xexp(\mathfrak{t}_{1}(F)\cap \omega))$ est \'egal \`a
$$\vert W(H_{x},T_{1})\vert \vert Z_{H^+}(x)(F)/H_{x}(F)\vert \vert Z_{H^+}(T_{1})(F)/Z_{H_{x}}(T_{1}) (F)\vert ^{-1}.$$

Soient $T$ et $t$ comme en (1). Supposons $c_{f}(t)\not=0$. Alors $\theta_{f}$ n'est nulle dans aucun voisinage de $t$.  Le support de $\theta_{f}$ est inclus dans la cl\^oture de $(Supp (f))^G$, donc dans $\Omega$. Donc $t\in \Omega$ et on peut fixer $y\in G(F)$ et $X\in \omega$ tels que $yty^{-1}=xexp(X)$.
Puisque $t\in T_{\natural}(F)$, le noyau de $t-1$ agissant dans $V$ est $V''_{T}$. Gr\^ace \`a la condition $(7)_{\rho}$ de 3.1, le noyau de $xexp(X)-1$ est contenu dans $V''$. Donc $y(V''_{T})\subset V''$. Comme dans la preuve de 8.2(2), on peut alors modifier $y$ et $X$ de telle sorte que $y$ conserve $D\oplus Z$. Cela entra\^{\i}ne que $xexp(X)$ agit sur cet espace par l'identit\'e, donc $X\in \mathfrak{h}_{x}(F)$. L'\'el\'ement $y$ conserve $W$. Notons $h$ sa restriction \`a cet espace, qui appartient \`a $H^+(F)$. Posons $T_{1}=hTh^{-1}$. On a $T\subset H_{t}$, donc $T_{1}\subset H_{xexp(X)}\subset H_{x}$. De plus, puisque $y(V''_{T})\subset V''$, on a $W'\subset h(W'_{T})$. Mais alors le tore $T_{1}$ appartient \`a $\underline{\cal T}_{x}$. Quitte \`a multiplier $h$ \`a gauche par un \'el\'ement de $H_{x}(F)$, on peut supposer $T_{1}\in {\cal T}_{x}$. En notant $w\in W^+(T_{1},T)$ l'isomorphisme induit par la conjugaison par $h^{-1}$, on a $t\in w(xexp(\mathfrak{t}_{1}(F)\cap \omega))$, ce qui prouve (1).
 
Les assertions (2) et (3) se prouvent comme (3) et (4) de 8.2.   Remarquons toutefois que le quotient $Z_{H^+}(T_{1})(F)/Z_{H_{x}}(T_{1}) (F)$ figurant dans (3) est fini car, puisque $x\in T_{1}(F)$, le groupe $Z_{H^+}(T_{1})$ est contenu dans $Z_{H^+}(x)$. On laisse les d\'etails au lecteur.

Les trois assertions pr\'ec\'edentes permettent d'\'ecrire
$$I(\theta,f)=\sum_{T_{1}\in {\cal T}_{x}}\sum_{T\in {\cal T}}\sum_{w_{1}\in W^+(T_{1},T)}w(T_{1})\vert W(H,T)\vert ^{-1}\nu(T)$$
$$\qquad\int_{\mathfrak{t}_{1}(F)\cap \omega}c_{\theta}(w_{1}(xexp(X)))c_{f}(w_{1}(xexp(X)))D^H(w_{1}(xexp(X)))\Delta(xexp(X))^rdX,$$
o\`u $w(T_{1})$ est l'inverse du nombre de couples calcul\'e en (3). Tous les termes contenant $w_{1}$ sont invariants par $H^+(F)$ et le $w_{1}$ dispara\^{\i}t. On  a aussi $\nu(T)=\nu(T_{1})$ et $\vert W(H,T)\vert =\vert W(H,T_{1})$ si $W^+(T_{1},T)$ n'est pas vide. On a les \'egalit\'es $c_{\theta}(xexp(X))=c_{\theta,x,\omega}(X)$, $c_{f}(xexp(X))=c_{f,x,\omega}(X)$ et, gr\^ace aux hypoth\`eses (7) et $(7)_{\rho}$ de 3.1, 
$$D^H(xexp(X))\Delta(xexp(X))^r=D^H(x)D^{H_{x}}(X)\Delta(x)^r\Delta''(X)^r.$$
 On obtient
$$(4) \qquad I(\theta,f)=D^H(x)\Delta(x)^r\sum_{T_{1}\in {\cal T}_{x}}w'(T_{1})\nu(T_{1})\int_{\mathfrak{t}_{1}(F)}c_{\theta,x,\omega}(X)c_{f,x,\omega}(X)D^{H_{x}}(X)\Delta''(X)^rdX,$$
o\`u 
$$w'(T_{1})=w(T_{1})\vert W(H,T_{1})\vert ^{-1} \vert \{(T,w_{1}); T\in {\cal T}, w_{1}\in W^+(T_{1},T)\}\vert .$$
On calcule ce terme comme dans la preuve du lemme 8.2. On obtient
$$D^H(x)\Delta(x)^rw'(T_{1})=C(x)\vert W(H_{x},T_{1})\vert ^{-1}.$$
Alors la formule (4) devient celle de l'\'enonc\'e. $\square$

\bigskip

\section{Utilisation de la transformation de Fourier}

\bigskip

\subsection{Position du probl\`eme }

Comme on l'a dit, on conserve la situation de 8.2. Posons $U''=U\cap G''$. Remarquons que $U''=U_{x}$. Soient $\theta''$ un quasi-caract\`ere de $\mathfrak{h}''(F)$ et $\varphi\in C_{c}^{\infty}(\mathfrak{g}''(F))$. Appliquant les d\'efinitions de 7.9 o\`u l'on remplace les espaces $V$ et $W$ par $V''$ et $W''$, on d\'efinit une fonction $\varphi^{\xi}$ sur $\mathfrak{h}''(F)$ et, pour $g\in G''(F)$, une int\'egrale $I(\theta'',\varphi,g)$.  Remarquons que, si le support de $\varphi$ est contenu dans $\omega''$, celui de $\varphi^{\xi}$ est contenu dans $\omega''\cap \mathfrak{h}''(F)$.   Soit $S\in \mathfrak{h}''(F)$. On suppose que $S$ est r\'egulier et que le noyau de $S$ agissant dans $W''$ est de dimension au plus $1$.  On suppose que, pour toute $\phi\in C_{c}^{\infty}( \mathfrak{h}''(F))$ \`a support dans $\omega''\cap \mathfrak{h}''(F)$, on a l'\'egalit\'e
$$\theta''(\phi)=J_{H''}(S,\hat{\phi}).$$
Soit enfin $\kappa''\in C_{c}^{\infty}(U''(F)H''(F)\backslash G''(F))$.  G\'en\'eralisant la d\'efinition de 7.8, on pose
$$I_{\kappa''}(\theta'',\varphi)=\int_{U''(F)H''(F)\backslash G''(F)}I(\theta'',\varphi,g)\kappa''(g)dg.$$
Cette int\'egrale est \`a support compact. Le but de la section est d'exprimer $I_{\kappa''}(\theta'',\varphi)$ \`a l'aide de la transform\'ee de Fourier $\hat{\varphi}$ de $\varphi$, quand $\varphi$ est \`a support dans $\omega''$.

\bigskip

\subsection{Premi\`ere transformation}

Soit $\Xi$ l'\'el\'ement de $\mathfrak{g}''(F)$ qui annule $W''$ et v\'erifie $\Xi v_{i+1}=\xi_{i}v_{i}$ pour tout $i=0,...,r-1$. Remarquons que l'on a $\Xi v_{0}=-2\nu_{0}\xi_{0}e_{-1}$, o\`u $\nu_{0}=q(v_{0})$, $\Xi v_{-i}=-\xi_{i}v_{-i-1}$ pour $i=1,...,r-1$ et $\Xi v_{-r}=0$. On a aussi $\xi(N)=<\Xi,N>$ pour tout $N\in \mathfrak{u}''(F)$. 

Posons $\Lambda_{0}=\{c(v_{0},v); v\in W''\}$. Cet espace est l'orthogonal de $\mathfrak{h}''(F)$ dans $\mathfrak{g}''_{0}(F)$.  La forme bilin\'eaire $<.,.>$ est non d\'eg\'en\'er\'ee sur $\Lambda_{0}$.   Posons $\Sigma=\mathfrak{a}(F) \oplus \Lambda_{0}\oplus \mathfrak{u}''(F)$. On munit les deux premiers espaces de la mesure autoduale. On a implicitement fix\'e une mesure sur $U''(F)$ dans le paragraphe pr\'ec\'edent, dont le choix n'importe pas. On en d\'eduit une mesure sur $\mathfrak{u}''(F)$, puis sur $\Sigma$.

\ass{Lemme}{Pour tout $\varphi\in C_{c}^{\infty}(\mathfrak{g}''(F))$ et tout $Y\in \mathfrak{h}''(F)$, on a l'\'egalit\'e
$$(\varphi^{\xi})\hat{}(Y)=\int_{\Sigma}\hat{\varphi}(\Xi+Y+X)dX.$$}

Preuve. Introduisons le groupe unipotent $\bar{\mathfrak{u}}''$ oppos\'e \`a $\mathfrak{u}''$. Les espaces $\bar{\mathfrak{u}}''(F)$ et $\mathfrak{u}''(F)$ sont en dualit\'e. La mesure sur le second espace se dualise en une mesure sur le premier et la transformation de Fourier \'echange $C_{c}^{\infty}(\mathfrak{u}''(F))$ et $C_{c}^{\infty}(\bar{\mathfrak{u}}''(F))$. On a l'\'egalit\'e
$$\mathfrak{g}''=\bar{\mathfrak{u}}''\oplus \mathfrak{a}\oplus \mathfrak{h}''\oplus \Lambda_{0}\oplus \mathfrak{u}''.$$
Par lin\'earit\'e, on peut supposer que
$$\varphi=\varphi_{\bar{\mathfrak{u}}''(F)}\otimes \varphi_{\mathfrak{a}(F)}\otimes \varphi_{\mathfrak{h}''(F)}\otimes \varphi_{\Lambda_{0}}\otimes\varphi_{\mathfrak{u}''(F)},$$
o\`u, pour chaque espace $E$ figurant en indice, $\varphi_{E}\in C_{c}^{\infty}(E)$. On a
$$\hat{\varphi}=\hat{\varphi}_{\mathfrak{u}''(F)}\otimes \hat{\varphi}_{\mathfrak{a}(F)}\otimes \hat{\varphi}_{\mathfrak{h}''(F)}\otimes \hat{\varphi}_{\Lambda_{0}}\otimes\hat{\varphi}_{\bar{\mathfrak{u}}''(F)}.$$
Pour $Y\in \mathfrak{h}''(F)$, on calcule
$$\varphi^{\xi}(Y)=\varphi_{\bar{\mathfrak{u}}''(F)}(0) \varphi_{\mathfrak{a}(F)}(0) \varphi_{\mathfrak{h}''(F)}(Y)\varphi_{\Lambda_{0}}(0)\hat{\varphi}_{\mathfrak{u}''(F)}(\Xi),$$
$$(\varphi^{\xi})\hat{}(Y)=\varphi_{\bar{\mathfrak{u}}''(F)}(0) \varphi_{\mathfrak{a}(F)}(0) \hat{\varphi}_{\mathfrak{h}''(F)}(Y) \varphi_{\Lambda_{0}}(0)\hat{\varphi}_{\mathfrak{u}''(F)}(\Xi),$$
$$\int_{\Sigma}\hat{\varphi}(\Xi+Y+X)dX=\hat{\varphi}_{\mathfrak{u}''(F)}(\Xi)\hat{\varphi}_{\mathfrak{h}''(F)}(Y)\int_{\Sigma}\hat{\varphi}_{\mathfrak{a}(F)}\otimes \hat{\varphi}_{\Lambda_{0}}\otimes \hat{\varphi}_{\bar{\mathfrak{u}}''(F)}(X)dX,$$
$$\qquad =\hat{\varphi}_{\mathfrak{u}''(F)}(\Xi)\hat{\varphi}_{\mathfrak{h}''(F)}(Y)\varphi_{\mathfrak{a}(F)}(0)\varphi_{\Lambda_{0}}(0)\varphi_{\bar{\mathfrak{u}}''(F)}(0).$$
Le lemme r\'esulte de la comparaison des \'egalit\'es ci-dessus. $\square$

\bigskip

\subsection{Description de l'espace affine $\Xi+S+\Sigma$}

Notons $\Lambda_{\mathfrak{u}''}$ le sous-espace de $\mathfrak{u}''(F)$ engendr\'e par les \'el\'ements $c(v_{i},v_{i+1})$ pour $i=0,...,r-1$. Si $d$ est impair ou si $r=0$, on pose $\Lambda=\Lambda_{0}\oplus \Lambda_{\mathfrak{u}''}$. Supposons $d$ pair, donc $ dim(W'')$ impair. Alors $S$, agissant dans $W''$, a un noyau de dimension $1$. On fixe un \'el\'ement non nul $w_{S}$ de ce noyau et on note $W''_{S}$ son orthogonal dans $W''$.  Supposons de plus $r>0$. On pose
$$\Lambda_{0,S}=\{c(v_{0},v); v\in W''_{S}\},$$
$$\Lambda=\Lambda_{0,S}\oplus Fc(w_{S},v_{r})\oplus \Lambda_{\mathfrak{u}''}.$$
Dans les deux cas, $\Lambda$ est un sous-espace de $\Sigma$. Puisque $\Sigma$ et $\Lambda$ sont des espaces vectoriels sur $F$,on peut les consid\'erer comme les ensembles de points sur $F$ de vari\'et\'es sur $\bar{F}$ que, dans ce paragraphe, on note encore $\Sigma$ et $\Lambda$.  

\ass{Lemme}{L'espace affine $\Xi+S+\Sigma$ est stable par conjugaison par $U''$. L'application
$$\begin{array}{ccc}U''\times (\Xi+S+\Lambda)&\to&\Xi+S+\Sigma\\ (u,X)&\mapsto&u^{-1}Xu\\ \end{array}$$
est un isomorphisme de vari\'et\'es alg\'ebriques.}

Preuve. L'annulateur de $\Sigma$ dans $\mathfrak{g}''$ est l'espace $\mathfrak{h}''\oplus \mathfrak{u}''$. Pour prouver la premi\`ere assertion, il suffit de prouver que, pour $u\in U''$, $X\in \Sigma$ et $Y\in \mathfrak{h}''\oplus \mathfrak{u}''$, on a l'\'egalit\'e
$$trace(u(\Xi+S+X)u^{-1}Y)=trace((\Xi+S)Y),$$
ou encore
$$trace((\Xi+S+X)u^{-1}Yu)=trace((\Xi+S)Y).$$
Posons $u^{-1}Yu=Y+N$. On a $N\in \mathfrak{u}''$ et
$$trace(u(\Xi+S+X)u^{-1}Y)=trace((\Xi+S)Y)+trace(\Xi N)+trace(XY)+trace((S+X)N).$$
Les deux derniers termes sont nuls: ce sont des traces d'\'el\'ements de $\mathfrak{u}''$. Il faut montrer que $trace(\Xi N)=0$, ou encore $\xi(N)=0$. Il suffit pour cela de prouver que $q(Nv_{i},v_{-i-1})=0$ pour $i=0,...,r-1$. Mais $u-1$ et $Y$ appartiennent \`a l'alg\`ebre de Lie du radical unipotent du sous-groupe parabolique de $GL(V'')$ qui conserve le drapeau 
$$Fv_{r}\subset Fv_{r}\oplus Fv_{r-1}\subset...\subset Fv_{r}\oplus...\oplus F_{v_{0}}.$$
  Donc $Nv_{i}$ appartient au sous-espace engendr\'e par les $v_{j}$ pour $j\geq i+2$. Donc $q(Nv_{i},v_{-i-1})=0$, ce qui prouve la premi\`ere assertion de l'\'enonc\'e.
  
   Si $r=0$, on a $\Lambda=\Sigma$, $U''=\{1\}$ et la seconde assertion est tautologique. Supposons $r>0$. Introduisons le sous-groupe parabolique $P_{2}$ de $G''$ qui conserve le sous-espace totalement isotrope $Z_{+}$, sa composante de L\'evi  $M_{2}$ qui conserve $Z_{+}$ et $Z_{-}$ et son radical unipotent $U_{2}$. Le groupe $M_{2}$ s'identifie \`a $GL(Z_{+})\times G''_{0}$. Notons $U_{4}$ le centre de $U_{2}$. Les groupes $U_{4}$ et $U_{2}/U_{4}$ sont ab\'eliens. Par l'application $(v,v')\mapsto exp( c(v,v'))$, ils  s'identifient respectivement \`a $\bigwedge^2(Z_{+})$ et $Hom(V''_{0},Z_{+})$. Ce dernier espace se d\'ecompose en $Hom(W'',Z^+)\oplus Hom(D,Z^+)$. On note $U_{3}$  le sous-groupe de $U_{2}$ tel que $U_{3}/U_{4}$ s'identifie \`a $Hom(W'',Z_{+})$ et $U_{D}$ celui tel que $U_{D}/U_{4}$ s'identifie \`a $Hom(D,Z_{+})$. On pose $U_{1}=U''$, $U_{5}=\{1\}$. Remarquons que l'on a les inclusions $U_{2}\subset U''\subset P_{2}$. On a donc la cha\^{\i}ne de sous-groupes
  $$U_{5}\subset U_{4}\subset U_{3}\subset U_{2}\subset U_{1},$$
  et chacun de ces sous-groupes est distingu\'e dans $U_{1}$. Posons $\mathfrak{r}=\{c(v_{-1},v); v\in Z_{+}\}$. C'est un sous-espace de $\mathfrak{u}''$. D\'efinissons les espaces
  $$\Sigma_{1}=\Sigma= \mathfrak{a}\oplus \Lambda_{0}\oplus \mathfrak{u}_{1};$$
  $$\Sigma_{2}=\Lambda_{0}\oplus \mathfrak{r}\oplus \mathfrak{u_{2}};$$
  $$\Sigma_{3}=\Lambda_{0}\oplus \mathfrak{u_{2}};$$
  $$\Sigma_{4}=\left\lbrace\begin{array}{cc}\Lambda_{0}\oplus \mathfrak{u}_{D},&\,\,{\rm si}\,\,d\,\,{\rm est\,\,impair},\\ \Lambda_{0,S}\oplus \bar{F}c(w_{S},v_{r})\oplus \mathfrak{u}_{D},&\,\,{\rm si}\,\,d\,\,{\rm est\,\,pair};\\ \end{array}\right.$$
  $$\Sigma_{5}=\Lambda.$$
  On a les inclusions 
  $$\Sigma_{5}\subset \Sigma_{4}\subset \Sigma_{3}\subset \Sigma_{2}\subset \Sigma_{1}.$$
  Pour $i=2,...,4$, $\Sigma_{i}$ est l'ensemble des \'el\'ements $X\in\Sigma_{i-1}$ qui v\'erifient les conditions suivantes:
  
  (1) si $i=2$, $Xv_{j}=0$ pour $j=2,...,r$;
  
  (2) si $i=3$, $Xv_{1}=0$;
  
  (3) si $i=4$ et $d$ est impair, $X(W'')\subset Z_{+}\oplus D$; si $i=4$ et $d$ est pair, $X(W''_{S})\subset Z_{+}\oplus D$ et $X(w_{S})\in\bar{F}v_{r}$.
  
  On a
  
  (4) pour $i=1,2,3$, les ensembles $\Sigma_{i}$ et $S+\Sigma_{i}$ sont stables par conjugaison par $U_{1}$; pour $i=4,5$, les ensembles $\Sigma_{i}$ et $S+\Sigma_{i}$ sont stables par conjugaison par $U_{4}$; l'ensemble $\Sigma_{4}$  est stable par conjugaison par $U_{2}$.
  
 Posons $M''=M\cap G''$.  En g\'en\'eral, si $E$ est un sous-ensemble de $\mathfrak{m}''$, $E\oplus \mathfrak{u}_{1}$ est invariant par conjugaison par $U_{1}$. Si $E$ est un sous-ensemble de $\mathfrak{g}''_{0}$, $E\oplus \mathfrak{u}_{2}$ est stable par conjugaison par $U_{1}$. Si $E$ est un sous-ensemble de $\mathfrak{g}''_{0}\oplus \mathfrak{u}_{1}$, $E$ est stable par conjugaison par $U_{4}$. On en d\'eduit que $\Sigma_{1}$, $S+\Sigma_{1}$, $\Sigma_{3}$ et $S+\Sigma_{3}$ 
 sont stables par conjugaison par $U_{1}$, et  $\Sigma_{4}$, $S+\Sigma_{4}$, $\Sigma_{5}$ et $S+\Sigma_{5}$ sont stables par conjugaison par $U_{4}$. On a $\Sigma_{2}=\Sigma_{3}\oplus \mathfrak{r}$. L'ensemble $\Sigma_{3}$ est stable par conjugaison par $U_{1}$. Pour prouver que $\Sigma_{2}$ l'est aussi, il suffit de prouver que, pour $u\in U_{1}$ et $X\in \mathfrak{r}$, on a $u^{-1}Xu\in \Sigma_{2}$. Il est clair que cet \'el\'ement appartient \`a $\mathfrak{u}_{1}$, donc \`a $\Sigma_{1}$. On doit montrer qu'il v\'erifie la condition (1). C'est clair puisque $u$ conserve le sous-espace de base $(v_{j})_{j=2,...,r}$ tandis que $X$ annule ce sous-espace. Le m\^eme raisonnement s'applique \`a l'ensemble $S+\Sigma_{2}$. Soient $u\in U_{2}$ et $X\in \Sigma_{4}$. Puisque $\Sigma_{3}$ est stable par conjugaison par $U_{1}$, on $u^{-1}Xu\in \Sigma_{3}$. Pour prouver que cet \'el\'ement appartient \`a $\Sigma_{4}$, on doit montrer qu'il v\'erifie (3). Soit $w\in W''$. On a $uw\in w+Z_{+}$, puis $Xuw=Xw$ car $X$ annule $Z_{+}$. On a $Xw\in Z_{+}\oplus D$ car $X\in \Sigma_{4}$. Or $u^{-1}$ conserve cet espace, donc $u^{-1}Xuw\in Z_{+}\oplus D$. Si $d$ est pair, on a $Xw_{S}\in \bar{F}v_{r}$ et $u^{-1}$ conserve cette droite, donc aussi $u^{-1}Xuw_{S}\in \bar{F}v_{r}$. Cela prouve (4).
  
  On va montrer
  
  (5) pour $i=1,...,4$, l'ensemble $\Xi+S+\Sigma_{i}$ est stable par conjugaison par $U_{i}$.
  
  Pour $i=1$, c'est la premi\`ere assertion de l'\'enonc\'e. Supposons $i\geq 2$.  On sait d\'ej\`a par (4) que  $S+\Sigma_{i}$ est stable par conjugaison par $U_{i}$. On doit donc prouver que, pour $u\in U_{i}$, on a $(u^{-1}\Xi u-\Xi)\in \Sigma_{i}$. En raisonnant par r\'ecurrence sur $i$, on peut supposer que l'on a en tout cas $(u^{-1}\Xi u-\Xi)\in \Sigma_{i-1}$ (pour $i=2$, cette hypoth\`ese r\'esulte de la premi\`ere assertion de l'\'enonc\'e). On doit montrer que cet \'el\'ement v\'erifie les conditions (1), resp. (2), (3), si $i=2$, resp. $i=3,4$. Supposons $i=2$.     Soit $j=2,...,r$. On a $uv_{j}=v_{j}$ et $u^{-1}v_{j-1}=v_{j-1}$ par d\'efinition de $U_{2}$. On a aussi $\Xi v_{j}=\xi_{j-1}v_{j-1}$ et on d\'eduit l'\'egalit\'e $(u^{-1}\Xi u-\Xi)v_{j}=0$ que l'on cherchait \`a prouver. Supposons $i=3$.   On a $uv_{1}=v_{1}$, $\Xi v_{1}=\xi_{0}v_{0}$, $u^{-1}v_{0}=v_{0}$ par d\'efinition de $U_{3}$, d'o\`u encore l'assertion. Supposons $i=4$. Pour $w\in W''$, on a $uw=w$ et $\Xi w=0$. Donc $(u^{-1}\Xi u-\Xi)w=0$ et $u^{-1}\Xi u$ v\'erifie la condition requise. Cela d\'emontre (5).
  
  Gr\^ace \`a (5),  pour $i=1,...,4$, on peut former le quotient $U_{i}\times_{U_{i+1}}\Sigma_{i+1}$ de $U_{i}\times \Sigma_{i+1}$ par la relation d'\'equivalence $(u,X)\equiv (u',X')$ si et seulement s'il existe $v\in U_{i+1}$ tel que $(u',X')=(uv,v^{-1}Xv)$. On va montrer que
  
  (6) l'application
    $$\begin{array}{ccc}U_{i}\times(\Xi+S+ \Sigma_{i+1})&\to&\Xi+S+\Sigma_{i}\\ (u,X)&\mapsto &u^{-1}Xu\\ \end{array}$$
    se descend en un isomorphisme de $U_{i}\times_{U_{i+1}}\Sigma_{i+1}$ sur $\Xi+S+\Sigma_{i}$.
  
  Supposons $i=1$. Posons $U_{B}=U_{1}\cap M_{2}$. Ce groupe s'identifie au radical unipotent du sous-groupe de Borel $B$ de $GL(Z_{+})$ qui conserve le drapeau 
  $$Fv_{r}\subset Fv_{r}\oplus Fv_{r-1}\subset...\subset Fv_{r}\oplus...\oplus Fv_{1}.$$
  L'application produit de $U_{B}\times U_{2}$ sur $U_{1}$ est un isomorphisme. Il suffit de prouver que l'application
  $$\begin{array}{ccc}U_{B}\times (\Xi+S+ \Sigma_{2})&\to&\Xi+S+\Sigma_{1}\\ (u,X)&\mapsto &u^{-1}Xu\\ \end{array}$$
  est un isomorphisme.   On a $\Sigma_{1}=\mathfrak{b}\oplus \Sigma_{3}$, $\Sigma_{2}=\mathfrak{r}\oplus \Sigma_{3}$ et $\mathfrak{r}$ est le sous-ensemble des \'el\'ements de $\mathfrak{b}$ dont seuls les termes de la derni\`ere colonne sont non nuls. D\'efinissons $\underline{\Xi}\in \mathfrak{g}''$ par $\underline{\Xi}v_{j}=\xi_{j-1}v_{j-1}$ pour $j=2,...,r$, $\underline{\Xi}v_{1}=0$ et $\underline{\Xi}$ annule $V''_{0}$. On a $\underline{\Xi}\in End(Z_{+})\subset \mathfrak{m}_{2}$. Pour $u\in U_{B}$, l'image de $u-1$ est contenu dans le sous-espace de $V''$ engendr\'e par les vecteurs $v_{j}$ pour $j=2,...,r$ et $v_{-j}$ pour $j=1,...,r-1$. L'\'el\'ement $\Xi-\underline{\Xi}$ annule cet espace. Son image est contenue dans le plan engendr\'e par $v_{0}$ et $v_{-1}$, lequel est annul\'e par $u^{-1}-1$. On en d\'eduit que
   $u^{-1}\Xi u-\Xi=u^{-1}\underline{\Xi}u-\underline{\Xi}$. On est ramen\'e \`a prouver que l'application
  $$\begin{array}{ccc}U_{B}\times (\underline{\Xi}+\mathfrak{r})&\to&\Xi+\mathfrak{b}\\ (u,X)&\mapsto &u^{-1}Xu\\ \end{array}$$
  est un isomorphisme. Tout se passe dans $End(Z_{+})$. L'assertion est bien connue et se prouve en filtrant $U_{B}$ de la fa\c{c}on habituelle.
  
  Supposons $i=2$. L'application
  $$\begin{array}{ccc}Hom(D,Z_{+})\times U_{3}&\to& U_{2}\\ (Y,u)&\mapsto &exp(Y)u\\ \end{array}$$
  est un isomorphisme. On est ramen\'e \`a prouver que l'application
  $$\begin{array}{ccc}Hom(D,Z_{+})\times(\Xi+S+ \Sigma_{3})&\to&\Xi+S+\Sigma_{2}\\ (Y,X)&\mapsto &exp(-Y)Xexp(Y)\\ \end{array}$$
  est un isomorphisme.  D'apr\`es (4), $S+\Sigma_{3}$ est stable par conjugaison  par $exp(Y)$ pour tout $Y\in Hom(D,Z_{+})$. Cela nous ram\`ene \`a prouver que l'application de $Hom(D,Z_{+})$ dans $\mathfrak{r}=\Sigma_{2}/\Sigma_{3}$ qui, \`a $Y\in Hom(D,Z_{+})$, associe l'image dans $\Sigma_{2}/\Sigma_{3}$ de $exp(-Y)\Xi exp(Y)-\Xi$, est un isomorphisme. L'espace $\mathfrak{r}$ s'identifie \`a $Z_{+}$ par $X\mapsto Xv_{1}$. Il s'agit donc de montrer que l'application
  $$\begin{array}{ccc}Hom(D,Z_{+})&\to &Z_{+}\\ Y&\mapsto& (exp(-Y)\Xi exp(Y)-\Xi)v_{1}\\ \end{array}$$
  est un isomorphisme. On a $exp(Y)v_{1}=v_{1}$, $\Xi v_{1}=\xi_{0}v_{0}$, $exp(-Y)v_{0}=-Yv_{0}+v_{0}$. L'application est donc $Y\mapsto -Yv_{0}$, qui est bien un isomorphisme.
  
  Supposons $i=3$. L'application
  $$\begin{array}{ccc}Hom(W'',Z_{+})\times U_{4}&\to& U_{3}\\ (Y,u)&\mapsto &exp(Y)u\\ \end{array}$$
  est un isomorphisme. D'apr\`es (4), l'ensemble $\Sigma_{4}$ est invariant par conjugaison par $U_{2}$. Comme dans le cas $i=2$, on est ramen\'e \`a prouver que l'application de $Hom(W'',Z_{+})$ dans $\Sigma_{3}/\Sigma_{4}$ qui, \`a $Y\in Hom(W'',Z_{+})$ associe l'image dans $\Sigma_{3}/\Sigma_{4}$ de $exp(-Y)(\Xi+S)exp(Y)-\Xi-S$, est un isomorphisme. Supposons $d$ impair. Notons $proj_{Z_{+}}$ la projection de $V''$ sur $Z_{+}$ de noyau $V''_{0}\oplus Z_{-}$. Alors $\Sigma_{3}/\Sigma_{4}$ s'identifie  \`a $Hom(W'',Z_{+})$ par l'application qui \`a $X\in \Sigma_{3}$ associe la restriction \`a $W''$ de $proj_{Z_{+}}\circ X$. Soit $w\in W''$. On a $exp(Y)w=w+Yw$, $Sexp(Y)w=Sw$, $exp(-Y)Sexp(Y)w=Sw-YSw$, $\Xi w=0$, $\Xi exp(Y)w=\Xi Yw$. Ce dernier \'el\'ement appartient \`a l'espace $Z_{+,0}$ de base $(v_{j})_{j=0,...,r-1}$. Puisque $Y$ annule $Z_{+,0}$, on a $exp(-Y)\Xi Yw=\Xi Yw$. Donc 
  $$proj_{Z_{+}}((exp(-Y)(\Xi+S)exp(Y)-\Xi-S)w)=proj_{Z_{+}}(\Xi Yw-YSw).$$
  On a introduit ci-dessus un \'el\'ement $\underline{\Xi}$. On a $proj_{Z_{+}}\circ \Xi=\underline{\Xi}$ sur $Z_{+}$. La formule ci-dessus devient
   $$proj_{Z_{+}}((exp(-Y)(\Xi+S)exp(Y)-\Xi-S)w)=(\underline{\Xi}Y-YS)w,$$
   et on est ramen\'e \`a prouver que l'application $Y\mapsto \underline{\Xi}Y-YS$ de $Hom(W'',Z_{+})$ dans lui-m\^eme est un isomorphisme. Pour $k=0,...,r$, introduisons le sous-espace $Z_{+}^k$  de $Z_{+}$ de base $(v_{j})_{j=1,...,k}$. L'espace $Hom(W'',Z_{+})$ est filtr\'e par les $Hom(W'',Z_{+}^k)$. L'application pr\'ec\'edente respecte cette filtration et l'application du gradu\'e qui s'en d\'eduit est la m\^eme que celle d\'eduite de $Y\mapsto YS$. Cette derni\`ere est un isomorphisme puisque les valeurs propres de $S$ agissant dans $W''$ sont non nulles. Supposons maintenant $d$ pair. Notons $proj_{Z_{+,0}}$ la projection de $V''$ sur $Z_{+,0}$ de noyau $Fv_{r}\oplus W''\oplus Z_{-}$. Alors $\Sigma_{3}/\Sigma_{4}$ s'identifie \`a $Hom(W''_{S},Z_{+})\oplus Hom(\bar{F}w_{S},Z_{+,0})$ par l'application qui, \`a $X\in \Sigma_{3}$ associe la somme de la restriction \`a $W''_{S}$ de $proj_{Z_{+}}\circ X$ et de la restriction \`a $Fw_{S}$ de $proj_{Z_{+,0}}\circ X$. Soit $Y\in Hom(W'',Z_{+})$, que l'on d\'ecompose en $Y=Y_{1}+Y_{2}$ avec $Y_{1}\in Hom(W''_{S},Z_{+})$ et $Y_{2}\in Hom(\bar{F}w_{S},Z_{+})$. On v\'erifie comme ci-dessus que l'image de $exp(-Y)(\Xi+S)exp(Y)-\Xi-S $ dans $\Sigma_{3}/\Sigma_{4}$ est la somme de la restriction \`a $W''_{S}$ de $\underline{\Xi}Y_{1}-Y_{1}S$ et de $\underline{\Xi}Y_{2}$. Parce que les valeurs propres de $S$ dans $W''_{S}$ sont non nulles, l'application $Y_{1}\mapsto \underline{\Xi}Y_{1}-Y_{1}$ est un isomorphisme pour la m\^eme raison que ci-dessus. L'application $Y_{2}\mapsto \underline{\Xi}Y_{1}$ est un isomorphisme car $\underline{\Xi}$ se restreint en un isomorphisme de $Z_{+}$ sur $Z_{+,0}$.
   
   Supposons $i=4$. Gr\^ace \`a (4), on est encore ramen\'e \`a prouver que l'application de $\mathfrak{u}_{4}$ dans $\Sigma_{4}/\Sigma_{5}$ qui, \`a $Y\in \mathfrak{u}_{4}$,  associe l'image de $exp(-Y)\Xi exp(Y)$ dans $\Sigma_{4}/\Sigma_{5}$, est un isomorphisme.  L'espace $\mathfrak{u}_{4}$, resp. $\mathfrak{u}_{D}$, $\Lambda_{\mathfrak{u}''}$,  a pour base les $c(v_{j},v_{k})$ pour $1\leq j<k\leq r$, resp. pour $0\leq j<k\leq r$, pour $0\leq j<k=j+1\leq r$.  L'injection de $\mathfrak{u}_{D}$ dans $\Sigma_{4}$ se quotiente en un isomorphisme de $\mathfrak{u}_{D}/\Lambda_{\mathfrak{u}''}$ sur $\Sigma_{4}/\Sigma_{5}$. Un calcul simple montre que l'application qui nous int\'eresse s'identifie \`a l'application  $\tau:\mathfrak{u}_{4}\to\mathfrak{u}_{D}/\Lambda_{\mathfrak{u}''}$ ainsi d\'efinie: pour $1\leq j<k\leq r$, $\tau(c(v_{j},v_{k}))$ est l'image dans $\mathfrak{u}_{D}/\Lambda_{\mathfrak{u}''}$ de $c(v_{j},v_{k-1})-c(v_{j-1},v_{k})$. Pour $l\in \{0,...,r\}$ notons $E_{l}$ le sous-espace de $\mathfrak{u}_{D}$ engendr\'e par les $c(v_{j},v_{k}$ tels que $0\leq j<k\leq l+j \leq r$. L'espace $\mathfrak{u}_{D}/\Lambda_{\mathfrak{u}''}$ est filtr\'e par les espaces $E_{l}/E_{1}$. L'espace $\mathfrak{u}_{4}$ est filtr\'e par les espaces $E_{l-1}\cap \mathfrak{u}_{4}$. On v\'erifie que $\tau$ est compatible avec ces filtrations et que l'application gradu\'ee qui s'en d\'eduit est un isomorphisme. Cela ach\`eve la preuve de (6).
   
   En appliquant (6) successivement pour $i=1,...,4$, on obtient la seconde assertion de l'\'enonc\'e. $\square$
   
   \bigskip
   
   \subsection{Polyn\^ome caract\'eristique}
   
    On introduit un syst\`eme hyperbolique  maximal $(w_{\pm j })_{j=1,...,m}$ de $W''\otimes_{F}\bar{F}$ form\'e de vecteurs propres pour $S$. On note $s_{j}$ la valeur propre de $S$ sur $w_{j}$, pour $j>0$. Si $d$ est impair, resp. pair,  $(w_{\pm j })_{j=1,...,m}$ est une base de $W''\otimes_{F}\bar{F}$, resp. $W''_{S}\otimes_{F}\bar{F}$. Si $d$ est pair, on pose $\nu_{S}=q(w_{S})$. On introduit des coordonn\'ees sur $\Lambda$ en \'ecrivant un \'el\'ement $X\in \Lambda$ sous la forme suivante:
 
 - si $d$ est impair,
 $$X=c(v_{0},\sum_{j=\pm 1,...,\pm m}z_{j}w_{j})+\sum_{i=0,...r-1}\lambda_{i}c(v_{i},v_{i+1});$$
 
 -si $d$ est pair et $r>0$,
 $$X=c(v_{0},\sum_{j=\pm 1,...,\pm m}z_{j}w_{j})+z_{0}c(w_{S},v_{r})+\sum_{i=0,...r-1}\lambda_{i}c(v_{i},v_{i+1});$$
 
 - si $d$ est pair et $r=0$,
 $$X=c(v_{0},z_{0 }w_{S}+\sum_{j=\pm 1,...,\pm m}z_{j}w_{j}) .$$

Notons $R_{S}$ le polyn\^ome caract\'eristique de $S$ agissant dans $W''$. On a donc
$$R_{S} (T)=\left\lbrace\begin{array}{cc}\prod_{j=1,...,m}(T^2-s_{j}^2),&\,\,{\rm si\,\,}d{\rm\,\,est\,\,impair,}\\ T \prod_{j=1,...,m}(T^2-s_{j}^2),&\,\,{\rm si\,\,}d{\rm\,\,est\,\,pair.}\\ \end{array}\right.$$
Pour $X\in \mathfrak{g}''$, on note $P_{X}$ le polyn\^ome caract\'eristique de $X$ agissant dans $V''$.  
  
  \ass{Lemme}{Soit $X\in \Lambda$, auquel on associe des coordonn\'ees comme ci-dessus. On a les \'egalit\'es suivantes:
  
  - si $d$ est impair,
  $$P_{\Xi+S+X}(T)=T^{2r+1}R_{S}(T)+  \sum_{j=1,...,m}4\nu_{0}z_{j}z_{-j} \frac {R_{S}(T)T^{2r+1}}{T^2-s_{j}^2}$$
  $$\qquad +\sum_{i=0,...,r-1}(-1)^{i+1}4\nu_{0}R_{S}(T)T^{2r-1-2i}\lambda_{i}\xi_{i}\prod_{i'=0,...,i-1}\xi_{i'}^2 ;$$
  
  - si $d$ est pair et $r>0$,
    $$P_{\Xi+S+X}(T)=T^{2r+1}R_{S}(T)+ \sum_{j=1,...,m}4\nu_{0}z_{j}z_{-j}\frac{R_{S}(T)T^{2r+1}}{T^2-s_{j}^2}$$
    $$\qquad  +(-1)^{r}4\nu_{S}\nu_{0}z_{0}^2\frac{R_{S}(T)}{T}(\prod_{i=0,...,r-1}\xi_{i}^2)+\sum_{i=0,...,r-1}(-1)^{i+1}4\nu_{0}R_{S}(T)T^{2r-1-2i}\lambda_{i}\xi_{i}\prod_{i'=0,...,i-1}\xi_{i'}^2    ;$$
  
  - si $d$ est pair et $r=0$,
  $$P_{\Xi+S+X}(T)=TR_{S}(T)+ \sum_{j=1,...,m}4\nu_{0}z_{j}z_{-j}\frac{R_{S}(T)T}{T^2-s_{j}^2} +4\nu_{S}\nu_{0}z_{0}^2\frac{R_{S}(T)}{T}.$$}
  
  Preuve. On \'ecrit  l'\'el\'ement $\Xi+S+X$ comme une matrice. Les m\'ethodes usuelles de d\'eveloppement selon les lignes ou les colonnes permettent d'exprimer son d\'eterminant comme une somme de termes ais\'es \`a calculer et d'un d\'eterminant analogue \`a celui de d\'epart mais associ\'e \`a des valeurs de $r$ ou $m$ strictement inf\'erieures. En raisonnant par r\'ecurrence, on obtient l'assertion. On renonce \`a r\'ediger davantage la preuve. Indiquons simplement la forme de la matrice dans deux exemples.
  
  Supposons $m=2$, $r=2$ et $d$ est impair. On choisit pour base ordonn\'ee de $V''$ la famille $v_{2},v_{1},w_{2},w_{1},v_{0},w_{-1},w_{-2},v_{-1},v_{-2}$. Dans cette base, la matrice de $\Xi+S+X$ est
  $$\left(\begin{array}{ccccccccc}0&0&0&0&0&0&0&\lambda_{1}&0\\ \xi_{1}&0&0&0& 2\nu_{0}\lambda_{0}&0&0&0&-\lambda_{1}\\ 0&0&s_{2}&0&2\nu_{0}z_{2}&0&0&0&0\\ 0&0&0&s_{1}&2\nu_{0}z_{1}&0&0&0&0\\ 0&\xi_{0}&-z_{-2}&-z_{-1}&0&-z_{1}&-z_{2}&-\lambda_{0}&0\\0&0&0&0&2\nu_{0}z_{-1}&-s_{1}&0&0&0\\0&0&0&0&2\nu_{0}z_{-2}&0&-s_{2}&0&0\\ 0&0&0&0&-2\nu_{0}\xi_{0}&0&0&0&0\\0&0&0&0&0&0&0&-\xi_{1}&0\\ \end{array}\right)$$
  Supposons $m=2$, $r=2$ et $d$ est pair. On choisit pour base ordonn\'ee de $V''$ la famille $v_{2},v_{1},w_{2},w_{1}, w_{S},v_{0},w_{-1},w_{-2},v_{-1},v_{-2}$. Dans cette base, la matrice de $\Xi+S+X$ est
  $$\left(\begin{array}{cccccccccc}0&0&0&0&2\nu_{S}z_{0}&0&0&0&\lambda_{1}&0\\ \xi_{1}&0&0&0& 0&2\nu_{0}\lambda_{0}&0&0&0&-\lambda_{1}\\ 0&0&s_{2}&0&0&2\nu_{0}z_{2}&0&0&0&0\\ 0&0&0&s_{1}&0&2\nu_{0}z_{1}&0&0&0&0\\ 0&0&0&0&0&0&0&0&0&-z_{0}\\
  0&\xi_{0}&-z_{-2}&-z_{-1}&0&0&-z_{1}&-z_{2}&-\lambda_{0}&0\\0&0&0&0&0&2\nu_{0}z_{-1}&-s_{1}&0&0&0\\0&0&0&0&0&2\nu_{0}z_{-2}&0&-s_{2}&0&0\\ 0&0&0&0&0&-2\nu_{0}\xi_{0}&0&0&0&0\\0&0&0&0&0&0&0&0&-\xi_{1}&0\\ \end{array}\right)$$
  $\square$
  
   Remarquons que les termes $z_{j}z_{-j}$, $\lambda_{i}$ et  $z_{0}^2$ dans le cas o\`u $d$ est pair sont d\'etermin\'es par $P_{\Xi+S+X}$. On a en particulier
   $$(1) \qquad z_{j}z_{-j}=\frac{P_{\Xi+S+X}(s_{j})}{4\nu_{0}s_{j}^{1+2r}R_{S,j}(s_{j})}$$
   pour $j=1,...m$, o\`u $R_{S,j}(T)=\frac{R_{S}(T)}{T^2-s_{j}^2}$,
    $$(2) \qquad z_{0}^2=\left\lbrace\begin{array}{cc}\frac{P_{\Xi+S+X}(0)}{(-1)^{r} 4\nu_{S}\nu_{0}R_{S,0}(0)\prod_{i=0,...,r-1}\xi_{i}^2} ,&\,\,{\rm si}\,\,r>0,\\ \frac{P_{\Xi+S+X}(0)}{4\nu_{S}\nu_{0}R_{S,0}(0)},&,\,\, {\rm si}\,\,r=0,\\ \end{array}\right.$$
   o\`u $R_{S,0}(T)=\frac{R_{S}(T)}{T}$. Posons $d''=dim(V'')$ et notons $Pol_{d''}$ l'espace des polyn\^omes de degr\'e $d''$, \`a coefficients dans $F$, de coefficient dominant \'egal \`a $1$ et ne contenant que des puissances de l'ind\'etermin\'ee $T$ de m\^eme parit\'e que $d''$. C'est exactement l'espace des polyn\^omes caract\'eristiques des \'el\'ements de $\mathfrak{g}''(F)$. Introduisons le sous-ensemble $Pol_{d''}^S$ form\'es des polyn\^omes $P$ tels que
   
   $P$ est le polyn\^ome caract\'eristique d'un \'el\'ement de $Y\in\mathfrak{g}''_{reg}(F)$;  
      
   $P(s_{j})\not=0$ pour tout $j=1,...,m$ et $P(0)\not=0$ si $d$ est pair.
   
   C'est un ouvert de Zariski non vide de $Pol_{d''}$. Notons $\Lambda^S$ le sous-ensemble des $X\in \Lambda$ tels que $\Xi+S+X\in \mathfrak{g}''_{reg}(F)$,   $z_{j}\not=0$ pour tout $j\in\{\pm 1,...,\pm m\}$ et de plus, si $d$ est pair,   $z_{0}\not=0$. C'est exactement l'image r\'eciproque de $Pol_{d''}^S$ dans $\Lambda$ par l'application $X\mapsto P_{\Xi+S+X}$. Donc $\Lambda^S$ est un ouvert de Zariski non vide de $\Lambda$. Les formules du lemme montrent que l'application pr\'ec\'edente restreinte \`a $\Lambda^S$ est une application $F$-analytique surjective et partout submersive de $\Lambda^S$ sur $Pol_{d''}^S$.

  \bigskip
  
  \subsection{Orbites dans $\Xi+S+\Lambda$}

   Notons $\Sigma^S$ le sous-ensemble de $\Sigma$ tel que  l'image de $U''(F)\times(\Xi+S+ \Lambda^S)$ par l'isomorphisme du lemme 9.3 soit $\Xi+S+\Sigma^S$ 
   
   \ass{Lemme}{Le groupe $H''_{S}(F)U''(F)$ agit par conjugaison dans $\Xi+S+\Sigma^S$ et cette action est libre. Deux \'el\'ements de $\Xi+S+\Sigma^S$ sont conjugu\'es par un \'el\'ement de $G''$ si et seulement s'ils le sont par un \'el\'ement de $H''_{S}(F)U''(F)$.}
   
   Preuve. Soient $Y\in \Xi+S+\Sigma^S$ et $g\in H''_{S}(F)U''(F)$ tels que $g^{-1}Yg=Y$.  Par d\'efinition de $\Sigma^S$, on peut \'ecrire $Y=u^{-1}Y'u$, avec $u\in U''(F)$ et $Y'\in \Xi+S+\Lambda^S$. Alors $ug^{-1}u^{-1}Y'ugu^{-1}=Y'$. Quitte \`a remplacer $Y$ par $Y'$ et $g$ par $ugu^{-1}$, on est ramen\'e au cas o\`u $Y\in \Xi+S+\Lambda^S$. On peut \'ecrire $g=tu$, avec $t\in H''_{S}(F)$ et $u\in U''(F)$. La conjugaison par $t$  fixe $\Xi+S$ et conserve $\Lambda$. Plus pr\'ecis\'ement, introduisons des coordonn\'ees sur $\Lambda$ comme en 9.4. L'\'el\'ement $t$ agit par homoth\'etie sur chaque droite $\bar{F}w_{j}$, pour $j=\pm 1,...\pm m$. Pour $j>0$, on note $t_{j}$ la valeur propre associ\'ee.  Alors la conjugaison par $t$ laisse inchang\'ees les coordonn\'ees $\lambda_{i}$ et $z_{0}$ dans le cas o\`u $d$ est pair. Elle agit sur les coordonn\'ees restantes par
   $$(1) \qquad (z_{r},...,z_{1},z_{-1},...,z_{-r})\mapsto (t_{r}z_{r},...,t_{1}z_{1},t_{1}^{-1}z_{-1},...,t_{r}^{-1}z_{-r}).$$
   Posons $Y'=t^{-1}Yt$. Alors $Y$ et $Y'$ sont deux \'el\'ements de $\Xi+S+\Lambda$ qui sont conjugu\'es par l'\'el\'ement $u\in U''(F)$. Le lemme 9.3 entra\^{\i}ne que $u=1$ et $Y=Y'$. Ecrivons $Y=\Xi+S+X$, avec $X\in \Lambda^S$. Les coordonn\'ees $z_{j}$ de $X$ sont toutes non nulles et la formule ci-dessus montre que $X$ ne peut \^etre fix\'e par $t$ que si tous les $t_{j}$ valent $1$, autrement dit $t=1$. Donc $g=tu=1$ et cela d\'emontre la premi\`ere assertion de l'\'enonc\'e.
   
   Comme ci-dessus, on peut remplacer dans la seconde assertion l'ensemble $\Xi+S+\Sigma^S$   par $\Xi+S+\Lambda^S$. Soient $X,\underline{X}\in \Lambda^S$, notons comme en 9.4 les coordonn\'ees de $X$ et notons par des lettres soulign\'ees celles de $\underline{X}$. Supposons $\Xi+S+X$ et $\Xi+S+\underline{X}$  conjugu\'es par un \'el\'ement de $G''$. Alors $P_{\Xi+S+X}=P_{\Xi+S+\underline{X}}$. D'apr\`es les remarques du paragraphe pr\'ec\'edent, on a $z_{j}z_{-j}=\underline{z}_{j}\underline{z}_{-j}$ pour tout $j=1,...,m$, $\lambda_{i}=\underline{\lambda}_{i}$ pour tout $i=0,...,r-1$ et $z_{0}^2=\underline{z}_{0}^2$ si $d$ est pair. Supposons d'abord $d$ impair. La formule (1) ci-dessus montre qu'il existe un unique $t\in H''_{S}(\bar{F})$ tel que $t^{-1}Xt=\underline{X}$. L'unicit\'e de $t$ et le fait que $X$ et $\underline{X}$ sont tous deux d\'efinis sur $F$ entra\^{\i}nent que $t\in H''_{S}(F)$. Alors $\Xi+S+X$ et $\Xi+S+\underline{X}$ sont conjugu\'es par un  \'el\'ement de $H''_{S}(F)$, ce que l'on voulait d\'emontrer. Supposons maintenant $d$ pair. On trouve comme dans le cas $d$ impair un unique \'el\'ement $t\in H''_{S}(F)$ tel que  $t^{-1}Xt=\underline{X}$ ou $\underline{X}'$, ce dernier \'el\'ement ayant les m\^emes coordonn\'ees que $\underline{X}$, \`a l'exception de $\underline{z}_{0}$ qui est chang\'e en $-\underline{z}_{0}$. On a alors soit $t^{-1}(\Xi+S+X)t=\Xi+S+\underline{X}$, soit $t^{-1}(\Xi+S+X)t=\Xi+S+\underline{X}'$. Il suffit pour conclure de prouver que cette deuxi\`eme possibilit\'e ne se produit pas. Consid\'erons l'\'el\'ement $\delta$ du groupe orthogonal $G^{_{''}+}(F)$ qui agit par  multiplication par $-1$ sur la droite $Fw_{S}$ et qui fixe tout \'el\'ement de l'orthogonal de cette droite. On v\'erifie que $\delta^{-1}(\Xi+S+\underline{X})\delta=\Xi+S+\underline{X}'$. On sait par hypoth\`ese que $\Xi+S+X$ est conjugu\'e \`a $\Xi+S+\underline{X}$ par un \'el\'ement de $G''$. S'il \'etait conjugu\'e par $t$ \`a  $\Xi+S+\underline{X}'$, les deux \'el\'ements $\Xi+S+\underline{X}$ et $\Xi+S+\underline{X}'$ seraient conjugu\'es par un \'el\'ement de $G''$ et l'ensemble $\delta G''$ couperait le centralisateur de $\Xi+S+\underline{X}$ dans $G^{_{''}+}$. Or ce centralisateur est contenu dans $G''$ parce que $\Xi+S+\underline{X}$ est r\'egulier et n'a pas de valeur propre nulle (cela parce que son polyn\^ome caract\'eristique n'est pas nul en $0$). Puisque $\delta\not\in G''$, on obtient une contradiction qui ach\`eve la preuve. $\square$
   
   \bigskip
   
   \subsection{Mesures autoduales}
   
  Consid\'erons l'application
  $$\mathfrak{g}''_{reg}(F)\to \bigsqcup_{T\in {\cal T}(G'')}(\mathfrak{t}(F)\cap \mathfrak{g}''_{reg}(F))/W(G'',T)$$
  qui, \`a un \'el\'ement de $\mathfrak{g}''_{reg}(F)$, associe l'unique \'el\'ement de l'ensemble d'arriv\'ee qui lui est conjugu\'e par un \'el\'ement de $G''(F)$. Elle est analytique. Pour tout sous-tore maximal $T$ de $G''$, on note $\mathfrak{t}(F)^S$ le sous-ensemble des \'el\'ements de $\mathfrak{t}(F)$ qui sont conjugu\'es \`a un \'el\'ement de $\Xi+S+\Sigma^S$ par un \'el\'ement de $G(F)$. L'application pr\'ec\'edente se restreint en une application analytique
  $$(1) \qquad \Xi+S+\Sigma^S\to \bigsqcup_{T\in {\cal T}(G'')}\mathfrak{t}(F)^S/W(G'',T).$$
 Elle est surjective.   Si on note $(\Xi+S+\Sigma^S)/H''_{S}(F)U''(F)$ l'ensemble des classes de conjugaison par $H''_{S}(F)U''(F)$ dans $\Xi+S+\Sigma^S$,  le lemme pr\'ec\'edent montre qu'elle se quotiente en une bijection 
   $$(2) \qquad (\Xi+S+\Sigma^S)/H''_{S}(F)U''(F)\to \bigsqcup_{T\in {\cal T}(G'')}\mathfrak{t}(F)^S/W(G'',T).$$
  On munit l'ensemble de d\'epart de la mesure quotient des mesures d\'ej\`a fix\'ees sur $\Xi+S+\Sigma^S$ et $H''_{S}(F)U''(F)$. Les remarques de la fin du paragraphe 9.4 montrent que l'application (1) est partout submersive. La mesure sur l'ensemble de d\'epart de (2) s'identifie donc \`a une mesure r\'eguli\`ere sur l'ensemble d'arriv\'ee. Pour tout $T\in {\cal T}(G'')$, l'ensemble $\mathfrak{t}(F)^S$ est ainsi muni d'une mesure que l'on note $d_{\Sigma}Y$. Rappelons que l'on note simplement $dY$ la mesure autoduale.
  
  \ass{Lemme}{Pour tout $T\in {\cal T}(G'')$, on a l'\'egalit\'e $d_{\Sigma}Y=D^{H''}(S)^{-1/2}D^{G''}(Y)^{1/2}dY$ en tout point $Y\in \mathfrak{t}(F)^S$.}
  
  Preuve. Fixons $T\in {\cal T}(G'')$. Un objet tel que $\mathfrak{t}(F)^S$ ou $\mathfrak{t}(F)^S/W(G'',T)$ n'a pas de structure alg\'ebrique naturelle. Commen\c{c}ons par alg\'ebriser la situation. On consid\`ere $\Sigma^S$ comme une vari\'et\'e alg\'ebrique (un ouvert d'un espace vectoriel). Notons $\bar{W}(G'',T)=Norm_{G''}(T)/T$, introduisons l'ensemble $\mathfrak{t}^S$ des \'el\'ements de $\mathfrak{t}$ qui sont conjugu\'es \`a un \'el\'ement de $\Xi+S+\Sigma^S$ puis le quotient $\mathfrak{t}/\bar{W}(G'',T)$. Ce sont des vari\'et\'es alg\'ebriques. Il y a une application alg\'ebrique
  $$(3) \qquad \tau:\Xi+S+\Sigma^S\to \mathfrak{t}^S/\bar{W}(G'',T)$$
  qui se quotiente en un isomorphisme
  $$(4) \qquad (\Xi+S+\Sigma^S)/H''_{S}U''\to \mathfrak{t}^S/\bar{W}(G'',T).$$
  La structure alg\'ebrique sur $\mathfrak{t}^S/\bar{W}(G'',T)$ d\'etermine une structure analytique sur $(\mathfrak{t}^S/\bar{W}(G'',T))(F)$. Il y a une application naturelle 
  $$\iota:\mathfrak{t}(F)^S \to 
  (\mathfrak{t}^S/\bar{W}(G'',T))(F),$$
  qui est localement un isomorphisme de vari\'et\'es analytiques. Cela va nous permettre de remplacer l'application (2) par son avatar alg\'ebrique (4).
  
   Rappelons que, une fois le corps $F$ muni de la mesure autoduale, pour toute vari\'et\'e alg\'ebrique lisse ${\cal X}$ d\'efinie sur $F$, une forme diff\'erentielle $\delta$ sur ${\cal X}$,   d\'efinie sur $F$  et de degr\'e maximal, d\'efinit une mesure $\vert \delta\vert _{F}$ sur ${\cal X}(F)$. Plus g\'en\'eralement, ne supposons plus $\delta$ d\'efinie sur $F$. Il existe une fonction alg\'ebrique $\alpha$ sur ${\cal X}$, non nulle et telle   que $\alpha \delta$ soit d\'efinie sur $F$. Etendons la valeur absolue de $F$ \`a $\bar{F}$. On d\'efinit une mesure $\vert \delta\vert _{F}$ sur ${\cal X}(F)$ par 
  $$\vert \delta\vert _{F}=\vert \alpha\vert _{F}^{-1}\vert \alpha\delta\vert _{F}.$$
  Cela ne d\'epend pas du choix de $\alpha$.En particulier, soit $E$ un sous-$F$-espace de $\mathfrak{g}''(F)$ sur lequel la forme $<.,.>$ est non d\'eg\'en\'er\'ee. Fixons une base $(e_{k})_{k=1,...l}$ de $E$ sur $F$, notons $Q$ la matrice $l\times l$ telle que $Q_{k,k'}=<e_{k},e_{k'}>$ et \'ecrivons tout \'el\'ement de $E$ sous la forme $e=\sum_{k=1,...,l}x_{k}e_{k}$. D\'efinissons la forme diff\'erentielle $\delta=\bigwedge_{k=1,...,l}dx_{k}$. On v\'erifie que la mesure autoduale sur $E$ est 
  $$(5) \qquad \vert det(Q)\vert _{F}^{-1/2}\vert \delta\vert _{F}.$$
    Supposons maintenant que $(e_{k})_{k=1,...,l}$ est une base de $E\otimes_{F}\bar{F}$. La forme diff\'erentielle $\delta=\bigwedge_{k=1,...,l} dx_{k}$ n'est pas, en g\'en\'eral, d\'efinie sur $F$ mais il existe $\alpha\in \bar{F}^{\times}$ tel que $\alpha\delta$ le soit et on peut d\'efinir $\vert \delta\vert _{F}$ comme plus haut. Un simple calcul de changement de bases montre que la mesure autoduale sur $E$ est encore donn\'ee par la formule (5).
    
     Choisissons des formes lin\'eaires $Y\mapsto y_{k}$, $k=1,...,l $,  sur $\mathfrak{t}$ de sorte que, pour un \'el\'ement  $Y\in \mathfrak{t}$ en position g\'en\'erale, l'action de $Y$ dans $V''$  ait pour valeurs propres non nulles $(\pm y_{k})_{k=1,...,l}$. Avec des notations \'evidentes, on a l'\'egalit\'e
    $$<Y,Y'>=\frac{1}{2}trace(YY')=\sum_{k=1,...,l}y_{k}y'_{k}.$$
   D\'efinissons la forme diff\'erentielle $\delta_{\mathfrak{t}}$ sur $\mathfrak{t}$ par $\delta_{\mathfrak{t}}=\bigwedge_{k=1,...,l}dy_{k}$. La formule (5) montre que la mesure autoduale sur $\mathfrak{t}(F)$ est $\vert \delta_{\mathfrak{t}}\vert _{F}$. Fixons un sous-ensemble positif de l'ensemble des racines de $T$ dans $\mathfrak{g}''$. Pour $Y\in \mathfrak{t}$, posons 
   $$d^{G''}(Y)=\prod_{\alpha>0}\alpha(Y),$$
   le produit \'etant pris sur cet ensemble de racines. On v\'erifie que la forme diff\'erentielle $d^{G''}\delta_{\mathfrak{t}}$ se descend en une forme diff\'erentielle sur $\mathfrak{t}/\bar{W}(G'',T)$, que l'on note $\delta_{\mathfrak{t}/W}$. Evidemment, la mesure autoduale sur $\mathfrak{t}(F)$ est
   $$(6) \qquad dY=\vert d^{G''}(Y)\vert _{F}^{-1} \vert \iota^*(\delta_{\mathfrak{t}/W})(Y)\vert _{F}=D^{G''}(Y)^{-1/2} \vert \iota^*(\delta_{\mathfrak{t}/W})(Y)\vert _{F}.$$
   
Introduisons comme en 9.4 un syst\`eme hyperbolique maximal $(w_{\pm j})_{j=1,...,m}$ de $W''\otimes_{F}\bar{F}$ form\'e de vecteurs propres pour $S$, donc aussi pour $H''_{S}$. Pour $t\in H''_{S}$ et $j=1,...,m$, notons $t_{j}$ la valeur propre de $t$ sur $w_{j}$. D\'efinissons $\delta_{H''_{S}}=(\prod_{j=1,...,m}t_{j})^{-1}\bigwedge_{j=1,...,m}dt_{j}$. La formule (5), remont\'ee au groupe par l'exponentielle, montre que $\vert \delta_{H''_{S}}\vert _{F}$ est la mesure que nous avons fix\'ee sur $H''_{S}(F)$.   Fixons une base de $\mathfrak{u}''(F)$ sur $F$ et prenons pour $\delta_{\mathfrak{u}''}$ le produit, dans un ordre fix\'e, des diff\'erentielles des coordonn\'ees relativement \`a cette base. On a implicitement fix\'e une mesure sur $\mathfrak{u}''(F)$, mais notre probl\`eme est insensible au choix de cette mesure. On peut donc supposer que cette mesure est $\vert \delta_{\mathfrak{u}''}\vert _{F}$. Via l'exponentielle, $\delta_{\mathfrak{u}''}$ d\'efinit une forme diff\'erentielle $\delta_{U''}$ sur $U''$ et la mesure de Haar sur $U''(F)$ n'est autre que $\vert \delta_{U''}\vert _{F}$. Introduisons des coordonn\'ees sur $\Lambda_{0}$  (qui est vu ici comme une vari\'et\'e alg\'ebrique sur $F$) en \'ecrivant tout \'el\'ement $X$ de cet ensemble sous la forme
   
   - si $d$ est impair, $X=c(v_{0},\sum_{j=\pm 1,...,\pm m}z_{j}w_{j})$;
   
   - si $d$ est pair, $X=c(v_{0},z_{0}w_{S}+\sum_{j=\pm 1,...,\pm m}z_{j}w_{j})$.
   
   Remarquons que l'\'eventuel terme $z_{0}$ n'est pas le m\^eme qu'en 9.4. Puisqu'on a ici \'etendu les scalaires, on peut supposer $q(v_{0})=1$ et, si $d$ est pair, $q(w_{S})=-1$. On v\'erifie alors que
   $$<X,X'>=\frac{1}{2}trace(XX')= [z_{0}z'_{0}]-\sum_{\pm 1,...,\pm m}z_{j}z'_{-j} ,$$
   o\`u, ici comme dans la suite, on indique symboliquement entre crochets les termes qui n'existent que dans le cas $d$ pair.   On pose 
   $$\delta_{\Lambda_{0}}=  \bigwedge_{j=[0],\pm 1,...,\pm m}dz_{j} .$$
   D'apr\`es (5), $\vert \delta_{\Lambda_{0}}\vert _{F}$ est la mesure autoduale sur $\Lambda_{0}$.  Notons $(a_{i})_{i=1,...,r}$ les valeurs propres sur les vecteurs $(v_{i})_{i=1,...,r}$ d'un \'el\'ement de $\mathfrak{a}$. D\'efinissons $\delta_{\mathfrak{a}}=\bigwedge_{i=1,...,r}da_{i}$. D'apr\`es (5), $\vert \delta_{\mathfrak{a}}\vert _{F}$ est la mesure autoduale sur $\mathfrak{a}_{F}$. Rappelons que $\Sigma=\mathfrak{a}\oplus \Lambda_{0}\oplus \mathfrak{u}''$.   On d\'efinit la forme diff\'erentielle $\underline{\delta}$ sur $\Xi+S+\Sigma^S$ qui, via la translation par $\Xi+S$, correspond \`a la forme diff\'erentielle $\delta_{\mathfrak{a}}\wedge \delta_{\Lambda_{0}}\wedge \delta_{\mathfrak{u}''}$ sur $\Sigma^S$.   Alors $\vert \underline{\delta}\vert _{F}$ est la mesure que nous avons fix\'ee sur $\Xi+S+\Sigma^S$. On v\'erifie que $\underline{\delta}$ est invariante par conjugaison par $H''_{S}U''$. Il y a alors une  forme diff\'erentielle $\bar{\delta}$ sur le quotient $ (\Xi+S+\Sigma^S)/H''_{S}U''$ de sorte que, via le choix de sections locales, on ait l'\'egalit\'e $\underline{\delta}=\delta_{H''_{S}}\wedge \delta_{U''}\wedge \bar{\delta}$. Par (4), $\bar{\delta}$ correspond \`a une forme $\beta \delta_{\mathfrak{t}/W}$  sur $\mathfrak{t}/\bar{W}(G'',T)$, o\`u $\beta$ est une fonction alg\'ebrique sur cette vari\'et\'e. Remontons $\beta$ en une fonction sur $\mathfrak{t}$.  En tenant compte de (6), on voit que l'on a l'\'egalit\'e
   $$(7) \qquad d_{\Sigma}(Y)=D^{G''}(Y)^{1/2}\vert \beta(Y)\vert _{F}dY$$
   pour tout $Y\in \mathfrak{t}(F)^S$.  
   
   Il s'agit de calculer la fonction $\beta$. Supposons d'abord $r=0$. Dans ce cas $\Xi=0$ et $U''=\{1\}$. Introduisons le sous-ensemble $\Lambda_{1}$ des \'el\'ements $X\in\Lambda_{0}$  \'ecrits comme plus haut, tels que $z_{-j}=1$ pour $j=1,...,m$. L'action $(t,S+X)\mapsto S+X'=t(S+X)t^{-1}$ de $H''_{S}$ sur $S+\Sigma$ s'\'ecrit, avec les syst\`emes de coordonn\'ees que l'on a introduits,
   $$((t_{j})_{j=1,...,m},(z_{j})_{j=[0],\pm 1,...,\pm m})\mapsto (z'_{j})_{j=[0],\pm 1,...,\pm m},$$
   o\`u $z'_{j}=t_{j}z_{j}$  et $z'_{-j}=t_{j}^{-1}z_{-j}$ pour $j=1,...,m$, et $z'_{0}=z_{0}$ dans le cas $d$ pair.  De cette action se d\'eduit un isomorphisme de $H''_{S}\times S+\Lambda_{1}$ sur l'ensemble des $S+X\in S+\Lambda_{0}$ dont toutes les coordonn\'ees $z_{-j}$ sont non nulles, lequel contient $S+\Sigma^S$. On peut identifier $(S+\Sigma^S)/H''_{S}$ avec un ouvert dense de $\Lambda_{1}$ et on v\'erifie que, modulo cette identification, $\bar{\delta}=\wedge_{j=[0],1,...,m}dz_{j}$. Soient $X\in \Lambda_{1}$ de coordonn\'ees $(z_{j})_{j=[0],1,...,m}$ et $Y\in \mathfrak{t}$ de coordonn\'ees $(y_{k})_{k=1,...,l}$. On suppose que l'image de $S+X$ par (4) est l'image de $Y$ dans $\mathfrak{t}/\bar{W}(G'',T)$. Supposons pour fixer les id\'ees $d$ pair. On a $l=m+1$ et 
   $$P_{S+X}(T)=P_{Y}(T)=\prod_{k=1,...,l}(T^2-y_{k}^2).$$
   Les formules 9.4(1) et 9.4(2) deviennent
   $$z_{j}=\frac{\prod_{k=1,...,l}(s_{j}^2-y_{k}^2)}{2s_{j}^2\prod_{j'=1,...,m; j'\not=j}(s_{j}^2-s_{j'}^2)} $$
   pour $j\not=0$ et
   $$z_{0}^2=\frac{\prod_{k=1,...,l}y_{k}^2}{\prod_{k=1,...,m}s_{j}^2}.$$
    Cette derni\`ere relation signifie qu'il existe $\epsilon\in \{\pm 1\}$ tel que 
    $$z_{0}=\epsilon\frac{\prod_{k=1,...,l}y_{k}}{\prod_{k=1,...,m}s_{j}}.$$
    En d\'erivant de fa\c{c}on usuelle, on obtient
    $$dz_{j}= -s_{j}^{-2}\prod_{j'=1,...,m; j'\not=j}(s_{j}^2-s_{j'}^2)^{-1}\sum_{k=1,...,l}A_{j,k}  dy_{k}$$
    pour $j\not=0$ et
    $$dz_{0}=\epsilon\prod_{k=1,...,m}s_{j}^{-1}\sum_{k=1,...,l}A_{0,k}dy_{k},$$
    o\`u
    $$A_{j,k}=\left\lbrace\begin{array}{cc}y_{k}\prod_{k'=1,...,l; k'\not=k}(s_{j}^2-y_{k'}^2),&\,\,{\rm si}\,\,j\not=0,\\ \prod_{k'=1,...,l; k'\not=k}y_{k'},&\,\,{\rm si}\,\,j=0.\\ \end{array}\right.$$
    D'o\`u
    $$\bigwedge_{j=0,...,m}dz_{j}=(-1)^m\epsilon\prod_{j=1,...,m}s_{j}^{-3}\prod_{j,j'=1,...,m; j\not=j'}(s_{j}^2-s_{j'}^2)^{-1}det(A)\bigwedge_{k=1,...,l}dy_{k},$$
    o\`u $A$ est la matrice carr\'ee de coefficients $A_{j,k}$. Posons $s_{0}=0$ et
    $$B_{j,k}=\prod_{k'=1,...,l; k'\not=k}(s_{j}^2-y_{k'}^2).$$
    On a
    $$A_{j,k}=\left\lbrace\begin{array}{cc}y_{k}B_{j,k},&\,\,{\rm si}\,\,j\not=0,\\ (-1)^my_{k}(\prod_{k'=1,...,l}y_{k'}^{-1})B_{0,k},&\,\,{\rm si}\,\,j=0.\\ \end{array}\right.$$
    Donc $det(A)=(-1)^mdet(B)$. On laisse au lecteur le calcul \'el\'ementaire  du d\'eterminant de $B$, qui vaut
    $$det(B)=(-1)^{l(l-1)/2}\prod_{0\leq j<j'\leq m}(s_{j}^2-s_{j'}^2)\prod_{1\leq k<k'\leq l}(y_{k}^2-y_{k'}^2)$$
  $$\qquad =(-1)^{m+l(l-1)/2}\prod_{j=1,...,m}s_{j}^2\prod_{1\leq j<j'\leq m}(s_{j}^2-s_{j'}^2) \prod_{1\leq k<k'\leq l}(y_{k}^2-y_{k'}^2).$$
  D'o\`u
  $$(8) \qquad \bigwedge_{j=0,...,m}dz_{j}=\epsilon' d^{H''}(S)^{-1}\prod_{1\leq k<k'\leq l}(y_{k}^2-y_{k'}^2) \bigwedge_{k=1,...,l}dy_{k},$$
  o\`u $\epsilon'=\pm 1$ et
  $$d^{H''}(S)=\prod_{j=1,...,m}s_{j}\prod_{1\leq j'<j\leq m}(s_{j}^2-s_{j'}^2).$$ Remarquons que le  produit intervenant dans (8) est \'egal \`a $\pm d^{G''}(Y)$. Alors (8) devient
  $$\bar{\delta}(S+X)=\epsilon''d^{H''}(S)^{-1}\delta_{\mathfrak{t}/W}(Y),$$
  avec $\epsilon''=\pm \epsilon$, d'o\`u $\beta(Y)=\epsilon''d^{H''}(S)^{-1}$. Puisque $\vert d^{H''}(S)\vert _{F}=D^{H''}(S)^{1/2}$, la formule (7) devient celle de l'\'enonc\'e.
   
     Passons au cas o\`u $r\not=0$. Quitte \`a conjuguer $T$ par un \'el\'ement de $G''$, on peut supposer $A\subset T$. On a alors $T=AT_{0}$ o\`u $T_{0}$ est un sous-tore maximal de $G''_{0}$. On peut supposer que les coordonn\'ees que l'on a introduites sur $\mathfrak{t}$ et $\mathfrak{a}$ sont compatibles. Pr\'ecis\'ement, soit $Y\in \mathfrak{t}$, \'ecrivons $Y=Y_{\mathfrak{a}}+Y_{0}$, avec $Y_{\mathfrak{a}}\in \mathfrak{a}$ et $Y_{0}\in \mathfrak{t}_{0}$, et introduisons les coordonn\'ees $(y_{k})_{k=1,...,l}$ de $Y$ et $(a_{i})_{i=1,...,r}$ de $Y_{\mathfrak{a}}$. On peut supposer $a_{i}=y_{i}$ pour tout $i=1,...,r$. En se pla\c{c}ant dans $G''_{0}$, on d\'efinit la vari\'et\'e $\mathfrak{t}_{0}/\bar{W}(G''_{0},T_{0})$ munie de sa forme diff\'erentielle $\delta_{\mathfrak{t}_{0}/W}$, la forme diff\'erentielle $\bar{\delta}_{0}$ sur $(S+\Lambda_{0})/H''_{S}$ et la fonction $\beta_{0}$ telle que l'application
     $$(S+\Lambda_{0})/H''_{S}\to \mathfrak{t}_{0}/W$$
 identifie $\bar{\delta}_{0}$ \`a $\beta_{0}\delta_{\mathfrak{t}_{0}/W}$, du moins sur un ouvert dense. Le calcul pr\'ec\'edent s'applique: $\beta_{0}$ est constante, de valeur $\epsilon_{0}d^{H''}(S)^{-1}$, o\`u $\epsilon_{0}\in \{\pm 1\}$. Consid\'erons le diagramme
 $$\begin{array}{ccc}U''\times(\Xi+S+\mathfrak{a}+\Lambda_{0})&\stackrel{f_{1}}{\to}&\Xi+S+\Sigma\\ f_{2}\downarrow\,\,&&\downarrow \\ (S+\Lambda_{0})\times \mathfrak{a}&\to&(\Xi+S+\Sigma)/U''\\ \downarrow&&\downarrow\\ (S+\Lambda_{0})/H''_{S}\times \mathfrak{a}&\to&(\Xi+S+\Sigma)/U''H''_{S}\\ f_{3}\downarrow\,\,&&\downarrow\\ (\mathfrak{t}_{0}/\bar{W}(G''_{0},T_{0}))\times \mathfrak{a}&\stackrel{f_{4}}{\to} &\mathfrak{t}/\bar{W}(G'',T)\\ \end{array}$$
 o\`u $f_{1}(u'',\Xi+S+X_{\mathfrak{a}}+X_{\Lambda_{0}})=u^{_{''}-1}(\Xi+S+X_{\mathfrak{a}}+X_{\Lambda_{0}})u''$ et $f_{2}(u'',\Xi+S+X_{\mathfrak{a}}+X_{\Lambda_{0}})=(S+X_{\Lambda_{0}},X_{\mathfrak{a}})$, les autres applications \'etant \'evidentes. Ce diagramme est commutatif. Pour le voir, soient $u''\in U''$, $X_{\mathfrak{a}}\in \mathfrak{a}$ et $X_{\Lambda_{0}}\in \Lambda_{0}$. Soit $Y_{0}\in \mathfrak{t}_{0}$ tel que la partie semi-simple de $S+X_{\Lambda_{0}}$ soit conjugu\'ee \`a $Y_{0}$ par un \'el\'ement de $G''_{0}$. Soit $Y\in \mathfrak{t}$ tel que la partie semi-simple de $\Xi+S+X_{\mathfrak{a}}+X_{\Lambda_{0}}$ soit conjugu\'ee \`a $Y$ par un \'el\'ement de $G''$. L'image de $(u'',\Xi+S+X_{\mathfrak{a}}+X_{\Lambda_{0}})$ par le chemin sud-ouest du diagramme est l'image de $Y_{0}+X_{\mathfrak{a}}$ dans $\mathfrak{t}/\bar{W}(G'',T)$.  Son image par le chemin nord-est est l'image de $Y$ dans cet ensemble. Mais $S+X_{\mathfrak{a}}+X_{\Lambda_{0}}$ appartient \`a $\mathfrak{m}''$ tandis que $\Xi$ appartient \`a $\bar{\mathfrak{u}}''$. Alors les parties semi-simples de $\Xi+S+X_{\mathfrak{a}}+X_{\Lambda_{0}}$ et de $S+X_{\mathfrak{a}}+X_{\Lambda_{0}}$ sont conjugu\'ees par un \'el\'ement de $G''$. Donc $Y$ est conjugu\'e \`a $Y_{0}+X_{\mathfrak{a}}$ et ces deux \'el\'ements ont m\^eme image dans $\mathfrak{t}/\bar{W}(G'',T)$. Cela d\'emontre la commutativit\'e du diagramme. Ce raisonnement et le lemme 9.5 montrent que les fl\`eches horizontales du diagramme sont des isomorphismes locaux, au moins si l'on se restreint \`a des ouverts denses de chaque vari\'et\'e, ce que l'on fait, ici et dans la suite. Du diagramme se d\'eduit une application
 $$(U''\times(\Xi+S+\mathfrak{a}+\Lambda_{0}))/U''H''_{S}\stackrel{g}{\to}\mathfrak{t}/\bar{W}(G'',T),$$
 qui est toujours un isomorphisme local. Munissons $\Xi+S+\mathfrak{a}+\Lambda_{0}\simeq \mathfrak{a}+\Lambda_{0}$ de la diff\'erentielle $\boldsymbol{\delta}=\delta_{\mathfrak{a}}\wedge \delta_{\Lambda_{0}}$. On en d\'eduit une forme diff\'erentielle $\bar{\boldsymbol{\delta}}$ sur l'espace de d\'epart de $g$. Calculons $g^*(\delta_{\mathfrak{t}/W})$ en utilisant le chemin sud-ouest du diagramme. D'apr\`es les d\'efinitions, $f_{4}^*(\delta_{\mathfrak{t}/W})=\pm d^{G''}(d^{G''_{0}})^{-1}\delta_{\mathfrak{t}_{0}/W}\wedge \delta_{\mathfrak{a}}$.  Puis $f_{3}^*f_{4}^*(\delta_{\mathfrak{t}/W})=\pm \beta_{0}^{-1}d^{G''}(d^{G''_{0}})^{-1}\bar{\delta}_{0}\wedge \delta_{\mathfrak{a}}$. Les deux applications verticales restantes identifient $\bar{\boldsymbol{\delta}}$ \`a $\bar{\delta}_{0}\wedge \delta_{\mathfrak{a}}$ et on obtient 
 $$(9) \qquad g^*(\delta_{\mathfrak{t}/W})=\pm \beta_{0}^{-1}d^{G''}(d^{G''_{0}})^{-1}\bar{\boldsymbol{\delta}}.$$
Utilisons le chemin nord-est. Par la suite d'applications verticales, $\delta_{\mathfrak{t}/W}$ se rel\`eve en la forme $\beta^{-1}\underline{\delta}$ sur $\Xi+S+\Sigma)$. Soit $\gamma$ la fonction telle que $f_{1}^*(\underline{\delta})=\gamma \delta_{U''}\wedge \boldsymbol{\delta}$.  Alors
$$(10) \qquad g^*(\delta_{\mathfrak{t}/W})=\pm \gamma \beta^{-1} \bar{\boldsymbol{\delta}}.$$
 Soient $X_{\mathfrak{a}}\in \mathfrak{a}$ et $X_{\Lambda_{0}}\in \Lambda_{0}$. Posons $X=S+X_{\mathfrak{a}}+X_{\Lambda_{0}}$. La diff\'erentielle de $f_{1}$ au point $(1,\Xi+X)$ se calcule ais\'ement. C'est l'application
 $$(11) \qquad \begin{array}{ccc}\mathfrak{u}''\times (\mathfrak{a}+\mathfrak{\Lambda_{0}})&\to& \Sigma=\mathfrak{a}+\mathfrak{\Lambda_{0}}+\mathfrak{u}''\\ (N,X'_{\mathfrak{a}}+X'_{\Lambda_{0}})&\mapsto&X'_{\mathfrak{a}}+X'_{\Lambda_{0}}-[N,\Xi+X]\\ \end{array}$$
 Parce que $X$ appartient \`a $\mathfrak{m}''$ et $\Xi$ \`a $\bar{\mathfrak{u}}''$, on peut trouver une base de $\mathfrak{u}''$ telle que l'application $N\mapsto [N,X]$ soit diagonale, tandis que l'application compos\'ee de $N\mapsto [N,\Xi]$ et de la projection sur $\mathfrak{u}''$ soit nilpotente sup\'erieure. Le d\'eterminant de l'application (11), c'est-\`a-dire $\gamma(X)$, est donc le m\^eme que celui de l'application $N\mapsto [N,X]$ de $\mathfrak{u}''$ dans lui-m\^eme. Celui-ci est \'egal \`a $\pm d^{G''}(X)d^{G''_{0}}(X)^{-1}$. En reportant cette valeur dans (10) et en comparant avec (9), on obtient $\beta=\pm \beta_{0}=\pm \epsilon_{0}d^{H''}(S)^{-1}$. Comme dans le cas $r=0$, la formule (7) devient celle de l'\'enonc\'e. $\square$
 
 \bigskip
 
 \subsection{Sections locales}
 
 L'application (1) de 9.6 est analytique. Le lemme 9.3 et la preuve du lemme 9.5 montrent qu'elle est partout submersive. Pour tout $T\in {\cal T}(G'')$, on peut donc fixer une application localement analytique
  $$\begin{array}{ccc} \mathfrak{t}(F)^S&\to& \Xi+S+\Sigma^S\\ Y&\mapsto&Y_{\Sigma}\\ \end{array}$$
   de sorte que le diagramme 
 $$\begin{array}{ccccc}\Xi+S+\Sigma^S&&\to&&\mathfrak{t}(F)^S/W(G'',T)\\&\nwarrow&&\nearrow&\\ &&\mathfrak{t}(F)^S&&\\ \end{array}$$
 soit commutatif. Il existe une application $Y\mapsto \gamma_{Y}$  de $\mathfrak{t}(F)^S$ dans $T(F)\backslash G''(F)$, localement analytique, de sorte que $Y_{\Sigma}=\gamma_{Y}^{-1}Y\gamma_{Y}$. Mais l'application $G''(F)\to T(F)\backslash G''(F)$ admet elle-m\^eme des sections localement analytiques. On peut donc supposer que l'application $Y\mapsto \gamma_{Y}$ est localement analytique \`a valeurs dans $G''(F)$. 
 
 On va montrer
 
 (1) soit $\omega_{T}$ un sous-ensemble compact de $\mathfrak{t}(F)$; on peut choisir l'application $Y\mapsto Y_{\Sigma}$ telle que l'image de $\mathfrak{t}(F)^S\cap\omega_{T}$ soit contenue dans un sous-ensemble compact de $\Xi+S+\Lambda$.
 
 D'apr\`es le lemme 9.3, on peut supposer que $Y_{\Sigma}\in\Xi+S+\Lambda^S$ pour tout $Y\in \mathfrak{t}(F)^S$.   Soit $Y\in \mathfrak{t}(F)^S\cap \omega_{T}$, posons $Y_{\Sigma}=\Xi+S+X$   et introduisons les coordonn\'ees de $X$ comme en 9.4. Les coordonn\'ees $\lambda_{i}$ sont lin\'eaires en les coefficients du polyn\^ome caract\'eristique $P_{Y}(T)$, et sont donc born\'ees. De m\^eme, l'\'eventuelle coordonn\'ee $z_{0}$ et les produits $z_{j}z_{-j}$, pour $j=1,...,m$, sont born\'es. Montrons que l'on peut supposer chaque $z_{\pm j}$ born\'e. Consid\'erons d'abord deux cas particuliers. Dans le premier, on suppose qu'il existe une extension $F_{1}$ de $F$ de degr\'e $m$ et une extension quadratique $F_{2}$ de $F_{1}$ telle que $H''_{S}(F)$ soit le noyau de la norme de $F_{2}^{\times}$ dans $F_{1}^{\times}$. Dans ce cas, on peut identifier $W''$ \`a $F_{2}$ et l'action de $H''_{S}(F)$ sur $W''$ \`a la multiplication. Posons $w=\sum_{j=\pm 1,...,\pm m}z_{j}w_{j}\in W''=F_{2}$. En normalisant convenablement les vecteurs $w_{j}$, les coordonn\'ees $z_{\pm j}$ sont les images de $w$ par les diff\'erents plongements de $F_{2}$ dans $\bar{F}$. Alors ces coordonn\'ees ont toutes la m\^eme valeur absolue. Puisque les produits $z_{j}z_{-j}$ sont born\'es, chaque terme $z_{\pm j}$ l'est aussi. Dans le deuxi\`eme cas particulier, on consid\`ere une extension $F_{1}$ comme ci-dessus et on suppose que $H_{S}(F)=F_{1}^{\times}$. Dans ce cas, on peut identifier $W''$ \`a $F_{1}\oplus F_{1}$ et l'action de $H''_{S}(F)$ sur $W''$ \`a l'application $(h,w_{+}\oplus w_{-})\mapsto hw_{+}\oplus h^{-1}w_{-}$. D\'efinissons $w$ comme ci-dessus. On peut supposer que ses deux composantes $w_{+}$ et $w_{-}$ sont respectivement \'egales \`a $\sum_{j=1,...,m}z_{j}w_{j}$ et $\sum_{j=1,...,m}z_{-j}w_{-j}$. Parce que $X$ appartient \`a $\Lambda^S$, $w_{-}$ est un \'el\'ement non nul de $F_{1}$. Posons $h=w_{-}\in H''_{S}(F)$. On peut remplacer $w$ par $hw$. Pour cet \'el\'ement, les coordonn\'ees $z_{-j}$ sont toutes \'egales \`a $1$ et on conclut encore que les autres coordonn\'ees $z_{j}$ sont born\'ees. Dans le cas g\'en\'eral, on peut d\'ecomposer $W''$ en somme directe de sous-espaces et d\'ecomposer conform\'ement $H''_{S}$ en produit de tores de sorte que chaque composante soit de l'un des deux cas particuliers que l'on vient de consid\'erer. On en d\'eduit la propri\'et\'e requise.
 
Supposons  (1) v\'erifi\'ee. En appliquant 2.3(1), on voit que l'on peut choisir l'application $Y\mapsto \gamma_{Y}$ de sorte qu'il existe  $c>0$ tels que
$$(2) \qquad  \sigma(\gamma_{Y})\leq c(1+\vert log\,D^{G''}(Y)\vert )$$
pour tout $Y\in  \mathfrak{t}(F)^S\cap \omega_{T}$.

   \bigskip
 
 \subsection{Calcul de $I_{\kappa''}(\theta'',\varphi)$}
 
 Revenons \`a la situation de 9.1 et supposons que $\varphi$ est \`a support dans $\omega''$. Soit $g\in G''(F)$. D'apr\`es l'hypoth\`ese sur $\theta''$, on a
 $$I(\theta'',\varphi,g)=J_{H''}(S,((^g\varphi)^{\xi})\hat{})=D^{H''}(S)^{1/2}\int_{H''_{S}(F)\backslash H''(F)}((^g\varphi)^{\xi})\hat{}(h^{-1}Sh)dh.$$
 En utilisant le lemme 9.2, on obtient
 $$I(\theta'',\varphi,g)=D^{H''}(S)^{1/2}\int_{H''_{S}(F)\backslash H''(F)}\int_{\Sigma}(^g\varphi)\hat{}(\Xi+h^{-1}Sh+X)dX\, dh,$$
 puis
 $$(1) \qquad I_{\kappa''}(\theta'',\varphi)=D^{H''}(S)^{1/2}\int_{H''(F)U''(F)\backslash G''(F)}\int_{H''_{S}(F)\backslash H''(F)}$$
 $$\qquad \int_{\Sigma}(^g\varphi)\hat{}(\Xi+h^{-1}Sh+X)dX\, dh\kappa''(g)dg.$$
 Remarquons que cette expression est absolument convergente: les trois int\'egrales sont \`a support compact. On transforme cette expression en
  $$I_{\kappa''}(\theta'',\varphi)=D^{H''}(S)^{1/2}\int_{H''(F)U''(F)\backslash G''(F)}\int_{H''_{S}(F)\backslash H''(F)}\int_{\Sigma}(^{hg}\varphi)\hat{}(\Xi+S+X)dX\, dh\kappa''(g)dg$$
 $$\qquad =D^{H''}(S)^{1/2}\int_{H''_{S}(F)U''(F)\backslash G''(F)}\int_{\Sigma}(^g\varphi)\hat{}(\Xi+S+X)dX\kappa''(g)dg.$$
 Le lemme 9.6 nous permet de remplacer l'int\'egrale int\'erieure par
 $$\sum_{T\in {\cal T}(G'')}\vert  W(G'',T)\vert ^{-1}\int_{H''_{S}(F)U''(F)}\int_{\mathfrak{t}(F)^S}(^g\varphi)\hat{}(y^{-1}\gamma_{Y}^{-1}Y\gamma_{Y}y)D^{H''}(S)^{-1/2}D^{G''}(Y)^{1/2}dY dy$$
 $$=\sum_{T\in {\cal T}(G'')}\vert  W(G'',T)\vert ^{-1}\int_{H''_{S}(F)U''(F)}\int_{\mathfrak{t}(F)^S}(^{\gamma_{Y}yg}\varphi)\hat{}(Y)D^{H''}(S)^{-1/2}D^{G''}(Y)^{1/2}dY dy.$$
 Un simple changement de variables conduit alors \`a l'\'egalit\'e
 $$I_{\kappa''}(\theta'',\varphi)=\sum_{T\in {\cal T}(G'')}\vert W(G'',T)\vert ^{-1}\int_{\mathfrak{t}(F)^S}\int_{G''(F)}\hat{\varphi}(g^{-1}Yg)\kappa''(\gamma_{Y}^{-1}g)dgD^{G''}(Y)^{1/2}dY.$$
 Pour $T\in {\cal T}(G'')$ et $Y\in \mathfrak{t}(F)^S$, d\'efinissons une fonction $\kappa''_{Y}$ sur $G''(F)$ par
 $$\kappa''_{Y}(g)=\nu(A_{T})\int_{A_{T}(F)}\kappa''(\gamma_{Y}^{-1}ag)da.$$
 Remarquons que cette expression ne d\'epend pas du choix de l'application $Y\mapsto \gamma_{Y}$: tout autre choix remplace $\gamma_{Y}$ par $\gamma_{Y}y$, avec $y\in H''_{S}(F)U''(F)$, mais $\kappa''$ est invariante \`a gauche par ce groupe. On obtient:
 $$(2) \qquad I_{\kappa''}(\theta'',\varphi)=\sum_{T\in {\cal T}(G'')}\nu(A_{T})^{-1}\vert W(G'',T)\vert ^{-1}$$
 $$\qquad \int_{\mathfrak{t}(F)^S}\int_{A_{T}(F)\backslash G''(F)}\hat{\varphi}(g^{-1}Yg)\kappa''_{Y}(g)dgD^{G''}(Y)^{1/2}dY.$$
 Les transformations que l'on a effectu\'ees sont justifi\'ees par la convergence absolue de l'expression (1) de d\'epart. 
 
 \bigskip
 
 \section{Calcul de la limite $lim_{N\to \infty}I_{x,\omega,N}(\theta,f)$}
 
 \bigskip
 
 \subsection{Convergence d'une premi\`ere expression}
 
 On se place dans la situation de 8.2. D'apr\`es la proposition 6.4 et le lemme 6.3(i), on peut fixer une famille finie $(Y_{i})_{i=1,...,n}$ d'\'el\'ements de $\mathfrak{h}_{x,reg}(F)$ et une famille finie $(c_{i})_{i=1,...,n}$ de nombres complexes de sorte que
 $$\theta_{x,\omega}(X)=\sum_{i=1,...,n}c_{i}\hat{j}^{H_{x}}(Y_{i},X)$$
 pour tout $X\in \omega\cap \mathfrak{h}_{x,reg}(F)$. Formulons cette propri\'et\'e diff\'eremment, en utilisant les notations introduites en 5.4: pour $X\in \mathfrak{g}_{x}(F)$, on note $X=X'+X''$ la d\'ecomposition de $X$ en somme d'un \'el\'ement $X'\in \mathfrak{g}'_{x}(F)$ et d'un \'el\'ement $X''\in \mathfrak{g}''(F)$. Il existe alors une famille finie ${\cal S}$ d'\'el\'ements de $\mathfrak{h}''_{reg}(F)$ et une famille finie $(\hat{j}_{S})_{S\in {\cal S}}$ de fonctions d\'efinies presque partout sur $\mathfrak{h}'_{x}(F)$ de sorte que
 $$\theta_{x,\omega}(X)=\sum_{S\in {\cal S}}\hat{j}_{S}(X')\hat{j}^{H''}(S,X'')$$
 pour tout $X\in \omega\cap \mathfrak{h}_{x,reg}(F)$. Les \'el\'ements $S$ sont les diff\'erentes projections $Y''_{i}$. La preuve du lemme 6.3(i) nous autorise \`a remplacer les $Y_{i}$ par des \'el\'ements assez voisins. On peut donc supposer que le noyau de chaque $S$ agissant dans $W''$ est de dimension au plus $1$. Les fonctions $\hat{j}_{S}$ sont combinaisons lin\'eaires de fonctions $X'\mapsto \hat{j}^{H'_{x}}(Y'_{i},X')$ et h\'eritent donc de leurs propri\'et\'es.
 
 Rappelons que, par construction, on a l'\'egalit\'e $H'_{x}=G'_{x}$. Pour $g\in G(F)$, on a
 $$I_{x,\omega}(\theta,f,g)=\int_{\mathfrak{g}'_{x}(F)\times\mathfrak{h}''(F)}\theta_{x,\omega}(X){^gf}^{\xi}_{x,\omega}(X)dX.$$
 En utilisant la formule de Weyl pour l'int\'egrale sur $\mathfrak{g}'_{x}(F)$, on obtient
 $$I_{x,\omega}(\theta,f,g)=\sum_{S\in{\cal S}}\sum_{T'\in {\cal T}(G'_{x})}\vert W(G'_{x},T')\vert ^{-1}\int_{\mathfrak{t}'(F)} \hat{j}_{S}(X')D^{G_{x}'}(X')$$
 $$\qquad \int_{T'(F)\backslash G'_{x}(F)} \int_{\mathfrak{h}''(F)}\hat{j}^{H''}(S,X''){^gf}^{\xi}_{x,\omega}(g^{_{'}-1}X'g'+X'')dX''\, dg'\, dX'.$$
 D'o\`u
 $$I_{x,\omega,N}(\theta,f)=\sum_{S\in{\cal S}}\sum_{T'\in {\cal T}(G'_{x})}\vert W(G'_{x},T')\vert ^{-1}\int_{\mathfrak{t}'(F)}  \hat{j}_{S}(X')D^{G_{x}'}(X')$$
 $$\qquad \int_{T'(F)H''(F)U_{x}(F)\backslash G(F)}\int_{\mathfrak{h}''(F)}\hat{j}^{H''}(S,X''){^gf}^{\xi}_{x,\omega}(X'+X'')dX''\kappa_{N}(g) dg\, dX'.$$
 On peut \'ecrire les deux derni\`eres int\'egrales ci-dessus sous la forme
$$ \int_{T'(F)G''(F)\backslash G(F)}\int_{H''(F)U_{x}(F)\backslash G''(F)}\int_{\mathfrak{h}''(F)}\hat{j}^{H''}(S,X''){^{g''g}f}^{\xi}_{x,\omega}(X'+X'')dX''\kappa_{N}(g''g)dg''\, dg.$$
Les deux int\'egrales int\'erieures sont \'egales \`a $I_{\kappa''}(\theta'',\varphi)$, o\`u $\theta''(X'')=\hat{j}^{H''}(S,X'')$, $\varphi(X'')={^gf}_{x,\omega}(X'+X'')$ et $\kappa''(g'')=\kappa_{N}(g''g)$. Utilisons la formule 9.8(2) qui calcule cette expression. La fonction $\hat{\varphi}$ qui y intervient est \'egale \`a ${^gf}_{x,\omega}^{\sharp}$, avec la notation de 5.4. Quelques remises en ordre conduisent alors \`a l'\'egalit\'e
$$(1) \qquad I_{x,\omega,N}(\theta,f)=\sum_{S\in{\cal S}}\sum_{T\in {\cal T}(G_{x})}\nu(A_{T''})^{-1}\vert W(G_{x},T)\vert ^{-1}\int_{\mathfrak{t}'(F)\times \mathfrak{t}''(F)^{S}}\hat{j}_{S}(X')D^{G_{x}'}(X')D^{G''}(X'')^{1/2}$$
$$\qquad \int_{T'(F)A_{T''}(F)\backslash G(F)}{^gf}_{x,\omega}^{\sharp}(X'+X'')\kappa_{N,X''}(g)dg\,dX''\,dX',$$
o\`u 
$$\kappa_{N,X''}(g)=\nu(A_{T''})\int_{A_{T''}(F)}\kappa_{N}(\gamma_{X''}^{-1}ag)da.$$
Ces manipulations formelles sont justifi\'ees par le lemme ci-dessous. Pour tout $S\in{\cal S}$ et tout $T\in {\cal T}(G_{x})$, fixons une famille finie ${\cal Q}_{S,T}$ de polyn\^omes non nuls sur $\mathfrak{t}(F)$. Pour tout $\epsilon>0$, notons $\mathfrak{t}(F)[S;\leq \epsilon]$ l'ensemble des $X\in  \mathfrak{t}(F)$ pour lesquels il existe $Q\in {\cal Q}_{S,T}$ tel que $\vert Q(X)\vert_{F} \leq \epsilon$, et notons $ \mathfrak{t}(F)[S;> \epsilon]$ l'ensemble des $X\in  \mathfrak{t}(F)$ pour lesquels $\vert Q(X)\vert _{F }>\epsilon$ pour tout $Q\in {\cal Q}_{S,T}$.
Notons $I_{N,\leq\epsilon}$, resp. $I_{N,>\epsilon}$, l'expression obtenue \`a partir de l'expression (1) en rempla\c{c}ant les int\'egrales sur $\mathfrak{t}'(F)\times \mathfrak{t}''(F)^{S}$ par les int\'egrales sur     $(\mathfrak{t}'(F)\times \mathfrak{t}''(F)^{S})\cap \mathfrak{t}(F)[S;\leq \epsilon]$, resp. $( \mathfrak{t}'(F)\times \mathfrak{t}''(F)^{S})\cap \mathfrak{t}(F)[S;> \epsilon]$. On a \'evidemment l'\'egalit\'e
$$I_{x,\omega,N}(\theta,f)=I_{N,\leq\epsilon}+I_{N,>\epsilon}.$$
Notons enfin $\vert I\vert _{x,\omega,N}(\theta,f)$ et $\vert I\vert _{N,\leq\epsilon}$ les expressions obtenues en rempla\c{c}ant dans $I_{x,\omega,N}(\theta,f)$ (ou plus exactement dans l'expression (1)) et $I_{N,\leq\epsilon}$ toutes les fonctions par leurs valeurs absolues.
 
 \ass{Lemme}{(i) Il existe $k\in {\mathbb N}$ et $c>0$ tel que $\vert I\vert _{x,\omega,N}(\theta,f)\leq cN^k$ pour tout $N\geq 1$.
 
 (ii) Il existe un entier $b\geq 1$ et $c>0$ tel que $\vert I\vert _{N,\leq N^{-b}}\leq cN^{-1}$ pour tout $N\geq 1$.}
 
 Preuve. Soit $S\in{\cal S}$. Notons $(\pm s_{j})_{j=1,...,m}$ les valeurs propres non nulles de l'action de $S$ sur $W''$. Pour $X''\in \mathfrak{g}''(F)$, posons 
 $$Q_{S}(X'')=\left\lbrace\begin{array}{cc}\prod_{j=1,...,m}s_{j}^{-1}P_{X''}(s_{j}),&\,\,{\rm si\,\,}d{\rm\,\,est\,\,impair}\\ \prod_{j=1,...,m}P_{X''}(s_{j}),&\,\,{\rm si\,\,}d{\rm\,\,est\,\,pair}.\\ \end{array}\right.$$
  Certainement, $Q_{S}$ est un polyn\^ome non nul sur l'alg\`ebre de Lie de tout sous-tore maximal de $G''$. Soient $T\in {\cal T}(G_{x})$ et $\omega_{T''}$ un sous-ensemble compact de $\mathfrak{t}''(F)$. On va montrer
 
 (2) il existe un entier $k\in {\mathbb N}$ et $c>0$ tels que
 $$\kappa_{N,X''}(g)\leq cN^k\sigma(g)^k(1+\vert log\vert Q_{S}(X'')\vert _{F}\vert )^k(1+\vert log\,D^{G''}(X'')\vert )^k$$
 pour tout $X''\in \mathfrak{t}''(F)^S\cap\omega_{T''}$, tout $g\in G(F)$ et tout $N\geq 1$.
 
 Commen\c{c}ons par d\'eduire l'\'enonc\'e de (2). On peut fixer $S\in{\cal S}$ et $T\in {\cal T}(G_{x})$ et consid\'erer l'int\'egrale
$$ \int_{\mathfrak{t}'(F)\times \mathfrak{t}''(F)^{S}}\vert \hat{j}_{S}(X')\vert D^{G_{x}'}(X')D^{G''}(X'')^{1/2} $$
$$\qquad \int_{T'(F)A_{T''}(F)\backslash G(F)}\vert {^gf}_{x,\omega}^{\sharp}(X'+X'')\vert \kappa_{N,X''}(g)dg\,dX''\,dX'.$$
 Introduisons une notation impr\'ecise mais commode. Soient deux nombres $a$ et $b$ d\'ependant de variables, ici $N$, $g$, $X'$ et $X''$. On \'ecrit $a<<b$ pour dire qu'il existe $c>0$ tel que, quelles que soient ces variables, on ait $a\leq cb$. D'apr\`es la d\'efinition de $\hat{j}_{S}$ et un r\'esultat de Harish-Chandra ([HCvD] th\'eor\`eme 13), on a $\vert \hat{j}_{S}(X')\vert <<D^{G'_{x}}(X')^{-1/2}$. D'apr\`es 3.1(5), on peut fixer un sous-ensemble compact $\Gamma\subset G(F)$ tel que ${^gf}_{x,\omega}^{\sharp}=0$ si $g\not\in G_{x}(F)\Gamma$. On peut donc fixer $\gamma\in \Gamma$ et remplacer l'int\'egrale sur $T'(F)A_{T''}(F)\backslash G(F)$ par celle sur $T'(F)A_{T''}(F)\backslash G_{x}(F)\gamma$. On peut majorer $\vert {^{\gamma}f}_{x,\omega}^{\sharp}\vert $ par une combinaison lin\'eaire de fonction $f'\otimes f''$ o\`u $f'\in C_{c}^{\infty}(\mathfrak{g}'_{x}(F))$, $f''\in C_{c}^{\infty}(\mathfrak{g}''(F))$ et $f'$ et $f''$ sont \`a valeurs positives ou nulles. On est ramen\'e \`a majorer
$$ \int_{\mathfrak{t}'(F)\times \mathfrak{t}''(F)^{S}}D^{G_{x}'}(X')^{1/2}D^{G''}(X'')^{1/2}\int_{T'(F)\backslash G'_{x}(F)}$$
$$\qquad \int_{A_{T''}(F)\backslash G''(F)}f'(g^{_{'}-1}X'g')f''(g^{_{''}-1}X''g'')\kappa_{N,X''}(g'g''\gamma)dg''\,dg'\,dX''\,dX'.$$
On peut fixer un sous-ensemble compact $\omega_{T''}\subset \mathfrak{t}''(F)$ tel que, pour tout $g''$, la fonction  $X''\mapsto f''(g^{_{''}-1}X''g'')$  sur $\mathfrak{t}''(F)$ soit \`a support dans $\omega_{T''}$. Gr\^ace \`a 2.3(1), on peut supposer que les $g''$ intervenant dans l'int\'egrale v\'erifient $\sigma(g'')<<1+\vert log\,D^{G''}(X'')\vert $. Puisque $G'_{x}=H'_{x}\subset H$, on a $\kappa_{N,X''}(g'g''\gamma)=\kappa_{N,X''}(g''\gamma)$. En appliquant (2), on obtient $\kappa_{N,X''}(g'g''\gamma)<<N^k\varphi(X'')$
o\`u 
$$\varphi(X'')=(1+\vert log\vert Q_{S}(X'')\vert _{F}\vert )^k(1+\vert log\,D^{G''}(X'')\vert )^{2k}.$$
L'expression \`a majorer devient
$$N^k \int_{\mathfrak{t}'(F)\times \mathfrak{t}''(F)^{S}}D^{G_{x}'}(X')^{1/2}D^{G''}(X'')^{1/2}\int_{T'(F)\backslash G'_{x}(F)}$$
$$\qquad \int_{A_{T''}(F)\backslash G''(F)}f'(g^{_{'}-1}X'g')f''(g^{_{''}-1}X''g'') \varphi(X'')dg''\,dg'\,dX''\,dX'.$$
Elle est major\'ee par
$$N^k\int_{\mathfrak{t}(F)}J_{G_{x}}(X'+X'',f'\otimes f'')\varphi(X'')dX''\,dX'.$$
D'apr\`es Harish-Chandra, l'int\'egrale orbitale est born\'ee. Elle est aussi \`a support compact, ce qui nous conduit \`a majorer
$$N^k\int_{\omega_{T}}\varphi(X'')dX''\,dX',$$
o\`u $\omega_{T}$ est un sous-ensemble compact de $\mathfrak{t}(F)$. Le lemme 2.4 nous dit que l'int\'egrale est convergente, ce qui entra\^{\i}ne la majoration du (i) de l'\'enonc\'e. Pour le (ii), on est de m\^eme conduit \`a majorer
$$N^k\int_{\omega_{T}\cap \mathfrak{t}(F)[S;\leq N^{-b}]} \varphi(X'')dX''\,dX'.$$
D'apr\`es l'in\'egalit\'e de Schwartz, il  suffit de majorer
$$N^k(\int_{\omega_{T}\cap \mathfrak{t}(F)[S;\leq N^{-b}]}dX)^{1/2}(\int_{\omega_{T}}\varphi(X'')^2dX''\,dX')^{1/2}.$$
Le dernier terme se majore comme ci-dessus. Pour $\epsilon>0$, on a
$$\int_{\omega_{T}\cap \mathfrak{t}(F)[S;\leq\epsilon]}dX\leq \sum_{Q\in {\cal Q}_{S,T}}mes(\{X\in \omega_{T}; \vert Q(X)\vert _{F}\leq\epsilon\}).$$
D'apr\`es [A5], lemme 7.1, il existe un r\'eel $r>0$ tel que chacun de ces termes soit $<< \epsilon^r$. On en d\'eduit une majoration
$$\vert I\vert _{N,\leq N^{-b}}<<N^{k-rb/2}.$$
En prenant $b>2(k+1)/r$, on obtient le (ii) de l'\'enonc\'e.

 Prouvons (2). Rempla\c{c}ons $V$ par $V''$ dans les d\'efinitions de 7.2. On fixe un r\'eseau sp\'ecial $R''$ de $V''$ de m\^eme que l'on a fix\'e $R$, on note $K''$ son stabilisateur dans $G''(F)$ et on d\'efinit une fonction $\kappa''_{N}$ sur $G''(F)$. Posons 
 $$\kappa''_{N,X''}(1)=\nu(A_{T''})\int_{A_{T''}(F)}\kappa_{N}(\gamma_{X''}^{-1}a)da.$$
 On va montrer
 
 (3) il existe un entier $k\in {\mathbb N}$ et $c>0$ tels que
 $$\kappa''_{N,X''}(1)\leq cN^k(1+\vert log(\vert Q_{S}(X'')\vert _{F})\vert )^k(1+\vert log\,D^{G''}(X'')\vert )^k$$
 pour tout $X''\in \mathfrak{t}''(F)^S\cap\omega_{T''}$ et tout $N\geq 1$.
 
 D\'eduisons d'abord (2) de (3). Pour un r\'eel $r>0$, posons $\kappa''_{r}=\kappa''_{N(r)}$, o\`u $N(r)$ est le plus petit entier sup\'erieur ou \'egal \`a $r$. On a
 
 (4) il existe $c>0$ tel que $\kappa_{N}(g''g)\leq \kappa''_{N+c\sigma(g)}(g'')$ pour tous $g\in G(F)$, $g''\in G''(F)$.
 
  Ecrivons $g''=m''u''k''$, avec $m''\in M''(F)$, $u''\in U''(F)$, $k''\in K''$, puis $k''g=muk$, avec $m\in M(F)$, $u\in U(F)$, $k\in K$. On a $\kappa_{N}(g''g)=\kappa_{N}(m''m)$. Supposons ce terme non nul (donc \'egal \`a $1$), d\'ecomposons $m''$ et $m$ en $m''=a''g''_{0}$ et $m=ag_{0}$, o\`u $a'',a\in A(F)$, $g''_{0}\in G''_{0}(F)$ et $g_{0}\in G_{0}(F)$. Alors $\vert val_{F}(a''_{i}a_{i})\vert \leq N$ pour tout $i=1,...,r$ et $g_{0}^{-1}g_{0}^{_{''}-1}v_{0}\in \varpi_{F}^{-N}R_{0}$.  On a $\sigma(m)<< \sigma(g)$. Donc  $\vert val_{F}(a_{i})\vert <<\sigma(g)$ pour tout $i=1,...,r$ et $\sigma(g_{0})<<\sigma(g)$. On en d\'eduit d'abord qu'il existe $c_{1}>0$ tel que  $\vert val_{F}(a''_{i})\vert \leq N+c_{1}\sigma(g)$ pour tout $i=1,...,r$. Il existe $c_{2}>0$ tel que $g_{0}R_{0}\subset \varpi_{F}^{-N(c_{2}\sigma(g_{0}))}R_{0}$. Il existe $c_{3}\in {\mathbb N}$ tel que $R_{0}\cap V''\subset \varpi_{F}^{-c_{3}}R''_{0}$. Alors $g_{0}^{_{''}-1}v_{0}\in \varpi_{F}^{-N'}R''_{0}$, o\`u $N'\leq N+c_{4}\sigma(g)$, pour $c_{4}>0$ convenable. En prenant $c>c_{1},c_{4}$, on voit que $g''$ v\'erifie les conditions requises pour que $\kappa''_{N+c\sigma(g)}(g'')=1$. Cela prouve (4).
 
 En utilisant (4), on  a  
 $$\kappa_{N,X''}(g)=\nu(A_{T''})\int_{A(T'')(F)}\kappa_{N}(\gamma_{X''}^{-1}ag)da$$
 $$\qquad \leq\nu(A_{T''})\int_{A(T'')(F)}\kappa''_{N+c\sigma(g)}(\gamma_{X''}^{-1}a)da\leq \kappa''_{N+c\sigma(g),X''}(1).$$
 La majoration (3) entra\^{\i}ne alors (2).

 Prouvons maintenant (3).  On suppose v\'erifi\'ees les conditions (1) et (2) de 9.7 pour le compact $\omega_{T''}$. Soit $a\in A_{T''}(F)$ tel que $\kappa''_{N}(\gamma_{X''}^{-1}a)=1$. Gr\^ace \`a 7.2(1), on peut \'ecrire $\gamma_{X''}^{-1}a=vhy$, avec $v\in U''(F)$, $h\in H''(F)$, $y\in G''(F)$ et $\sigma(y)\leq cN$. On a 
 $$(5) \qquad yX''y^{-1}=h^{-1}v^{-1}X''_{\Sigma}vh,$$
  o\`u $X''_{\Sigma}=\gamma_{X''}^{-1}X''\gamma_{X''}$.  La condition impos\'ee \`a $y$ implique que $\sigma(yX''y^{-1})<<N$ (cf. 1.1 pour la d\'efinition de la fonction $\sigma$ sur $\mathfrak{g}''(F)$). On a $h^{-1}v^{-1}X''_{\Sigma}vh\in \Xi+h^{-1}Sh+\Sigma$ et $h^{-1}Sh\in \mathfrak{h}''(F)$. Or $\mathfrak{h}''(F)$ et $\Sigma$ sont en somme directe. Alors (5) entra\^{\i}ne   que $\sigma(h^{-1}Sh)<<N$. Donc il existe un entier $k>0$ tel que $\sigma(h^{-1}\varpi_{F}^{kN}Sh)<< 1$. En appliquant 2.3(1), on peut \'ecrire $h=tz$, avec $t\in H''_{S}(F)$, $z\in H''(F)$ et 
  $$\sigma(z)<<(1+\vert log\,D^{H''}(\varpi_{F}^{kN}S)\vert )<<N.$$
  On peut r\'ecrire $vhy=tug$, avec $u\in U''(F)$ et $g=zy$, donc $\sigma(g)<<N$. L'\'egalit\'e (5) se r\'ecrit $gX''g^{-1}=u^{-1}Yu$, o\`u $Y=t^{-1}X''_{\Sigma}t$. On a $\sigma(gX''g^{-1})<<N$. D'apr\`es la condition (1) de 9.7, $Y$ appartient \`a $\Xi+S+\Lambda$. Le lemme 9.3 nous dit que $u$ et $Y$ d\'ependent alg\'ebriquement de $u^{-1}Yu$. Donc  $\sigma(u)<<N$ et $\sigma(Y)<<N$. Posons $X''_{\Sigma}=\Xi+S+X$,  $Y=\Xi+S+X^*$, introduisons les coordonn\'ees de $X$ et $X^*$ comme en 9.4 (on affecte celles de $X^*$ d'un exposant $^*$) et, pour tout $j=1,...,m$, notons $t_{j}$ la valeur propre de $t$ sur $w_{j}$. On a $z_{j}^*=t_{j}^{-1}z_{j}$ et $z_{-j}^*=t_{j}z_{-j}$ pour tout $j=1,...,m$. La condition (1) de 9.7 et celle ci-dessus portant sur $Y$ nous disent qu'il existe $c>0$ tel que
  $$val_{F}(z_{j}^*)\geq -cN,\,\,val_{F}(z_{-j}^*)\geq -cN,\,\,val_{F}(z_{j})\geq -c,\,\,val_{F}(z_{-j})\geq -c$$
  pour tout $j$. On en d\'eduit
  $$\vert val_{F}(t_{j})\vert \leq c(N+1)+val_{F}(z_{j}z_{-j})\leq c(N+2m-1)+val_{F}(\prod_{j'=1,...,m}z_{j'}z_{-j'}).$$
  La formule 9.4(1) montre qu'il existe $c'>0$ tel que le dernier terme soit major\'e par $c'(1+\vert log\vert Q_{S}(X'')\vert _{F}\vert )$. On en d\'eduit $\sigma(t)<< N+\vert log\vert Q_{S}(X'')\vert _{F}\vert $. Alors 
  $$\sigma(\gamma_{X''}^{-1}a)=\sigma(tug)<<N+\vert log\vert Q_{S}(X'')\vert _{F}\vert .$$
   En appliquant 9.7(2), on en d\'eduit 
   $$\sigma(a)<<N+\vert log(\vert Q_{S}(X'')\vert _{F})\vert +\vert log\,D^{G''}(X'')\vert .$$
  Le terme $\kappa''_{N,X''}(1) $ est born\'e par la mesure de l'ensemble des $a$ v\'erifiant cette condition.  Il est facile de montrer que, pour tout r\'eel $r\geq 1$,
   $$mes(\{a\in A_{T''}(F); \sigma(a)\leq r\})<<r^k,$$
   o\`u $k=dim(A_{T''})$. On en d\'eduit
   $$\kappa''_{N,X''}(1)<<N^k(1+\vert log\vert Q_{S}(X'')\vert _{F}\vert)^k(1+\vert log\,D^{G''}(X'')\vert )^k,$$
   ce qui prouve (3) et ach\`eve la d\'emonstration. $\square$

   Pour tout $S\in{\cal S}$ et tout $T\in {\cal T}(G_{x})$, on note ${\cal Q}_{S,T}$ la famille des trois polyn\^omes sur $\mathfrak{t}(F)$ suivants
   $$X\mapsto det(ad(X')_{\vert \mathfrak{g}'_{x}/\mathfrak{t}'}),\,\,X\mapsto det(ad(X'')_{\vert \mathfrak{g}''/\mathfrak{t}''}),\,\,X''\mapsto Q_{S}(X'')$$
   o\`u $Q_{S}$ a \'et\'e d\'efini ci-dessus.  Appliqu\'e \`a ces donn\'ees, le (ii) du lemme  nous fournit un entier $b$ que l'on fixe. On pose
   $$I^*_{x,\omega,N}(\theta,f)=I_{N,>N^{-b}}.$$
   Le lemme entra\^{\i}ne que
   $$lim_{N\to \infty}(I_{x,\omega,N}(\theta,f)-I^*_{x,\omega,N}(\theta,f))=0.$$
 
   \bigskip
   
   \subsection{Commutant d'un tore}

On fixe  $T\in {\cal T}(G_{x})$. Notons $M_{\natural}$ le commutant de $A_{T''}$ dans $G$. C'est un L\'evi de $G$, qui contient $G'$. Notons $V_{2}''$ l'intersection des noyaux des actions de $a$ dans $V''$, pour $a\in A_{T''}$.  

\ass{Lemme}{On a l'\'egalit\'e $A_{T''}=A_{M_{\natural}}$ sauf dans le cas o\`u les  conditions suivantes sont satisfaites: $d$ est pair, $W'$ est hyperbolique et de dimension $2$, $V_{2}''=\{0\}$. Dans ce cas, on a $A_{M_{\natural}}=T'A_{T''}$.}

Preuve.   Posons $V_{2}=W'\oplus V_{2}''$, notons $V_{1}$ son orthogonal dans $V$ et $G_{1}$, resp. $G_{2}$ les groupes sp\'eciaux orthogonaux de $V_{1}$, resp. $V_{2}$. L'espace $V_{2}$ est l'intersection des noyaux des actions de $a$ dans $V$, pour $a\in A_{T''}$. Donc $M_{\natural}$ conserve $V_{2}$. Par cons\'equent, $M_{\natural}$ conserve aussi $V_{1}$, donc $M_{\natural}\subset G_{2}\times G_{1}$. Puisque $A_{T''}$ agit trivialement dans $V_{2}$, on a $A_{T''}\subset G_{1}$. Donc $M_{\natural}=G_{2}\times M_{1,\natural}$, o\`u $M_{1,\natural}$ est le commutant de $A_{T''}$ dans $G_{1}$, puis $A_{M_{\natural}}=A_{G_{2}}\times A_{M_{1,\natural}}$.  On a $V_{1}\subset V''$, donc $G_{1}\subset G''$. D'autre part $T''$ commute \`a $A_{T''}$, donc $T''\subset M_{\natural}$, donc $A_{M_{\natural}}$ commute \`a $T''$. Alors $A_{M_{1,\natural}}$ est contenu dans le commutant de $T''$ dans $G''$. Puisque $T''$ est un sous-tore maximal de $G''$, ce commutant est \'egal \`a $T''$. Donc $A_{M_{1,\natural}}\subset T''$, ce qui entra\^{\i}ne $A_{M_{1,\natural}}\subset A_{T''}$. L'inclusion oppos\'ee est imm\'ediate puisque $A_{T''}$ est \'evidemment un tore d\'eploy\'e central dans $M_{1,\natural}$. Donc $A_{M_{1,\natural}}=A_{T''}$. On a $A_{G_{2}}=\{1\}$ sauf dans le cas o\`u $V_{2}$ est hyperbolique de dimension $2$. Supposons cette condition v\'erifi\'ee. Puisque $G_{1}$ contient un sous-tore $A_{T''}$ qui agit sans point fixe non nul dans $V$, $dim(V_{1})$ est paire et $d$ aussi. Si $W'=\{0\}$, on a $G=G''$ et $T''$ est un tore maximal de $G$. Le m\^eme raisonnement que ci-dessus montre que $A_{M_{\natural}}\subset A_{T''}$ contrairement \`a l'hypoth\`ese $A_{G_{2}}\not=\{1\}$. Donc $W'\not=\{0\}$. Puisque $dim(W')$ est paire et $V_{2}=W'\oplus V_{2}''$, on a $W'=V_{2}$ et $ V_{2}''=\{0\}$. Inversement, si $W'$ est hyperbolique de dimension $2$ et $V_{2}''=\{0\}$, on a $G_{2}=G'$ et ce groupe est un tore d\'eploy\'e. Puisque $T'$ est un sous-tore maximal de $G'$, on a $T'=G'=A_{G_{2}}$ et la conclusion du lemme s'ensuit. $\square$

\bigskip

\subsection{D\'efinitions combinatoires}

Appelons cas exceptionnel celui de l'\'enonc\'e pr\'ec\'edent. Supposons tout d'abord que l'on n'est pas dans ce cas et rappelons quelques d\'efinitions d'Arthur. Soit ${\cal Y}=(Y_{P_{\natural}})_{P_{\natural}\in {\cal P}(M_{\natural})}$ une famille d'\'el\'ements de ${\cal A}_{M_{\natural}}$, $(G,M_{\natural})$-orthogonale et positive. Pour $Q=LU_{Q}\in {\cal F}(M_{\natural})$, on note $\zeta\mapsto\sigma_{M_{\natural}}^Q(\zeta,{\cal Y})$ la fonction caract\'eristique  dans ${\cal A}_{M_{\natural}}$ de la somme de ${\cal A}_{L}$ et de l'enveloppe convexe de la famille $(Y_{P_{\natural}})_{P_{\natural}\in {\cal P}(M_{\natural}); P_{\natural}\subset Q}$.   On note $\tau_{Q}$ la fonction caract\'eristique dans ${\cal A}_{M_{\natural}}$ de la somme   ${\cal A}_{M_{\natural}}^L+{\cal A}_{Q}^+$. Rappelons que ${\cal A}^+_{Q}$ est la chambre positive ouverte de ${\cal A}_{L}$ associ\'ee \`a $Q$. On a

(1) la fonction
$$ \zeta\mapsto\sigma_{M_{\natural}}^Q(\zeta,{\cal Y})\tau_{Q}(\zeta-Y_{Q})$$
sur ${\cal A}_{M_{\natural}}$ est la fonction caract\'eristique de la somme de ${\cal A}_{Q}^+$ et de l'enveloppe convexe de la famille $(Y_{P_{\natural}})_{P_{\natural}\in {\cal P}(M_{\natural}); P_{\natural}\subset Q}$;  

$$ (2) \qquad \sum_{Q\in {\cal F}(M_{\natural})} \sigma_{M_{\natural}}^Q(\zeta,{\cal Y})\tau_{Q}(\zeta-Y_{Q})=1$$
pour tout $\zeta\in {\cal A}_{M_{\natural}}$.

L'assertion (1) est imm\'ediate et (2) est l'assertion 3.9 de [A3].

Consid\'erons maintenant le cas exceptionnel. Alors ${\cal A}_{M_{\natural}}={\cal A}_{T'}\oplus {\cal A}_{T''}$ et ${\cal A}_{T'}$ est une droite. Conform\'ement \`a cette d\'ecomposition, on d\'efinit la projection $\zeta\mapsto \zeta_{T''}$ de ${\cal A}_{M_{\natural}}$ sur ${\cal A}_{T''}$. On note $\tilde{{\cal F}}(M_{\natural})$ l'ensemble des $Q\in {\cal F}(M_{\natural})$ tels que ${\cal A}_{Q}^+\cap {\cal A}_{T''}\not=\emptyset$. On note $\tilde{{\cal P}}(M_{\natural})={\cal P}(M_{\natural})\cap\tilde{{\cal F}}(M_{\natural})$. Pour $Q=LU_{Q}\in \tilde{{\cal F}}(M_{\natural})$, on note $\zeta\mapsto\tilde{\sigma}_{M_{\natural}}^Q(\zeta,{\cal Y})$ la fonction caract\'eristique dans ${\cal A}_{T''}$ de la somme de ${\cal A}_{L}\cap {\cal A}_{T''}$ et de l'enveloppe convexe de la famille $(Y_{P_{\natural},T''})_{P_{\natural}\in \tilde{{\cal P}}(M_{\natural}); P_{\natural}\subset Q}$.  On note $\tilde{\tau}_{Q}$ la fonction caract\'eristique dans ${\cal A}_{T''}$ de la somme   $({\cal A}_{T''}\cap{\cal A}_{M_{\natural}}^L)+({\cal A}_{T''}\cap{\cal A}_{Q}^+)$. On a encore

(3) la fonction
$$ \zeta\mapsto\tilde{\sigma}_{M_{\natural}}^Q(\zeta,{\cal Y})\tilde{\tau}_{Q}(\zeta-Y_{Q,T''})$$
sur ${\cal A}_{T''}$ est la fonction caract\'eristique de la somme de ${\cal A}_{T''}\cap {\cal A}_{Q}^+$ et de l'enveloppe convexe de la famille $(Y_{P_{\natural},T''})_{P_{\natural}\in \tilde{{\cal P}}(M_{\natural}); P_{\natural}\subset Q}$;

$$(4)\qquad \sum_{Q\in \tilde{\cal F}(M_{\natural})} \tilde{\sigma}_{M_{\natural}}^Q(\zeta,{\cal Y})\tilde{\tau}_{Q}(\zeta-Y_{Q,T''})=1$$
pour tout $\zeta\in {\cal A}_{T''}$.

Pour d\'emontrer ces propri\'et\'es et \`a des fins ult\'erieures, introduisons un espace $\tilde{V}$ somme directe de $V''$ et d'une droite $\tilde{D}$, muni de la somme orthogonale $\tilde{q}$  de la restriction de $q$ \`a $V''$ et d'une forme quadratique non d\'eg\'en\'er\'ee sur $\tilde{D}$. On note $\tilde{G}$ son groupe sp\'ecial orthogonal. Le tore $T''$ est inclus dans $G''$, donc dans $\tilde{G}$. Notons $\tilde{M}$ son commutant dans $\tilde{G}$. Le lemme 10.2 s'applique \`a $\tilde{G}$ et montre que $A_{T''}=A_{\tilde{M}}$.

\ass{Lemme}{(i) Il existe une unique bijection $Q=LU_{Q}\mapsto \tilde{Q}=\tilde{L}U_{\tilde{Q}}$ de $\tilde{{\cal F}}(M_{\natural})$ sur ${\cal F}(\tilde{M})$ telle que, pour tout $Q\in \tilde{{\cal F}}(M_{\natural})$, on ait ${\cal A}_{T''}\cap {\cal A}_{M_{\natural}}^L={\cal A}_{\tilde{M}}^{\tilde{L}}$, ${\cal A}_{T''}\cap {\cal A}_{L}={\cal A}_{\tilde{L}}$ et ${\cal A}_{T''}\cap {\cal A}_{Q}^+={\cal A}_{\tilde{Q}}^+$. Cette bijection conserve la relation d'inclusion et envoie $\tilde{{\cal P}}(M_{\natural})$ sur ${\cal P}(\tilde{M})$.

(ii) Pour $Q\in \tilde{{\cal F}}(M_{\natural})$, posons $Y_{\tilde{Q}}=Y_{Q,T''}$. La famille $\tilde{{\cal Y}}=(Y_{\tilde{P}})_{\tilde{P}\in {\cal P}(\tilde{M})}$ est $(\tilde{G},\tilde{M})$-orthogonale et positive. Pour tout $\tilde{Q}\in {\cal F}(\tilde{M})$, $Y_{\tilde{Q}}$ est associ\'e \`a cette famille comme en 2.1.

(iii) Pour tout $Q\in \tilde{{\cal F}}(M_{\natural})$ et tout $\zeta\in {\cal A}_{T''}={\cal A}_{\tilde{M}}$, on a les \'egalit\'es $\tilde{\sigma}_{M_{\natural}}^Q(\zeta,{\cal Y})=\sigma_{\tilde{M}}^{\tilde{Q}}(\zeta,\tilde{{\cal Y}})$ et $\tilde{\tau}_{Q}(\zeta)=\tau_{\tilde{Q}}(\zeta)$.}

Preuve. Les conditions impos\'ees \`a $A_{T''}$ impliquent que $V''$ est hyperbolique et qu'il existe un syst\`eme hyperbolique maximal $(e_{k})_{k=\pm 2,...,\pm d/2}$ dans $V''$ et une suite d'entiers $(d_{i})_{i=2,...,l}$ v\'erifiant les propri\'et\'es suivantes. On a $d_{i}\geq 1$ pour tout $i$ et $1+d_{2}+...+d_{l}=d/2$. Pour $\epsilon=\pm 1$ et $i=2,...,l$, notons $E_{\epsilon i}$ le sous-espace de $V''$ engendr\'e par les $e_{\epsilon k}$ pour $1+d_{2}+...+d_{i-1}<k\leq 1+d_{2}+...+d_{i}$. Alors $A_{T''}$ est le sous-groupe des \'el\'ements de $G''$ qui conservent chaque $E_{\pm i}$ et y agissent par homoth\'etie. Fixons une base hyperbolique $\{e_{1},e_{-1}\}$ de $W'$. Notons $E_{\pm 1}$ la droite port\'ee par $e_{\pm 1}$. Posons $I=\{\pm 1,...,\pm l\}$. Alors $A_{M_{\natural}}$ est le sous-groupe des \'el\'ements de $G$ qui conservent chaque $E_{i}$ pour $i\in I$ et y agissent par homoth\'etie. L'espace ${\cal A}_{M_{\natural}}$ est celui des familles $\zeta=(\zeta_{i})_{i\in I}$ de nombres r\'eels telles que $\zeta_{-i}=-\zeta_{i}$ pour tout $i$. L'espace ${\cal A}_{T''}$ est le sous-espace des $\zeta$ tels que $\zeta_{1}=\zeta_{-1}=0$, ou encore, en posant $ \tilde{I}=\{\pm 2,...,\pm l\}$, ${\cal A}_{T''}$ est l'espace des familles $\zeta=(\zeta_{i})_{i\in \tilde{I}}$ de nombres r\'eels telles que $\zeta_{-i}=-\zeta_{i}$ pour tout $i$.

  On d\'ecrit l'ensemble ${\cal F}(\tilde{M})$ de la fa\c{c}on habituelle suivante.  Notons $\Phi(\tilde{M})$ l'ensemble des applications $\varphi: \tilde{I}\to \{0,\pm 1,...,\pm j_{\varphi}\}$, o\`u $j_{\varphi}$ est un \'el\'ement quelconque de ${\mathbb N}$, telles que  $\varphi(-i)=-\varphi(i)$ pour tout $i$ et $\varphi^{-1}(j)\not=\emptyset$ pour tout $j=\pm 1,...,\pm j_{\varphi}$.  Pour une telle $\varphi$ et pour $j\in \{1,...,j_{\varphi}\}$, on pose $E_{\varphi,j}=\oplus_{\varphi^{-1}(j)}E_{i}$. On note $\tilde{Q}_{\varphi}=\tilde{L}_{\varphi}\tilde{U}_{\varphi}$ le sous-groupe parabolique de $\tilde{G}$ form\'e des \'el\'ements qui conservent le drapeau
$$E_{\varphi,j_{\varphi}}\subset E_{\varphi,j_{\varphi}}\oplus E_{\varphi,j_{\varphi}-1}\subset...\subset E_{\varphi,j_{\varphi}}\oplus ...\oplus E_{\varphi,1}.$$
Alors $\varphi\mapsto \tilde{Q}_{\varphi}$ est une bijection de $\Phi(\tilde{M})$ sur ${\cal F}(\tilde{M})$. L'espace ${\cal A}_{\tilde{M}}^{\tilde{L}_{\varphi}}$ est form\'e des $\zeta=(\zeta_{i})_{i\in \tilde{I}}$  tels que, pour tout $j\in \{\pm 1,...,\pm j_{max}\}$, $\sum_{i\in \varphi^{-1}(j) }\zeta_{i}=0$. L'espace ${\cal A}_{\tilde{L}_{\varphi}}$ est  form\'e des $\zeta=(\zeta_{i})_{i\in \tilde{I}}$  tels que, pour tout $j\in \{0,\pm 1,...,\pm j_{max}\}$, $\zeta_{i}$ est constant pour $i\in \varphi^{-1}(j)$. Notons $\zeta_{\varphi,j}$ cette valeur constante. Alors ${\cal A}_{\tilde{Q}_{\varphi}}^+$ est form\'e des $\zeta\in {\cal A}_{\tilde{L}_{\varphi}}$ tels que
$$\zeta_{\varphi,j_{\varphi}}>...>\zeta_{\varphi,1}>0.$$
 
 On d\'ecrit l'ensemble ${\cal F}(M_{\natural})$ de fa\c{c}on analogue. On d\'efinit $\Phi(M_{\natural})$ en rempla\c{c}ant $\tilde{I}$ par $I$ dans la d\'efinition de $\Phi(\tilde{M})$.  Pour $\varphi\in \Phi(M_{\natural})$, on d\'efinit de m\^eme un sous-groupe parabolique $Q_{\varphi}=L_{\varphi}U_{\varphi}$ de $G$. L'application $\varphi\mapsto Q_{\varphi}$ est une surjection de $\Phi(M_{\natural})$ sur ${\cal F}(M_{\natural})$. Les fibres ont $1$ ou $3$ \'el\'ements. Une fibre a $3$ \'el\'ements si et seulement si elle contient un \'el\'ement $\varphi_{0}$ pour lequel $\varphi_{0}^{-1}(0)$ a deux \'el\'ements $\{i_{h},-i_{h}\}$ et $dim(E_{i_{h}})=1$. Alors les deux autres \'el\'ements de la fibre sont les applications $\varphi_{1}$ et $\varphi_{-1}$ ainsi d\'efinies. Pour $\epsilon=\pm 1$, $\varphi_{\epsilon}(i_{h})=\epsilon$, $\varphi_{\epsilon}(-i_{h})=-\epsilon$. Pour $i\in I\setminus \{i_{h},-i_{h}\}$, si $\varphi_{0}(i)=\pm j$, avec $j\in \{1,...,j_{\varphi_{0}}\}$, on a $\varphi_{\epsilon}(i)=\pm(j+1)$. Si $\varphi$ appartient \`a une fibre \`a $1$ \'el\'ement, les espaces ${\cal A}_{M_{\natural}}^{L_{\varphi}}$, ${\cal A}_{L_{\varphi}}$ et l'ensemble ${\cal A}_{Q_{\varphi}}^+$ se d\'ecrivent comme ci-dessus, en rempla\c{c}ant $\tilde{I}$ par $I$. Si $\varphi$ est l'un des \'el\'ements $\varphi_{\pm 1}$ d'une fibre \`a $3$ \'el\'ements, les espaces ${\cal A}_{M_{\natural}}^{L_{\varphi}}$ et ${\cal A}_{L_{\varphi}}$ se d\'ecrivent de m\^eme. L'ensemble ${\cal A}_{Q_{\varphi}}^+$ est form\'e des $\zeta\in {\cal A}_{\tilde{L}_{\varphi}}$ tels que
$$\zeta_{\varphi,j_{\varphi}}>...>\zeta_{\varphi,2}, \,\,\zeta_{\varphi,2}+\zeta_{\varphi,1}>0,\,\,\zeta_{\varphi,2}-\zeta_{\varphi,1}>0.$$
Le sous-groupe parabolique $Q_{\varphi}$ appartient \`a $\tilde{{\cal F}}(M_{\natural})$ si et seulement si l'ensemble ${\cal A}_{Q_{\varphi}}^+$ contient un \'el\'ement de ${\cal A}_{T''}$, c'est-\`a-dire un \'el\'ement $\zeta$ pour lequel $\zeta_{1}=\zeta_{-1}=0$. La description ci-dessus montre que, si $\varphi$ appartient \`a une fibre \`a $1$ \'el\'ement, cette condition \'equivaut \`a $\varphi(1)=\varphi(-1)=0$. Si $\varphi$ est un \'el\'ement $\varphi_{\pm 1}$ d'une fibre \`a $3$ \'el\'ements, elle \'equivaut \`a $\varphi(1),\varphi(-1)\in \{\pm 1\}$. Il revient au m\^eme de dire que, si $\varphi$ est l'\'el\'ement $\varphi_{0}$ d'une fibre \`a $3$ \'el\'ements, la condition \'equivaut \`a $\varphi(1)=\varphi(-1)=0$. Notons $\tilde{\Phi}(M_{\natural})$ le sous-ensemble des $\varphi\in \Phi(M_{\natural})$ tels que $\varphi(1)=\varphi(-1)=0$. Alors l'application $\varphi\mapsto Q_{\varphi}$ se restreint en une bijection de $\tilde{\Phi}(M_{\natural})$ sur $\tilde{{\cal F}}(M_{\natural})$. L'application qui \`a $\varphi$ associe sa restriction \`a $\tilde{I}$ est une bijection de $\tilde{\Phi}(M_{\natural})$ sur $\Phi(\tilde{M})$. Elle param\`etre une bijection de $\tilde{{\cal F}}(M_{\natural})$ sur ${\cal F}(\tilde{M})$. En utilisant les descriptions ci-dessus, on v\'erifie qu'elle poss\`ede toutes les propri\'et\'es indiqu\'ees dans l'\'enonc\'e et c'est bien s\^ur la seule possible. $\square$

Gr\^ace \`a ce lemme, les propri\'et\'es (3) et (4) ne sont autres que (1) et (2) pour le groupe $\tilde{G}$.

Dans le cas non exceptionnel,  on pourra affecter d'un $\tilde{}$ les notations afin de les unifier avec celles du cas exceptionnel. Par exemple, on notera $\tilde{{\cal F}}(M_{\natural})={\cal F}(M_{\natural})$ et, pour un \'el\'ement $Q$ de cet ensemble, on  notera $\tilde{\tau}_{Q}=\tau_{Q}$. De m\^eme, on notera dans ce cas $\zeta\mapsto \zeta_{T''}$ l'application identit\'e de ${\cal A}_{M_{\natural}}$.

\bigskip

\subsection{Changement de fonction de troncature}

On fixe $S\in {\cal S}$ et $T\in {\cal T}(G_{x})$. On utilise les notations de 10.2 et 10.3.  Fixons un L\'evi minimal $M_{min}$ de $G$ contenu dans $M_{\natural}$ et contenant $A_{M_{\natural}}$. On fixe un sous-groupe compact sp\'ecial $K_{min}$ de $G(F)$ en bonne position relativement \`a $M_{min}$. Il nous sert \`a d\'efinir les fonctions $H_{Q}$ sur $G(F)$ pour $Q\in {\cal F}(M_{min})$. Fixons $P_{min}=M_{min}U_{min}\in {\cal P}(M_{min})$ et notons $\Delta_{min}$ l'ensemble des racines simples de $A_{M_{min}}$ dans $\mathfrak{u}_{min}$. Soit $Y_{P_{min}}\in {\cal A}_{P_{min}}^+$. Pour tout $P'\in {\cal P}(M_{min})$, il y a un unique $w\in W(G,M_{min})$ tel que $wP_{min}w^{-1}=P'$. On pose $Y_{P'}=wY_{P_{min}}$. La famille $(Y_{P'})_{P'\in {\cal P}(M_{min})}$ est $(G,M_{min})$-orthogonale et positive.  Pour $g\in G(F)$, d\'efinissons la famille ${\cal Y}(g)=(Y(g)_{Q})_{Q\in {\cal P}(M_{\natural})}$ par
$$Y(g)_{Q}=Y_{Q}-H_{\bar{Q}}(g).$$
Il est clair qu'il existe $c_{1}>0$ tel que

(1) pour tout $g\in G(F)$ tel que $\sigma(g)<c_{1}inf\{\alpha(Y_{P_{min}});\alpha\in \Delta_{min}\}$, la famille ${\cal Y}(g)$ est $(G,M_{\natural})$-orthogonale et positive; de plus $Y(g)_{Q}\in {\cal A}_{Q}^+$ pour tout $Q\in {\cal F}(M_{\natural})$.

On fixe un tel $c_{1}$. Remarquons que, pour $m\in M_{\natural}(F)$, la famille ${\cal Y}(mg)$ se d\'eduit de ${\cal Y}(g)$ par translations. Il en r\'esulte que la famille ${\cal Y}(g)$ est $(G,M_{\natural})$-orthogonale et positive pour tout  
$$g\in M_{\natural}(F)\{g'\in G(F); \sigma(g')<c_{1}inf\{\alpha(Y_{P_{min}});\alpha\in \Delta_{min}\}\}.$$
Pour un tel $g$, on pose
$$\tilde{v}(g)=\nu(A_{T''})\int_{A_{T''}(F)}\tilde{\sigma}_{M_{\natural}}^G(H_{M_{\natural}}(a),{\cal Y}(g))da.$$

On a d\'efini en 2.3 la fonction $\sigma_{T}$. Montrons que

(2) il existe $c_{2}>0$ et un sous-ensemble compact $\omega_{T}$ de $\mathfrak{t}(F)$ tels que les propri\'et\'es suivantes soient v\'erifi\'ees; soient $g\in G(F) $ et $X\in \mathfrak{t}(F)[S;>N^{-b}]\cap(\mathfrak{t}'(F)\times \mathfrak{t}''(F)^{S})$ tels que ${^gf}^{\sharp}_{x,\omega}(X)\not=0$; alors $X\in \omega_{T}$ et $\sigma_{T}(g)<c_{2}log(N)$.  

Preuve. Il suffit de reprendre la preuve du lemme 10.1. Les \'el\'ements $X'$ et $X''$ restent dans des sous-ensembles compacts de $\mathfrak{t}'(F)$ et $\mathfrak{t}''(F)$, ces sous-ensembles compacts, ainsi que les suivants \'etant bien s\^ur ind\'ependants de $X$ et $g$. On peut \'ecrire $g=g'g''\gamma$, o\`u $\gamma$ appartient \`a un sous-ensemble compact  de $G(F)$, $g'\in G'_{x}(F)$ et $g''\in G''(F)$ sont tels que $g^{_{'}-1}X'g'$ et $g^{_{''}-1}X''g''$ appartiennent \`a des sous-ensembles compacts de $\mathfrak{g}'_{x}(F)$ et $\mathfrak{g}''(F)$. D'apr\`es 2.3(1), quitte \`a multiplier $g$ \`a gauche par un \'el\'ement de $T(F)$, on a
$$\sigma(g)<<(1+\vert log\,D^{G'_{x}}(X')\vert )(1+\vert log\,D^{G''}(X'')\vert ).$$
Pour $X\in \mathfrak{t}(F)[S;>N^{-b}]$, on a
$$\vert log\,D^{G'_{x}}(X')\vert << log(N),\,\,\vert log\,D^{G''}(X'')\vert <<log(N).$$
Donc $\sigma(g)<<log(N)$ et l'assertion s'ensuit. $\square$

On fixe de tels  $\omega_{T}$ et $c_{2}$ et on suppose d\'esormais 
$$c_{2}log(N)<c_{1}inf\{\alpha(Y_{P_{min}});\alpha\in \Delta_{min}\}.$$
Puisque $T\subset M_{\natural}$, $\tilde{v}(g)$ est d\'efini pour tout $g$ satisfaisant la condition de (2). Le membre de droite de l'\'egalit\'e de la proposition ci-dessous est donc bien d\'efini.

\ass{Proposition}{Il existe $c>0$ et un entier $N_{0}\geq1$ tels que, si $N\geq N_{0}$ et
$$clog(N)<inf\{\alpha(Y_{P_{min}});\alpha\in \Delta_{min}\} ,$$
on ait l'\'egalit\'e
$$\int_{T'(F)A_{T''}(F)\backslash G(F)}{^gf}_{x,\omega}^{\sharp}(X)\kappa_{N,X''}(g)dg=\int_{T'(F)A_{T''}(F)\backslash G(F)}{^gf}_{x,\omega}^{\sharp}(X)\tilde{v}(g)dg$$
pour tout $X\in \mathfrak{t}(F)[S;>N^{-b}]\cap(\mathfrak{t}'(F)\times \mathfrak{t}''(F)^{S})$.}

Preuve. Soit $Z_{P_{min}}\in {\cal A}_{P_{min}}^+$. En rempla\c{c}ant $Y_{P_{min}}$ par cet \'el\'ement, on construit une famille ${\cal Z}(g)$ pour tout $g\in G(F)$. On impose
$$(3) \qquad c_{2}log(N)<c_{1}inf\{\alpha(Z_{P_{min}});\alpha\in \Delta_{min}\}.$$
Soit $g\in G(F)$ tel que $\sigma_{T}(g)<c_{2}log(N)$. Pour $a\in A_{T''}$, on a l'\'egalit\'e
$$\sum_{Q\in \tilde{{\cal F}}(M_{\natural})}\tilde{\sigma}_{M_{\natural}}^Q(H_{M_{\natural}}(a),{\cal Z}(g))\tilde{\tau}_{Q}(H_{M_{\natural}}(a)-Z(g)_{Q,T''})=1,$$
 cf. 10.3(4). En se rappelant la d\'efinition ci-dessus de $\tilde{v}(g)$, on peut \'ecrire
 $$\tilde{v}(g)=\nu(A_{T''})\sum_{Q\in\tilde{{\cal F}}(M_{\natural})}\tilde{v}(Q,g),$$
 o\`u
 $$\tilde{v}(Q,g)=\int_{A_{T''}(F)}\tilde{\sigma}_{M_{\natural}}^G(H_{M_{\natural}}(a),{\cal Y}(g))\tilde{\sigma}_{M_{\natural}}^Q(H_{M_{\natural}}(a),{\cal Z}(g))\tilde{\tau}_{Q}(H_{M_{\natural}}(a)-Z(g)_{Q,T''})da.$$
 De m\^eme, soit $X\in \mathfrak{t}'(F)\times \mathfrak{t}''(F)^{S}$. On a
 $$\kappa_{N,X''}(g)=\nu(A_{T''})\sum_{Q\in \tilde{{\cal F}}(M_{\natural})}\kappa_{N,X''}(Q,g),$$
 o\`u
 $$\kappa_{N,X''}(Q,g)=\int_{A_{T''}(F)} \kappa(\gamma_{X''}^{-1}ag)\tilde{\sigma}_{M_{\natural}}^Q(H_{M_{\natural}}(a),{\cal Z}(g))\tilde{\tau}_{Q}(H_{M_{\natural}}(a)-Z(g)_{Q,T''})da.$$
 On a
 
 (4) les fonctions $g\mapsto \tilde{v}(Q,g)$ et $g\mapsto \kappa_{N,X''}(Q,g)$ sont invariantes \`a gauche par $T'(F)A_{T''}(F)$.
 
 Soit $t\in T'(F) A_{T''}(F)$. On a $H_{P'}(tg)=H_{M_{\natural}}(t)+H_{P'}(g)$ pour tout $P'\in {\cal P}(M_{\natural})$. Supposons $t\in A_{T''}(F)$. Alors remplacer $g$ par $tg$ dans la d\'efinition de $\tilde{v}(Q,g)$ ou de $\kappa_{N,X''}(Q,g)$ revient \`a changer la variable d'int\'egration $a$ en $at$. Cela ne change pas l'int\'egrale. Soit maintenant $t\in T'(F)$. On a $T'\subset H$, donc $\kappa(\gamma_{X''}^{-1}atg)=\kappa(\gamma_{X''}^{-1}ag)$ pour tout $a$. On a $T'\subset {G_{2}}\subset M_{\natural}$, avec la notation de la preuve du lemme 10.2. Si l'on n'est pas dans le cas exceptionnel, $G_{2}$ est semi-simple et $H_{M_{\natural}}(t)=0$. On ne change donc rien en rempla\c{c}ant $g$ par $tg$. Si l'on est dans le cas exceptionnel, ce ne sont pas les termes $H_{P'}(g)$ qui interviennent dans les d\'efinitions, mais leurs projections $H_{P'}(g)_{T''}$. Or $H_{M_{\natural}}(t)_{T''}=0$ et, de nouveau, remplacer $g$ par $tg$ ne change rien. Cela prouve (4).
 
 Soit $X\in \mathfrak{t}'(F)\times \mathfrak{t}''(F)^{S}$. Gr\^ace \`a (4), on peut \'ecrire  
 $$(5) \qquad \int_{T'(F)A_{T''}(F)\backslash G(F)}{^gf}_{x,\omega}^{\sharp}(X)\kappa_{N,X''}(g)dg=\nu(A_{T''})\sum_{Q\in\tilde{{\cal F}}(M_{\natural})}I(Q,X),$$  
  $$(6) \qquad \int_{T'(F)A_{T''}(F)\backslash G(F)}{^gf}_{x,\omega}^{\sharp}(X)\tilde{v}(g)dg =\nu(A_{T''})\sum_{Q\in\tilde{{\cal F}}(M_{\natural})}J(Q,X),$$
  o\`u
  $$I(Q,X)=\int_{T'(F)A_{T''}(F)\backslash G(F)}{^gf}_{x,\omega}^{\sharp}(X)\kappa_{N,X''}(Q,g)dg,$$
  $$J(Q,X)=\int_{T'(F)A_{T''}(F)\backslash G(F)}{^gf}_{x,\omega}^{\sharp}(X)\tilde{v}(Q,g)dg.$$
 
Consid\'erons d'abord les termes index\'es par $Q=G$. Supposons
$$(7)\qquad sup\{\alpha(Z_{min});\alpha\in \Delta_{min}\}\leq \left\lbrace\begin{array}{c}inf\{\alpha(Y_{min}); \alpha\in \Delta_{min}\},\\ log(N)^2.\\ \end{array}\right.$$
Fixons un sous-ensemble compact $\omega_{T''}$ de $\mathfrak{t}''(F)$ tel que $X''\in \omega_{T''}$ pour tout $X\in \omega_{T}$. L'ensemble $\mathfrak{t}(F)[S;>N^{-b}]$ a \'et\'e d\'efini comme celui des $X\in \mathfrak{t}(F)$ tels que $X'$ et $X''$ satisfassent certaines minorations. On d\'efinit $\mathfrak{t}''(F)^S[>N^{-b}]$ comme celui des $X''\in \mathfrak{t}''(F)^S$ v\'erifiant celles de ces minorations qui portent sur $X''$. On a

(8) il existe un entier $N_{1}\geq 1$ tel que, pour tout  $N\geq N_{1}$, pour tout $g\in G(F)$ tel que $\sigma_{T}(g)\leq c_{2}log(N)$ et tout $X''\in \omega_{T''}\cap\mathfrak{t}''(F)^S[>N^{-b}] $, on ait l'\'egalit\'e $\kappa_{N,X''}(G,g)=\tilde{v}(G,g)$.

Il suffit de prouver que pour tout $a\in A_{T''}(F)$ tel que $\tilde{\sigma}_{M_{\natural}}^G(H_{M_{\natural}}(a),{\cal Z}(g))=1$, on a l'\'egalit\'e $ \tilde{\sigma}_{M_{\natural}}^G(H_{M_{\natural}}(a),{\cal Y}(g))=  \kappa_{N}(\gamma_{X''}^{-1}ag)$. La premi\`ere in\'egalit\'e de (5) entra\^{\i}ne que l'enveloppe convexe de la famille $(Z(g)_{P',T''})_{P'\in \tilde{{\cal P}}(M_{\natural})}$ est incluse dans celle de la famille $(Y(g)_{P',T''})_{P'\in \tilde{{\cal P}}(M_{\natural})}$. Alors l'hypoth\`ese $\tilde{\sigma}_{M_{\natural}}^G(H_{M_{\natural}}(a),{\cal Z}(g))=1$ entra\^{\i}ne $\tilde{\sigma}_{M_{\natural}}^G(H_{M_{\natural}}(a),{\cal Y}(g))=1$. Munissons ${\cal A}_{M_{\natural}}$ d'une norme $ \vert . \vert $. La deuxi\`eme in\'egalit\'e de (5) et l'hypoth\`ese sur $g$ entra\^{\i}nent une majoration $\vert Z(g)_{P'} \vert << log(N)^2$ pour tout $P'\in {\cal P}(M_{\natural})$. L'hypoth\`ese $\tilde{\sigma}_{M_{\natural}}^G(H_{M_{\natural}}(a),{\cal Z}(g))=1$ entra\^{\i}ne alors $\sigma(a)<< log(N)^2$. D'apr\`es 9.7(2), on a $\sigma(\gamma_{X''})<<1+\vert log\,D^{G''}(X'')\vert $. Puisque $X''\in \mathfrak{t}''(F)^S[>N^{-b}]$, cela entra\^{\i}ne $\sigma(\gamma_{X''})<<log(N)$. D'o\`u $\sigma(\gamma_{X''}^{-1}ag)<<log(N)^2$. Mais on voit facilement qu'il existe $c_{3}>0$ tel que, pour tout $g'\in G(F)$ tel que $\sigma(g')<c_{3}N$, on a $\kappa_{N}(g')=1$. Si $N$ est assez grand, on a $\sigma(\gamma_{X''}^{-1}ag)<c_{3}N$, donc $\kappa_{N}(\gamma_{X''}^{-1}ag)=1$. Cela prouve (8).

De (2) et (8) r\'esulte l'\'egalit\'e 
$$(9) \qquad I(G,X)=J(G,X)$$
 pour tout  $N\geq N_{1}$ et tout $X\in \mathfrak{t}(F)[S;>N^{-b}]\cap(\mathfrak{t}'(F)\times \mathfrak{t}''(F)^{S})$.
 
 Soit maintenant $Q=LU_{Q}\in \tilde{{\cal F}}(M_{\natural})$ avec $Q\not=G$. On d\'ecompose les int\'egrales
 $$I(Q,X)=\int_{K_{min}}\int_{T'(F)A_{T''}(F)\backslash L(F)}\int_{U_{\bar{Q}}(F)}{^{\bar{u}lk}f}_{x,\omega}^{\sharp}(X)\kappa_{N,X''}(Q,\bar{u}lk)d\bar{u}\delta_{Q}(l)dl\,dk,$$
  $$J(Q,X)=\int_{K_{min}}\int_{T'(F)A_{T''}(F)\backslash L(F)}\int_{U_{\bar{Q}}(F)}{^{\bar{u}lk}f}_{x,\omega}^{\sharp}(X)\tilde{v}(Q,\bar{u}lk)d\bar{u}\delta_{Q}(l)dl\,dk.$$
  Nous montrerons aux paragraphes 10.5 et 10.8  les propri\'et\'es suivantes.
  
  (10) Soient $g\in G(F)$ et $\bar{u}\in U_{\bar{Q}}(F)$ tels que $\sigma(g),\sigma(\bar{u}g)<c_{1}inf\{\alpha(Z_{P_{min}}); \alpha\in \Delta_{min}\}$. Alors on a l'\'egalit\'e $\tilde{v}(Q,\bar{u}g)=\tilde{v}(Q,g)$.
  
  (11) Soit $c_{4}>0$. Il existe $c_{5}>0$ tel que si $c_{5}log(N)<inf\{\alpha(Z_{P_{min}}); \alpha\in \Delta_{min}\}$, les conditions suivantes soient v\'erifi\'ees. Soient $X''\in  \omega_{T''}\cap\mathfrak{t}''(F)^S[>N^{-b}] $, $g\in G(F)$ et $\bar{u}\in U_{\bar{Q}}(F)$. Supposons $\sigma(g), \sigma(\bar{u}), \sigma(\bar{u}g)<c_{4}log(N)$. Alors on a l'\'egalit\'e $\kappa_{N,X''}(Q,\bar{u}g)=\kappa_{N,X''}(Q,g)$;

  Admettons ces propri\'et\'es. Montrons
  
  (12) il existe $c_{5}>0$ tel que, si $c_{5}log(N)<inf\{\alpha(Z_{P_{min}}); \alpha\in \Delta_{min}\}$, on a les \'egalit\'es $I(Q,X)=J(Q,X)=0$ pour tout $X\in \mathfrak{t}(F)[S;>N^{-b}]\cap(\mathfrak{t}'(F)\times \mathfrak{t}''(F)^{S})$.
  
  Gr\^ace \`a (2), on peut supposer $X\in \omega_{T}$. Consid\'erons l'int\'egrale $I(Q,X)$. D'apr\`es (2), on peut limiter l'int\'egrale sur $T'(F)A_{T''}(F)\backslash L(F)$ aux \'el\'ements $l$ pour lesquels il existe $\bar{u}\in U_{\bar{Q}}(F)$ et $k\in K_{min}$ tels que $\sigma_{T}(\bar{u}lk)< c_{2}log(N)$. Un tel $l$ est repr\'esent\'e par un \'el\'ement de $L(F)$ tel que $\sigma(l)<c_{6}log(N)$, pour une constante $c_{6}$ convenable. Il existe $c_{7}>0$ tel que, pour $l$ v\'erifiant l'in\'egalit\'e pr\'ec\'edente, pour $k\in K_{min}$ et pour $\bar{u}\in U_{\bar{Q}}(F)$, l'in\'egalit\'e $\sigma(\bar{u}lk)<c_{2}log(N)$ entra\^{\i}ne $\sigma(\bar{u})<c_{7}log(N)$. Soit $c_{4}=c_{2}+c_{7}$ et prenons pour $c_{5}$ le nombre issu de (11). Fixons $k\in K_{min}$ et $l\in L(F)$ tel que $\sigma(l)<c_{6}log(N)$. On a alors
  $${^{\bar{u}lk}f}_{x,\omega}^{\sharp}(X)\kappa_{N,X''}(Q,\bar{u}lk)={^{\bar{u}lk}f}_{x,\omega}^{\sharp}(X)\kappa_{N,X''}(Q,lk)$$
  pour tout $\bar{u}\in U_{\bar{Q}}(F)$. En effet, si   $\sigma(\bar{u}lk)\geq c_{2}log(N)$, les deux termes sont nuls d'apr\`es (2). Si $\sigma(\bar{u}lk)<c_{2}log(N)$, on a aussi $\sigma(\bar{u})<c_{7}log(N)$, puis $\sigma(lk)<c_{4}log(N)$. La relation (11) s'applique \`a $g=lk$ et $u$. Donc $\kappa_{N,X''}(Q,\bar{u}lk)=\kappa_{N,X''}(Q,lk)$ et l'\'egalit\'e affirm\'ee s'ensuit. Il r\'esulte de cette \'egalit\'e que, dans $I(Q,X)$, l'int\'egrale int\'erieure est simplement
  $$\int_{U_{\bar{Q}}(F)}{^{\bar{u}lk}f}_{x,\omega}^{\sharp}(X)d\bar{u}.$$ 
  Or cette int\'egrale est nulle d'apr\`es le lemme 5.5(i), puisque $Q\not=G$. Donc $I(Q,X)=0$. On prouve de m\^eme que $J(Q,X)=0$.
  
  Le terme $Z_{P_{min}}$ est un terme auxiliaire. Il est clair qu'il existe  $c>0$ et un entier $N_{2}\geq1$ tels que, si $N\geq N_{2}$ et
$$clog(N)<inf\{\alpha(Y_{P_{min}});\alpha\in \Delta_{min}\} ,$$
on peut choisir $Z_{P_{min}}$ satisfaisant les hypoth\`eses (3) et (7) et celle de la relation (12). Si on suppose de plus $N\geq N_{1}$, les conclusions de (9) et (12) s'appliquent. Alors les \'egalit\'es (5) et (6) entra\^{\i}nent la conclusion de l'\'enonc\'e. $\square$
  
  \bigskip
  
  \subsection{Preuve de la propri\'et\'e 10.4(10)}
  
  Soient $g$ et $\bar{u}$ comme dans cette relation. Rappelons que
  $$\tilde{v}(Q,g)=\int_{A_{T''}(F)}\tilde{\sigma}_{M_{\natural}}^G(H_{M_{\natural}}(a),{\cal Y}(g))\tilde{\sigma}_{M_{\natural}}^Q(H_{M_{\natural}}(a),{\cal Z}(g))\tilde{\tau}_{Q}(H_{M_{\natural}}(a)-Z(g)_{Q,T''})da.$$
 Les fonctions $\zeta\mapsto \tilde{\sigma}_{M_{\natural}}^Q(\zeta,{\cal Y}(g))$ et $\zeta\mapsto \tilde{\tau}_{Q}(\zeta-Z(g)_{Q,T''})$ ne d\'ependent de $g$ que par l'interm\'ediaire des termes $H_{\bar{P}'}(g)$ pour des sous-groupes paraboliques $P'\in \tilde{{\cal F}}(M_{\natural})$ tels que $P'\subset Q$.  Elles ne changent donc pas quand on remplace $g$ par $\bar{u}g$. On peut alors fixer $a\in A_{T''}(F)$ tel que  
 $$\tilde{\sigma}_{M_{\natural}}^Q(H_{M_{\natural}}(a),{\cal Z}(g))\tilde{\tau}_{Q}(H_{M_{\natural}}(a)-Z(g)_{Q,T''})\not=0$$
 et prouver que
 $$(1) \qquad \tilde{\sigma}_{M_{\natural}}^G(H_{M_{\natural}}(a),{\cal Y}(g))=\tilde{\sigma}_{M_{\natural}}^G(H_{M_{\natural}}(a),{\cal Y}(\bar{u}g)).$$
 Supposons que l'on n'est pas dans le cas exceptionnel. Tout sous-groupe parabolique $P'\in {\cal P}(M_{\natural})$ tel que $P'\subset Q$ d\'etermine une chambre ${\cal A}_{P'}^{L,+}$ dans ${\cal A}^L$.  Posons $\zeta=H_{M_{\natural}}(a)$, fixons un tel $P'$ de sorte que $proj_{M_{\natural}}^L(\zeta)\in Cl({\cal A}_{P'}^{L,+})$, o\`u, pour tout sous-ensemble $E$ de ${\cal A}_{M_{\natural}}$, $Cl(E)$ d\'esigne sa cl\^oture. Montrons que
 
 (2) $\zeta\in Cl({\cal A}_{P'}^+)$.
 
 D'apr\`es 10.3(1), l'hypoth\`ese sur $a$ signifie que $\zeta$ est la somme d'un \'el\'ement $\zeta'\in {\cal A}_{Q}^+$ et d'un \'el\'ement $\zeta''$ dans l'enveloppe convexe des $Z(g)_{P''}$, pour $P''\in {\cal P}(M_{\natural})$ tel que $P''\subset Q$. Soit $\alpha$ une racine de $A_{M_{\natural}}$ dans $\mathfrak{g}$, positive pour $P'$. Si $\alpha$ intervient dans $\mathfrak{u}_{Q}$, $\alpha$ est positive pour tous les $P''$ ci-dessus. Or $Z(g)_{P''}\in {\cal A}_{P''}^+$ d'apr\`es 10.4(1), donc $\alpha(Z(g)_{P''})>0$. Il en r\'esulte que $\alpha(\zeta'')>0$. On a aussi $\alpha(\zeta')>0$, donc $\alpha(\zeta)>0$. Si maintenant $\alpha$ intervient dans $\mathfrak{u}_{P'}\cap \mathfrak{l}$, on a $\alpha(\zeta)=\alpha(proj_{M_{\natural}}^L(\zeta))\geq 0$ d'apr\`es le choix de $P'$. Cela prouve (2).
 
 D'apr\`es [A3] lemme 3.1, pour $\zeta\in Cl({\cal A}_{P'}^+)$, la condition $\tilde{\sigma}_{M_{\natural}}^G(\zeta,{\cal Y}(g))=1$ \'equivaut \`a certaines in\'egalit\'es portant sur $\zeta-Y(g)_{P'}$. Cette condition ne d\'epend de $g$ que par l'interm\'ediaire de $H_{\bar{P}'}(g)$. Comme ci-dessus, elle est donc insensible au changement de $g$ en $\bar{u}g$. Cela d\'emontre (1).
 
 Dans le cas exceptionnel, on utilise le lemme 10.3 pour interpr\'eter nos fonctions comme leurs analogues pour le groupe $\tilde{G}$. dans ce groupe, on peut faire le m\^eme raisonnement que ci-dessus et on obtient la m\^eme conclusion. $\square$
 
 \bigskip
 
 \subsection{Calcul d'un polyn\^ome}

On aura besoin du lemme ci-dessous.  Soient ${\mathbb F}$ un corps alg\'ebriquement clos, $l$ un entier tel que  $l \geq 1$ et $R=R(T)$ un polyn\^ome en une indetermin\'ee, \`a coefficients dans ${\mathbb F}$, de degr\'e $l-1$ et de coefficient dominant $1$. Introduisons le polyn\^ome en $l+1$ ind\'etermin\'ees
$$Q=Q(T,S_{1},...,S_{l})=\prod_{j=1,...,l}(T-S_{j}),$$
et la fraction rationnelle
$$P=P(T,S_{1},...,S_{l})=1+\sum_{j=1,...,l}\frac{(T+S_{j})R(S_{j})}{(T-S_{j})\prod_{j'=1,...,l; j'\not=j}(S_{j}-S_{j'})}.$$

\ass{Lemme}{On a l'\'egalit\'e
$$P(T,S_{1},...,S_{l})=\frac{2TR(T)}{Q(T,S_{1},...,S_{l})}.$$}

Preuve. Introduisons le polyn\^ome
$$\Delta=\Delta(S_{1},...,S_{l})=\prod_{j,j'=1,...,l; j<j'}(S_{j}-S_{j'}).$$
Par r\'eduction au m\^eme d\'enominateur,
$$P=\frac{P_{\flat}}{\Delta Q}$$
o\`u
$$P_{\flat}=P_{\flat}(T,S_{1},...,S_{l})=\Delta Q+$$
$$\qquad\sum_{j=1,...,l}(-1)^{j-1}(T+S_{j})R(S_{j})\prod_{j'=1,...,l; j'\not=j}(T-S_{j'})\prod_{j',j''=1,...,l;j'<j''; j',j''\not=j}(S_{j'}-S_{j''}).$$
La fraction rationnelle $P$ est sym\'etrique en les $S_{j}$. Alors le polyn\^ome $P_{\flat}$ est antisym\'etrique et donc divisible par $\Delta$: il existe un polyn\^ome $P_{\natural}=P_{\natural}(T,S_{1},...,S_{l})$ tel que $P_{\flat}=P_{\natural}\Delta$. Le polyn\^ome $P_{\flat}$ est de degr\'e au plus $l$ en $T$. Le coefficient de $T^l$ est
$$\Delta+\sum_{j=1,...,l}(-1)^{j-1}R(S_{j})\prod_{j',j''=1,...,l; j'<j''; j',j''\not=j}(S_{j'}-S_{j''}).$$
Ce polyn\^ome en les $S_{j}$ est divisible par $\Delta$. Or son degr\'e total est inf\'erieur ou \'egal \`a celui de $\Delta$. Il est donc proportionnel \`a $\Delta$. On calcule le coefficient de proportionnalit\'e en calculant le coefficient de $S_{1}^{l-1}S_{2}^{l-2}...S_{l-1}$. On obtient  que ce coefficient est $2$. On en d\'eduit que le polyn\^ome $P_{\natural}$ est de degr\'e $l$ en $T$ et que son coefficient dominant est $2$. Pour $T=S_{j}$, on calcule
$$P_{\flat}(S_{j},S_{1},...,S_{l})=2S_{j}R(S_{j})\Delta.$$
Donc 
$$P_{\natural}(S_{j},S_{1},...,S_{l})=2S_{j}R(S_{j}).$$
Alors $P_{\natural}(T,S_{1},...,S_{j})$ et $2TR(T)$ sont des polyn\^omes de degr\'e $l$ en $T$, de m\^eme coefficient dominant, et prenant les m\^emes valeurs aux $l$ points $T=S_{j}$. Ils sont donc \'egaux. On obtient $P_{\flat}=2TR(T)\Delta$ et la formule de l'\'enonc\'e s'ensuit. $\square$

\bigskip

\subsection{R\'eseaux sp\'eciaux et extension de corps de base}

Pour ce paragraphe, on oublie les d\'efinitions de l'espace $Z$ et du r\'eseau $R$. Pour toute extension finie $F'$ de $F$, on note $V_{F'}=V\otimes_{F}F'$.  La forme $q$ se prolonge en une forme $F'$-bilin\'eaire $q_{F'}$ sur $V_{F'}$. Si $R$ est un $\mathfrak{o}_{F}$-r\'eseau de $V$, on note $R_{F'}=R\otimes_{\mathfrak{o}_{F}}\mathfrak{o}_{F'}$. Soient $R$ un $\mathfrak{o}_{F}$-r\'eseau de $V$ et $R'$ un $\mathfrak{o}_{F'}$-r\'eseau de $V_{F'}$. Notons $K$, resp. $K'$, le stabilisateur de $R$ dans $G(F)$, resp. de $R'$ dans $G(F')$. Disons que $R$ et $R'$ sont compatibles s'ils v\'erifient les deux conditions

(1) $R'\cap V=R$;

(2) $K'\cap G(F)=K$.

Remarquons que la premi\`ere condition entra\^{\i}ne $K'\cap G(F)\subset K$: un \'el\'ement de $K'\cap G(F)$ conserve $R'$ et $V$, donc aussi leur intersection $R$. On peut aussi bien remplacer (2) par

(2') $K\subset K'$.

\ass{Lemme}{Soit $R$ un r\'eseau sp\'ecial de $V$. Il existe une extension finie $E$ de $F$ telle que, pour toute extension finie $F'$ de $E$, il existe un r\'eseau sp\'ecial $R'$ de $V_{F'}$ qui soit compatible avec $R$.}

Preuve. On imagine qu'il y a une d\'emonstration immobili\`ere g\'en\'erale. Donnons une d\'emonstration d'alg\`ebre lin\'eaire \'el\'ementaire. On a la propri\'et\'e \'evidente

(3) soient $F'$ une extension finie de $F$, $F''$ une extension finie de $F'$, $R'$ un r\'eseau de $V_{F'}$ et $R''$ un r\'eseau de $V_{F''}$; supposons $R'$ compatible avec $R$  et $R''$ compatible avec   $R'$; alors $R''$ est compatible avec  $R$.

On a d'autre part

(4) si $d_{an}(V)\leq 1$, le lemme est v\'erifi\'e pour $E=F$.

En effet, pour toute extension finie $F'$ de $F$, le r\'eseau $R_{F'}$ est compatible avec $R$. Si $d_{an}(V)\leq 1$, $R_{F'}$ est sp\'ecial, d'o\`u (4).

A l'aide de ces deux propri\'et\'es, un raisonnement par r\'ecurrence descendante sur $d_{an}(V)$ montre qu'il suffit de prouver l'assertion suivante

(5) si $d_{an}(V)\geq2$, il existe une extension finie $F'$ de $F$ et un r\'eseau sp\'ecial $R'$ de $V_{F'}$ tels que $d_{an}(V_{F'})<d_{an}(V)$ et $R'$ soit compatible avec $R$.

On choisit comme en 7.1 une d\'ecomposition orthogonale $V=Z\oplus V_{an}$, un entier $c$ tel qu'il existe $v\in V_{an}$ de sorte que $val_{F}(q(v))=c$,  et une base hyperbolique $(v_{i})_{i=\pm1,...,\pm n}$ de sorte que $R$ soit la somme du r\'eseau $R_{an}$ form\'e des \'el\'ements $v\in V_{an}$ tels que $val_{F}(q(v))\geq c$ et du r\'eseau $R_{Z}$ engendr\'e par les $v_{i}$  et  les $\varpi_{F}^cv_{-i}$ pour $i=1,...,n$. Supposons $d_{an}(V)=2$. Il existe une extension quadratique $E$ de $F$ et un \'el\'ement $\lambda\in F^{\times}$ tel que $val_{F}(\lambda)=c$, de sorte que l'on puisse identifier $V_{an}$ \`a $E$ et la restriction $q_{an}$ de $q$ \`a $V_{an}$ \`a la forme quadratique $(v,v')\mapsto \lambda Trace_{E/F}(\tau(v)v')$, o\`u $\tau$ l'\'el\'ement non trivial de $Gal(E/F)$. L'espace $V_{an,E}$ s'identifie \`a un espace de dimension $2$ sur $E$, muni d'une base $(w_{+},w_{-})$, et $q_{an,E}$ \`a la forme
$$(x_{+}w_{+}+x_{-}w_{-},y_{+}w_{+}+y_{-}w_{-})=\lambda(x_{+}y_{-}+x_{-}y_{+}).$$
L'espace $V_{an}$ est form\'e des $x_{+}w_{+}+x_{-}w_{-}\in V_{an,E}$ tels que $x_{-}=\tau(x_{+})$. Posons $v_{n+1}=w_{+}$, $v_{-n-1}=\lambda^{-1}w_{-}$. Alors $(v_{i})_{i=\pm 1,...,\pm (n+1)}$ est une base hyperbolique de $V_{E}$. Notons $R'$ le $\mathfrak{o}_{E}$-r\'eseau de $V_{E}$ engendr\'e par les $v_{i}$ et les $\varpi_{F}^cv_{-i}$ pour $i=1,...,n+1$. Il est sp\'ecial. Montrons que $R'$ est compatible avec $R$. Remarquons que $R'$ est la somme de $R_{Z,E}$ et du r\'eseau $R'_{an}$ engendr\'e par $w_{+}$ et $w_{-}$.  Pour prouver que $R'\cap V=R$, il suffit de prouver que $R'_{an}\cap V_{an}=R_{an}$. Le r\'eseau $R'_{an}\cap V_{an}$, resp. $R_{an}$, est form\'e des $xw_{+}+\tau(x)w_{-}$, avec $x\in E$, tels que $x\in \mathfrak{o}_{E}$, resp. $val_{F}(x\tau(x))\geq0$. Ces deux derni\`eres conditions sont \'equivalentes, d'o\`u l'assertion. Montrons que $K\subset K'$. Soit $k\in K$. Il suffit de prouver que, pour tout  \'el\'ement $v$ de la base de $R'$, on a $kv\in R'$. Si $v=v_{i}$ ou $\varpi_{F}^cv_{-i}$ avec $i= 1,...,n$, on a $v\in R$, donc $kv\in R\subset R'$. On peut remplacer les deux \'el\'ements de base restants par $w_{+}$ et $w_{-}$. Ecrivons
$$kw_{+}=x_{+}w_{+}+x_{-}w_{-}+\sum_{i=1,...,n}(x_{i}v_{i}+\varpi_{F}^cx_{-i}v_{-i}),$$
$$kw_{+}=y_{+}w_{+}+y_{-}w_{-}+\sum_{i=1,...,n}(y_{i}v_{i}+\varpi_{F}^cy_{-i}v_{-i}).$$
Soit $i=1,...,n$. On a $\varpi_{F}^cx_{-i}=q_{E}(kw_{+},v_{i})=q_{E}(w_{+},k^{-1}v_{i})$. On a d\'ej\`a prouv\'e que $k^{-1}v_{i}$ appartenait \`a $R'$, donc $q_{E}(w_{+},k^{-1}v_{i})\in \varpi_{F}^c\mathfrak{o}_{E}$, puis $x_{-i}\in \mathfrak{o}_{E}$. De m\^eme $x_{i}$, $y_{-i}$ et $y_{i}$ appartiennent \`a $\mathfrak{o}_{E}$. On a $q_{E}(w_{+})=0$, donc $q_{E}(kw_{+})=0$, ce qui s'\'ecrit
$$\lambda x_{+}x_{-}+\sum_{i=1,...,n}\varpi_{F}^cx_{i}x_{-i}=0.$$
D'apr\`es ce que l'on vient de d\'emontrer, cela entra\^{\i}ne $x_{+}x_{-}\in \mathfrak{o}_{E}$. De m\^eme, $y_{+}y_{-}\in \mathfrak{o}_{E}$. L'automorphisme galoisien $\tau$ de $E$ induit un automorphisme antilin\'eaire de $V_{E}$ que l'on note aussi $\tau$. On a $\tau(v_{\pm i})=v_{\pm i}$ pour $i=1,...,n$, $\tau(w_{+})=w_{-}$ et $\tau(w_{-})=w_{+}$. Puisque $k\in G(F)$, il commute \`a $\tau$, donc $kw_{-}=\tau(kw_{+})$, d'o\`u $y_{-}=\tau(x_{+})$ et $y_{+}=\tau(x_{-})$. On a $q_{E}(w_{+},w_{-})=\lambda$, donc
$q_{E}(kw_{+},kw_{-})=\lambda$, ce qui s'\'ecrit
$$\lambda (x_{+}\tau(x_{+})+x_{-}\tau(x_{-}))+\varpi_{F}^c\sum_{i=1,...,n}(x_{i}y_{-i}+x_{-i}y_{-})=\lambda.$$
Cela entra\^{\i}ne $x_{+}\tau(x_{+})+x_{-}\tau(x_{-})\in \mathfrak{o}_{E}$. Si par exemple $x_{+}\not\in \mathfrak{o}_{E}$, cette relation implique que $x_{-}$ n'appartient pas non plus \`a $\mathfrak{o}_{E}$. Alors, la relation $x_{+}x_{-}\in \mathfrak{o}_{E}$ n'est pas v\'erifi\'ee, contrairement \`a ce que l'on a d\'ej\`a prouv\'e. Cette contradiction prouve que $x_{+}\in \mathfrak{o}_{E}$. De m\^eme, $x_{-}$, $y_{+}$ et $y_{-}$ appartiennent \`a $\mathfrak{o}_{E}$. Alors $kw_{+}$ et $kw_{-}$ appartiennent \`a $R'$ comme on le voulait. Cela prouve (4) sous l'hypoth\`ese $d_{an}(V)=2$.

Supposons $d_{an}(V)=3$. Soit $E$ l'extension quadratique non ramifi\'ee de $F$ et $\tau$ l'\'el\'ement non trivial de $Gal(E/F)$. On peut identifier $V_{an}$ \`a $E\oplus F$ et $q_{an}$ \`a une forme
$$(w\oplus z,w'\oplus z')\mapsto \lambda\mu Trace_{E/F}(\tau(w)w')+2\lambda\nu zz',$$
o\`u $\lambda,\mu,\nu\in F^{\times}$, $val_{F}(\lambda)=c$ et, ou bien $\mu=1$ et $val_{F}(\nu)=1$, ou bien $\nu=1$ et $val_{F}(\mu)=1$. Le r\'eseau $R_{an}$ s'identifie \`a $\mathfrak{o}_{E}\oplus \mathfrak{o}_{F}$. L'espace $V_{an,E}$ s'identifie \`a un espace de dimension $3$ sur $E$ muni d'une base $(w_{+},w_{-},w_{0})$ de sorte que la forme $q_{an,E}$ s'\'ecrive
$$q(x_{+}w_{+}+x_{-}w_{-}+x_{0}w_{0},y_{+}w_{+}+y_{-}w_{-}+y_{0}w_{0})=\lambda\mu(x_{+}y_{-}+x_{-}y_{+})+2\lambda\nu x_{0}y_{0}.$$
Parce que $E$ est non ramifi\'ee sur $F$, $R_{an,E}$ s'identifie au $\mathfrak{o}_{E}$-r\'eseau engendr\'e par les \'el\'ements de base. Posons $v_{n+1}=w_{+}$, $v_{-n-1}=\lambda^{-1}\mu^{-1}w_{-}$. La famille $(v_{i})_{i=\pm 1,...,\pm (n+1)}$ est un syst\`eme hyperbolique maximal de $V_{E}$ et $w_{0}$ est une base de son orthogonal. Le r\'eseau $R_{E}$ est engendr\'e sur $\mathfrak{o}_{E}$ par les $v_{i}$ et les $\varpi_{F}^cv_{-i}$ pour $i=1,...,n$, les \'el\'ements $v_{n+1}$ et $w_{0}$ et l'\'el\'ement $\varpi_{F}^c v_{-n-1}$, resp. $\varpi_{F}^{c+1}v_{-n-1}$, si $\mu=1$, resp. $\nu=1$. Il est compatible avec $R$ mais n'est pas sp\'ecial (dans le cas o\`u $\mu=1$, la droite $Fw_{0}$ ne repr\'esente aucun \'el\'ement de valuation $c$). N\'eanmoins, gr\^ace \`a (3), il suffit de trouver une extension $F'$ de $E$ et un r\'eseau sp\'ecial $R'$ de $V_{F'}$ qui soit compatible \`a $R_{E}$. Introduisons une racine carr\'ee $z$ de $\mu\nu$, soit $F'=E(z)$ et $R'$ le r\'eseau de $V_{F'}$  engendr\'e par les $v_{i}$ et les $\varpi_{F}^cv_{-i}$ pour $i=1,...,n$ et par

- $v_{n+1}$, $\varpi_{F}^cv_{-n-1}$ et $z^{-1}w_{0}$ si $\mu=1$;

- $z^{-1}v_{n+1}$, $\varpi_{F}^czv_{-n-1}$ et $w_{0}$ si $\nu=1$.

Ce r\'eseau est sp\'ecial. Il est clair que $R'\cap V_{E}=R_{E}$. Notons $K_{E}$ le stabilisateur de $R_{E}$ dans $G(E)$ et $K'$ celui de $R'$ dans $G(F')$. Il reste \`a prouver que $K_{E}\subset K'$.  Introduisons le r\'eseau $R_{F'}$, qui est aussi \'egal \`a $(R_{E})_{F'}$, et son dual $R_{F'}^*=\{v\in V_{F'}; \forall w\in R_{F'}, q_{F'}(v,w)\in \mathfrak{o}_{F'}\}$. On v\'erifie que
$$R'=z^{-1}R_{F'}\cap \varpi_{F}^cR_{F'}^*\cap\{v\in V_{F'}; val_{F'}(q_{F'}(v))\geq val_{F'}(\lambda)\}.$$
Un \'el\'ement de $K_{E}$ stabilise forc\'ement $R_{F'}$, donc aussi son dual, et il stabilise aussi le dernier ensemble ci-dessus. Donc il stabilise $R'$ et appartient \`a $K'$. Cela prouve (4) sous l'hypoth\`ese $d_{an}(V)=3$. 

Supposons enfin $d_{an}(V)=4$. Avec les m\^emes notations que dans le cas pr\'ec\'edent, on peut identifier $V_{an}$ \`a $E\oplus E$ et $q_{an}$ \`a la forme
$$(w_{1}\oplus w_{2},w'_{1}\oplus w'_{2})=\lambda Trace_{E/F}(\tau(w_{1})w'_{1})+\varpi_{F}\lambda Trace_{E/F}(\tau(w_{2})w'_{2}).$$
Introduisons une racine carr\'ee $z$ de $\varpi_{F}$, posons $F'=E(z)$. On v\'erifie comme dans le cas pr\'ec\'edent que le r\'eseau $R'=z^{-1}R_{F'}\cap \lambda R_{F'}^*$ de $V_{F'}$ satisfait les conditions de (4). Cela ach\`eve la preuve. $\square$

Revenons au r\'eseau $R$ que l'on a fix\'e en 7.2. On lui a impos\'e d'\^etre somme d'un r\'eseau de $V_{0}$ et d'un r\'eseau de $Z$ engendr\'e par des \'el\'ements proportionnels aux \'el\'ements $v_{i}$ pour $i=\pm 1,...,\pm r$. La preuve du lemme montre que l'on peut imposer aux r\'eseaux $R'$ de l'\'enonc\'e de v\'erifier les m\^emes conditions.

 \subsection{Preuve de la relation 10.4(11)}
 
 On va \'elargir les hypoth\`eses sur la fonction $X''\mapsto \gamma_{X''}$. Celles que l'on impose dans ce paragraphe sont
 
 (1) il existe un sous-ensemble compact $\Omega$ de $\Xi+S+\Sigma$ tel que $X''_{\Sigma}=\gamma_{X''}^{-1}X''\gamma_{X''}\in \Omega$ pour tout $X''\in \omega_{T''}\cap\mathfrak{t}''(F)^S$;
 
 (2) il existe $c_{1}>0$ tel que $\sigma(\gamma_{X''})<c_{1}log(N)$ pour tout $X''\in \omega_{T''}\cap\mathfrak{t}''(F)^S[>N^{-b}]$.

La fonction que nous avons utilis\'ee jusque-l\`a v\'erifie ces hypoth\`eses: (1) r\'esulte de 9.7(1) et on a vu dans la preuve de 10.4(8) que 9.7(2) entra\^{\i}nait (2).

Soit $Q=LU_{Q}\in \tilde{{\cal F}}(M_{\natural})$. On note $\Sigma_{Q}^+$ l'ensemble des racines de $A_{M_{\natural}}$ dans $\mathfrak{u}_{Q}$.

\ass{Lemme}{Soit $c>0$. Il existe $c'>0$ tel que la propri\'et\'e suivante soit v\'erifi\'ee. Soient $a\in A_{T''}(F)$, $g\in G(F)$, $\bar{u}\in U_{\bar{Q}}(F)$ et $X\in \omega_{T''}\cap\mathfrak{t}''(F)^S[>N^{-b}]$. On suppose $\sigma(g),\,\sigma(\bar{u}),\,\sigma(\bar{u}g)<clog(N)$ et $\alpha(H_{M_{\natural}}(a))>c'log(N)$ pour tout $\alpha\in \Sigma_{Q}^+$. Alors on a l'\'egalit\'e $\kappa_{N}(\gamma_{X''}^{-1}a\bar{u}g)=\kappa_{N}(\gamma_{X''}^{-1}ag)$.}

Preuve. Soient $E$ une extension finie de $F$  v\'erifiant la condition du lemme pr\'ec\'edent et $F'$ une extension finie de $E$. Fixons un r\'eseau $R'$ de $V_{F'}$ v\'erifiant la conclusion de ce lemme. Comme on l'a dit, on peut supposer qu'il v\'erifie des propri\'et\'es analogues \`a celles que l'on a impos\'ees \`a $R$. A l'aide de ce r\'eseau, on construit la fonction $\kappa_{N}^{F'}$ sur $G(F')$ analogue \`a $\kappa_{N}$. On v\'erifie que la restriction \`a $G(F)$ de la fonction $\kappa_{Nval_{F'}(\varpi_{F})}^{F'}$ est \'egale \`a $\kappa_{N}$. Le lemme se d\'eduit alors du m\^eme lemme o\`u le corps des scalaires a \'et\'e \'etendu \`a $F'$. On va maintenant oublier ces constructions, en retenant que l'on a le droit d'\'etendre le corps $F$. En particulier, on peut supposer les tores $T$ et $H_{S}''$ d\'eploy\'es.

Dans tout ce qui suit, les \'el\'ements $X''$ sont implicitement suppos\'es appartenir \`a $ \omega_{T''}\cap  \mathfrak{t}''(F)^{S}$. Montrons  que l'on peut changer de fonction $X''\mapsto \gamma_{X''}$. Consid\'erons une fonction $X''\mapsto \gamma_{X''}$  soumise aux conditions (1) et (2) ci-dessus, et notons ici $X''\mapsto \underline{\gamma}_{X''}$ la fonction initiale, soumise aux conditions de 9.7.  On pose $X''_{\Sigma}=\gamma_{X''}^{-1}X''\gamma_{X''}$ et $\underline{X}''_{\Sigma}=\underline{\gamma}_{X''}^{-1}X''\underline{\gamma}_{X''}$. D'apr\`es le lemme 9.5, pour tout $X''\in \mathfrak{t}''(F)^{S}$, il y a d'uniques \'el\'ements $u(X'')\in U''(F)$ et $t(X'')\in H''_{S}(F)$ tels que $X''_{\Sigma}=u(X'')^{-1}t(X'')^{-1}\underline{X''}_{\Sigma}t(X'')u(X'')$. On a $t(X'')^{-1}\underline{X''}_{\Sigma}t(X'')\in \Xi+S+\Lambda$. D'apr\`es le lemme 9.3, les coefficients de $u(X'')$ et $t(X'')^{-1}\underline{X''}_{\Sigma}t(X'')$ sont polynomiaux en ceux de $X''_{\Sigma}$. Ce dernier terme  reste dans un compact, donc $u(X'')$ et $t(X'')^{-1}\underline{X''}_{\Sigma}t(X'')$ sont born\'es. Reprenons la preuve du lemme 10.1, pr\'ecis\'ement celle de l'assertion (3). On voit que le fait que $t(X'')^{-1}\underline{X''}_{\Sigma}t(X'')$ soit born\'e entra\^{\i}ne une majoration $\sigma(t(X''))<<1+\vert log\vert Q_{S}(X'')\vert _{F}\vert $. Pour $X''\in\mathfrak{t}''(F)^S[>N^{-b}]$, on a donc $\sigma(t(X''))<<log(N)$. Puisque les images de $X''$ par les conjugaisons par $\gamma_{X''}$ et par $\underline{\gamma}_{X''}t(X'')u(X'')$ sont \'egales, il existe $ y(X'')\in T''(F)$ tel que $\gamma_{X''}=y(X'')\underline{\gamma}_{X''}t(X'')u(X'')$. Les majorations (2) pour nos deux fonctions et celles que nous venons de d\'emontrer impliquent   $\sigma(y(X''))<<log(N)$ pour $X\in\mathfrak{t}(F)[S;>N^{-b}]$. Soient $c$, $a$, $g$, $\bar{u}$ et $X''$ comme dans le lemme. Puisque $\kappa_{N}$ est invariante \`a gauche par $U(F)H(F)$, on a
$$\kappa_{N}(\gamma_{X''}^{-1}a\bar{u}g)=\kappa_{N}(\underline{\gamma}_{X''}^{-1}a\bar{u}'g'),$$
$$\kappa_{N}(\gamma_{X''}^{-1}ag)=\kappa_{N}(\underline{\gamma}_{X''}^{-1}ag'),$$
   o\`u $g'=y(X'')^{-1}g$, $\bar{u}'=y(X'')^{-1}\bar{u}y(X'')$. Il existe $\underline{c}>0$ (ind\'ependant des variables) tel que $\sigma(g'),\,\sigma(\bar{u}')\,\sigma(u'g')<\underline{c}log(N)$. Supposons le lemme d\'emontr\'e pour la fonction $X''\mapsto \underline{\gamma}_{X''}$. A $\underline{c}$, ce lemme associe une constante $\underline{c}'$. Les \'egalit\'es ci-dessus montrent que, pour la fonction $X''\mapsto \gamma_{X''}$, la conclusion du lemme est valide pour la constante $c'=\underline{c}'$. Le m\^eme calcul s'applique dans l'autre sens: la validit\'e du lemme pour la fonction $X''\mapsto \gamma_{X''}$ entra\^{\i}ne sa validit\'e pour la fonction $X''\mapsto \underline{\gamma}_{X''}$.

 D\'emontrons maintenant le lemme dans le cas o\`u $r=0$.  Comme en 9.4, on introduit des coordonn\'ees en posant 
 $$X''_{\Sigma}=S+c(v_{0},[z_{0}w_{S}]+\sum_{j=\pm 1,...,\pm m}z_{j}w_{j}),$$
  o\`u, comme en 9.6, les termes entre crochets n'existent que si $d$ est pair.  Le tore $T''$ \'etant d\'eploy\'e, on peut introduire un syst\`eme hyperbolique maximal $(\epsilon_{\pm k})_{k=1,...,l}$ de $V''_{0}$ form\'e de vecteurs propres pour $T''$. On note $x_{k}$ la valeur propre de $X''$ sur $\epsilon_{k}$. Soit $k=1,...,l$, posons
  $$\gamma_{X''}^{-1}\epsilon_{k}=Y(2\nu_{0})^{-1}v_{0}+[y_{0}w_{S}]+\sum_{j=\pm 1,...\pm m}y_{j}w_{j},$$ 
    $$\gamma_{X''}^{-1}\epsilon_{-k}=Y'(2\nu_{0})^{-1}v_{0}+[y'_{0}w_{S}]+\sum_{j=\pm 1,...\pm m}y'_{j}w_{j}.$$ 
Ces \'el\'ements sont vecteurs propres de $X''_{\Sigma}$, de valeurs propres $x_{k}$ et $-x_{-k}$. Cela se traduit par les \'egalit\'es
$$s_{j}y_{j}+z_{j}Y=x_{k}y_{j},\,\,-s_{j}y_{-j}+z_{-j}Y=x_{k}y_{-j},\,\,[z_{0}Y=x_{k}y_{0}],$$
$$s_{j}y'_{j}+z_{j}Y'=-x_{k}y'_{j},\,\,-s_{j}y'_{-j}+z_{-j}Y'=-x_{k}y'_{-j},\,\,[z_{0}Y'=-x_{k}y'_{0}],$$
o\`u $j=1,...,m$ et les $s_{j}$ sont les valeurs propres de $S$. Supposons $X\in \mathfrak{t}(F)[S;>N^{-b}]$ Alors $Q_{S}(X'')\not=0$, a fortiori $x_{k}\pm s_{j}\not=0$ pour tous $j,k$ et les \'egalit\'es ci-dessus impliquent
$$y_{ j}=\frac{z_{j}Y}{x_{k}-s_{j}},\,\,y_{-j}=\frac{z_{-j}Y}{x_{k}+s_{j}}\,\,[y_{0}=\frac{z_{0}Y}{x_{k}}],$$  
$$y'_{ j}=\frac{z_{j}Y'}{-x_{k}-s_{j}},\,\,y'_{-j}=\frac{z_{-j}Y'}{-x_{k}+s_{j}}\,\,[y'_{0}=\frac{z_{0}Y'}{-x_{k}}],$$    
On a d'autre part l'\'egalit\'e $q(\gamma_{X''}^{-1}\epsilon_{k},\gamma_{X''}^{-1}\epsilon_{-k})=1$. Les formules pr\'ec\'edentes traduisent cette \'egalit\'e par $YY'P(X'')=2\nu_{0}$, o\`u
$$P(X'')=1-[4\nu_{0}\nu_{S}\frac{z_{0}^2}{x_{k}^2}]-2\nu_{0}\sum_{j=1,...,m}z_{j}z_{-j}(\frac{1}{(x_{k}-s_{j})^2}+\frac{1}{(x_{k}+s_{j})^2}).$$ 
En utilisant les \'egalit\'es  9.4(1) et (2), on obtient
$$P(X'')=1+[ \frac{\prod_{k'=1,...,l; k'\not=k}x_{k'}^2}{\prod_{j=1,...,m}s_{j}^2}]-\sum_{j=1,...,m}\frac{(x_{k}^2+s_{j}^2)\prod_{k'=1,...,l; k'\not=k}(s_{j}^2-x_{k'}^2)}{(s_{j}^2-x_{k}^2)[s_{j}^2]\prod_{j'=1,...,m;j'\not=j}(s_{j}^2-s_{j'}^2)}.$$
Cette expression est calcul\'ee par le lemme 10.6. En effet, si $d$ est impair, on a $l=m$. On prend pour polyn\^ome $R(T)=\prod_{k'=1,...,l; k'\not=k}(T- x_{k'}^2)$ et on remplace les ind\'etermin\'ees $T$, $S_{1}$, ...,$S_{l}$ de 10.6 par $x_{k}^2$, $s_{1}^2$,...,$s_{m}^2$. Si $d$ est pair, on a $l=m+1$. On prend le m\^eme polyn\^ome $R(T)$ et on remplace les ind\'etermin\'ees par $x_{k}^2$, $0$, $s_{1}^2$,...,$s_{m}^2$. Le lemme 10.6 conduit ainsi \`a l'\'egalit\'e
$$(3) \qquad YY'=2\nu_{0}\frac{R_{1}(X'')}{R_{2}(X'')},$$
o\`u
$$(4)\qquad R_{1}(X'')=\prod_{j=1,...,m}(s_{j}^2-x_{k}^2),$$
 $$R_{2}(X'')=\left\lbrace\begin{array}{cc}-2x_{k}^2\prod_{k'=1,...,l;k'\not=k}(x_{k'}^2-x_{k}^2),&\,\,\text{si}\,\,d\,\,\text{est impair,}\\ 2\prod_{k'=1,...,l;k'\not=k}(x_{k'}^2-x_{k}^2),&\,\,\text{si}\,\,d\,\,\text{est pair.}\\ \end{array}\right.$$
Rappelons que
$$Q_{S}(X'')=\prod_{j=1,...,m; k'=1,...,l}(s_{j}^2-x_{k'}^2).$$
On peut donc aussi \'ecrire $YY'Q(X'')=Q_{S}(X'')$, o\`u $Q$ est un polyn\^ome  sur $\mathfrak{t}''(F)$. 
 Puisque $X''$ reste dans un compact, $val_{F}(Q(X''))$ est minor\'e, donc $val_{F}(YY')<<val_{F}(Q_{S}(X''))$. Supposons $X''\in \mathfrak{t}''(F)^S[>N^{-b}]$. Alors $val_{F}(Q_{S}(X''))<<log(N)$. On a aussi $\sigma(\gamma_{X''})<<log(N)$, donc $val_{F}(Y), val_{F}(Y')>>-log(N)$. On en d\'eduit 
 $$val_{F}(Y), val_{F}(Y')<<log(N).$$
  Remarquons que $Y=q(\gamma_{X''}^{-1}\epsilon_{k},v_{0})=q(\epsilon_{k},\gamma_{X''}v_{0})$ et de m\^eme $Y'=q(\epsilon_{-k},\gamma_{X''}v_{0})$. Cela d\'emontre que

(5) on a $val_{F}(q(\epsilon_{\pm k},\gamma_{X''}v_{0}))<<log(N)$  pour tout $k=1,...,l$ et tout $X\in \omega_{T''}\cap\mathfrak{t}''(F)^S[>N^{-b}]$.

Consid\'erons l'ensemble $E_{N}$ des $a\in A_{T''}(F)$ tels qu'il existe $g\in G(F)$ et $X''\in \omega_{T''}\cap\mathfrak{t}''(F)^S[>N^{-b}]$ de sorte que $\sigma(g)<c\,log(N)$ et $\kappa_{N}(\gamma_{X''}^{-1}ag)=1$. Puisque $T''$ est d\'eploy\'e, on a $A_{T''}=T''$ et tout $a\in A_{T''}(F)$ est d\'etermin\'e par ses valeurs propres $a_{k}$ sur les vecteurs $\epsilon_{k}$. Montrons que

(6) on peut fixer $c_{2}>0$ tel que

(i) $val_{R}(g^{-1}v)-val_{R}(v)>-c_{2}log(N)$ pour tout $v\in V$, $v\not=0$, et tout $g\in G(F)$ tel que $\sigma(g)<c\,log(N)$;

(ii) $val_{R}(\gamma_{X''}v_{0})>-c_{2}log(N)$ pour tout $X''\in \omega_{T''}\cap \mathfrak{t}''(F)^S[>N^{-b}]$;

(iii) $val_{R}(a^{-1}v)-val_{R}(v)>-N-c_{2}log(N)$ pour tout $v\in V$, $v\not=0$, et tout $a\in E_{N}$.

 Les deux premi\`eres majorations sont ais\'ees. Montrons la troisi\`eme. Soient $a\in E_{N}$, $g$ tel que $\sigma(g)<c\,log(N)$  et $X''\in \omega_{T''}\cap \mathfrak{t}''(F)^S[>N^{-b}]$ tels que $\kappa_{N}(\gamma_{X''}^{-1}ag)=1$. Cette condition signifie que $val_{R}(g^{-1}a^{-1}\gamma_{X''}v_{0})\geq -N$. Gr\^ace \`a (6)(i), elle entra\^{\i}ne $val_{R}(a^{-1}\gamma_{X''}v_{0})+N>> -log(N)$. Donc $$val_{F}(q(a^{-1}\gamma_{X''}v_{0},\epsilon_{\pm k}))+N>>-log(N)$$
pour tout $k$. On a 
$$q(a^{-1}\gamma_{X''}v_{0},\epsilon_{\pm k})=q(\gamma_{X''}v_{0},a\epsilon_{\pm k})=a_{k}^{\pm 1}q(\gamma_{X''}v_{0},\epsilon_{\pm k}).$$ 
La majoration pr\'ec\'edente  et (5) entra\^{\i}nent
$$\pm val_{F}(a_{k})+N >>-log(N),$$
et on en d\'eduit la majoration (6)(iii).

Enfin

(7) il existe $c'>0$ v\'erifiant la propri\'et\'e suivante; pour tout $a\in A_{T''}(F)$ tel que $\alpha(H_{M_{\natural}}(a))>c'log(N)$ pour tout $\alpha\in \Sigma_{Q}^+$, pour tout $\bar{u}\in U_{\bar{Q}}(F)$ tel que $\sigma(\bar{u})<c\,log(N)$ et tout $v\in V$, $v\not=0$, on a $val_{R}(a\bar{u}a^{-1}v-v)>3c_{2}log(N)$.

En effet, les valuations des coefficients de $\bar{u}-1$, disons dans une base de $R$, sont minor\'es par $-c_{3}log(N)$, pour une constante $c_{3}$ convenable. On en d\'eduit que les valuations des  coeffients de $a\bar{u}a^{-1}-1$ sont minor\'es par $-c_{3}log(N)+inf\{\alpha(H_{M_{\natural}}(a)); \alpha\in \Sigma_{Q}^+\}$. L'assertion s'ensuit.

La constante $c'$ \'etant maintenant fix\'ee, soient $a$, $g$, $\bar{u}$ et $X''$ comme dans l'\'enonc\'e. Si $a\not\in E_{N}$, on a $\kappa_{N}(\gamma_{X''}^{-1}a\bar{u}g)=\kappa_{N}(\gamma_{X''}^{-1}ag)=0$. Supposons $a\in E_{N}$. D'apr\`es la d\'efinition de $\kappa_{N}$, pour prouver l'\'egalit\'e $\kappa_{N}(\gamma_{X''}^{-1}a\bar{u}g)=\kappa_{N}(\gamma_{X''}^{-1}ag)$, il suffit de prouver que $val_{R}(v)\geq -N$, o\`u
$$v=g^{-1}\bar{u}^{-1}a^{-1}\gamma_{X''}v_{0}-g^{-1}a^{-1}\gamma_{X''}v_{0}.$$
On pose $v_{1}=gv$, $v_{2}=av_{1}$, $v_{3}=\gamma_{X''}v_{0}$. On a $v_{2}=a\bar{u}^{-1}a^{-1}v_{3}-v_{3}$ et
$$val_{R}(v)=val_{R}(v)-val_{R}(v_{1})+val_{R}(v_{1})-val_{R}(v_{2})+val_{R}(v_{2})-val_{R}(v_{3})+val_{R}(v_{3}).$$
Les minorations (6) et (7) entra\^{\i}nent la minoration $val_{R}(v)\geq -N$ cherch\'ee.

Passons au cas g\'en\'eral o\`u on ne suppose plus $r=0$. On peut fixer un \'el\'ement $P_{\natural}=M_{\natural}U_{\natural}\in \tilde{{\cal P}}(M_{\natural})$ et se borner \`a consid\'erer des $a\in A_{T''}(F)$ tels que $H_{M_{\natural}}(a)\in Cl({\cal A}_{P_{\natural}}^+)$. Si $P_{\natural}$ n'est pas inclus dans $Q$, l'assertion \`a prouver est vide car il n'y a pas de tels $a$ pour lesquels $\alpha(H_{M_{\natural}}(a))>0$ pour tout $\alpha\in \Sigma_{Q}^+$. On suppose donc $P_{\natural}\subset Q$. Montrons que

(8) il existe $\delta\in G''(F)$ tel que $ \delta P_{\natural}\delta^{-1}\subset \bar{P}$ et $A\subset \delta A_{T''}\delta^{-1}$.

 Puisque $A\subset G''$ et $A_{T''}=T''$ est un sous-tore maximal de $G''$, on peut en tout cas trouver $\delta\in G''(F)$ tel que $\delta^{-1}A\delta \subset  A_{T''}$. On a alors $\delta^{-1}\bar{P}\delta\in{\cal F}(M_{\natural})$. Supposons que l'on n'est pas dans le cas exceptionnel. Fixons un \'el\'ement $P'\in {\cal P}(M_{\natural})$ tel que $P'\subset \delta^{-1}\bar{P}\delta$. Le L\'evi $M_{\natural}$ a une forme particuli\`ere: il est produit d'un groupe sp\'ecial orthogonal et de groupes $GL(1)$. On sait qu'alors deux \'el\'ements de ${\cal P}(M_{\natural})$ sont conjugu\'es par un \'el\'ement de $Norm_{G(F)}(M_{\natural})$. Si  $d$ est impair ou si $W'=\{0\}$, l'application naturelle
$$Norm_{G''(F)}(A_{T''})\to Norm_{G(F)}(M_{\natural})/M_{\natural}(F)$$
est surjective. Donc $P_{\natural}$ et $P'$ sont conjugu\'es par un \'el\'ement de $Norm_{G''(F)}(A_{T''})$. Quitte \`a multiplier $\delta$ \`a droite par un \'el\'ement de cet ensemble, on peut supposer $P'=P_{\natural}$ et la conclusion de (6) est v\'erifi\'ee. Supposons $d$ pair et $W'\not=\{0\}$. Fixons $w'\in W'$ tel que $q(w')\not=0$. Identifions $G^{_{''}+}$ \`a un sous-groupe de $G$ en faisant agir un \'el\'ement
$g\in G^{_{''}+}$ par $det(g)$ sur $Fw'$ et par l'identit\'e sur l'orthogonal de $w'$ dans $W'$. Alors l'application naturelle
 $$Norm_{G^{_{''}+}(F)}(A_{T''})\to Norm_{G(F)}(M_{\natural})/M_{\natural}(F)$$
 est surjective, donc $P_{\natural}=g^{-1}P'g$ pour un \'el\'ement $g\in Norm_{G^{_{''}+}(F)}(A_{T''})$. Si $g\in G''(F)$, on conclut comme ci-dessus. Sinon, on remarque que $A$ n'est pas un sous-tore maximal de $G''$, car $A$ fixe le vecteur $v_{0}\in V''$. L'inclusion $\delta^{-1}A\delta \subset  A_{T''}$ est stricte et on voit qu'il existe un \'el\'ement $y\in Norm_{G^{_{''}+}(F)}(A_{T''})$ tel que $det(y)=-1$ et la conjugaison par $y$ fixe tout point de $\delta^{-1}A\delta$. Cette conjugaison conserve donc $P'$ et on peut remplacer $g$ par $yg$, ce qui nous ram\`ene au cas pr\'ec\'edent.
 
 Dans le cas exceptionnel, on v\'erifie que $\delta^{-1}\bar{P}\delta$ appartient \`a $\tilde{{\cal F}}(M_{\natural})$. On remplace dans la preuve ci-dessus les ensembles ${\cal P}(M_{\natural})$ et $Norm_{G(F)}(M_{\natural})$ par $\tilde{{\cal P}}(M_{\natural})$ et $Norm_{G(F)}(A_{T''})$. En utilisant le lemme 10.3, on voit que le raisonnement reste valable. Cela prouve (8).
 
 Soient $\delta$ v\'erifiant (8) et $a$, $g$, $\bar{u}$ et $X''$ comme dans l'\'enonc\'e. On a
 $$\kappa_{N}(\gamma_{X''}^{-1}a\bar{u}g)=\kappa_{N}(\underline{\gamma}^{-1}_{\underline{X}''}\underline{a}\underline{\bar{u}}\underline{g}),$$
 o\`u $\underline{X}''=\delta X''\delta^{-1}$, $\underline{\gamma}_{\underline{X}''}=\delta\gamma_{X''}$, $\underline{a}=\delta a\delta^{-1}$, $\underline{\bar{u}}=\delta \bar{u}\delta^{-1}$, $\underline{g}=\delta g$. On remarque que ces termes v\'erifient les m\^emes hypoth\`eses que les termes de d\'epart (avec une autre constante $c$), mais avec le tore $T''$ et le sous-groupe parabolique $P_{\natural}$ remplac\'es par $\delta T''\delta^{-1}$ et $\delta P_{\natural} \delta^{-1}$. Cela nous ram\`ene \`a d\'emontrer le lemme pour ces nouveaux termes. 
 
 On va plut\^ot oublier ces constructions, mais supposer que nos objets de d\'epart v\'erifient les m\^emes hypoth\`eses que ceux que l'on vient de construire. C'est-\`a-dire que l'on suppose d\'esormais $P_{\natural}\subset \bar{P}$ et $A\subset A_{T''}$. Le tore $T''$ se d\'ecompose en $T''=AT''_{0}$, o\`u $T''_{0}$ est un sous-tore maximal de $G''_{0}$. En travaillant dans ce groupe $G''_{0}$, on d\'efinit l'ensemble $\mathfrak{t}''_{0}(F)^S$ et une fonction $X''_{0}\mapsto \gamma_{0,X''_{0}}\in G''_{0}(F)$ sur cet ensemble, v\'erifiant les analogues de (1) et (2). Pour $X''\in \mathfrak{t}''(F)^S$, \'ecrivons $X''=X_{a}''+X''_{0}$, avec $X''_{a}\in \mathfrak{a}(F)$ et $X''_{0}\in \mathfrak{t}''_{0}(F)$. On v\'erifie que $X''_{0}\in \mathfrak{t}''_{0}(F)^S$. Posons $X''_{\Sigma}=\Xi+X''_{a}+\gamma_{0,X''_{0}}^{-1}X''_{0}\gamma_{0,X''_{0}}$. On a $X''_{\Sigma}\in \Xi+S+\Sigma$ et, comme dans la preuve du lemme 9.6, on montre que $X''_{\Sigma}$ est conjugu\'e \`a $X''$ par un \'el\'ement de $G''(F)$. D'autre part, consid\'erons l'\'el\'ement $\gamma_{0,X''_{0}}\Xi\gamma_{0,X''_{0}}^{-1}$. Il appartient \`a $ \bar{\mathfrak{u}}(F)$. Puisque $X''$ est un \'el\'ement r\'egulier de $\mathfrak{m}''(F)$,   il existe $v_{X''}\in \bar{U}(F)$ tel que $\gamma_{0,X''_{0}}\Xi\gamma_{0,X''_{0}}^{-1}=v_{X''}^{-1}X''v_{X''}$. Cet \'el\'ement $v_{X''}$ est unique. Ainsi qu'il est bien connu, ses coefficients sont  des fractions rationnelles en les coefficients de $\gamma_{0,X''_{0}}\Xi\gamma_{0,X''_{0}}^{-1}$ et $X''$  et les d\'enominateurs de ces fractions rationnelles divisent  le polyn\^ome $det(X''\vert \mathfrak{g}''/\mathfrak{t}'')$. Pour $X''\in \omega_{T''}\cap\mathfrak{t}''(F)^S[>N^{-b}]$, on a donc une majoration $\sigma(v_{X''})<<log(N)$. Posons $\gamma_{X''}=v_{X''}\gamma_{0,X''_{0}}$. On a alors $\gamma_{X''}^{-1}X''\gamma_{X''}=X''_{\Sigma}$ et on voit que l'application $X''\mapsto \gamma_{X''}$ v\'erifie les propri\'et\'es (1) et (2). On peut travailler avec cette application. On a $v_{X''}\in U_{\natural}(F)$. D\'ecomposons cet \'el\'ement en $v_{X''}=n_{X''}\nu_{X''}$, avec $n_{X''}\in U_{\natural}(F)\cap L(F)$ et $\nu_{X''}\in U_{Q}(F)$. Soient $a$, $X''$, $g$ et $\bar{u}$ comme dans l'\'enonc\'e. On a $\gamma_{X''}^{-1}a\bar{u}g=\gamma_{0,X_{0}''}^{-1}a\bar{u}'g'k)$, o\`u $\bar{u}'=(a^{-1}n_{X''}a)^{-1}\bar{u}(a^{-1}n_{X''}a)$, $g'=a^{-1}n_{X''}^{-1}ag$, $k=g^{-1}\bar{u}^{-1}a^{-1}\nu_{X''}a\bar{u}g$. Puisque $\nu_{X''}\in U_{Q}(F)$ et ainsi qu'on l'a vu au cours de la preuve du cas $r=0$, on peut fixer $c_{4}>0$ tel que la condition $inf\{\alpha(H_{M_{\natural}}(a)); \alpha\in \Sigma_{Q}^+\}>c_{4}log(N)$ entra\^{\i}ne que les valuations des coefficients de $k-1$ soient $>>log(N)$. On impose cette condition sur $a$. Alors $k\in K$. La conjugaison par $a^{-1}$ contracte $U_{\natural}(F)$, donc $\sigma(a^{-1}n_{X''}a)<<log(N)$ et aussi $\sigma(\bar{u}')<<log(N)$, $\sigma(g')<<log(N)$. On a donc $\kappa_{N}(\gamma_{X''}^{-1}a\bar{u}g)=\kappa_{N}(\gamma_{0,X''_{0}}^{-1}a\bar{u}'g')$ et de m\^eme $\kappa_{N}(\gamma_{X''}^{-1}ag)=\kappa_{N}(\gamma_{0,X''_{0}}^{-1}ag')$, o\`u les \'el\'ements $\bar{u}'$ et $g'$ v\'erifient des conditions analogues \`a celles de d\'epart. Soient $\bar{u}_{0}\in U_{\bar{Q}}(F)\cap G_{0}(F)$ tel que $\bar{u}'\in (U_{\bar{Q}}(F)\cap U(F))\bar{u}_{0}$, $y'\in A(F)$ et $g_{0}\in G_{0}(F)$ tels que $g'\in U(F)y'g_{0}K$ et enfin $y\in A(F)$ et $a_{0}\in A_{T''_{0}}(F)$ tels que $a=y a_{0}$. Alors $\gamma_{0,X''_{0}}^{-1}a\bar{u}'g'\in U(F)yy'\gamma_{0,X''_{0}}^{-1}a_{0}\bar{u}_{0}g_{0}K$, donc
 $$\kappa_{N}(\gamma_{X''}^{-1}a\bar{u}g)=\kappa_{N}(\gamma_{0,X''_{0}}^{-1}a\bar{u}'g')=\kappa_{N}(yy'\gamma_{0,X''_{0}}^{-1}a_{0}\bar{u}_{0}g_{0}).$$
  De m\^eme
$$\kappa_{N}(\gamma_{X''}^{-1}ag)=\kappa_{N}(\gamma_{0,X''_{0}}^{-1}ag')=\kappa_{N}(yy'\gamma_{0,X''_{0}}^{-1}a_{0}g_{0}).$$
 Par d\'efinition de $\kappa_{N}$, ces expressions se r\'ecrivent 
 $$\kappa_{A,N}(yy')\kappa_{0,N}(\gamma_{0,X''_{0}}^{-1}a_{0}\bar{u}_{0}g_{0})\text{, resp.} \,\,\kappa_{A,N}(yy')\kappa_{0,N}(\gamma_{0,X''_{0}}^{-1}a_{0}g_{0}),$$
 o\`u $\kappa_{A,N}$ est une certaine fonction sur $A(F)$ et $\kappa_{0,N}$ est l'analogue de $\kappa_{N}$ pour le groupe $G_{0}$. Mais les donn\'ees affect\'ees d'un indice $0$ v\'erifient des conditions similaires \`a celles de d\'epart. Cela nous ram\`ene au cas du groupe $G_{0}$, autrement dit au cas $r=0$ que nous avons d\'ej\`a trait\'e. Cela ach\`eve la d\'emonstration. $\square$

  D\'emontrons 10.4(11). Soit $c_{4}>0$. On impose \`a $Z_{P_{min}}$ la minoration $c_{4}log(N)<c_{1}inf\{\alpha(Z_{P_{min}}); \alpha\in \Delta_{min}\}$ pour que tous les termes ci-dessous soient d\'efinis. Comme en 10.5, pour $g$ et $\bar{u}$ comme en 10.4(11), la fonction
 $$\zeta\mapsto  \tilde{\sigma}_{M_{\natural}}^Q(\zeta,{\cal Z}(g))\tilde{\tau}_{Q}(\zeta-Z(g)_{Q,T''})$$
 est insensible au changement de $g$ en $\bar{u}g$. Alors
 $$\kappa_{N,X''}(Q,\bar{u}g)-\kappa_{N,X''}(Q,g)=\int_{A_{T''}(F)}\tilde{\sigma}_{M_{\natural}}^Q(H_{M_{\natural}}(a),{\cal Z}(g))\tilde{\tau}_{Q}(H_{M_{\natural}}(a)-Z(g)_{Q,T''})$$
 $$\qquad (\kappa_{N}(\gamma_{X''}^{-1}a\bar{u}g)-\kappa_{N}(\gamma_{X''}^{-1}ag))da.$$
 Il nous suffit que la condition  $c_{5}log(N)<inf\{\alpha(Z_{P_{min}}); \alpha\in \Delta_{min}\}$ entra\^{\i}ne la propri\'et\'e suivante. Soit $a\in A_{T''}(F)$ tel que 
 $$(9) \qquad \tilde{\sigma}_{M_{\natural}}^Q(H_{M_{\natural}}(a),{\cal Z}(g))\tilde{\tau}_{Q}(H_{M_{\natural}}(a)-Z(g)_{Q,T''})\not=0.$$
 Alors $\kappa_{N}(\gamma_{X''}^{-1}a\bar{u}g)=\kappa_{N}(\gamma_{X''}^{-1}ag)$.   Prenons $c=c_{4}$ dans le lemme pr\'ec\'edent. On en d\'eduit une constante $c'$. Le m\^eme calcul qu'en 10.5 montre que  (9) implique 
 $$   inf\{\alpha(H_{M_{\natural}}(a)); \alpha\in \Sigma_{Q}^+\}-inf\{\alpha(Z_{P_{min}}); \alpha\in \Delta_{min}\}>>-log(N).$$
 Il existe donc $c_{5}>0$ (et $c_{5}>c_{4}/c_{1}$) tel que la condition $c_{5}log(N)<inf\{\alpha(Z_{P_{min}}); \alpha\in \Delta_{min}\}$ entra\^{\i}ne $inf\{\alpha(H_{M_{\natural}}(a)); \alpha\in \Sigma_{Q}^+\}>c'log(N)$. Alors l'\'egalit\'e cherch\'ee est la conclusion du lemme ci-dessus. $\square$
   
\bigskip

\subsection{Apparition des int\'egrales orbitales pond\'er\'ees}

On fixe $S\in {\cal S}$ et $T\in {\cal T}(G_{x})$. Soit $N_{0}$ l'entier d\'etermin\'e par la proposition 10.4.

\ass{Proposition}{Pour tout  $N\geq N_{0}$ et tout $X\in \mathfrak{t}(F)[S;>N^{-b}]\cap(\mathfrak{t}'(F)\times \mathfrak{t}''(F)^S)$, on a les \'egalit\'es
$$ \int_{T'(F)A_{T''}(F)\backslash G(F)}{^gf}^{\sharp}_{x,\omega}(X)\kappa_{N,X''}(g)dg=0$$
si $A_{T'}\not=\{1\}$;  
$$ \int_{T'(F)A_{T''}(F)\backslash G(F)}{^gf}^{\sharp}_{x,\omega}(X)\kappa_{N,X''}(g)dg=\nu(T') \nu(A_{T''})\theta^{\sharp}_{f,x,\omega}(X)$$
si $A_{T'}=\{0\}$.}

Preuve. En 10.4, on avait fix\'e $Y_{P_{min}}$ et construit une fonction $g\mapsto \tilde{v}(g)$.  Il convient maintenant de la noter plus pr\'ecis\'ement  $g\mapsto \tilde{v}(g,Y_{P_{min}})$.  Soit $X\in \mathfrak{t}(F)[S;>N^{-b}]\cap(\mathfrak{t}'(F)\times \mathfrak{t}''(F)^S)$. Dans les int\'egrales ci-dessus, on peut remplacer $\kappa_{N,X''}(g)$ par $\tilde{v}(g,Y_{P_{min}})$, pourvu que $Y_{P_{min}}$ v\'erifie la minoration de la proposition 10.4. Supposons que l'on n'est pas dans le cas exceptionnel. Alors $\tilde{v}(g,Y_{P_{min}})$ est la fonction introduite par Arthur dans [A3] p.30: avec les notations de cette r\'ef\'erence, c'est $v_{M_{\natural}}(1,g,Y_{P_{min}})$ (il n'y a pas de $\nu(A_{T''})$ dans la d\'efinition d'Arthur, car sa mesure sur $A_{T''}(F)$ n'est pas la m\^eme que la n\^otre). Remarquons que les int\'egrales de l'\'enonc\'e sont \`a support compact. On peut faire tendre $Y_{P_{min}}$ vers l'infini. Alors $\tilde{v}(g,Y_{P_{min}})$ est calcul\'e en [A3] p.46: pour $Y_{P_{min}}$ dans un r\'eseau convenable ${\cal R}\subset {\cal A}_{M_{min}}$, c'est une somme de fonctions $Y_{P_{min}}\mapsto q_{\zeta}(Y_{P_{min}})exp(\zeta(Y_{P_{min}}))$, o\`u   $q_{\zeta}$ est un polyn\^ome et $\zeta\in Hom({\cal R},2\pi i{\mathbb Q}/2\pi i{\mathbb Z})$. De telles fonctions sont lin\'eairement ind\'ependantes. Puisque l'expression que l'on calcule est ind\'ependante de $Y_{P_{min}}$, on peut aussi bien remplacer $\tilde{v}(g,Y_{P_{min}})$ par $q_{0}(0)$. Avec les notations d'Arthur, on a 
$$q_{0}(0)=\tilde{v}_{M_{\natural}}(1,g)=(-1)^{a_{M_{\natural}}}\sum_{Q\in {\cal F}(M_{\natural})}c'_{Q}v_{M_{\natural}}^Q(g),$$
cf. [A3] (6.6) et p.92. Les $c'_{Q}$ sont des constantes et on a $c'_{G}=1$. On a obtenu l'\'egalit\'e
   $$(1) \qquad \int_{T'(F)A_{T''}(F)\backslash G(F)}{^gf}^{\sharp}_{x,\omega}(X)\kappa_{N,X''}(g)dg=(-1)^{a_{M_{\natural}}}\sum_{Q\in {\cal F}(M_{\natural})}c'_{Q}I(Q),$$
   o\`u
   $$I(Q)=\int_{T'(F)A_{T''}(F)\backslash G(F)}{^gf}^{\sharp}_{x,\omega}(X)v_{M_{\natural}}^Q(g)dg.$$
   Pour $Q=LU_{Q}\not=G$, on d\'ecompose l'int\'egrale en produit d'int\'egrales sur $T'(F)A_{T''}(F)\backslash L(F)$, $K_{min}$ et $U_{Q}(F)$. On voit appara\^{\i}tre une int\'egrale
   $$\int_{U_{Q}(F)}{^{ulk}f}^{\sharp}_{x,\omega}(X)du.$$
   Or cette int\'egrale est nulle d'apr\`es le lemme 5.5(i). Donc $I(Q)=0$. Pour $Q=G$, on peut remplacer l'int\'egration sur $T'(F)A_{T''}(F)\backslash G(F)$ par l'int\'egration sur $T(F)\backslash G(F)$, \`a condition de multiplier par $mes(T(F)/T'(F)A_{T''}(F))$. On obtient
   $$I(G)=mes(T(F)/T'(F)A_{T''}(F))D^{G_{x}}(X)^{-1/2}J^{\sharp}_{M_{\natural},x,\omega}(X,f)$$
   avec la notation de 5.4. Si $A_{T'}\not=\{1\}$, on a $A_{M_{\natural}}=A_{T''}\subsetneq A_{G_{x,X}}=A_{T'}A_{T''}$. Donc $J^{\sharp}_{M_{\natural},x,\omega}(X,f)=0$ d'apr\`es le lemme 5.5(ii). Supposons $A_{T'}=\{1\}$. Alors $M_{\natural}$ est le L\'evi not\'e ${\bf M}(X)$ en 5.6 et, en appliquant les d\'efinitions de ce paragraphe, on obtient
   $$I(G)=(-1)^{a_{M_{\natural}}}\nu(T)mes(T(F)/T'(F)A_{T''}(F))\theta^{\sharp}_{f,x,\omega}(X).$$
   On v\'erifie que $\nu(T)mes(T(F)/T'(F)A_{T''}(F))=\nu(T')\nu(A_{T''})$ et la formule (1) devient celle de l'\'enonc\'e.
   
   Supposons maintenant que l'on est dans le cas exceptionnel. Alors $T'=G'$ est un tore d\'eploy\'e de dimension $1$ et on peut supposer $T'\subset M_{min}$.  Il faut remarquer que les fonctions d'Arthur que nous avons utilis\'ees ci-dessus ne d\'ependent de $g$ et $Y_{P_{min}}$ que par l'interm\'ediaire  des familles de points $(H_{P_{\natural}}(g))_{P_{\natural}\in{\cal P}(M_{\natural})}$ et $(Y_{P_{\natural}})_{P_{\natural}\in {\cal P}(M_{\natural})}$. On peut en fait associer de telles fonctions \`a  deux familles $(G,M_{\natural})$-orthogonales de points de ${\cal A}_{M_{\natural}}$, la seconde \'etant "assez positive".  Introduisons le groupe $\tilde{G}$ de 10.3 et rappelons que le lemme de ce paragraphe nous permet d'identifier $\tilde{{\cal P}}(M_{\natural})$ et $\tilde{{\cal F}}(M_{\natural})$ \`a ${\cal P}(\tilde{M})$ et ${\cal F}(\tilde{M})$. En rempla\c{c}ant $G$ par  $\tilde{G}$, $M_{\natural}$ par $\tilde{M}$, la famille $(H_{P_{\natural}}(g))_{P_{\natural}\in{\cal P}(M_{\natural})}$ par $(H_{P_{\natural}}(g)_{T''})_{P_{\natural}\in \tilde{{\cal P }}(M_{\natural})}$ et la famille $(Y_{P_{\natural}})_{P_{\natural}\in {\cal P}(M_{\natural})}$ par $(Y_{P_{\natural},T''})_{P_{\natural}\in \tilde{{\cal P}}(M_{\natural})}$, on remplace la fonction $v_{M_{\natural}}(1,g,Y_{P_{min}})$ utilis\'ee ci-dessus par une autre fonction, notons-la $v_{\tilde{M}}(1,g,Y_{P_{min}})$. Elle est \'egale \`a notre fonction $\tilde{v}(g,Y_{P_{min}})$. Les calculs d'Arthur restent valables pour cette fonction ainsi que les arguments ci-dessus. On obtient la formule (1) modifi\'ee de la fa\c{c}on suivante: la somme est limit\'ee aux $Q\in \tilde{{\cal F}}(M_{\natural})$; les fonctions $v_{M_{\natural}}^Q(g)$ sont remplac\'ees par des fonctions, notons-les $v_{\tilde{M}}^{\tilde{Q}}(g)$  Ce terme est la constante associ\'ee \`a la $(\tilde{G},\tilde{M})$-famille de points $(H_{P_{\natural}}(g)_{T''})_{P_{\natural}\in\tilde{{\cal P}}(M_{\natural}); P_{\natural}\subset Q}$. Le sous-espace ${\cal A}_{T''}$ de ${\cal A}_{M_{\natural}}$ v\'erifie les conditions de [A2] paragraphe 7. Pour $Q=LU_{Q}\in \tilde{{\cal F}}(M_{\natural})$, on peut appliquer  le corollaire 7.2 de [A2] et on obtient une \'egalit\'e
   $$v_{\tilde{M}}^{\tilde{Q}}(g)=\sum_{L'\in {\cal F}^L(M_{\natural})}d(L')v_{M_{\natural}}^{Q'}(g).$$
   Comme en 2.2(3), $Q'$ est un \'el\'ement de ${\cal P}(L')$. La constante $d(L')$ est non nulle si et seulement si on a l'\'egalit\'e
   $${\cal A}_{M_{\natural}}^L=proj_{M_{\natural}}^L({\cal A}_{T''})\oplus {\cal A}_{L'}^L.$$
   Les $L'$ qui interviennent sont tous diff\'erents de $G$: c'est \'evident si $Q\not=G$ puisque $L'\subset L$; si $Q=G$, cela r\'esulte de l'\'egalit\'e ci-dessus et du fait que ${\cal A}_{T''}$ est strictement inclus dans ${\cal A}_{M_{\natural}}$. Le m\^eme argument que dans le cas non exceptionnel montre alors que tous les termes de la formule (1) sont nuls. $\square$
   
   \bigskip
   
   \subsection{La proposition principale}
   
Si $A_{G'_{x}}=\{1\}$, posons
  $$(1) \qquad J_{x,\omega}(\theta,f)= \sum_{S\in {\cal S}}\sum_{T=T'T''\in {\cal T}_{ell}(G'_{x})\times {\cal T}(G'')}\nu(T')\vert W(G_{x},T)\vert ^{-1}$$
  $$\qquad \int_{\mathfrak{t}'(F)\times \mathfrak{t}''(F)^S}\hat{j}_{S}(X')D^{G'_{x}}(X')D^{G''}(X'')^{1/2}\theta^{\sharp}_{f,x,\omega}(X) dX.$$
  Si $A_{G'_{x}}\not=\{1\}$, posons
  $$J_{x,\omega}(\theta,f)=0.$$
  
  \ass{Proposition}{(i) L'expression (1) est absolument convergente.
  
  (ii) On a l'\'egalit\'e
  $$lim_{N\to \infty}I_{x,\omega,N}(\theta,f)=J_{x,\omega}(\theta,f).$$}
  
  Preuve. Les int\'egrales de l'expression (1) sont \`a support compact. D'apr\`es le lemme 5.4(iii), on a une majoration
  $$D^{G_{x}}(X)^{1/2}\vert \theta^{\sharp}_{f,x,\omega}(X)\vert <<(1+\vert log(D^{G_{x}}(X))\vert )^k.$$
  D'apr\`es [HCvD] th\'eor\`eme 13, la fonction
  $$X'\mapsto D^{G'_{x}}(X')^{1/2}\vert \hat{j}_{S}(X')\vert $$
  est born\'ee.  Le lemme 2.4 entra\^{\i}ne le (i) de l'\'enonc\'e.
  
  Pour le (ii), on utilise la derni\`ere formule de 10.1 qui nous ram\`ene \`a prouver que $lim_{N\to \infty}I^*_{x,\omega,N}(\theta,f)=J_{x,\omega}(\theta,f)$. Le terme $I^*_{x,\omega,N}(\theta,f)$ est d\'efini par la formule 10.4(1) o\`u on limite les int\'egrales en $X$ aux ensembles $\mathfrak{t}(F)[S;>N^{-b}]\cap(\mathfrak{t}'(F)\times \mathfrak{t}''(F)^S)$. Pour $N$ assez grand, la proposition pr\'ec\'edente nous permet de remplacer les int\'egrales int\'erieures de cette expression par
  
$$\left\lbrace\begin{array}{c}  0 \text{  si  } A_{T'}\not=\{1\};\\
  
  \nu(T')\nu(A_{T''})\theta^{\sharp}_{f,x,\omega}(X)\text{ si  } A_{T'}=\{1\}.\\ \end{array}\right.$$
  
  Si $A_{G'_{x}}\not=\{1\}$, il n'y a aucun $T'$ tel que $A_{T'}=\{1\}$ donc $I^*_{x,\omega,N}(\theta,f)=0$. Supposons $A_{G'_{x}}=\{1\}$. Alors $I^*_{x,\omega,N}(\theta,f)$ est \'egale \`a l'expression obtenue \`a partir de (1) ci-dessus en limitant les int\'egrales aux ensembles $\mathfrak{t}(F)[S;>N^{-b}]\cap(\mathfrak{t}'(F)\times \mathfrak{t}''(F)^S)$. Quand $N$ tend vers l'infini, cette expression tend vers $J_{x,\omega}(\theta,f)$. $\square$
  
  \bigskip
  
  \section{Etude au voisinage de l'origine}
  
  \bigskip
  
  \subsection{Enonc\'e de la proposition}
  
  Consid\'erons l'hypoth\`ese
  
  \ass{Hypoth\`ese}{Pour tout quasi-caract\`ere $\theta$ sur $\mathfrak{h}(F)$ et toute fonction tr\`es cuspidale $f\in C_{c}^{\infty}(\mathfrak{g}(F))$ dont le support ne contient aucun \'el\'ement nilpotent, on a l'\'egalit\'e
  $$lim_{N\to \infty}I_{N}(\theta,f)=I(\theta,f).$$}
  
  Le but de la  section est de prouver l'assertion suivante.
  
  \ass{Proposition}{ Sous cette hypoth\`ese, on a l'\'egalit\'e
  $$lim_{N\to \infty}I_{N}(\theta,f)=I(\theta,f).$$
   pour tout quasi-caract\`ere $\theta$ sur $\mathfrak{h}(F)$ et toute fonction tr\`es cuspidale $f\in C_{c}^{\infty}(\mathfrak{g}(F))$.}
  
  \bigskip
  
  \subsection{Calcul de $lim_{N\to \infty}I_{N}(\theta,f)$}
  
 Soient $\theta$ un quasi-caract\`ere sur $\mathfrak{h}(F)$ et $f\in C_{c}^{\infty}(\mathfrak{g}(F))$ une fonction tr\`es cuspidale. On peut fixer un ensemble fini ${\cal S}$ d'\'el\'ements de $\mathfrak{h}_{reg}(F)$ et une famille finie $(c_{S})_{S\in {\cal S}}$ de nombres complexes de sorte que
 $$\theta(Y)=\sum_{S\in {\cal S}}c_{S}\hat{j}^H(S,Y)$$
 pour tout $Y\in Supp(f)^G\cap \mathfrak{h}(F)$. On peut supposer que le noyau de l'action de chaque $S$ agissant sur $W$ est de dimension au plus $1$. Pour tout $T\in {\cal T}(G)$, on d\'efinit l'ensemble $\mathfrak{t}(F)^S$ en appliquant la d\'efinition de 9.6 au cas o\`u $V''=V$. On pose
 $$J(\theta,f)=\sum_{S\in {\cal S}}\sum_{T\in {\cal T}(G)}c_{S}\vert W(G,T)\vert ^{-1}\int_{\mathfrak{t}(F)^S}D^G(X)^{1/2}\hat{\theta}_{f}(X)dX.$$
 On a not\'e $\hat{\theta}_{f}$ la transform\'ee de Fourier de $\theta_{f}$, cf. lemme 6.1. A priori,  $J(\theta,f)$ d\'epend des choix des familles ${\cal S}$ et $(c_{S})_{S\in {\cal S}}$. Le lemme suivant montre que ce n'est pas le cas.
 
 \ass{Lemme}{(i) Cette expression est absolument convergente.
 
 (ii) On a l'\'egalit\'e $lim_{N\to \infty}I_{N}(\theta,f)=J(\theta,f)$.}
 
 Preuve. Cet \'enonc\'e n'est qu'une version relative aux alg\`ebres de Lie de la proposition 10.10 . On peut arguer qu'une d\'emonstration analogue \`a celle de cette proposition s'applique. Comme nous aurons besoin plus loin de la construction qui suit, expliquons plut\^ot comment on peut d\'eduire le lemme de cette proposition. Soit $\omega\subset \mathfrak{g}(F)$ un bon voisinage de $0$ (au sens de 3.1 appliqu\'e au cas $x=1$). Supposons d'abord $Supp(f)\subset \omega$.  Posons $\theta_{\omega}=\theta{\bf 1}_{\omega\cap \mathfrak{h}(F)}$. On a  $I_{N}(\theta,f)=I_{N}(\theta_{\omega},f)$.  Par l'exponentielle, on rel\`eve $\theta_{\omega}$ et $f$ en des fonctions $\boldsymbol{\theta}_{\omega}$ sur $H(F)$ et ${\bf f}$ sur $G(F)$, \`a support dans $exp(\omega\cap\mathfrak{h}(F))$, resp. $exp(\omega)$. On a l'\'egalit\'e $I_{N}(\boldsymbol{\theta}_{\omega},{\bf f})=I_{N}(\theta,f)$. On peut appliquer la proposition 10.10  aux fonctions $\boldsymbol{\theta}$ et ${\bf f}$ et au point $x=1$. Donc $lim_{N\to \infty}I_{N}(\theta,f)=J_{1,\omega}(\boldsymbol{\theta}_{\omega},{\bf f})$. Dans cette derni\`ere expression figure une fonction $\theta^{\sharp}_{{\bf f},x,\omega}$ qui n'est autre que $\hat{\theta}_{f}$. En effet, d'apr\`es la proposition 5.8, $\theta^{\sharp}_{{\bf f},x,\omega}$ est la transform\'ee de Fourier partielle de $\theta_{{\bf f},x,\omega}=\theta_{f}$. Mais, pour $x=1$, on a $V''=V$ et la transformation de Fourier partielle est simplement la transform\'ee de Fourier. Alors $J_{1,\omega}(\boldsymbol{\theta}_{\omega},{\bf f})$ co\"{\i}ncide avec $J(\theta,f)$, ce qui prouve le lemme sous l'hypoth\`ese $Supp(f)\subset \omega$.
 
 En g\'en\'eral, rappelons-nous que l'on a fix\'e au d\'epart une famille $(\xi_{i})_{i=0,...,r-1}$ d'\'el\'ements  de $F^{\times}$, dont d\'ependent le caract\`ere $\xi$ et nos constructions.   Soit $\lambda\in F^{\times }$.   
  A la famille $(\lambda\xi_{i})_{i=0,...,r-1}$ est associ\'ee un caract\`ere $\xi'$. Il convient d'affecter d'indices $\xi$, resp.  $\xi'$, les objets construits \`a l'aide du caract\`ere $\xi$, resp. $\xi'$. Posons $\theta'=\theta^{\lambda}$ et $f'=f^{\lambda}$. Comparons les termes $I_{\xi,N}(\theta,f)$ et $J_{\xi}(\theta,f)$ avec leurs analogues $I_{\xi',N}(\theta',f')$ et $J_{\xi'}(\theta',f')$.  Pour $Y\in \mathfrak{h}(F)$, on a
 $$f^{_{'}\xi'}(Y)=\vert \lambda\vert _{F}^{-dim(U)}f^{\xi}(\lambda Y).$$
 On en d\'eduit
 $$(1) \qquad I_{\xi',N}(\theta',f')=\vert \lambda\vert _{F}^{-dim(U)-dim(H)}I_{\xi,N}(\theta,f).$$
  Gr\^ace \`a 2.6(1), on a
  $$\theta'(Y)=\sum_{S\in {\cal S}}c_{S}\hat{j}^H(S,\lambda Y)=\sum_{S\in {\cal S}}\vert \lambda\vert _{F}^{-\delta(H)/2}c_{S}\hat{j}^H(\lambda S,Y)$$
   pour tout $Y\in Supp(f')^G\cap \mathfrak{h}(F)$. Pour d\'efinir $J_{\xi'}(\theta',f')$, on peut donc prendre pour famille ${\cal S}'$ la famille $(\lambda S)_{S\in {\cal S}}$ et pour constantes les $c_{\lambda S}=\vert \lambda\vert _{F}^{-\delta(H)/2}c_{S}$. Soient $T\in {\cal T}(G)$ et $S\in {\cal S}$.  On v\'erifie que, quand on remplace $\xi$ par $\xi'$ et $S$ par $\lambda S$, l'ensemble $\mathfrak{t}(F)^S$ est remplac\'e par $\lambda \mathfrak{t}(F)^S$. Donc
    $$J_{\xi'}(\theta',f')=\sum_{S\in {\cal S}}\sum_{T\in {\cal T}(G)}c_{S}\vert \lambda\vert _{F}^{-\delta(H)/2}\vert W(G,T)\vert ^{-1}\int_{\lambda\mathfrak{t}(F)^S}D^G(X)^{1/2}\hat{\theta}_{f'}(X)dX.$$
  On a  $\hat{f}'(X)=\vert \lambda\vert _{F}^{-dim(G)}\hat{f}(\lambda^{-1}X)$, donc  $\hat{\theta}_{f'} (X)=\vert \lambda\vert _{F}^{-dim(G)}\hat{\theta}_{f}(\lambda^{-1}X)$ gr\^ace au lemme 6.1. On a aussi $D^G(\lambda X)^{1/2}=\vert \lambda\vert _{F}^{\delta(G)/2}D^G(X)^{1/2}$. Par changement de variable, on obtient
  $$(2) \qquad J_{\xi'}(\theta',f')=\vert \lambda\vert ^{-dim(G)+dim(T)+\delta(G)/2-\delta(H)/2}J_{\xi}(\theta,f).$$
  On a 
  $$-dim(G)+dim(T)+\delta(G)/2=-\delta(G)/2=-dim(U)-\delta(G_{0})/2.$$
   A l'aide de 7.7(2), on v\'erifie que $\delta(G_{0})+\delta(H)=2dim(H)$ et l'\'egalit\'e pr\'ec\'edente devient
  $$J_{\xi'}(\theta',f')=\vert \lambda\vert _{F}^{-dim(U)-dim(H)}J_{\xi}(\theta,f).$$
En comparant avec (1), on voit que la relation $lim_{N\to \infty}I_{\xi,N}(\theta,f)=J_{\xi}(\theta,f)$ est \'equivalente \`a $lim_{N\to \infty}I_{\xi',N}(\theta',f')=J_{\xi'}(\theta',f')$. Prenons $\lambda$ tel que $Supp(f')\subset \omega$. Alors la deuxi\`eme relation a d\'ej\`a \'et\'e d\'emontr\'ee (la d\'emonstration  est insensible au changement de $\xi$ en $\xi'$). La premi\`ere s'en d\'eduit, ce qui ach\`eve la preuve. $\square$

\subsection{Une premi\`ere expression du terme d'erreur}

Dans ce paragraphe, on suppose  v\'erifi\'ee l'hypoth\`ese du paragraphe 11.1. Consid\'erons l'application
$$(\theta,f)\mapsto E(\theta,f)=lim_{N\to \infty}I_{N}(\theta,f)-I(\theta,f)=J(\theta,f)-I(\theta,f),$$
d\'efinie sur l'espace des couples $(\theta,f)$ form\'es d'un quasi-caract\`ere $\theta$ sur $\mathfrak{h}(F)$ et d'une fonction tr\`es cuspidale $f\in C_{c}^{\infty}(\mathfrak{g}(F))$. Elle est bilin\'eaire.

\ass{Lemme}{L'application $E$ est combinaison lin\'eaire des applications $(\theta,f)\mapsto c_{\theta,{\cal O}^H}c_{\theta_{f},{\cal O}}$, o\`u ${\cal O}^H$, resp. ${\cal O}$, parcourt l'ensemble des orbites nilpotentes  r\'eguli\`eres de $\mathfrak{h}(F)$, resp. $\mathfrak{g}(F)$.}

Preuve. On commence par prouver

(1) on a $E(\theta,f)=0$ si $c_{\theta,{\cal O}^H}=0$ pour tout  ${\cal O}^H\in Nil(\mathfrak{h}(F))$ ou si $c_{\theta_{f},{\cal O}}=0$ pour tout ${\cal O}\in Nil(\mathfrak{g}(F))$.

Supposons $c_{\theta_{f},{\cal O}}=0$ pour tout ${\cal O}\in Nil(\mathfrak{g}(F))$. On peut alors fixer un $G$-domaine $\omega$ dans $\mathfrak{g}(F)$, compact modulo conjugaison et contenant $0$, tel que $\theta_{f}(X)=0$ pour tout $X\in \omega$. Posons $f'=f{\bf 1}_{\omega}$ et $f''=f-f'$. Ces fonctions sont tr\`es cuspidales. Le support de $f''$ ne contient pas de nilpotent donc $E(\theta,f'')=0$ d'apr\`es l'hypoth\`ese de 11.1. On a $\theta_{f'}=0$ donc aussi $\hat{\theta}_{f'}=0$. Des d\'efinitions r\'esultent les \'egalit\'es $J(\theta,f')=0=I(\theta,f')$. D'o\`u $E(\theta,f')=0$, puis $E(\theta,f)=0$. Supposons maintenant $c_{\theta,{\cal O}^H}=0$ pour tout  ${\cal O}^H\in Nil(\mathfrak{h}(F))$. On peut fixer un $G$-domaine $\omega$ dans $\mathfrak{g}(F)$, compact modulo conjugaison et contenant $0$, tel que $\theta(X)=0$ pour tout $X\in \omega\cap \mathfrak{h}(F)$. D\'efinissons $f'$ et $f''$ comme pr\'ec\'edemment.   Puisque $\theta$ est nul sur $Supp(f')^G\cap \mathfrak{h}(F)$, il r\'esulte des d\'efinitions que $J(\theta,f')=0=I(\theta,f')$. On conclut comme pr\'ec\'edemment. Cela prouve (1).

Soit $\lambda\in F^{\times 2}$, posons $\theta'=\theta^{\lambda}$, $f'=f^{\lambda}$. On a $\theta_{f'}=(\theta_{f})^{\lambda}$. Pour  ${\cal O}^H\in Nil(\mathfrak{h}(F))$ et ${\cal O}\in Nil(\mathfrak{g}(F))$, on a l'\'egalit\'e

$$(2) \qquad c_{\theta',{\cal O}^H}c_{\theta_{f'},{\cal O}}=\vert \lambda\vert ^{-dim({\cal O}^H)/2-dim({\cal O})/2}c_{\theta,{\cal O}^H}c_{\theta_{f},{\cal O}}$$
d'apr\`es 4.2(2).

Montrons que

(3) $E(\theta',f')=\vert \lambda\vert _{F}^{-\delta(G)/2-\delta(H)/2}E(\theta,f)$.
 
Reportons-nous \`a l'\'egalit\'e 11.2(2). Il y figure un terme $J_{\xi'}(\theta',f')$. En fait, il est \'egal \`a $J_{\xi}(\theta',f')$. En effet, ce terme ne d\'epend de $\xi$ que par la d\'efinition des ensembles $\mathfrak{t}(F)^S$, pour $T\in {\cal T}(G)$ et $S\in {\cal S}'$ (l'ensemble associ\'e \`a $\theta'$). Cet ensemble est celui des \'el\'ements de $\mathfrak{t}(F)$ v\'erifiant certaines conditions de r\'egularit\'e et conjugu\'es \`a un \'el\'ement de $\Xi+S+\Sigma$, si le caract\`ere utilis\'e est $\xi$, \`a un \'el\'ement de $\lambda \Xi+S+\Sigma$, si le caract\`ere utilis\'e est $\xi'$. Or, soit $a\in A(F)$ tel que $a_{i}=\lambda^{-i}$ pour tout $i=1,...,r$. Alors $a(\Xi+S+\Sigma)a^{-1}=\lambda \Xi+S+\Sigma$. Donc l'ensemble $\mathfrak{t}(F)^S$ est insensible au changement de $\xi$ en $\xi'$. L'\'egalit\'e 11.2(2) nous dit  alors que
$$(4)\qquad J(\theta',f')=\vert \lambda\vert _{F}^{-\delta(G)/2-\delta(H)/2}J(\theta,f).$$

Soit $T\in {\cal T}$. Introduisons comme en 7.3 la d\'ecomposition $W=W'\oplus W''$ relative \`a $T$ et les notations aff\'erentes. Soit $X\in \mathfrak{t}_{\natural}(F)$. D'apr\`es 4.2(2), on a les \'egalit\'es
$$c_{\theta'}(X)=\vert \lambda\vert _{F}^{-\delta(H'')/2}c_{\theta}(\lambda X),\,\,c_{f'}(X)=\vert \lambda\vert _{F}^{-\delta(G'')/2}c_{f}(\lambda X).$$
On calcule comme en 7.5
$$D^H(\lambda^{-1}X)=\vert \lambda\vert _{F}^{\delta(H'')-\delta(H)} D^H( X),\,\, \Delta(\lambda^{-1}X)=\vert \lambda\vert _{F}^{-dim(W')}\Delta(X).$$
Par changement de variable, on obtient
$$\int_{\mathfrak{t}(F)}c_{\theta'}(X)c_{f'}(X)D^H(X)\Delta(X)^rdX=\vert \lambda\vert _{F}^{b}\int_{\mathfrak{t}(F)}c_{\theta}(X)c_{f}(X)D^H(X)\Delta(X)^rdX,$$
o\`u
$$b=-\delta(G'')/2-dim(T)+\delta(H'')/2-\delta(H)-rdim(W').$$
En utilisant 7.7(2), on v\'erifie que $b=-\delta(G)/2-\delta(H)/2$. De l'\'egalit\'e pr\'ec\'edente r\'esulte alors l'\'egalit\'e
$$I(\theta',f')=\vert \lambda\vert ^{-\delta(G)/2-\delta(H)/2}I(\theta,f).$$
Jointe \`a (4), cette \'egalit\'e d\'emontre (3).

La relation (1) entra\^{\i}ne que la forme bilin\'eaire $E$ est combinaison lin\'eaire des applications $(\theta,f)\mapsto c_{\theta,{\cal O}^H}c_{\theta_{f},{\cal O}}$, o\`u ${\cal O}^H$, resp. ${\cal O}$, parcourt $Nil(\mathfrak{h}(F))$, resp. $Nil(\mathfrak{g}(F))$. Les relations (2) et (3) nous disent que $E$, ainsi que toutes ces applications, sont homog\`enes pour la transformation $(\theta,f)\mapsto (\theta^{\lambda},f^{\lambda})$. Il en r\'esulte que $E$ est combinaison lin\'eaire de celles des applications ci-dessus qui sont de m\^eme degr\'e que $E$. On a toujours $dim({\cal O}^H)\leq \delta(H)$, $dim({\cal O})\leq \delta(G)$. L'\'egalit\'e $dim({\cal O}^H)+dim({\cal O})=\delta(H)+\delta(G)$ est donc \'equivalente \`a la r\'eunion des deux \'egalit\'es $dim({\cal O}^H)= \delta(H)$ et $dim({\cal O})= \delta(G)$. Celles-ci sont v\'erifi\'ees si et seulement si ${\cal O}^H$ et ${\cal O}$ sont r\'eguli\`eres. $\square$

\bigskip

\subsection{Calcul de germes de Shalika}

Dans ce paragraphe, on suppose $G$ quasi-d\'eploy\'e. Soient $B$ un sous-groupe de Borel de $G$ et $T_{qd}$ un sous-tore maximal de $B$. Soit $X_ {qd}$ un \'el\'ement de $\mathfrak{t}_{qd}(F)\cap \mathfrak{g}_{reg}(F)$.

Supposons $d$ pair et $d\geq4$. Pour toute extension quadratique $E$ de $F$, notons $\tau_{E}$ l'\'el\'ement non trivial de $Gal(E/F)$ et $\chi_{E}$ le caract\`ere quadratique de $F^{\times}$ associ\'e \`a $E$. Si $d_{an}(V)=0$ (ou encore, si $G$ est d\'eploy\'e), on note $\chi_{V}$ le caract\`ere trivial de $F^{\times}$ et on pose $\eta=1$. Si $d_{an}(V)=2$, il y a une extension quadratique $E$ de $F$ et un \'el\'ement $\eta\in F^{\times}$ tel que le noyau anisotrope de $q$ soit la forme $\eta Norm_{E/F}$. L'\'el\'ement $\eta$ n'est pas unique, on le fixe. On pose $\chi_{V}=\chi_{E}$. Soient $F_{1}$ et $F_{2}$ deux extensions quadratiques de $F$ telles que $\chi_{F_{1}}\chi_{F_{2}}=\chi_{V}$. Pour $i=1,2$, soient $a_{i}\in F_{i}^{\times}$ tel que $\tau_{F_{i}}(a_{i})=-a_{i}$. On suppose $a_{1}\not=\pm a_{2}$. Soit $c\in F^{\times}$ tel que $\chi_{V}(\eta cNorm_{F_{1}/F}(a_{1}))=1$. On peut identifier $V$, comme espace quadratique, \`a la somme orthogonale $F_{1}\oplus F_{2}\oplus \tilde{Z}$, o\`u $F_{1}$ est muni de la forme $cNorm_{F_{1}/F}$, $F_{2}$ est muni de la forme $-cNorm_{F_{2}/F}$ et $\tilde{Z}$ est un espace hyperbolique de dimension $d-4$. Fixons une telle identification et un sous-tore d\'eploy\'e maximal $\tilde{T}$ du groupe sp\'ecial orthogonal de $\tilde{Z}$. Pour $\tilde{S}\in \tilde{\mathfrak{t}}(F)$, consid\'erons l'\'el\'ement $X_ {F_{1}}\in \mathfrak{g}(F)$ qui agit par multiplication par $a_{1}$, resp. $a_{2}$, sur $F_{1}$, resp. $F_{2}$, et par $\tilde{S}$ sur $\tilde{Z}$. Les \'el\'ements $a_{1}$ et $a_{2}$ \'etant fix\'es, on peut choisir $\tilde{S}$ tel que $X_ {F_{1}}$ soit r\'egulier. On fixe un tel $\tilde{S}$. Les \'el\'ements $a_{1}$, $a_{2}$ et $\tilde{S}$ \'etant fix\'es, on peut faire varier $c$. La classe de conjugaison de $X_ {F_{1}}$ ne d\'epend que de $\chi_{F_{1}}(c)$.  On note $X_ {F_{1}}^+$ l'\'el\'ement correspondant \`a un $c=c^+$ tel que $\chi_{F_{1}}(c^+)=\chi_{F_{1}}(\eta)\chi_{F_{1}}(Norm_{F_{1}/F}(a_{1})-Norm_{F_{2}/F}(a_{2}))$ et $X_ {F_{1}}^-$ celui qui correspond \`a un $c=c^-$ tel que $\chi_{F_{1}}(c^-)=-\chi_{F_{1}}(c^+ )$.

 On se rappelle que l'on a classifi\'e les orbites nilpotentes r\'eguli\`eres en 7.1.
 
 \ass{Lemme}{Soit ${\cal O}\in Nil(\mathfrak{g}(F))$.
 
 (i)   On a les \'egalit\'es
 $$\Gamma_{{\cal O}}(X_ {qd})=\left\lbrace\begin{array}{cc}0,&\,\,{\rm si}\,\,{\cal O}\,\,\text{ n'est pas\,\,r\'eguli\`ere;}\\ 1,&\,\,{\rm si}\,\,{\cal O}\,\,\text{ est\,\,r\'eguli\`ere.}\\ \end{array}\right.$$
 
 (ii) Supposons $d$ pair et $d\geq 4$. On a les \'egalit\'es
 $$\Gamma_{{\cal O}}(X_ {F_{1}}^+)-\Gamma_{{\cal O}}(X_ {F_{1}}^-)=\left\lbrace\begin{array}{cc}0,&\,\,{\rm si}\,\,{\cal O}\,\,{\rm n'est\,\,pas\,\,r\acute{e}guli\grave{e}re;}\\ \chi_{F_{1}}(\nu \eta),&\,\,{\rm si}\,\,{\cal O}={\cal O}_{\nu}\,\,{\rm avec}\,\,\nu\in {\cal N}^V.\\ \end{array}\right.$$}

 Preuve.  Le tore $T_{qd}$ est un L\'evi de $G$ et  la distribution $f\mapsto J_{G}(X_ {qd},f)$ est induite de la distribution $f\mapsto f(X_ {qd})$ sur $\mathfrak{t}_{qd}(F)$. On a \'evidemment $\Gamma_{\{0\}}^{T_{qd}}(X_ {qd})=1$. Alors $\Gamma_{{\cal O}} (X_ {qd})$ est non nul si et seulement si ${\cal O}$ intervient dans l'orbite induite de l'orbite $\{0\}$ de $\mathfrak{t}_{qd}(F)$. Cette condition \'equivaut \`a ce que ${\cal O}$ soit r\'eguli\`ere. On en d\'eduit la premi\`ere \'egalit\'e de (i).
 
 Supposons $d$ pair et $d\geq4$ et reprenons les constructions qui pr\'ec\`edent l'\'enonc\'e. Notons $G_{1}$ le groupe sp\'ecial orthogonal de $F_{1}\oplus \tilde{Z}$ et $G_{2}$ celui de $F_{2}$. Ils sont quasi-d\'eploy\'es. Pour $i=1,2$, on fixe un sous-tore  maximal $T_{i,qd}$ de $G_{i}$ inclus dans un sous-groupe de Borel (on a $T_{2,qd}=G_{2}$). Le groupe $G_{1}\times G_{2}$ est un groupe endoscopique de $G$. La distribution
 $$f\mapsto J_{G}(X^+_{F_{1}},f)-J_{G}(X^-_{F_{1}},f)$$
 est le transfert endoscopique d'une distribution
 $$(f_{1},f_{2})\mapsto J_{G_{1}}(X_ {1},f_{1})J_{G_{2}}(X_ {2},f_{2})$$
 sur $\mathfrak{g}_{1}(F)\times \mathfrak{g}_{2}(F)$, o\`u, pour $i=1,2$, $X_ {j}$ est un certain \'el\'ement de $\mathfrak{t}_{j,qd}(F)$. Il en r\'esulte  que le d\'eveloppement en germes de la premi\`ere distribution s'obtient en transf\'erant celui de la seconde distribution. Comme on vient de le voir, ce dernier ne contient que des orbites nilpotentes r\'eguli\`eres de $\mathfrak{g}_{1}(F)\times \mathfrak{g}_{2}(F)$. Le transfert endoscopique d'une int\'egrale nilpotente r\'eguli\`ere est combinaison lin\'eaire de telles int\'egrales. On en d\'eduit la premi\`ere \'egalit\'e de (ii).
 
 Il ne reste plus qu'\`a calculer des germes relatifs \`a des orbites nilpotentes r\'eguli\`eres. Ceux-ci ont \'et\'e calcul\'es par Shelstad ([S]). Il faut d'abord voir que les mesures utilis\'ees par Shelstad sont compatibles avec les n\^otres.  Shelstad suppose les mesures sur les tores maximaux "alg\'ebiques" au sens suivant. On  fixe une forme diff\'erentielle $\delta_{T_{qd}}$ de degr\'e maximal sur $T_{qd}$, invariante par translations, et un r\'eel $\lambda>0$. Pour sous-tore maximal $T$ de $G$, l'isomorphisme $T\simeq T_{qd}$ sur $\bar{F}$ permet de transf\'erer $\delta_{T_{qd}}$ en une forme diff\'erentielle $\delta_{T}$ sur $T$. On prend alors pour mesure sur $T(F)$ la mesure $\lambda \vert \delta_{T}\vert _{F}$, cf. 9.6 pour la notation. Mais on a vu dans la preuve du lemme 9.6 que nos mesures autoduales s'obtenaient par ce proc\'ed\'e, pour $\delta_{T_{qd}}$ et $\lambda$ convenables. Soit ${\cal O}$ une orbite nilpotente r\'eguli\`ere. La mesure sur ${\cal O}$ utilis\'ee par Shelstad est d\'efinie de la fa\c{c}on suivante. Soit $N\in {\cal O}$. Consid\'erons une suite $(Y_{j})_{j\in {\mathbb N}}$ d'\'el\'ements de $\mathfrak{g}_{reg}(F)$ telle que $lim_{j\to \infty}Y_{j}=N$. Alors, l'espace tangent $Tang_{Y_{j}}$ en $Y_{j}$ \`a la classe de conjugaison de $Y_{j}$ tend, en un sens que l'on va pr\'eciser, vers l'espace tangent $Tang_{N}$ en $N$ \`a ${\cal O}$. Notons $T_{j}=G_{Y_{j}}$. L'espace $Tang_{Y_{j}}$ est \'egal \`a $\mathfrak{g}(F)/\mathfrak{t}_{j}(F)$, sur lequel on a une mesure $m_{j}$. Alors les mesures $D^G (Y_{j})^{1/2}m_{j}$ tendent vers une mesure $m_{N}$ sur $Tang_{N}$. La mesure sur ${\cal O}$ est celle qui, en $N$, co\"{\i}ncide infinit\'esimalement avec $m_{N}$. Fixons un suppl\'ementaire $\mathfrak{r}$ du 
   noyau de $ad(N)$ dans $\mathfrak{g}(F)$. Pour $j$ assez grand, $\mathfrak{r}$ est encore un suppl\'ementaire de $\mathfrak{t}_{j}(F)$ dans $\mathfrak{g}(F)$ et on peut identifier $Tang_{Y_{j}}=Tang_{N}=\mathfrak{r}$. Soit $\mathfrak{l}_{j}$ l'orthogonal de $\mathfrak{t}_{j}$ dans $\mathfrak{g}$. On a aussi $Tang_{Y_{j}}\simeq \mathfrak{l}_{j}(F)$. Modulo cette identification, la mesure $m_{j}$ est la mesure autoduale associ\'ee \`a la forme quadratique $(Y,Z)\mapsto \frac{1}{2}trace(YZ)$ sur $\mathfrak{l}_{j}(F)$. Mais le jacobien de $ad(Y_{j})$ agissant dans $\mathfrak{l}(F)$ est $D^G(Y_{j})$. Donc $D^G(Y_{j})^{1/2}m_{j}$ est aussi la mesure associ\'ee \`a la forme symplectique $(Y,Z)\mapsto \frac{1}{2}trace([Y_{j},Y]Z)$ sur $\mathfrak{l}_{j}(F)$. La m\^eme formule d\'efinit une forme antisym\'etrique sur tout $\mathfrak{g}(F)$, de noyau $\mathfrak{t}_{j}(F)$. On peut donc remplacer $\mathfrak{l}_{j}(F)$ par $\mathfrak{r}$ et la mesure $D^G(Y_{j})^{1/2}m_{j}$ est la mesure associ\'ee \`a la forme symplectique $(Y,Z)\mapsto \frac{1}{2}trace([Y_{j},Y]Z)$ sur $\mathfrak{r}$. Quand $Y_{j}$ tend vers $N$, cette forme tend vers $(Y,Z)\mapsto \frac{1}{2}trace([N,Y]Z)$ et  $D^G(Y_{j})^{1/2}m_{j}$  tend vers la mesure associ\'ee \`a cette forme. Mais c'est pr\'ecis\'ement la fa\c{c}on dont nous avons d\'efini notre mesure sur ${\cal O}$ en 1.2.
 
 Cela \'etant, Shelstad montre qu'un germe $\Gamma_{{\cal O}}(S)$ associ\'e \`a une orbite nilpotente r\'eguli\`ere ${\cal O}$ vaut $1$ ou $0$, selon qu'un certain invariant  est \'egal ou non \`a $1$. Pour l'\'el\'ement $X_ {qd}$, il est facile de voir que l'invariant est $1$ et on en d\'eduit la seconde \'egalit\'e de (i).
 Consid\'erons la situation de (ii). Pour un signe $\zeta=\pm $, notons $T^{\zeta}$ le sous-tore maximal de $G$ tel que $X_ {F_{1}}^{\zeta}\in \mathfrak{t}^{\zeta}(F)$. Soient $\nu\in {\cal N}^V$ et $N\in {\cal O}_{\nu}$. Shelstad note l'invariant $inv(X_ {F_{1}}^{\zeta})inv(T^{\zeta})/inv_{T^{\zeta}}(N)$. Tous ces \'el\'ements appartiennent au groupe de cohomologie $H^1(T^{\zeta})=H^1(Gal(\bar{F}/F),T^{\zeta})$. On a ici $H^1(T^{\zeta})=\{\pm 1\}\times \{\pm 1\}$. Les invariants d\'ependent du choix d'un \'epinglage. On effectue ce choix comme en [W3] X.3 en prenant pour \'el\'ement $\eta$ de cette r\'ef\'erence notre \'el\'ement $\eta$ multipli\'e par $2(-1)^{d/2-1}$. Dans le lemme X.7 de [W3], nous avons calcul\'e le produit $inv(X_ {F_{1}}^{\zeta})inv(T^{\zeta})$.
 On a
 $$(1) \qquad inv(X_ {F_{1}}^{\zeta})inv(T^{\zeta})=$$
 $$\qquad (\chi_{F_{1}}(2(-1)^{d/2-1}\eta(c^{\zeta})^{-1}a_{1}^{-1}P'(a_{1})),\chi_{F_{2}}(2(-1)^{d/2}\eta(c^{\zeta})^{-1}a_{2}^{-1}P'(a_{2}))),$$
 o\`u $P$ est le polyn\^ome caract\'eristique de $X_ {F_{1}}$ agissant dans $V$ et $P'$ est le polyn\^ome d\'eriv\'e. Notons $(\pm \tilde{s}_{j})_{j=3,...,d/2}$ les valeurs propres de l'action de $\tilde{S}$ dans $\tilde{Z}$. Elles appartiennent \`a $F^{\times}$ puisque $\tilde{T}$ est d\'eploy\'e. On a
 $$P(T)=(T^2+Norm_{F_{1}/F}(a_{1}))(T^2+Norm_{F_{2}/F}(a_{2}))\prod_{j=3,...,d/2}(T^2-\tilde{s}_{j}^2).$$
 Donc
 $$a_{1}^{-1}P'(a_{1})=2(-Norm_{F_{1}/F}(a_{1})+Norm_{F_{2}/F}(a_{2}))\prod_{j=3,...,d/2}(-Norm_{F_{1}/F}(a_{1})-\tilde{s}^2_{j})$$
 $$\qquad =2(-1)^{d/2-1}(Norm_{F_{1}/F}(a_{1})-Norm_{F_{2}/F}(a_{2}))\prod_{j=3,...,d/2}(Norm_{F_{1}/F}(a_{1})+\tilde{s}^2_{j}).$$
  On a $\tilde{s}_{j}^2+Norm_{F_{1}/F}(a_{1})=Norm_{F_{1}/F}(\tilde{s}_{j}+a_{1})$, donc $\chi_{F_{1}}(\tilde{s}_{j}^2+Norm_{F_{1}/F}(a_{1}))=1$. Un calcul similaire vaut en \'echangeant les r\^oles de $F_{1}$ et $F_{2}$. La formule (1) se simplifie en
 $$inv(X_ {F_{1}}^{\zeta})inv(T^{\zeta})=$$
 $$\qquad
 (\chi_{F_{1}}(\eta(c^{\zeta})^{-1}(Norm_{F_{1}/F}(a_{1})-Norm_{F_{2}/F}(a_{2}))),\chi_{F_{2}}(\eta(c^{\zeta})^{-1}(Norm_{F_{1}/F}(a_{1})-Norm_{F_{2}/F}(a_{2}))).$$
 D'apr\`es la d\'efinition de $c^{\zeta}$, on obtient
 $$inv(X_ {F_{1}}^{\zeta})inv(T^{\zeta})=(\zeta,\zeta),$$
 o\`u on identifie $\zeta$ \`a un \'el\'ement de $\{\pm 1\}$.
 L'\'epinglage d\'etermine un \'el\'ement nilpotent r\'egulier $N^*$: avec les notations de [W3] page 313, $N^*=\sum_{j=1,...,d/2}X_{\alpha_{j}}$. Notons  $\nu^*$ l'\'el\'ement de ${\cal N}^V$ tel que $N^*\in {\cal O}_{\nu^*}$. On d\'efinit un cocycle $d_{N}$ de $Gal(\bar{F}/F)$ dans $T^{\zeta}$ de la fa\c{c}on suivante. Si $\nu=\nu^*$, $d_{N}=1$. Si $\nu\not=\nu^*$, on pose $E_{N}=F(\sqrt{\nu/\nu^*})$. Alors, pour $\sigma\in Gal(\bar{F}/F)$, on a $d_{N}(\sigma)=1$ si $\sigma\in Gal(\bar{F}/E_{N})$ et $d_{N}(\sigma)=-1$ sinon. L'invariant $inv_{T^{\zeta}}(N)$ est l'image dans $H^1(T)$ du cocycle  $d_{N}$. On  calcule cette image
 $$inv_{T^{\zeta}}(N)=(\chi_{F_{1}}(\nu/\nu^*),\chi_{F_{2}}(\nu/\nu^*)).$$
 En fait, on a forc\'ement $\chi_{V}(\nu/\nu^*)=1$, donc
 $$inv_{T^{\zeta}}(N)=\chi_{F_{1}}(\nu/\nu^*),\chi_{F_{1}}(\nu/\nu^*)).$$
 L'\'el\'ement $N^*$ laisse stable l'hyperplan engendr\'e par $e_{1},...,e_{d/2-1},e_{d/2}+e_{d/2+1},e_{d/2+2},...,e_{d}$, avec les notations de [W]. Le noyau anisotrope de la restriction de la forme $q$ \`a cet hyperplan est la restriction de $q$ \`a la droite port\'ee par $e_{d/2}+e_{d/2+1}$. Or $q(e_{d/2}+e_{d/2+1})=\eta$, donc $\nu^*=\eta$. Finalement
 $$inv(X_ {F_{1}}^{\zeta})inv(T^{\zeta})/inv_{T^{\zeta}}(N)=(\zeta\chi_{F_{1}}(\nu\eta),\zeta\chi_{F_{1}}(\nu\eta)).$$
 D'apr\`es Shelstad, on a donc
 $$\Gamma_{{\cal O}_{\nu}}(X_ {F_{1}}^{\zeta})=\left\lbrace\begin{array}{cc}1,&\,\,\text{si}\,\,\chi_{F_{1}}(\nu \eta)=\zeta,\\ 0,&\,\,\text{si}\,\,\chi_{F_{1}}(\nu\eta)=-\zeta.\\ \end{array}\right.$$
 Cela entra\^{\i}ne la seconde \'egalit\'e du (ii) de l'\'enonc\'e. $\square$

 \bigskip

 \subsection{Preuve de la proposition 11.1 dans le cas $d$ impair}
 
 On suppose v\'erifi\'ee l'hypoth\`ese du paragraphe 11.1 et on suppose $d$ impair. On veut prouver $E(\theta,f)=0$ pour tout quasi-caract\`ere $\theta$ sur $\mathfrak{h}(F)$ et toute fonction tr\`es cuspidale $f\in C_{c}^{\infty}(\mathfrak{g}(F))$. Si $G$ n'est pas d\'eploy\'e ou si $H$ n'est pas quasi-d\'eploy\'e, l'une des alg\`ebres de Lie de ces groupes n'a pas d'\'el\'ement nilpotent r\'egulier et la conclusion r\'esulte du lemme 11.3. Supposons $G$ d\'eploy\'e et $H$ quasi-d\'eploy\'e. Le m\^eme lemme nous dit que, si $d_{W}\leq2$, resp. $d_{W}\geq4$, il existe un nombre complexe $c_{reg}$, resp. une famille de nombres complexes $(c_{\nu})_{\nu\in {\cal N}^W}$, de sorte que
 $$E(\theta,f)=\left\lbrace\begin{array}{cc}c_{reg}c_{\theta,{\cal O}^H_{reg}}c_{\theta_{f},{\cal O}_{reg}},&\,\,\text{si}\,\,d_{W}\leq2,\\ \sum_{\nu\in {\cal N}^W}c_{\nu}c_{\theta,{\cal O}_{\nu}^H}c_{\theta_{f},{\cal O}_{reg}},&\,\,\text{si}\,\,d_{W}\geq4,\\ \end{array}\right.$$
 pour tous $\theta$, $f$ (on a introduit des exposants $H$ pour pr\'eciser la notation). Les constantes sont uniquement d\'etermin\'ees d'apr\`es le lemme 6.3(iii). 
 
  Soient $T\in {\cal T}(G)$ et $X\in \mathfrak{t}(F)\cap \mathfrak{g}_{reg}(F)$. D'apr\`es 6.3(3), on peut construire  un voisinage $\omega_{X}$ de $X$ dans $\mathfrak{t}(F)$ et une fonction tr\`es cuspidale $f[X]\in C_{c}^{\infty}(\mathfrak{g}(F))$ v\'erifiant les propri\'et\'es suivantes:
 
 (1) pour $T'\in {\cal T}(G)$ et $T'\not=T$, la restriction de $\hat{\theta}_{f[X]}$ \`a $\mathfrak{t}'(F)$ est nulle;
 
 (2) pour toute fonction localement int\'egrable $\varphi$ sur $\mathfrak{t}(F)$, invariante par $W(G,T_{d})$, 
 $$\int_{\mathfrak{t}(F)}\varphi(X')D^G(X')^{1/2}\hat{\theta}_{f[X]}(X')dX'=mes(\omega_{X} )^{-1}\int_{\omega_{X}}\varphi(X')dX';$$
 
 (3) pour tout $Y\in \mathfrak{g}_{reg}(F)$, 
 $$\theta_{f[X]}(Y)=mes(\omega_{X})^{-1}\int_{\omega_{X}}\hat{j}^G(X',Y)dX'$$
 \noindent
  (avec les notations de 6.3(3), $f[X]=\theta_{f'}(X)^{-1}D^G(X)^{-1/2}mes(\omega_{X})^{-1}\hat{f}'$). Au voisinage de $0$, l'\'egalit\'e   (3) se simplifie en
 $$\theta_{f[X]}(Y)=\hat{j}^G(X,Y)=\sum_{{\cal O}\in Nil(\mathfrak{g})}\Gamma_{{\cal O}}(X)\hat{j}^G({\cal O},Y).$$
 Donc 
 $$(4) \qquad c_{\theta_{f[X]},{\cal O}}=\Gamma_{{\cal O}}(X)$$
 pour tout ${\cal O}\in Nil(\mathfrak{g})$. Notons $T_{d}$ l'unique tore d\'eploy\'e dans ${\cal T}(G)$ et fixons $X_{d}\in \mathfrak{t}_{d}(F)\cap \mathfrak{g}_{reg}(F)$. On fixe $\omega_{X_{d}}$ et $f[X_{d}]$ comme ci-dessus et on pose $f=f[X_{d}]$. Gr\^ace \`a (4) et au lemme 11.4(i), on a $c_{\theta_{f},{\cal O}_{reg}}=1$.  Fixons un $G$-domaine $\omega$ dans $\mathfrak{g}(F)$, compact modulo conjugaison, contenant $0$ et $Supp(f)$. Soit $S\in \mathfrak{h}_{reg}(F)$. On suppose que l'action de $S$ dans $W$ est de noyau nul.  D'apr\`es le lemme 6.3(ii), on peut choisir un quasi-caract\`ere $\theta[S]$ sur $\mathfrak{h}(F)$ tel que $\theta[S](Y) =\hat{j}^H(S,Y)$ pour tout $  Y\in \omega$. Comme ci-dessus, on a $c_{\theta[S],{\cal O}^H}=\Gamma_{{\cal O}^H}(S)$ pour tout ${\cal O}^H\in Nil(\mathfrak{h})$. Rempla\c{c}ons $G$ et $V$ par $H$ et $W$ dans les d\'efinitions de 11.4. On d\'efinit dans $\mathfrak{h}(F)$ un \'el\'ement $X_{qd}$ et, si $d_{W}\geq4$, des \'el\'ements $X_{F_{1}}^{\pm}$. 
 
 Si $d_{W}\leq2$, posons $\theta=\theta[X_{qd}]$. On pose aussi ${\cal S}=\{X_{qd}\}$ et $c_{X_{qd}}=1$. On a $c_{\theta,{\cal O}_{reg}^H}=1$, donc $c_{reg}=E(\theta,f)$.

 Supposons $d_{W}\geq4$ et fixons $\nu\in {\cal N}^W$. Si $H$ est d\'eploy\'e, notons ${\cal F}^V$ l'ensemble des extensions quadratiques de $F$.  Si $H$ n'est pas d\'eploy\'e, soient $E$ et $\eta$ comme en 11.4. Les extensions quadratiques de $F$ distinctes de $E$ vont par paire: \`a $F_{1}$ est associ\'e $F_{2}$ tel que $\chi_{F_{1}}\chi_{F_{2}}=\chi_{E}$. On fixe un sous-ensemble ${\cal F}_{V}$ qui contient un \'el\'ement de chaque paire. Remarquons que, dans les deux cas, on a l'\'egalit\'e
 $$\vert {\cal N}^V\vert =1+\vert {\cal F}^V\vert .$$
 Posons
 $$\theta= \vert {\cal N}^V\vert ^{-1}(\theta[X_{qd}]+\sum_{F_{1}\in {\cal F}^V}\chi_{F_{1}}(\nu\eta)(\theta[X_{F_{1}}^+]-\theta[X_{F_{1}}^-])).$$
  On note ${\cal S}$ l'ensemble des \'el\'ements $S$ tels que $\theta[S]$ intervienne dans ces formules et, pour $S\in {\cal S}$, $c_{S}$ le coefficient dont $\theta[S]$ y est affect\'e. Le lemme 11.4 entra\^{\i}ne que, pour $\nu'\in {\cal N}^W$, on a l'\'egalit\'e
 $$c_{\theta,{\cal O}_{\nu'}^H}= \delta_{\nu,\nu'},$$
 o\`u ce dernier terme est le symbole de Kronecker. Donc $c_{\nu}=E(\theta,f)$.
 
 Pour $X\in \omega_{X_{d}}$, la distribution $\varphi\mapsto J_{G}(X,\hat{\varphi})$ est induite d'une distribution sur $\mathfrak{t}_{d}(F)$. Cela entra\^{\i}ne que la fonction $Y\mapsto \hat{j}^G(X,Y)$ est \`a support dans l'ensemble des \'el\'ements appartenant \`a une sous-alg\`ebre de Borel de $G$. D'apr\`es (3), c'est aussi le cas de la fonction $\theta_{f}$. Comme dans la preuve du lemme 7.6, cela entra\^{\i}ne que, si $T$ est un \'el\'ement de ${\cal T}$ diff\'erent du tore $\{1\}$, la fonction $c_{f}$ est nulle sur $\mathfrak{t}(F)$. Donc $I(\theta,f)$ se r\'eduit \`a la contribution de l'unique tore $T=\{1\}\in {\cal T}$. Par d\'efinition, celle-ci est $c_{\theta,{\cal O}_{reg}^H}c_{\theta_{f},{\cal O}_{reg}}$ si $d_{W}\leq2$, $c_{\theta,{\cal O}_{-\nu_{0}}^H}c_{\theta_{f},{\cal O}_{reg}}$ si $d_{W}\geq4$. Avec les calculs ci-dessus, on obtient
 $$I(\theta,f)=\left\lbrace\begin{array}{cc}1,&\,\,\text{si}\,\,d_{W}\leq2,\\ \delta_{\nu,-\nu_{0}},&\,\,\text{si}\,\,d_{W}\geq4.\\ \end{array}\right.$$
 
 L'ensemble ${\cal S}$ et la famille $(c_{S})_{S\in {\cal S}}$ permettent de calculer $J(\theta,f)$, cf. 5.2. Pour tout $S\in {\cal S}$, posons 
 $$m(S)=mes(\omega_{X})^{-1}mes(\omega_{X}\cap \mathfrak{t}_{d}(F)^S).$$
 En utilisant les propri\'et\'es (1) et (2), on obtient
 $$(5) \qquad J(\theta,f)=\sum_{S\in {\cal S}}c_{S}m(S).$$
 Soit $S\in {\cal S}$. On utilise les d\'efinitions et notations de 9.4, appliqu\'ees au cas $V''=V$. Soit $X\in \mathfrak{t}_{d}(F)\cap \mathfrak{g}_{reg}(F)$. On a
 
 (6) $X\in \mathfrak{t}_{d}(F)^S$ si et seulement s'il existe une famille $(z_{\pm j})_{j=1,...,d_{W}/2}$ d'\'el\'ements de $\bar{F}^{\times}$ telle que
 $$\left\lbrace\begin{array}{c}
   z_{j}z_{-j}=\frac{P_{X}(s_{j})}{4\nu_{0}s_{j}^{1+2r}R_{S,j}(s_{j})}\text{ pour tout }j=1,...,d_{W}/2;\\ \sum_{j=\pm 1,...\pm d_{W}/2}z_{j}w_{j}\in W.\\ \end{array}\right.$$

  En effet, $X$ appartient \`a $\mathfrak{t}_{d}(F)^S$ si et seulement s'il est conjugu\'e par un \'el\'ement de $G(F)$ \`a un \'el\'ement de $\Xi+S+\Lambda^S$.  Il n'y a pas d'indiscernabilit\'e pour le tore d\'eploy\'e $T_{d}$.  La condition \'equivaut donc \`a ce qu'il existe $Y\in \Xi+S+\Lambda^S$ tel que l'on ait l'\'egalit\'e des polyn\^omes caract\'eristiques $P_{X}=P_{Y}$. Cela \'equivaut \`a ce qu'il existe  une famille $(\lambda_{i})_{i=0,...,r-1}$ d'\'el\'ements de $F$ et une famille $(z_{\pm j})_{j=1,...,d_{W}/2}$ d'\'el\'ements de $\bar{F}^{\times}$ telles que, d'une part, le polyn\^ome $P_{X}$ soit \'egal \`a celui figurant dans l'\'enonc\'e du lemme 9.4, d'autre part, l'\'el\'ement $\sum_{j=\pm 1,...\pm d_{W}/2}z_{j}w_{j}$ appartienne \`a $W$. D'apr\`es   9.4(1), les conditions sur les $z_{j}$ sont celles de (6). Les $\lambda_{i}$ sont ensuite d\'etermin\'es par un syst\`eme  inversible d'\'equations lin\'eaires \`a coefficients dans $F$.  Cela prouve (6). 
 
La condition (6) impose que $P_{X}(s_{j})\not=0$ pour tout $j$. On suppose cette condition v\'erifi\'ee. Il existe
 
- une d\'ecomposition de $W$ en somme directe 
 $$\oplus_{j=1,...,h}F_{j}\oplus \tilde{Z},$$
 o\`u $h\leq2$ et les $F_{j}$ sont des extensions quadratiques de $F$;

- des \'el\'ements $c_{j}\in F^{\times}$, pour $j=1,...,h$, de sorte que $q_{W}$ soit la somme directe orthogonale des formes $c_{j}Norm_{F_{j}/F}$ sur $F_{j}$ et d'une forme hyperbolique sur $\tilde{Z}$;

- des \'el\'ements $a_{j }\in F_{j}^{\times}$, pour $j=1,...,h$ tels que $\tau_{F_{j}}(a_{j})=-a_{j}$ et un \'el\'ement $\tilde{S}$ appartenant \`a l'alg\`ebre de Lie d'un sous-tore d\'eploy\'e maximal du groupe sp\'ecial orthogonal de $\tilde{Z}$, de sorte que $S$ agisse par multiplication par $a_{j}$ dans $F_{j}$ et par $\tilde{S}$ dans $\tilde{Z}$.

Pour $j=1,...,h$, soit $(e_{j},e_{-j})$ la base de $F_{j}\otimes_{F}\bar{F}$ telle que tout \'el\'ement $x\in F_{j}$ soit \'egal \`a $xe_{j}+\tau_{F_{j}}(x)e_{-j}$. Elle est d\'efinie sur $F_{j}$ et on a l'\'egalit\'e $\tau_{F_{j}}(e_{j})=e_{-j}$. On a l'\'egalit\'e $q(e_{j},e_{-j})=c_{j}$. On peut donc  supposer que $w_{j}=e_{j}$, $w_{-j}=c_{j}^{-1}e_{-j}$ pour $j=1,...,h$ et $(w_{j})_{j=\pm(h+1),...,\pm d_{W}/2}$ est une base hyperbolique de $\tilde{Z}$. On a $s_{j}=a_{j}$ pour $j\leq h$. La condition (6) se d\'ecompose en $d/2$ conditions $(6)_{j}$ portant sur les couples $(z_{j},z_{-j})$. Pour $j>h$, on satisfait $(6)_{j}$ en prenant $z_{-j}=1$ et $z_{j}$ \'egal au membre de droite de la premi\`ere relation de (6). Pour $j\leq h$, la seconde relation de (6) \'equivaut \`a  $z_{j}\in F_{j}^{\times}$ et $\tau_{F_{j}}(z_{j})=c_{j}^{-1}z_{-j}$. La condition $(6)_{j}$ \'equivaut donc \`a
$$\chi_{F_{j}}(\frac{P_{X}(a_{j})}{c_{j}\nu_{0}a_{j}^{1+2r}R_{S,j}(a_{j})})=1.$$
On a 
$$R_{S,j}(a_{j})=\prod_{j'=1,...,d_{W}/2; j'\not=j}(a_{j}^2-s_{j'}^2)=(-1)^{d_{W}/2-1}(\prod_{j'=1,...,h; j'\not=j}(Norm_{F_{j}/F}(a_{j})-Norm_{F_{j'}/F}(a_{j'})))$$
$$\qquad \prod_{j'=h+1,...,d_{W}/2}(s_{j'}^2-a_{j}^2).$$
Notons $(\pm x_{k})_{k=1,...,(d-1)/2}$ les valeurs propres non nulles de $X$. On a
$$P_{X}(a_{j})=a_{j}\prod_{k=1,...,(d-1)/2}(a_{j}^2-x_{k}^2)=(-1)^{(d-1)/2}a_{j}\prod_{k=1,...,(d-1)/2}(x_{k}^2-a_{j}^2).$$
Pour $j'=h+1,...,d_{W}/2$, $s_{j'}$ appartient \`a $F^{\times}$, donc $s_{j'}^2-a_{j}^2$ est la norme d'un \'el\'ement de $F_{j}^{\times}$. De m\^eme, pour tout $k=1,...,(d-1)/2$, $x_{k}^2-a_{j}^2$ est une norme. Enfin $(-1)^ra_{j}^{2r}$ est aussi une norme. Ces termes disparaissent de notre calcul et la condition $(6)_{j}$ \'equivaut \`a
$$(7)\qquad \chi_{F_{j}}(-c_{j}\nu_{0}\prod_{j'=1,...,h; j'\not=j}(Norm_{F_{j}/F}(a_{j})-Norm_{F_{j'}/F}(a_{j'})))=1.$$
Cette relation est ind\'ependante de $X$. On a impos\'e \`a $X$ des conditions de non nullit\'e qui sont v\'erifi\'ees sur un ouvert dont le compl\'ementaire est de mesure nulle. Cela d\'emontre que $m(S)=1$ si la condition (7) est v\'erifi\'ee pour tout $j=1,...,h$, $m(S)=0$ sinon.

Supposons $H$ d\'eploy\'e et $S=X_{qd}$. Alors $h=0$, donc $m(S)=1$. Supposons $H$ non d\'eploy\'e et $S=X_{qd}$. Alors $h=1$, $F_{1}=E$ et $c_{1}=\eta$ avec les notations de 11.4 appliqu\'ees \`a $W$. Le noyau anisotrope de $q$ est le m\^eme que celui de la forme quadratique sur $E\oplus D$  
 $$e\oplus xv_{0}\mapsto \eta Norm_{E/F}(e)+\nu_{0}x^2.$$
 Puisqu'on a suppos\'e $G$ d\'eploy\'e, cette forme n'est pas anisotrope donc $-\eta\nu_{0}$ est la norme d'un \'el\'ement de $E^{\times}$. La condition (7) est v\'erifi\'ee et $m(S)=1$. Supposons maintenant $d_{W}\geq4$ et $S=X_{F_{1}}^{\zeta}$, o\`u $\zeta=\pm$, $F_{1}\in {\cal F}^V$.  Alors $h=2$,  et les termes $F_{1}$, $F_{2}$, $a_{1}$ et $a_{2}$ co\"{\i}ncident avec ceux de 11.4. On a $c_{1}=c^{\zeta}$ et $c_{2}=-c^{\zeta}$.  D'apr\`es la d\'efinition de ces termes, la condition (7) \'equivaut \`a $\chi_{F_{1}}(-\eta \nu_{0})=\zeta$ pour $j=1$, resp. $\chi_{F_{2}}(-\eta \nu_{0})=\zeta$ pour $j=2$. Mais le calcul ci-dessus montre que $\chi_{V}(-\eta \nu_{0})=1$. Les conditions pour $j=1$ et $j=2$ sont donc \'equivalentes. On obtient que $m(S)=1$ si $\chi_{F_{1}}(-\eta \nu_{0})=\zeta$, $m(S)=0$ sinon. Reportons ces valeurs de $m(S)$ dans l'\'egalit\'e (5). Dans le cas $d_{W}\leq2$, on obtient imm\'ediatement la formule ci-dessous. Dans le cas $d_{W}\geq4$, celle-ci r\'esulte d'une inversion de Fourier sur le groupe $F^{\times}/F^{\times 2}$ si $H$ est d\'eploy\'e, sur le groupe $Norm_{E/F}(E^{\times})/F^{\times 2}$ sinon. La formule est
$$J(\theta,f)=\left\lbrace\begin{array}{cc}1,&\,\,\text{si}\,\,d_{W}\leq2,\\ \delta_{\nu,-\nu_{0}},&\,\,\text{si}\,\,d_{W}\geq4.\\ \end{array}\right.$$
Alors $J(\theta,f)=I(\theta,f)$ et $E(\theta,f)=0$. Donc $c_{reg}=0$ dans le cas $d_{W}\leq2$ et $c_{\nu}=0$ dans le cas $d_{W}\geq4$. Mais alors l'application bilin\'eaire $E$ est identiquement nulle, ce que l'on voulait d\'emontrer.

\bigskip

 \subsection{Preuve de la proposition 11.1 dans le cas $d$ pair}
 
 On suppose v\'erifi\'ee l'hypoth\`ese du paragraphe 11.1 et on suppose $d$ pair. Eliminons le cas $d=2$. Dans ce cas, $G$ est un tore de dimension $1$. Il r\'esulte des d\'efinitions que
 $$J(\theta,f)=\theta(0)\int_{\mathfrak{g}(F)}\hat{f}(X)dX$$
 et 
 $$I(\theta,f)=\theta(0)f(0).$$
 Ces deux expressions sont \'egales par inversion de Fourier. On suppose maintenant $d\geq4$.  En imitant ce que l'on a fait au paragraphe pr\'ec\'edent, on peut supposer $G$ quasi-d\'eploy\'e et $H$ d\'eploy\'e. Il existe une unique famille de nombres complexes $(c_{\nu})_{\nu\in {\cal N}^V}$ de sorte que
 $$E(\theta,f)=\sum_{\nu\in {\cal N}^V}c_{\nu}c_{\theta,{\cal O}_{reg}^H}c_{\theta_{f},{\cal O}_{\nu}}$$
 pour tous $\theta$, $f$. Fixons $\nu\in {\cal N}^V$. On introduit des \'el\'ements $X_{qd}$ et $X_{F_{1}}^{\pm}$ de $\mathfrak{g}(F)$ comme en 11.4. Soit $X$ un de ces \'el\'ements. On peut supposer qu'il appartient \`a un tore appartenant \`a ${\cal T}(G)$, que l'on note $T_{X}$. On introduit un voisinage $\omega_{X}$ de $X$ dans $\mathfrak{t}_{X}(F)$ et une fonction $f[X]$ v\'erifiant les conditions (1), (2) et (3) de 11.5. On pose
 $$f= \vert {\cal N}^V\vert ^{-1}(f[X_{qd}]+\sum_{F_{1}\in {\cal F}^V}\chi_{F_{1}}(\nu\eta)(f[X_{F_{1}}^+]-f[X_{F_{1}}^-])).$$
  Les notations sont celles introduites dans le paragraphe 11.5, l'extension $E$ \'etant maintenant associ\'ee \`a $V$ et non plus \`a $W$. On note ${\cal X}$ l'ensemble des \'el\'ements $X$ tels que $f[X]$ apparaisse dans ces formules et, pour $X\in {\cal X}$, $c_{X}$ le coefficient dont il est affect\'e. On fixe un $G$-domaine $\omega$ dans $\mathfrak{g}(F)$, compact modulo conjugaison, contenant $0$ et $Supp(f)$. Notons $T_{d}$ l'unique tore d\'eploy\'e dans ${\cal T}(H)$et fixons un \'el\'ement $S\in \mathfrak{t}_{d}(F)\cap \mathfrak{h}_{reg}(F)$. On choisit un quasi-caract\`ere $\theta$ sur $\mathfrak{h}(F)$ tel que $\theta(Y)=\hat{j}^H(S,Y)$ pour tout $Y\in \omega\cap \mathfrak{h}_{reg}(F)$. Pour $T\in {\cal T}$, $T\not=\{1\}$, la fonction $c_{f}$ est nulle hors de $\omega$ tandis que $c_{\theta}$ est nulle sur $\mathfrak{t}(F)\cap \omega$ pour la m\^eme raison que $c_{f}$ l'\'etait en 11.5. Donc $c_{\theta}c_{f}$ est nulle sur $\mathfrak{t}(F)$. Comme en 11.5, on calcule alors
 $$I(\theta,f)=\delta_{\nu,\nu_{0}}.$$
 
 Pour $X\in {\cal X}$, on pose
 $$m(X)=mes(\omega_{X})^{-1}mes(\omega_{X}\cap \mathfrak{t}_{X}(F)^S).$$
 On a
 $$(1) \qquad J(\theta,f)=\sum_{X\in {\cal X}}c_{X}m(X),$$
 et on est ramen\'e \`a calculer ces termes $m(X)$.
 
 On utilise les d\'efinitions et notations de 9.4 appliq\'ees au cas $W''=W$. Puisque $T_{d}$ est d\'eploy\'e, on peut supposer que les vecteurs $w_{\pm j}$ appartiennent \`a $W$. Soit $X\in \mathfrak{g}_{reg}(F)$. Supposons $P_{X}(s_{j})\not=0$ pour tout $j=1,...,m$ et $P_{X}(0)\not=0$. Soit $X_{1},...,X_{l}$ un ensemble de repr\'esentants des classes de conjugaison par $G(F)$ dans la classe de conjugaison stable de $X$. Montrons que
 
 (2) il existe un unique $i\in \{1,...,l\}$ tel que $X_{i}$ soit conjugu\'e \`a un \'el\'ement de $\Xi+S+\Sigma^S$ par un \'el\'ement de $G(F)$.
 
 La forme bilin\'eaire $(v,v')\mapsto q(v,Xv')$ est symplectique. Son d\'eterminant est donc un carr\'e dans $F^{\times}$. Ce d\'eterminant est $det(q)det(X)$. On a $det(X)=P_{X}(0)$ et $det(q)=(-1)^{d/2-1}4\nu_{0}\nu_{S}$. On a $R_{S,0}(0)=(-1)^{m}\prod_{j=1,...,m}s_{j}^2$ et $m=(d_{W}-1)/2$. On en d\'eduit que
 $$(3) \qquad \frac{P_{X}(0)}{(-1)^r\nu_{S}\nu_{0}R_{S,0}(0)}\in F^{\times 2}.$$
  On peut choisir des coordonn\'ees $(\lambda_{i})_{i=1,...,r}$, $(z_{\pm j})_{j=1,...,m}$ et $z_{0}$, avec $\lambda_{i}\in F$ et $z_{j}\in F^{\times}$ pour $j=0,\pm 1,...,\pm m$, de sorte que $P_{X}$ soit le polyn\^ome du lemme 3.4. En effet, on pose $z_{-j}=1$ pour $j=1,...,m$ et, gr\^ace \`a (3), on peut choisir les $z_{j}$, pour $j=0,...,m$, de sorte que les \'egalit\'es   9.4(1) et 9.4(2) soient v\'erifi\'ees. Ensuite, les $\lambda_{i}$ sont d\'etermin\'es par un syst\`eme inversible d'\'equations lin\'eaires \`a coefficients dans $F$. Cela montre qu'il existe $Y\in \Xi+S+\Lambda^S$ tel que $P_{Y}=P_{X}$. Un tel $Y$ est conjugu\'e \`a $X$ par un \'el\'ement de $G^+$. Comme on l'a dit dans la preuve du lemme 9.5, quitte \`a changer $z_{0}$ en $-z_{0}$, on peut assurer que $Y$ est conjugu\'e \`a $X$ par un \'el\'ement de $G$. Donc $Y$ appartient \`a la classe de conjugaison stable de $X$ et est conjugu\'e \`a l'un des $X_{i}$ par un \'el\'ement de $G(F)$. L'unicit\'e de cet \'el\'ement $X_{i}$ est assur\'ee par le lemme 9.5: deux \'el\'ements de $\Xi+S+\Sigma^S$ ne peuvent \^etre conjugu\'es par un \'el\'ement de $G$ que s'ils le sont par un \'el\'ement de $G(F)$. D'o\`u (2).
  
 On a choisi un \'el\'ement $S$. On peut supposer que les hypoth\`eses de non-nullit\'e impos\'ees ci-dessus \`a $X$ sont v\'erifi\'ees pour tout  \'el\'ement de $ \bigcup_{X\in {\cal X}}\omega_{X}$. Pour $X\in \omega_{X_{qd}}$, la classe de conjugaison stable de $X$ se r\'eduit \`a sa classe de conjugaison par $G(F)$. Donc $X\in \mathfrak{t}_{X_{qd}}(F)^S$, puis $m(X_{qd})=1$. Soit $F_{1}\in {\cal F}^V$. Pour $\zeta=\pm$, posons $X^{\zeta}=X_{F_{1}}^{\zeta}$. Alors $X^+$ et $X^-$ sont stablement conjugu\'es et ces deux \'el\'ements sont un ensemble de repr\'esentants des classes de conjugaison par $G(F)$ dans leur classe de conjugaison stable. L'assertion (3) nous dit qu'il y a un unique $\zeta$ tel que $X^{\zeta}$ appartienne \`a $\mathfrak{t}_{X^{\zeta}}(F)^S$. Notons $\boldsymbol{\zeta}$ cet $\zeta$. On va le d\'eterminer. Dans ce qui suit, on fixe $F_{1}$ et $\zeta=\pm$, on pose $X=X^{\zeta}$. Comme dans le paragraphe pr\'ec\'edent (en changeant $W$ en $V$), on d\'ecompose $V$ en somme directe orthogonale
 $$V=F_{1}\oplus F_{2}\oplus \tilde{Z}$$
 et, pour $j=1,2$, on introduit la base $(e_{j},e_{-j})$ de $F_{j}$.
 
 Supposons d'abord $d_{W}\geq3$ et $r=0$. Pour $j=1,2$, posons $\epsilon_{j}=e_{j}$ et $\epsilon_{-j}=c_{j}^{-1}e_{-j}$. Choisissons une base hyperbolique $(\epsilon_{\pm k})_{k=3,...,d/2}$ de $\tilde{Z}$. Supposons  $\zeta=\boldsymbol{\zeta}$ et choisissons un \'el\'ement $\gamma\in G(F)$ tel que $\gamma^{-1}X\gamma\in S+\Lambda^S$. Pour $k=1,...,d/2$, on a \'etudi\'e dans la preuve de 10.8  les coordonn\'ees de $\gamma^{-1}\epsilon_{\pm k}$ dans la base $\{v_{0},w_{S},w_{\pm 1},...,w_{\pm m}\}$ de $V$. En notant $Y_{\pm k}$ sa coordonn\'ee sur $v_{0}$, les formules 10.8(3) et 10.8(4) nous disent que
 $$Y_{k}Y_{-k}=\nu_{0}\frac{\prod_{j=1,...,m}(s_{j}^2-x_{k}^2)}{\prod_{k'=1,...,d/2;k'\not=k}(x_{k'}^2-x_{k}^2)}.$$
 Appliquons cela \`a $k=1$. On a $\tau_{F_{1}}(\epsilon_{1})=c^{\boldsymbol{\zeta}}\epsilon_{-1}$, donc aussi $\tau_{F_{1}}(\gamma^{-1}\epsilon_{1})=c^{\boldsymbol{\zeta}}\gamma^{-1}\epsilon_{-1}$ puis $\tau_{F_{1}}(Y_{1})=c^{\boldsymbol{\zeta}}Y_{-1}$. Alors $\chi_{F_{1}}(c^{\boldsymbol{\zeta}}Y_{1}Y_{-1})=1$. On a $x_{1}=a_{1}$ et on d\'ej\`a utilis\'e plusieurs fois que $y^2-a_{1}^2$ \'etait la norme d'un \'el\'ement de $F_{1}^{-1}$ pour tout $y\in F$. On d\'eduit de la formule ci-dessus que
 $$\chi_{F_{1}}(Y_{1}Y_{-1})=\chi_{F_{1}}(\nu_{0}(x_{2}^2-x_{1}^2))=\chi_{F_{1}}(\nu_{0}(Norm_{F_{1}/F}(a_{1})-Norm_{F_{2}/F}(a_{2}))).$$
 D'apr\`es la d\'efinition de $c^{\boldsymbol{\zeta}}$, on en d\'eduit $\boldsymbol{\zeta}=\chi_{F_{1}}(\eta\nu_{0})$.
 
 Passons au cas $d_{W}\geq3$ et $r>0$. On peut supposer que $Z\subset \tilde{Z}$ et que $\epsilon_{\pm k}=v_{\pm(k+r-d/2)}$ pour $k=d/2+1-r,...,d/2$. On peut \'ecrire $X=X_{0}+X_{a}$, avec $X_{0}\in \mathfrak{g}_{0}(F)$ et $X_{a}\in \mathfrak{a}(F)$. L'\'el\'ement $X_{0}$ v\'erifie les m\^emes conditions que $X$, relativement \`a l'espace $V_{0}$. Supposons $\zeta=\chi_{F_{1}}(\eta\nu_{0})$. On vient de voir que $X_{0}$ est conjugu\'e \`a un \'el\'ement de $S+\Sigma_{0,\flat}$ par un \'el\'ement de $G_{0}(F)$. Par un argument que l'on a d\'ej\`a utilis\'e plusieurs fois, $X$ est conjugu\'e \`a un \'el\'ement de $\Xi+S+\Sigma^S$ par un \'el\'ement de $G(F)$. Donc $\boldsymbol{\zeta}=\zeta=\chi_{F_{1}}(\eta\nu_{0})$. 

Consid\'erons enfin le cas $d_{W}=1$. Dans ce cas, $S=0$, $\Xi$ est un \'el\'ement nilpotent r\'egulier et $\Xi+\Lambda$ est une section de Kostant relative  \`a cet \'el\'ement. D'apr\`es [Kot] th\'eor\`eme 5.1, si $X$ est conjugu\'e \`a un \'el\'ement de $\Xi+\Lambda$ par un \'el\'ement de $G(F)$, on a l'\'egalit\'e
$$inv(X)inv(T_{X})=inv_{T}(\Xi),$$
avec les notations de 11.4. L'\'el\'ement $\Xi$ laisse stable l'hyperplan  $D\oplus Z$. Le noyau anisotrope de la restriction de $q$ \`a cet hyperplan est la restriction de $q$ \`a   $D$. Donc $\Xi\in {\cal O}_{\nu_{0}}$. L'\'egalit\'e ci-dessus jointe \`a 11.4(2) entra\^{\i}ne $\boldsymbol{\zeta}=\chi_{F_{1}}(\eta\nu_{0})$.

Le raisonnement ci-dessus s'\'etend \`a tout \'el\'ement de $\omega_{X^+}\cup \omega_{X^-}$. En effet, tout \'el\'ement de cet ensemble est du m\^eme type que $X^{\pm}$, avec des valeurs propres diff\'erentes. On obtient alors
$$m(X^{\zeta})=\left\lbrace\begin{array}{cc}1,&\,\,\text{si}\,\,\zeta=\chi_{F_{1}}(\eta\nu_{0}),\\ 0,&\,\,\text{si}\,\,\zeta=-\chi_{F_{1}}(\eta\nu_{0}).\\ \end{array}\right.$$
En reportant ces valeurs dans l'\'egalit\'e (1), on obtient $J(\theta,f)=\delta_{\nu,\nu_{0}}$. Donc $J(\theta,f)=I(\theta,f)$ et $E(\theta,f)=0$. Comme dans le paragraphe pr\'ec\'edent, cela implique que $E$ est identiquement nulle. $\square$

\bigskip

\section{Preuve des th\'eor\`emes 7.8 et 7.9}

\bigskip

\subsection{Du groupe \`a l'alg\`ebre de Lie}

Consid\'erons les assertions suivantes
\ass{$(th)_{G}$}{Pour tout quasi-caract\`ere $\theta$ sur $H(F)$ et toute fonction tr\`es cuspidale $f\in C_{c}^{\infty}(G(F))$, on a l'\'egalit\'e $lim_{N\to \infty}I_{N}(\theta,f)=I(\theta,f)$.}
\ass{$(th')_{G}$}{Pour tout quasi-caract\`ere $\theta$ sur $H(F)$ et toute fonction tr\`es cuspidale $f\in C_{c}^{\infty}(G(F))$,  dont le support ne contient pas d'\'el\'ement unipotent, on a l'\'egalit\'e $lim_{N\to \infty}I_{N}(\theta,f)=I(\theta,f)$.}
\ass{$(\mathfrak{th})_{\mathfrak{g}}$}{Pour tout quasi-caract\`ere $\theta$ sur $\mathfrak{h}(F)$ et toute fonction tr\`es cuspidale $f\in C_{c}^{\infty}(\mathfrak{g}(F))$, on a l'\'egalit\'e $lim_{N\to \infty}I_{N}(\theta,f)=I(\theta,f)$.}
\ass{$(\mathfrak{th}')_{\mathfrak{g}}$}{Pour tout quasi-caract\`ere $\theta$ sur $\mathfrak{h}(F)$ et toute fonction tr\`es cuspidale $f\in C_{c}^{\infty}(\mathfrak{g}(F))$, dont le support ne contient pas d'\'el\'ement nilpotent, on a l'\'egalit\'e $lim_{N\to \infty}I_{N}(\theta,f)=I(\theta,f)$.}
\ass{Lemme}{ L'assertion $(th)_{G}$ entra\^{\i}ne $(\mathfrak{th})_{\mathfrak{g}}$. L'assertion $(th')_{G}$ entra\^{\i}ne $(\mathfrak{th}')_{\mathfrak{g}}$. }

Preuve. Supposons v\'erifi\'ee $(th)_{G}$. Soient $\theta$ un quasi-caract\`ere sur $\mathfrak{h}(F)$ et $f\in C_{c}^{\infty}(\mathfrak{g}(F))$ une fonction tr\`es cuspidale. On veut montrer que  $E(\theta,f)=0$, avec la notation de 11.3. Dans la preuve de ce paragraphe, on a vu que $E$ \'etait homog\`ene pour la transformation $(\theta,f)\mapsto (\theta^{\lambda},f^{\lambda})$. Cela permet de supposer que le support de $f$  est contenu dans un bon voisinage $\omega$ de $0$ dans $\mathfrak{g}(F)$. Comme dans la preuve du lemme 11.2, on introduit  un quasi-caract\`ere $\boldsymbol{\theta}_{\omega}$ sur $H(F)$ et une fonction ${\bf f}\in C_{c}^{\infty}(G(F))$ tr\`es cuspidale. On v\'erifie que $J(\theta,f)=J_{1,\omega}(\boldsymbol{\theta}_{\omega},{\bf f})$ et $I(\theta,f)=I_{1,\omega}(\boldsymbol{\theta}_{\omega},{\bf f})$. D'apr\`es les lemmes 8.2 et 8.3 et la proposition 10.9, on a $J_{1,\omega}(\boldsymbol{\theta}_{\omega},{\bf f})=lim_{N\to \infty}I_{N}(\boldsymbol{\theta}_{\omega},{\bf f})$, $I_{1,\omega}(\boldsymbol{\theta}_{\omega},{\bf f})=I(\boldsymbol{\theta}_{\omega},{\bf f})$. D'apr\`es  $(th)_{G}$, ces deux termes sont \'egaux. La conclusion $E(\theta,f)=0$ s'ensuit. 

La preuve de la seconde assertion est identique: si le support de $f$ ne contient pas d'\'el\'ement nilpotent, celui de ${\bf f}$ ne contient pas d'\'el\'ement unipotent. $\square$

\bigskip

\subsection{ La r\'ecurrence}

On va raisonner par r\'ecurrence sur $d$. Consid\'erons des donn\'ees $V''$, $W''$, $d''$ analogues \`a $V$, $W$, $d$ (avec le m\^eme $r$). Pour ces donn\'ees, il y a une assertion $(th)_{G''}$  analogue \`a $(th)_{G}$. Consid\'erons l'assertion

\ass{$(th)_{<d}$}{L'assertion $(th)_{G''}$ est v\'erifi\'ee si $d''<d$.}

\ass{Lemme}{L'assertion $(th)_{<d}$ entra\^{\i}ne $(th')_{G}$. Les assertions $(th)_{<d}$ et $(\mathfrak{th})_{\mathfrak{g}}$ entra\^{\i}nent $(th)_{G}$. }

Preuve. Supposons v\'erifi\'ee $(th)_{<d}$, soient $\theta$ un quasi-caract\`ere sur $G(F)$ et $f\in C_{c}^{\infty}(G(F))$ une fonction tr\`es cuspidale, dont le support ne contient pas d'\'el\'ement unipotent. Soient $x\in G_{ss}(F)$  et $\omega_{x}$  un bon voisinage de $0$ dans $\mathfrak{g}_{x}(F)$. Posons $\Omega_{x}=(xexp(\omega_{x}))^G$. On impose \`a $\omega_{x}$ les conditions suivantes. Si $x=1$, on suppose que $\Omega_{1}\cap Supp(f)=\emptyset$ (c'est loisible d'apr\`es l'hypoth\`ese sur $f$). Si $x$ n'est  conjugu\'e \`a aucun \'el\'ement de $H_{ss}(F)$, on impose \`a $\omega_{x}$ les conditions de 8.1. Si $x$ est conjugu\'e \`a un \'el\'ement de $H_{ss}(F)$, on fixe un tel \'el\'ement $x'$, on choisit un bon voisinage $\omega_{x'}$ de $0$ dans $\mathfrak{g}_{x'}(F)$ v\'erifiant les conditions de 8.2 et on d\'efinit $\omega_{x}$ comme l'image par conjugaison de $\omega_{x'}$. Le m\^eme proc\'ed\'e qu'\`a la fin de la preuve de la proposition 6.4 permet de choisir un sous-ensemble fini ${\cal X}\subset G_{ss}(F)$ et de d\'ecomposer $f$ en somme finie  $f=\sum_{x\in {\cal X}}f_{x}$,  de sorte que, pour tout $x\in {\cal X}$, $f_{x}$ soit le produit de $f$ et de la fonction caract\'eristique d'un $G$-domaine inclus dans $\Omega_{x}$.  Par lin\'earit\'e, on peut aussi bien fixer $x\in {\cal X}$ et supposer $f=f_{x}$. Si $x=1$, cette fonction est nulle d'apr\`es le choix de $\omega_{1}$. L'assertion \`a prouver est triviale. D'apr\`es 8.1, c'est aussi le cas si $x$ n'est  conjugu\'e \`a aucun \'el\'ement de $H(F)$. Supposons que $x\not=1$ et $x$ est conjugu\'e \`a un \'el\'ement de $H(F)$. On peut aussi bien supposer $x\in H_{ss}(F)$. Les lemmes 8.2 et 8.3 nous ram\`enent \`a prouver l'\'egalit\'e
$$(1) \qquad lim_{N\to \infty}I_{x,\omega,N}(\theta,f)=I_{x,\omega}(\theta,f).$$
On utilise les notations des sections 8 \`a 10. Si $A_{G'_{x}}\not=\{1\}$, le membre de gauche est nul d'apr\`es la proposition 10.9.  L'ensemble $\underline{\cal T}_{x}$ est vide et le membre de droite est nul lui-aussi. Supposons $A_{G'_{x}}=\{1\}$. En 10.1, on a d\'ecompos\'e $\theta_{x,\omega}$ en une somme finie
$$\theta_{x,\omega}(X)=\sum_{S\in {\cal S}}\hat{j}_{S}(X')\hat{j}^{H''}(S,X'')$$
pour $X\in \omega\cap \mathfrak{h}(F)$. On peut aussi bien \'ecrire
$$\theta_{x,\omega}(X)=\sum_{S\in {\cal S}}\theta'_{S}(X')\theta''_{S}(X'')$$
pour tout $X\in \mathfrak{h}_{x}(F)$, o\`u $\theta'_{S}(X')={\bf 1}_{\omega'}\hat{j}_{S}(X')$ et $\theta''_{S}(X'')={\bf 1}_{\omega''\cap \mathfrak{h}''(F)}\hat{j}^{H''}(S,X'')$.
Remarquons que  $\theta'_{S}$, resp. $\theta''_{S}$,  est un quasi-caract\`ere sur $\mathfrak{g}'_{x}(F)=\mathfrak{h}'_{x}(F)$, resp. $\mathfrak{h}''(F)$,   \`a support compact modulo conjugaison. On peut de m\^eme d\'ecomposer $\theta_{f,x,\omega}$ en somme 
$$\theta_{f,x,\omega}(X)=\sum_{b\in B}\theta'_{f,b}(X')\theta''_{f,b}(X'')$$
o\`u $B$ est un ensemble fini d'indices et, pour tout $b\in B$, $\theta'_{f,b}$, resp. $\theta''_{f,b}$, est un quasi-caract\`ere sur $\mathfrak{g}'_{x}(F)$, resp. $\mathfrak{g}''(F)$, \`a support compact modulo conjugaison. D'apr\`es la proposition 6.4, pour tout $b\in B$, on peut choisir une fonction $f''_{b}\in C_{c}^{\infty}(\mathfrak{g}''(F))$, tr\`es cuspidale et telle que $\theta''_{f,b}=\theta_{f''_{b}}$. En comparant les formules de 7.9 et 8.3, on obtient l'\'egalit\'e
$$I_{x,\omega}(\theta,f)=\sum_{S\in {\cal S}, b\in B} I'(S,b)I(\theta''_{S},f''_{b})$$
o\`u
$$I'(S,b)=\sum_{T'\in {\cal T}_{ell}(G'_{x})}\vert W(G'_{x},T')\vert ^{-1}\nu(T')\int_{\mathfrak{t'}(F)}\theta'_{S}(X')\theta'_{f,b}(X')D^{G'_{x}}(X')dX'.$$
De m\^eme, en comparant les formules de 10.9 et 11.2, on obtient
$$J_{x,\omega}(\theta,f)=\sum_{S\in{\cal S},b\in B}I'(S,b)J(\theta''_{S},f''_{b}).$$
Puisque $x\not=1$, on a $dim(W'')<d$. L'assertion $(th)_{<d}$ et le lemme 12.1 entra\^{\i}nent que $(\mathfrak{th})_{\mathfrak{g''}}$ est v\'erifi\'ee.  Donc $J(\theta''_{S},f''_{b})=I(\theta''_{S},f''_{b})$ gr\^ace au lemme 11.2. Alors les formules ci-dessus et la proposition 10.9 impliquent l'\'egalit\'e (1) qu'il fallait prouver.

La seconde assertion de l'\'enonc\'e se d\'emontre de m\^eme. On ne peut plus \'eliminer le point $x=1$. Mais, pour ce point, l'\'egalit\'e $J(\theta''_{S},f''_{b})=I(\theta''_{S},f''_{b})$ provient de l'assertion $(\mathfrak{th})_{\mathfrak{g}}$ que l'on suppose v\'erifi\'ee. $\square$

\bigskip

\subsection{Finale}

Pour prouver le th\'eor\`eme 7.8, et aussi le th\'eor\`eme 7.9 d'apr\`es le lemme 12.1, il suffit de prouver que l'assertion $(th)_{<d}$ entra\^{\i}ne $(th)_{G}$. Or $(th)_{<d}$ entra\^{\i}ne $(th')_{G}$ (lemme 12.2), qui entra\^{\i}ne $(\mathfrak{th}')_{\mathfrak{g}}$ (lemme 12.1). Cette derni\`ere assertion entra\^{\i}ne $(\mathfrak{th})_{\mathfrak{g}}$: c'est une reformulation de la proposition 11.1. Jointe \`a $(th)_{<d}$, cette derni\`ere assertion entra\^{\i}ne $(th)_{G}$ (lemme 12.2). Cela ach\`eve la preuve.

\bigskip

\section{Un cas de la version faible de la conjecture locale de Gross-Prasad}

\bigskip

\subsection{D\'efinition des multiplicit\'es}

 Soit $G$ un groupe r\'eductif connexe d\'efini sur $F$. Une repr\'esentation  de $G(F)$ est un couple $(\pi,E_{\pi})$, o\`u $E_{\pi}$ est l'espace de la repr\'esentation (pour nous, un espace vectoriel complexe), et $\pi$ un homomorphisme de $G(F)$ dans $GL_{{\mathbb C}}(E_{\pi})$. On oubliera souvent l'un des termes en la notant simplement $\pi$ ou $E_{\pi}$. On notera aussi une classe d'isomorphie de repr\'esentations  comme une repr\'esentation dans cette classe. On note $\check{\pi}$ la repr\'esentation contragr\'ediente de $\pi$. On note $Irr(G)$ l'ensemble des  classes d'isomorphie de repr\'esentations admissibles irr\'eductibles de $G(F)$ et $Temp(G)$ le sous-ensemble des repr\'esentations
qui sont aussi temp\'er\'ees. A tout \'el\'ement $\pi\in Irr(G)$ est associ\'e un caract\`ere $\theta_{\pi}$, que l'on peut consid\'erer comme une distribution sur $G(F)$ ou comme une fonction d\'efinie presque partout. Dans cette derni\`ere interpr\'etation, $\theta_{\pi}$ est un quasi-caract\`ere sur $G(F)$ d'apr\`es un th\'eor\`eme de Harish-Chandra ([HCDS] th\'eor\`eme 16.2).

Consid\'erons la situation de 7.2, soient $(\pi,E_{\pi})\in Irr(G)$ et $(\sigma,E_{\sigma})\in Irr(H)$. Notons $Hom_{H,\xi}(\pi,\sigma)$ l'espace des applications lin\'eaires $\varphi:E_{\pi}\to E_{\sigma}$ telles que
$$\varphi(\pi(hu)e)=\xi(u)\sigma(h)\varphi(e)$$
pour tous $h\in H(F)$, $u\in U(F)$, $e\in E_{\pi}$. Cet espace d\'epend des constantes $\xi_{i}$ qui nous ont permis de d\'efinir $\xi$, mais de fa\c{c}on inessentielle. En effet, si on associe un caract\`ere $\xi'$ \`a d'autres constantes, on v\'erifie qu'il existe $a\in A(F)$ tel que l'application $\varphi\mapsto \varphi\circ\pi(a)$ soit un isomorphisme de $Hom_{H,\xi}(\pi,\sigma)$ sur $Hom_{H,\xi'}(\pi,\sigma)$. Tout ce qui suit est donc essentiellement ind\'ependant des constantes $\xi_{i}$.

On a
$$(1) \qquad dim_{{\mathbb C}}(Hom_{H,\xi}(\pi,\sigma))\leq1.$$
Si $r=0$, c'est le th\'eor\`eme 1' de [AGRS]. Ce r\'esultat est g\'en\'eralis\'e \`a tout $r$ dans [GGP] corollaire 20.4 (pour les groupes unitaires, mais la preuve est la m\^eme pour les groupes sp\'eciaux orthogonaux). On pose
$$m(\sigma,\pi)=dim_{{\mathbb C}}(Hom_{H,\xi}(\pi,\sigma)).$$
Soit $T\in {\cal T}$. En appliquant les d\'efinitions de 7.3 aux quasi-caract\`eres $\theta=\theta_{\check{\sigma}}$ et $\tau=\theta_{\pi}$, on d\'efinit les fonctions $c_{\theta_{\check{\sigma}}}$ et $c_{\theta_{\pi}}$ sur $T(F)$, que l'on note simplement $c_{\check{\sigma}}$ et $c_{\pi}$. On pose
$$m_{geom}(\sigma,\pi)=\sum_{T\in {\cal T}}\vert W(H,T)\vert ^{-1}\nu(T)\int_{T(F)}c_{\check{\sigma}}(t)c_{\pi}(t)D^H(t)\Delta(t)^rdt.$$
Cette expression est absolument convergente d'apr\`es la proposition 7.3.

\ass{Proposition}{Supposons $\pi\in Irr(G)$, $\sigma\in Irr(H)$ et $\pi$ supercuspidale. Alors on a l'\'egalit\'e $m(\sigma,\pi)=m_{geom}(\sigma,\pi)$.}

Preuve. Si $V$ est hyperbolique de dimension $2$, on v\'erifie directement que les deux termes de l'\'egalit\'e valent $1$. On exclut ce cas. Introduisons les repr\'esentations de $G(F)$: $\boldsymbol{\rho}=Ind_{H(F)U(F)}^{G(F)}(\sigma\otimes \bar{\xi})$ et $\rho=ind_{H(F)U(F)}^{G(F)}(\check{\sigma}\otimes \xi)$. La premi\`ere est l'induite lisse de la repr\'esentation $\sigma\otimes \bar{\xi}$ de $H(F)U(F)$, o\`u $\bar{\xi}$ est le conjugu\'e complexe de $\xi$, la seconde est l'induite \`a supports compacts de la repr\'esentation $\check{\sigma}\otimes \xi$. Par r\'eciprocit\'e de Frobenius pour la premi\`ere \'egalit\'e et d'apr\`es un r\'esultat standard pour la seconde, on a
$$Hom_{H,\bar{\xi}}(\pi,\sigma)=Hom_{G(F)}(\pi,\boldsymbol{\rho})=Hom_{G(F)}( \rho,\check{\pi}).$$
Puisque $\check{\pi}$ est supercuspidale et le centre de $G(F)$ est fini, la th\'eorie de Bernstein nous dit que $ \rho$ se d\'ecompose en somme d'une repr\'esentation $\tau$ dont aucun sous-quotient n'est isomorphe \`a $\check{\pi}$ et d'un certain nombre de facteurs tous isomorphes \`a $\check{\pi}$. Le nombre de ces facteurs est pr\'ecis\'ement $m(\sigma,\pi)$. Soit $f$ un coefficient de $\pi$. On a l'\'egalit\'e
$$trace(\check{\pi}(f)\vert E_{\check{\pi}})=f(1)d(\pi)^{-1}$$
o\`u $d(\pi)$ est le degr\'e formel de $\pi$. L'op\'erateur $\tau(f)$ est nul. Donc $\rho(f)$ est de rang fini et 
$$(2) \qquad trace(\rho(f))=m(\sigma,\pi)f(1)d(\pi)^{-1}.$$
Montrons que, si $N$ est un entier assez grand,
$$(3)\qquad trace(\rho(f))=I_{N}(\theta_{\check{\sigma}},f).$$
Fixons un sous-groupe ouvert compact $K'\subset K$ tel que $f$ soit biinvariante par $K'$. Notons $\Omega_{N}$ le support de la fonction $\kappa_{N}$ et $E_{\rho,N}$ le sous-espace de $E_{\rho}$ form\'e des fonctions \`a support dans $\Omega_{N}$. Selon l'usage, notons $E_{\rho,N}^{K'}$ le sous-espace des \'el\'ements invariants par l'action de $K'$. Puisque l'image de $\rho(f)$ est de dimension finie, elle est contenue dans $E_{\rho,N}^{K'}$ si $N$ est assez grand. Alors $trace(\rho(f))$ est la trace de la restriction de $\rho(f)$ \`a $E_{\rho,N}^{K'}$. Fixons un ensemble de repr\'esentants $\Gamma_{N}$ des doubles classes $H(F)U(F)\backslash \Omega_{N}/K'$. Notons $\Gamma'_{N}$ le sous-ensemble des $\gamma\in \Gamma_{N}$ tels que $\xi$ soit trivial sur $\gamma K' \gamma^{-1}\cap U(F)$. Pour $\gamma\in \Gamma_{N}$, posons $K^H[\gamma]=\gamma K'\gamma^{-1}\cap H(F)$ et fixons une base ${\cal B}[\gamma]$ de l'espace $E_{{\sigma}}^{K^H[\gamma]}$. Notons $\{\check{b}; b\in {\cal B}[\gamma]\}$ la base duale de $E_{\check{\sigma}}^{K^H[\gamma]}$. Pour $\gamma\in \Gamma'_{N}$ et $b\in {\cal B}[\gamma]$, il existe un unique \'el\'ement $\varphi[b,\gamma]\in E_{\rho}$, \`a support dans $H(F)U(F)\gamma K'$, invariant \`a droite par $K'$ et tel que $\varphi[b,\gamma](\gamma)=\check{b}$. L'ensemble $\{\varphi[b,\gamma]; \gamma\in \Gamma'_{N},b\in {\cal B}[\gamma]\}$ est une base de $E_{\rho,N}^{K'}$. Alors, la trace de $\rho(f)$ agissant dans $E_{\rho,N}^{K'}$ est
$$\sum_{\gamma\in \Gamma'_{N}, b\in {\cal B}[\gamma]}<(\rho(f)\varphi[b,\gamma])(\gamma),b>.$$
On a
$$<(\rho(f)\varphi[b,\gamma])(\gamma),b>=\int_{G(F)}<\varphi[b,\gamma](\gamma g),b>f(g)dg$$
$$\qquad =mes(H(F)U(F)\backslash H(F)U(F)\gamma K')\int_{H(F)U(F)}<\check{\sigma}(h)\check{b},b>\xi(u)f(\gamma^{-1}hu\gamma)du\,dh$$
$$\qquad =mes(H(F)U(F)\backslash H(F)U(F)\gamma K')\int_{H(F)}<\check{\sigma}(h)\check{b},b>{^{\gamma}f}^{\xi}(h)dh.$$
La somme sur $b\in {\cal B}[\gamma]$ de  cette int\'egrale est $trace(\check{\sigma}(^{\gamma}f^{\xi}))$, ou encore $I(\theta_{\check{\sigma}},f,\gamma)$. On v\'erifie que ce terme est nul pour $\gamma\in \Gamma_{N}$, $\gamma\not\in \Gamma'_{N}$ et on obtient
$$trace(\rho(f))=\sum_{\gamma\in \Gamma_{N}}mes(H(F)U(F)\backslash H(F)U(F)\gamma K')I(\theta_{ \check{\sigma}},f,\gamma)$$
$$\qquad =\int_{H(F)U(F)\backslash G(F)}I(\theta_{\check{\sigma}},f,g)\kappa_{N}(g)dg,$$
d'o\`u (3).

 La fonction $f$ est tr\`es cuspidale. On calcule son quasi-caract\`ere associ\'e
 $$(4) \qquad \theta_{f}=f(1)d(\pi)^{-1}\theta_{\pi}.$$
 En effet, cela r\'esulte de notre d\'efinition 5.3 de $\theta_{f}$ et de [A6] th\'eor\`eme p.3 (il y a un probl\`eme de passage \`a la contragr\'ediente dans ce th\'eor\`eme).
 
 Gr\^ace \`a (3) et au th\'eor\`eme 7.8, on a $trace(\rho(f))=I(\theta_{\check{\sigma}},f)$. 
 Gr\^ace \`a (4), $I(\theta_{\check{\sigma}},f)=f(1)d(\pi)^{-1}m_{geom}(\sigma,\pi)$. Alors la proposition r\'esulte de (2). $\square$
 
 \bigskip
 
 \subsection{Les $L$-paquets}

Consid\'erons de nouveau la situation de 7.2. On suppose d\'esormais que $G$ et $H$ sont quasi-d\'eploy\'es, autrement dit que $d_{an}(V)\leq2$, $d_{an}(W)\leq2$. On affecte les notations d'un indice $i$: $V_{i}$, $W_{i}$, $G_{i}$, $H_{i}$ etc...

A \'equivalence pr\`es,  il existe au plus un couple $(V',q')$ analogue \`a $(V_{i},q_{i})$, tel que $dim(V')=d$, le discriminant de $q'$ soit \'egal \`a celui de $q_{i}$, mais $(V',q')$ ne soit pas \'equivalent \`a $(V_{i},q_{i})$. Un tel couple v\'erifie $d_{an}(V')+d_{an}(V_{i})=4$.  Le couple $(V',q')$ existe si et seulement si $d+d_{an}(V_{i})\geq 4$. S'il existe, on le note $(V_{a},q_{a})$ et on introduit pour ce couple les m\^emes objets que pour $(V_{i},q_{i})$, affect\'es d'un indice $a$. On introduit de m\^eme un couple $(W_{a},q_{W_{a}})$, s'il existe, c'est-\`a-dire si $d_{W_{i}}+d_{an}(W_{i})\geq4$. Quand ces deux couples existent, on a l'\'egalit\'e $V_{a}=W_{a}\oplus D$ et $q_{a}$ est la somme orthogonale de $q_{W_{a}}$ et de la forme d\'ej\`a fix\'ee sur $D$.

{\bf Remarque.} Les indices $i$ et $a$ signifient "isotrope" et "anisotrope", les donn\'ees affect\'ees de l'indice $i$ ayant tendance \`a \^etre "plus isotropes" que celles affect\'ees de l'indice $a$. Pr\'ecis\'ement, on a $d_{an}(V_{i})+d_{an}(W_{i})<d_{an}(V_{a})+d_{an}(W_{a})$.
\bigskip

 Pour simplifier, si le couple $(V_{a},q_{a})$, resp. $(W_{a},q_{W_{a}})$, n'existe pas, on pose $Temp(G_{a})=\emptyset$, resp. $Temp(H_{a})=\emptyset$. Selon une conjecture due essentiellement \`a Langlands, les ensembles de repr\'esentations $Temp(G_{i})$ et $Temp(G_{a})$  se d\'ecomposent en r\'eunions disjointes de $L$-paquets, qui poss\`edent les propri\'et\'es (1), (2) et (3) ci-dessous.
 
 (1) Soit $\Pi$ un $L$-paquet de $Temp(G_{i})$ ou $Temp(G_{a})$. L'ensemble $\Pi$ est fini. Posons $\theta_{\Pi}=\sum_{\pi\in \Pi}\theta_{\pi}$. Alors $\theta_{\Pi}$ est une distribution stable. 
 
 (2) Il existe une application   qui, \`a un $L$-paquet dans $Temp(G_{a})$, associe un $L$-paquet dans $Temp(G_{i})$, et satisfait les conditions suivantes. Elle est injective. Soient $\Pi_{a}$ un $L$-paquet dans $Temp(G_{a})$ et $\Pi_{i}$ son image. Alors $(-1)^d\theta_{\Pi_{a}}$ est le transfert \`a $G_{a}(F)$ de la distribution $\theta_{\Pi_{i}}$ sur $G_{i}(F)$. Soit $\Pi_{i}$ un $L$-paquet dans $Temp(G_{i})$ qui n'est pas dans l'image de l'application. Alors le transfert \`a $G_{a}(F)$ de $\theta_{\Pi_{i}}$ est nul (si le groupe $G_{a}$ existe).
 
 {\bf Remarque.} Le signe $(-1)^d$ s'interpr\`ete comme $(-1)^{rang_{F}(G_{a})-rang_{F}(G_{i})}$, o\`u $rang_{F}(G_{a})$, resp. $rang_{F}(G_{i})$, est la dimension d'un sous-tore d\'eploy\'e maximal de $G_{a}$, resp. $G_{i}$.
 \bigskip
 
On sait d\'efinir la notion de mod\`ele de Whittaker d'une repr\'esentation dans $Irr(G_{i})$. Plus pr\'ecis\'ement, une telle notion est associ\'ee \`a chaque orbite nilpotente r\'eguli\`ere ${\cal O}$ de $\mathfrak{g}_{i}(F)$. Soient ${\cal O}$ une telle orbite et $\bar{N}\in {\cal O}$. On peut compl\'eter $\bar{N}$ en un $\mathfrak{sl}_{2}$-triplet qui d\'etermine un sous-tore maximal $T$  et un sous-groupe de Borel $B $ de $G_{i}$ de sorte que $T\subset B$ et $\bar{N}\in \bar{\mathfrak{b}}(F)$. Notons $U_{B}$ le radical unipotent de $B$ et d\'efinissons une fonction $\xi_{\bar{N}}$ sur $U_{B}(F)$ par $\xi_{\bar{N}}(exp(N))=\psi(<\bar{N},N>)$. C'est un caract\`ere. Pour $\pi\in Irr(G_{i})$, on dit que $\pi$ admet un mod\`ele de Whittaker relatif \`a ${\cal O}$ s'il existe une forme lin\'eaire $l$ sur $E_{\pi}$, non nulle et telle que $l(\pi(u)e)=\xi_{\bar{N}}(u)l(e)$ pour tous $u\in U_{B}(F)$, $e\in E_{\pi}$. La derni\`ere propri\'et\'e des $L$-paquets est

(3) pour tout $L$-paquet $\Pi$ dans $Temp(G_{i})$ et toute orbite nilpotente r\'eguli\`ere ${\cal O}$ dans $\mathfrak{g}_{i}(F)$, il existe un et un seul \'el\'ement de $\Pi$ qui admet un mod\`ele de Whittaker relatif \`a ${\cal O}$.

{\bf Remarques.} La propri\'et\'e (1) pour le groupe $G_{i}$ est annonc\'ee par Arthur. Il n'y a gu\`ere de doute que, dans un avenir proche, les travaux d'Arthur d\'emontreront \'egalement cette propri\'et\'e pour le groupe $G_{a}$ et la propri\'et\'e (2). Konno a montr\'e que des r\'esultats \'egalement annonc\'es par Arthur entra\^{\i}naient la propri\'et\'e (3) ([Konno] th\'eor\`eme 3.4). Signalons que cette propri\'et\'e (3) est une conjecture de Shahidi.
\bigskip

{\bf Dans la suite de l'article, on admet l'existence de $L$-paquets poss\'edant ces propri\'et\'es}. Bien \'evidemment, on les admet aussi pour les groupes $H_{i}$ et $H_{a}$ (l'entier $d$ \'etant chang\'e en $d_{W}$).

\bigskip

\subsection{Un r\'esultat en direction de la conjecture locale  de Gross-Prasad}

Soient $\Pi_{i}$ un $L$-paquet dans $Temp(G_{i})$ et $\Sigma_{i}$ un $L$-paquet dans $Temp(H_{i})$. Si $\Pi_{i}$ est dans l'image de l'application 13.2(2), on note $\Pi_{a}$ le $L$-paquet dans $Temp(G_{a})$ dont il est l'image. Sinon, on pose $\Pi_{a}=\emptyset$. On d\'efinit de m\^eme $\Sigma_{a}$. Pour $(\sigma,\pi)\in (\Sigma_{i}\times\Pi_{i})\cup(\Sigma_{a}\times\Pi_{a})$, on d\'efinit la multiplicit\'e $m(\sigma,\pi)$.

\ass{Th\'eor\`eme}{Supposons que tout \'el\'ement de $\Pi_{i}\cup \Pi_{a}$ soit supercuspidal. Alors il existe un unique couple $(\sigma,\pi)\in (\Sigma_{i}\times\Pi_{i})\cup(\Sigma_{a}\times\Pi_{a})$ tel que $m(\sigma,\pi)=1$.}

Sous les hypoth\`eses indiqu\'ees, c'est une partie de la conjecture 6.9 de [GP]. La d\'emonstration sera donn\'ee au paragraphe 13.6.

\bigskip

\subsection{Calcul de fonctions $\hat{j}$}

On consid\`ere la situation du paragraphe 11.4 dont on reprend les notations. On reprend aussi la notation ${\cal F}_{V}$ introduite en 11.5. Dans le cas o\`u $d$ est pair et $d\geq4$, on a d\'efini des \'el\'ements $X_{F_{1}}^{\pm}$. On pose $\epsilon_{F_{1}}=\chi_{F_{1}}(-Norm_{F_{2}/F}(a_{2}))$. En fait, ce terme ne d\'epend pas de $a_{2}$ puisque l'image de $Norm_{F_{2}/F}(a_{2})$ dans $F^{\times}/F^{\times 2}$ est uniquement d\'etermin\'ee. On  note $T_{F_{1}}$ le commutant de $X_{F_{1}}^+$ dans $G$.

\ass{Lemme}{(i) Supposons $d$ impair ou $d\leq 2$. On a l'\'egalit\'e
$$\hat{j}({\cal O}_{reg},X_{qd})=\vert W^G\vert D^G(X_{qd})^{-1/2}.$$

(ii) Supposons $d$ pair, $d\geq4$. Soit $\nu\in {\cal N}^V$. On a l'\'egalit\'e
$$\hat{j}({\cal O}_{\nu},X_{qd}) =\vert {\cal N}^V\vert ^{-1}\vert W^G\vert D^G(X_{qd})^{-1/2}.$$
Pour $F_{1}\in {\cal F}^V$, on a l'\'egalit\'e
$$\hat{j}({\cal O}_{\nu},X_{F_{1}}^+)=-\hat{j}({\cal O}_{\nu},X_{F_{1}}^-)=\epsilon_{F_{1}}\chi_{F_{1}}(\nu\eta)\frac{\vert W(G,T_{F_{1}})\vert }{2\vert {\cal N}^V\vert }D^G(X_{F_{1}}^+)^{-1/2}.$$}

Preuve. Supposons $d$ impair ou $d\leq2$. D'apr\`es 2.6(4) et le lemme 11.4, on a l'\'egalit\'e
$$\hat{j}({\cal O}_{reg},X_{qd})=\hat{j}(\lambda X_{qd},X_{qd})$$
pour tout $\lambda\in F^{\times 2}$ assez voisin de $0$. Le commutant $T_{qd}$ de $X_{qd}$ dans $G$ est un L\'evi minimal. La formule 2.6(5) exprime $\hat{j}^G(\lambda X_{qd},X_{qd})$ \`a l'aide de la fonction $\hat{j}^{T_{qd}}$. Mais, pour tout tore $T$, la fonction $(X,Y)\mapsto \hat{j}^T(X,Y)$ est constante de valeur $1$, ainsi qu'il r\'esulte de sa d\'efinition. D'autre part, un \'el\'ement de $\mathfrak{t}_{qd}(F)$ est conjugu\'e \`a $X_{qd}$ par un \'el\'ement de $G(F)$ si et seulement s'il l'est par un \'el\'ement de $Norm_{G(F)}(T_{qd})$. Il y a $\vert W^G\vert $ tels \'el\'ements. Alors, la formule 2.6(5) conduit \`a l'\'egalit\'e du (i) de l'\'enonc\'e.

Dans la situation du (ii), notons ${\cal X}$ l'ensemble form\'e des \'el\'ements $X_{qd}$ et $X_{F_{1}}^{\pm}$, pour $F_{1}\in {\cal F}^V$. La formule 2.6(4) et le lemme 11.4 entra\^{\i}nent que pour $\lambda\in F^{\times 2}$ assez voisin de $0$, on a l'\'egalit\'e
$$(1) \qquad \hat{j}({\cal O}_{\nu},Y)=\vert {\cal N}^V\vert ^{-1}(\hat{j}(\lambda X_{qd},Y)+\sum_{F_{1}\in{\cal F}^V}\chi_{F_{1}}(\nu \eta)(\hat{j}(\lambda X_{F_{1}}^+,Y)-\hat{j}(\lambda X_{F_{1}}^-,Y)))$$
pour tout $Y\in {\cal X}$. Fixons un tel $\lambda$. 

Comme ci-dessus, on a
$$ \hat{j}(\lambda X_{qd},X_{qd})=\vert W^G\vert D^G(X_{qd})^{-1/2}.$$
D'autre part,  pour $F_{1}\in {\cal F}^V$, l'\'el\'ement $X_{F_{1}}^{\pm}$ n'est conjugu\'e \`a aucun \'el\'ement de $\mathfrak{t}_{qd}(F)$ et la formule 2.6(5) entra\^{\i}ne
$$\hat{j}(\lambda X_{qd},X_{F_{1}}^{\pm})=0.$$

Soit $F_{1}\in {\cal F}^V$. Il nous faut calculer $\hat{j}(\lambda X_{F_{1}}^{+},Y)-\hat{j}(\lambda X_{F_{1}}^{-},Y)$ pour $Y\in {\cal X}$. On va d'abord supposer $d=4$. Pour $j=1,2$, notons $G_{i}$ le groupe sp\'ecial orthogonal de $F_{j}$ muni de la forme $Norm_{F_{j}/F}$. C'est un tore de dimension $1$. Le groupe $G'=G_{1}\times G_{2}$ est un groupe endoscopique elliptique de $G$. La classe de conjugaison stable de $X_{F_{1}}^{\pm}$ dans $\mathfrak{g}(F)$ est l'image de la classe de conjugaison stable d'un \'el\'ement $X'\in \mathfrak{g}'(F)$, cette derni\`ere classe se r\'eduisant \`a $\{X'\}$ puisque $G'$ est un tore. On peut normaliser le facteur de transfert $\Delta_{G,G'}$ de sorte que $\Delta_{G,G'}(X',X_{F_{1}}^{\zeta})=\zeta$ pour $\zeta=\pm=\pm 1$. Gr\^ace \`a Ngo Bao Chau, la conjecture 1.2 de [W4] est maintenant d\'emontr\'ee. La fonction $\hat{i}^{G}(X,Y)$ de [W4] est \'egale \`a $\hat{j}^G(X,Y)D^G(Y)^{1/2}$. La fonction $\hat{i}^{G'}$ est constante de valeur $1$ puisque $G'$ est un tore. Avec les notations de cette r\'ef\'erence, on a donc l'\'egalit\'e
 $$\gamma_{\psi}(\mathfrak{g})(\hat{j}^G(\lambda X_{F_{1}}^+,Y)-\hat{j}^G(\lambda X_{F_{1}}^-,Y))D^G(Y)^{1/2}=\gamma_{\psi}(\mathfrak{g}')\sum_{Z\in \mathfrak{g}'(F)}\Delta_{G,G'}(Z,Y)$$
 pour tout $Y\in \mathfrak{g}_{reg}(F)$. La classe de conjugaison stable de l'\'el\'ement $X_{qd}$ n'est l'image d'aucun \'el\'ement de $\mathfrak{g}'(F)$.  Il n'y a donc aucun $Z$ pour lequel $\Delta_{G,G'}(Z,X_{qd})\not=0$. On obtient
$$ \hat{j}^G(\lambda X_{F_{1}}^+,X_{qd})-\hat{j}^G(\lambda X_{F_{1}}^-,X_{qd})=0.$$
Soit $F'_{1}\in {\cal F}^V$. Si $F'_{1}\not=F_{1}$, la classe de conjugaison stable d'un \'el\'ement $X_{F'_{1}}^{\zeta}$, pour $\zeta=\pm$, n'est l'image d'aucun \'el\'ement de $\mathfrak{g}'(F)$ et on obtient de m\^eme
$$\hat{j}^G(\lambda X_{F_{1}}^+,X_{F'_{1}}^{\zeta})-\hat{j}^G(\lambda X_{F_{1}}^-,X_{F'_{1}}^{\zeta})=0.$$
Si $F_{1}\not=F_{2}$,  la classe de conjugaison stable de $X_{F_{1}}^{\zeta}$ est l'image d'un unique \'el\'ement de $\mathfrak{g}'(F)$, \`a savoir $X'$ et on conna\^{\i}t la valeur $\Delta_{G,G'}(X',X_{F_{1}}^{\zeta})=\zeta$ (en identifiant $\zeta$ \`a un \'el\'ement de $\{\pm 1\}$). Si $F_{1}=F_{2}$, la classe de conjugaison stable de $X_{F_{1}}^{\zeta}$ est l'image de deux \'el\'ements de $\mathfrak{g}'(F)$: $X'$ et l'\'el\'ement $X''$ obtenu en \'echangeant les deux facteurs de $X'$. Un argument g\'en\'eral nous dit que, parce que $G$ est quasi-d\'eploy\'e, le facteur de transfert est insensible \`a l'action d'un automorphisme du groupe endoscopique $G'$. Donc $\Delta_{G,G'}(X'',X_{F_{1}}^{\zeta})=\Delta_{G,G'}(X',X_{F_{1}}^{\zeta})=\zeta$. On obtient
$$\hat{j}^G(\lambda X_{F_{1}}^+,X_{F_{1}}^{\zeta})-\hat{j}^G(\lambda X_{F_{1}}^-,X_{F_{1}}^{\zeta})=\zeta (1+\delta_{F_{1},F_{2}})\gamma_{\psi}(\mathfrak{g}')\gamma_{\psi}(\mathfrak{g})^{-1}D^G(X_{F_{1}}^{\pm})^{-1/2}.$$
Remarquons que $D^G(X_{F_{1}}^+)=D^G(X_{F_{1}}^-)$. Si $F_{1}\not= F_{2}$, le groupe $W(G,T_{F_{1}})$ a deux \'el\'ements: l'action de l'\'el\'ement non trivial de ce groupe envoie $X_{F_{1}}^{\zeta}$ sur un \'el\'ement analogue o\`u les valeurs propres $a_{1}$ et $a_{2}$ sont chang\'ees en $-a_{1}$ et $-a_{2}$. Si $F_{1}=F_{2}$, le groupe $W(G,T_{F_{1}})$ a $4$ \'el\'ements: on peut de plus permuter les deux facteurs. Donc
$$1+\delta_{F_{1},F_{2}}=\frac{\vert W(G,T_{F_{1}})\vert }{2}.$$
Il reste \`a calculer les facteurs $\gamma_{\psi}$. Ces facteurs sont les "constantes de Weil" associ\'ees \`a $\psi$ et aux formes quadratiques $(X,Y)\mapsto trace(XY)/2$ sur $\mathfrak{g}(F)$ et $\mathfrak{g}'(F)$. Fixons $\xi_{1},\xi_{2}\in F^{\times}$ tels que $F_{j}=F(\sqrt{\xi_{j}})$ pour $j=1,2$. Si $G$ est d\'eploy\'e, posons $\xi=1$. Si $G$ n'est pas d\'eploy\'e, soit $\xi\in F^{\times}$ tel que $E=F(\sqrt{\xi})$. La forme quadratique sur $\mathfrak{g}(F)$ a m\^eme noyau anisotrope que la forme  $x^2+\xi y^2$ de dimension $2$. La forme quadratique sur $\mathfrak{g}'(F)$ est \'equivalente \`a la forme $\xi_{1}x^2+\xi_{2}y^2$. Ces deux formes ont m\^eme d\'eterminant. Elles sont \'equivalentes si et seulement si $\chi_{F_{1}}(\xi_{2})=1$, ou, ce qui revient au m\^eme, si $\epsilon_{F_{1}}=1$. On en d\'eduit l'\'egalit\'e
$$\gamma_{\psi}(\mathfrak{g}')\gamma_{\psi}(\mathfrak{g})^{-1}=\epsilon_{F_{1}},$$
puis
$$(2)\qquad \hat{j}^G(\lambda X_{F_{1}}^+,X_{F_{1}}^{\zeta})-\hat{j}^G(\lambda X_{F_{1}}^-,X_{F_{1}}^{\zeta})=\zeta\epsilon_{F_{1}}\frac{\vert W(G,T_{F_{1}})\vert }{2}D^G(X_{F_{1}}^+)^{-1/2}.$$

{\bf Remarque.} Le calcul serait le m\^eme si l'on rempla\c{c}ait $X_{F_{1}}^{\zeta}$ par un \'el\'ement similaire, mais de valeurs propres diff\'erentes. Par exemple si l'on rempla\c{c}ait $X_{F_{1}}^{\zeta}$ par un  conjugu\'e par un \'el\'ement du groupe $G^+(F)$.
\bigskip
 
Levons l'hypoth\`ese $d=4$. Dans la construction de 11.4 des \'el\'ements $X_{F_{1}}^{\pm}$, on a fix\'e un espace hyperbolique $\tilde{Z}$ et un sous-tore d\'eploy\'e maximal $\tilde{T}$ du groupe sp\'ecial orthogonal de cet espace. Notons $M$ le commutant de $\tilde{T}$ dans $G$. On a $M=\tilde{T}G'$, o\`u $G'$ est le groupe sp\'ecial orthogonal de l'espace quadratique $F_{1}\oplus F_{2}$ de dimension $4$. Soit $X\in {\cal X}$ et $Y$ un \'el\'ement de $ \mathfrak{m}(F)$ conjugu\'e \`a $X$ par un \'el\'ement de $G(F)$. Il est clair que $Y=Y_{\tilde{T}}+Y'$, o\`u $Y_{\tilde{T}}\in \tilde{\mathfrak{t}}(F)$ et $Y'$ est un \'el\'ement de $\mathfrak{g}'(F)$ construit de la m\^eme fa\c{c}on que $X$, \'eventuellement conjugu\'e par un \'el\'ement du groupe $G^{_{'}+}(F)$. Comme on l'a dit ci-dessus, la fonction $\hat{j}^{\tilde{T}}$ est constante de valeur $1$. La formule 2.6(4)  et nos r\'esultats ci-dessus appliqu\'es au groupe $G'$ conduisent \`a l'\'egalit\'e
$$ \hat{j}^G(\lambda X_{F_{1}}^+,X)-\hat{j}^G(\lambda X_{F_{1}}^-,X)=0$$
 pour $X\in {\cal X}$, $X\not=X_{F_{1}}^{\pm}$. Pour $X=X_{F_{1}}^{\pm}$, on doit calculer le nombre de classes de conjugaison par $M(F)$ contenues dans l'intersection de $\mathfrak{m}(F)$ et de la classe de conjugaison par $G(F)$ de $X$.   Parce que $X$ est elliptique dans $\mathfrak{m}(F)$, tout \'el\'ement $g\in G(F)$ tel que $gXg^{-1}\in\mathfrak{m}(F)$ normalise $M$. Le nombre cherch\'e est donc le nombre d'\'el\'ements du groupe $Norm_{G(F)}(M)/M(F) $. On v\'erifie que  tout \'el\'ement de ce quotient a un repr\'esentant dans $Norm_{G(F)}(T_{F_{1}})$. Le nombre cherch\'e est donc $\vert W(G,T_{F_{1}})\vert \vert W(M,T_{F_{1}})\vert ^{-1}$. Le deuxi\`eme facteur est l'inverse de celui qui appara\^{\i}t dans  l'\'egalit\'e (2) appliqu\'ee au groupe $G'$. On obtient alors la m\^eme \'egalit\'e (2) pour notre groupe $G$.
 
 On a maintenant calcul\'e tous les termes intervenant dans la formule (1). Cette formule conduit \`a l'\'egalit\'e du (ii) de l'\'enonc\'e. $\square$

 \bigskip
 
 \subsection{Classes de conjugaison stable de tores}
 
 Dans ce paragraphe, on travaille soit avec l'une des s\'eries de donn\'ees $V_{i}$, $W_{i}$, $G_{i}$ etc.. ou $V_{a}$, $W_{a}$, $G_{a}$ etc..., soit avec les deux. Dans le premier cas, pour simplifier,  on note $\flat$ l'indice $i$ ou $a$. On a d\'efini l'ensemble $\underline{\cal T}_{\flat}$ en 7.3. A $T\in \underline{\cal T}_{\flat}$, on a associ\'e des espaces $W'_{\flat}$ etc... et des groupes $H'_{\flat}$ etc... On pr\'ecise la notation en les notant plut\^ot $W'_{\flat,T}$, $H'_{\flat,T}$ etc... Pour $T\in \underline{{\cal T}}_{\flat}$, on introduit le groupe de cohomologie $H^1(T)=H^1(Gal(\bar{F}/F),T)$. Puisque $A_{T}=\{1\}$, ce groupe est un produit de facteurs ${\mathbb Z}/2{\mathbb Z}$ et il est non trivial si $T\not=\{1\}$. On pose
 $$h(T)=\left\lbrace\begin{array}{cc}\vert H^1(T)\vert /2,&\,\,\text{si}\,\,T\not=\{1\},\\ 1,&\,\,\text{si}\,\,T=\{1\}.\\ \end{array}\right.$$
 On introduit le groupe $\bar{W}(H_{\flat},T)=Norm_{H_{\flat}}(T)/Z_{H_{\flat}}(T)$. Le groupe de Galois $Gal(\bar{F}/F)$ agit sur $\bar{W}(H_{\flat},T)$, on note $\bar{W}_{F}(H_{\flat},T)$ le sous-groupe des points fixes.
  
  Fixons un isomorphisme $\Phi:W_{a}\otimes_{F}\bar{F}\to W_{i}\otimes_{F}\bar{F}$ tel que $q_{W_{i}}(\Phi(w),\Phi(w'))=q_{W_{a}}(w,w')$ pour tous $w,w'\in W_{a}\otimes_{F}\bar{F}$. Pour $h\in H_{a}$, notons $\phi(h)=\Phi\circ h\circ\Phi^{-1}$. C'est un \'el\'ement de $H_{i}$ et l'isomorphisme $\phi$ ainsi d\'efini de $H_{a}$ sur $H_{i}$ est un torseur int\'erieur. Soient $T, T'\in \underline{\cal T}_{a}\cup \underline{\cal T}_{i}$. On dit que $T$ et $T'$ sont stablement conjugu\'es si et seulement si l'une des conditions suivantes (1) ou (2) est v\'erifi\'ee.
  
  (1) Il existe un indice $\flat$ tel que $T,T'\in \underline{\cal T}_{\flat}$. Il existe $h\in H_{\flat}$ tel que $hTh^{-1}=T'$ et l'homomorphisme $t\mapsto hth^{-1}$ de $T$ sur $T'$ est d\'efini sur $F$. 
  
  (2) Quitte \`a \'echanger $T$ et $T'$, on a $T\in \underline{\cal T}_{a}$ et $T'\in \underline{\cal T}_{i}$. Il existe $h\in H_{i}$ de sorte que $h\phi(T_{a})h^{-1}=T_{i}$ et l'homomorphisme $t\mapsto h\phi(t)h^{-1}$ de $T_{a}$ sur $T_{i}$ est d\'efini sur $F$. 
  
  On note $T\sim_{st}T'$ cette relation. On v\'erifie que c'est une relation d'\'equivalence.

\ass{Lemme}{(i) Soient $T$ et $T'$ deux \'el\'ements de $ \underline{\cal T}_{a}\cup\underline{\cal T}_{i}$ stablement conjugu\'es.  Dans la situation de (1), l'espace quadratique $W'_{\flat,T}$, resp. $W''_{\flat,T}$, est isomorphe \`a $W'_{\flat,T'}$, resp. $W''_{\flat,T'}$ et on peut choisir $h$ v\'erifiant (1) de sorte que la restriction de $h$ \`a $W''_{\flat,T}$ soit un isomorphisme d\'efini sur $F$ de $W''_{\flat,T}$ sur $W''_{\flat,T'}$. Dans la situation de (2), les espaces quadratiques $W''_{a,T}$ et $W''_{i,T'}$ sont isomorphes et on peut choisir $h$ v\'erifiant (1) de sorte que la restriction de $h\circ \Phi$ \`a $W''_{a,T}$ soit un isomorphisme d\'efini sur $F$ de $W''_{a,T}$ sur $W''_{i,T'}$.

(ii) Soit $T\in \underline{\cal T}_{\flat}$. On a l'\'egalit\'e
$$\sum_{T'\in {\cal T}_{\flat}; T'\sim_{st}T}\vert W(H_{\flat},T')\vert ^{-1}=h(T)\vert \bar{W}_{F}(H_{\flat},T)\vert ^{-1}.$$
 Les nombres $h(T) $ et $\vert \bar{W}_{F}(H_{\flat},T)\vert $  ne d\'ependent que de la classe de conjugaison stable de $T$.

(iii) Toute classe de conjugaison stable dans $\underline{\cal T}_{a}\cup\underline{\cal T}_{i}$ coupe $\underline{\cal T}_{i}$. La seule classe de conjugaison stable qui ne coupe pas $\underline{\cal T}_{a}$   est la classe r\'eduite au tore $\{1\}\in \underline{\cal T}_{i}$.}

Preuve. Soient $T$ et $T'$ deux \'el\'ements de $ \underline{\cal T}_{\flat}$ stablement conjugu\'es et $h$ v\'erifiant (1). Supposons $d_{W_{\flat}}$ impair. On a $h(W''_{\flat,T}\otimes_{F}\bar{F})=W''_{\flat,T'}\otimes_{F}\bar{F}$. Notons $h'':W''_{\flat,T}\otimes_{F}\bar{F}\to W''_{\flat,T'}\otimes_{F}\bar{F}$ la restriction de $h$. L'application $\phi'':x\mapsto h''xh^{_{''}-1}$ est un isomorphisme de $H''_{\flat,T}$ sur $H'_{\flat,T}$. Ces deux groupes \'etant quasi-d\'eploy\'es, on peut fixer dans chacun d'eux un sous-groupe de Borel, un sous-tore maximal de ce groupe et un \'epinglage, ces donn\'ees \'etant d\'efinies sur $F$. Quitte \`a multiplier $h''$ \`a droite par un \'el\'ement de $H''_{\flat,T}$, ce qui est loisible, on peut supposer que $\phi''$ envoie ces donn\'ees du groupe $H''_{\flat,T}$ sur celles du groupe $H''_{\flat,T'}$. Soit $\sigma\in Gal(\bar{F}/F)$. La condition (1) entra\^{\i}ne que $\sigma(h)^{-1}h$ appartient au commutant de $T$ dans $H_{\flat}$, c'est-\`a-dire \`a $T\times H''_{\flat,T}$. Donc $\sigma(h'')^{-1}h''\in H''_{\flat,T}$. Cet \'el\'ement conserve le sous-groupe de Borel, le tore maximal et l'\'epinglage de $H''_{\flat,T}$. Donc il appartient au centre de $H''_{\flat,T}$, qui est r\'eduit \`a $\{1\}$ puisque $dim(W''_{\flat})$ est impaire. Donc $h''$ est d\'efini sur $F$ et c'est un isomorphisme de $W''_{\flat,T}$ sur $W''_{\flat,T'}$. Supposons $d_{W_{\flat }}$ pair. On consid\`ere $h$ comme un \'el\'ement de $G$ et on remplace les espaces $W''_{\flat,T}$ et $W''_{\flat,T'}$ par $V''_{\flat,T}$ et $V''_{\flat,T'}$ dans le raisonnement pr\'ec\'edent. On conclut de m\^eme que ces deux derniers espaces sont isomorphes. Puisque $W''_{\flat,T}$ et $W''_{\flat,T'}$ sont les orthogonaux dans ces espaces de l'espace commun $D\oplus Z$, ils sont eux-aussi isomorphes. De m\^eme, maintenant que l'on a prouv\'e que $W''_{\flat,T}$ et $W''_{\flat,T'}$ \'etaient isomorphes, $W'_{\flat,T}$ et $W'_{\flat,T'}$ le sont aussi.
On peut remplacer $h$ par un \'el\'ement qui a m\^eme restriction \`a $W'_{\flat,T}\otimes_{F}\bar{F}$ et qui est un isomorphisme d\'efini sur $F$ de $W''_{\flat,T}$ sur $W''_{\flat,T'}$ (on peut choisir ce dernier tel que $h$ soit dans $H_{\flat}$ et non seulement dans $H_{\flat}^+$). Un tel $h$ v\'erifie encore (1) et la condition du (i) de l'\'enonc\'e.

Soient $T\in \underline{\cal T}_{a}$ et $T'\in \underline{\cal T}_{i}$ deux \'el\'ements stablement conjugu\'es et $h$ v\'erifiant (2). Dans le cas $d_{W_{i}}$ impair, le m\^eme raisonnement s'applique en rempla\c{c}ant l'application $h''$ par la restriction de $h\circ\Phi$ \`a $W''_{a,T}\otimes_{F}\bar{F}$ et $\phi''$ par la restriction \`a $H''_{a,T}$ de $x\mapsto h\phi(x)h^{-1}$. Dans le cas $d_{W_{i}}$ pair, on \'etend $\Phi$ en un isomorphisme de $V_{a}\otimes_{F}\bar{F}$ sur $V_{i}\otimes_{F}\bar{F}$ qui est l'identit\'e sur $D\oplus Z$. On en d\'eduit un prolongement de $\phi$ en un torseur int\'erieur de $G_{a}$ sur $G_{i}$. On remplace alors $h''$ par la restriction de $h\circ\Phi$ \`a $V''_{a,T}\otimes_{F}\bar{F}$ et $\phi''$ par la restriction \`a $G''_{a,T}$ de $x\mapsto h\phi(x)h^{-1}$. Cela prouve (i)

 Soit $T\in \underline{\cal T}_{\flat}$. Notons ${\cal H}_{T}$ l'ensemble des $h\in H_{\flat}$ tels que $\sigma(h)^{-1}h\in  T$ pour tout $\sigma\in Gal(\bar{F}/F)$.  On vient de voir que, pour tout \'el\'ement $T'\in \underline{\cal T}_{\flat}$ stablement conjugu\'e \`a $T$, il existait $h\in {\cal H}_{T}$ tel que $hTh^{-1}=T'$. Inversement, pour $h\in {\cal H}_{T}$, le tore $T'=hTh^{-1}$ est d\'efini sur $F$. La restriction de $h$ \`a $W''_{\flat,T}$ est d\'efinie sur $F$ et le tore $T'$ appartient \`a $\underline{\cal T}_{\flat}$: on a $W''_{\flat,T'}=h(W''_{\flat,T})$ et $W'_{\flat,T'}$ est l'orthogonal de cet espace dans $W_{\flat}$. A tout $h\in {\cal H}_{T}$, associons l'unique \'el\'ement $T_{h}\in {\cal T}_{\flat}$ tel que $hTh^{-1}$ soit conjugu\'e \`a $T_{h}$ par un \'el\'ement de $H_{\flat}(F)$. L'application $h\mapsto T_{h}$ se quotiente en une surjection de $H_{\flat}(F)\backslash{\cal H}_{T}/T$ sur l'ensemble des $T'\in {\cal T}_{\flat}$ tels que $T'\sim_{st}T$. Pour $h\in H_{\flat}(F)\backslash{\cal H}_{T}/T$, notons $n(h)$ le nombre d'\'el\'ements de la fibre de cette application au-dessus de $T_{h}$. Le membre de gauche de l'\'egalit\'e du (ii) de l'\'enonc\'e est \'egal \`a
 $$\sum_{h\in H_{\flat}(F)\backslash{\cal H}_{T}/T}n(h)^{-1}\vert W(H_{\flat},T_{h})\vert ^{-1}.$$
 Soit $h\in H_{\flat}(F)\backslash{\cal H}_{T}/T$, que l'on rel\`eve en un \'el\'ement de ${\cal H}_{T}$ tel que $hTh^{-1}=T_{h}$.  Notons ${\cal H}_{h,T}$ l'ensemble des $x\in {\cal H}_{T}$ tels que $xTx^{-1}=T_{h}$. L'entier $n(h)$ est le nombre d'\'el\'ements de l'image de ${\cal H}_{h,T}$ dans $H_{\flat}(F)\backslash{\cal H}_{T}/T$. On v\'erifie que ${\cal H}_{h,T}=h({\cal H}_{T}\cap Norm_{H_{\flat}}(T))$. Donc $n(h)$ est le nombre d'\'el\'ements de l'image de l'application
 $$\begin{array}{ccc}{\cal H}_{T}\cap Norm_{H_{\flat}}(T)&\to&H_{\flat}(F)\backslash{\cal H}_{T}/T\\ n&\mapsto &H_{\flat}(F)hnT.\\ \end{array}$$
 On v\'erifie que cette application se quotiente en une bijection de $Norm_{H_{\flat}(F)}(T_{h})\backslash ({\cal H}_{T}\cap Norm_{H_{\flat}}(T))/TH''_{\flat,T}(F)$ sur son image. Le groupe $TH''_{\flat,T}(F)$ est un sous-groupe distingu\'e de ${\cal H}_{T}\cap Norm_{H_{\flat}}(T)$ et son intersection avec $Norm_{H_{\flat}(F)}(T_{h})$ est $Z_{H_{\flat}(F)}(T_{h})$. Donc $n(h)=n_{1}\vert W(H_{\flat},T_{h})\vert ^{-1}$, o\`u
 $n_{1}$ est le nombre d'\'el\'ements de $({\cal H}_{T}\cap Norm_{H_{\flat}}(T))/TH''_{\flat,T}(F)$. Le membre de gauche de l'\'egalit\'e du (ii) de l'\'enonc\'e est donc \'egal \`a $n_{2}n_{1}^{-1}$, o\`u $n_{2}$ est le nombre d'\'el\'ements de $H_{\flat}(F)\backslash{\cal H}_{T}/T$. On v\'erifie que l'application qui, \`a $h\in {\cal H}_{T}$, associe le cocycle $\sigma\mapsto \sigma(h)^{-1}h$, se quotiente en une bijection de $H_{\flat}(F)\backslash{\cal H}_{T}/T$ sur le noyau de l'application $H^1(T)\to H^1(H_{\flat})$. On sait bien que ce noyau a $h(T)$ \'el\'ements. Donc $n_{2}=h(T)$. Consid\'erons l'application naturelle
 $$(3) \qquad ({\cal H}_{T}\cap Norm_{H_{\flat}}(T))/TH''_{\flat,T}(F)\to \bar{W}(H_{\flat},T).$$
 On v\'erifie qu'elle est injective et que son image est contenue dans $\bar{W}_{F}(H_{\flat},T)$. Inversement, soit $n\in Norm_{H_{\flat}}(T)$ dont l'image dans $\bar{W}(H_{\flat},T)$ appartient \`a $\bar{W}_{F}(H_{\flat},T)$. Tout \'el\'ement de $Norm_{H_{\flat}}(T)$ conserve les espaces $W'_{\flat,T}\otimes_{F}\bar{F}$ et $W'_{\flat,T}\otimes_{F}\bar{F}$. Notons $n'$ et $n''$ les restrictions de $n$ \`a ces espaces. Choisissons $n''_{0}\in H^{_{''}+}_{\flat,T}(F)$ de m\^eme d\'eterminant que $n''$. Consid\'erons l'\'el\'ement $n_{0}\in H_{\flat}$ de restriction $n'$ \`a $W'_{\flat,T}\otimes_{F}\bar{F}$ et de restriction $n''_{0}$ \`a $W''_{\flat,T}\otimes_{F}\bar{F}$. Il appartient \`a $Norm_{H_{\flat}}(T)$ et a m\^eme image que $n$ dans $\bar{W}(H_{\flat},T)$. Puisque cette image appartient \`a $\bar{W}_{F}(H_{\flat},T)$, on a $\sigma(n_{0})^{-1}n_{0}\in T\times H''_{\flat,T}$ pour tout $\sigma\in Gal(\bar{F}/F)$. Par construction, la restriction de $\sigma(n_{0})^{-1}n_{0}$ \`a $W''_{\flat,T}\otimes_{F}\bar{F}$ est l'identit\'e. Donc cet \'el\'ement appartient \`a $T$ et $n_{0}$ appartient \`a ${\cal H}_{T}\cap Norm_{H_{\flat}}(T)$. L'image de l'injection (3) est donc \'egale \`a $\bar{W}_{F}(H_{\flat},T)$ et $n_{1}$ est le nombre d'\'el\'ements de cet ensemble. Cela prouve l'\'egalit\'e du (ii) de l'\'enonc\'e.
 
 Deux tores $T$ et $T'$ stablement conjugu\'es sont isomorphes sur $F$, donc $h(T)=h(T')$. Dans la situation de (2), on v\'erifie que l'application $n\mapsto h\phi(n)h^{-1}$ est un isomorphisme de $\bar{W}_{F}(H_{a},T)$ sur $\bar{W}_{F}(H_{i},T')$ d'o\`u l'\'egalit\'e du nombre d'\'el\'ements de ces ensembles. Un r\'esultat analogue vaut dans la situation de (1). Cela prouve (ii).
 
 Soit $T\in \underline{\cal T}_{a}$.  Puisqu'au moins l'un des groupes $H_{a}$ ou $G_{a}$ n'est pas  quasi-d\'eploy\'e, on a $W'_{a,T}\not=\{0\}$. Cet espace est de dimension paire. S'il est de dimension $2$, on a $T=H'_{a,T}$ et, puisque $A_{T}=\{1\}$, $W'_{a,T}$ n'est pas hyperbolique. Il existe donc un espace quadratique, notons-le $W'_{i,T}$, qui a m\^eme dimension que $W'_{a,T}$, dont la forme quadratique a m\^eme d\'eterminant que celle de $W'_{a,T}$, mais qui n'est pas isomorphe \`a $W'_{a,T}$. Il est clair que la somme orthogonale $W'_{i,T}\oplus W''_{a,T}$ est isomorphe \`a $W_{i}$. Puisque $T$ est un sous-tore maximal elliptique de $H'_{a,T}$, on peut le transf\'erer en un sous-tore maximal elliptique $T'$ du groupe sp\'ecial orthogonal $H'_{i,T}$ de l'espace $W'_{i,T}$. Par l'isomorphisme pr\'ec\'edent, $T'$ devient un sous-tore de $H_{i}$. Ce tore $T'$ appartient \`a $\underline{\cal T}_{i}$ et est stablement conjugu\'e \`a $T$. Donc la classe de conjugaison stable de $T$ coupe $\underline{\cal T}_{i}$. Si $T\in \underline{\cal T}_{i}$ et $T\not=\{1\}$, un  raisonnement analogue montre que la classe de conjugaison stable de $T$ coupe $\underline{\cal T}_{a}$. Par contre, si $T=\{1\}$, sa classe de conjugaison stable se r\'eduit \`a $T$ lui-m\^eme. Le seul sous-tore de $H_{a}$ qui pourrait lui correspondre est le sous-tore $\{1\}$ de $H_{i}$. Mais celui-ci n'appartient pas \`a $\underline{\cal T}_{a}$ puisque  l'un des groupes $H_{a}$ ou $G_{a}$ n'est pas quasi-d\'eploy\'e. $\square$

 \bigskip
 
 \subsection{D\'emonstration du th\'eor\`eme 13.3}
 
 Pour un indice $\flat=i$ ou $a$, consid\'erons la somme
 $$m(\Sigma_{\flat},\Pi_{\flat})=\sum_{\sigma\in \Sigma_{\flat},\pi\in \Pi_{\flat}}m(\sigma,\pi).$$
 D'apr\`es l'hypoth\`ese du th\'eor\`eme, toutes les repr\'esentations $\pi$ qui interviennent sont supercuspidales. D'apr\`es la proposition 13.1,on peut remplacer les multiplicit\'es $m(\sigma,\pi)$ par $m_{geom}(\sigma,\pi)$.  On obtient
 $$(1) \qquad m(\Sigma_{\flat},\Pi_{\flat})=\sum_{T\in {\cal T}_{\flat}}\vert W(H_{\flat},T)\vert ^{-1}\nu(T)\int_{T(F)}c_{\check{\Sigma}_{\flat}}(t)c_{\Pi_{\flat}}(t)D^H(t)\Delta(t)^rdt,$$
 o\`u 
 $$c_{\Pi_{\flat}}(t)=\sum_{\pi\in \Pi_{\flat}}c_{\pi}(t),$$
 et $c_{\check{\Sigma}_{\flat}}(t)$ est d\'efini de fa\c{c}on analogue. Fixons $T\in {\cal T}_{\flat}$, introduisons les espaces $W''_{T}$ et $V''_{T}$ et les groupes $H''_{\flat,T}$ et $G''_{\flat,T}$ qui lui sont associ\'es.   Pour $\pi\in \Pi_{\flat}$, $t\in T_{\sharp}(F)$ et $X\in \mathfrak{g}''_{\flat,T,reg}(F)$, avec $X$ assez proche de $0$, on a un d\'eveloppement
 $$\theta_{\pi}(texp(X))D^{G''_{\flat,T}}(X)^{1/2}=\sum_{{\cal O}\in Nil(\mathfrak{g}''_{\flat,T})}c_{\theta_{\pi},{\cal O}}(t)\hat{j}^{G''_{\flat,T}}({\cal O},X)D^{G''_{\flat,T}}(X)^{1/2}.$$
 Le terme $c_{\pi}(t)$ est \'egal \`a $c_{\theta_{\pi},{\cal O}_{reg}}(t)$ si $d$ est impair ou si $dim(V''_{\flat,T})\leq2$, \`a $c_{\theta_{\pi},{\cal O}_{\nu_{0}}}(t)$ si $d$ est pair et $dim(V''_{\flat,T})\geq4$.
 Rempla\c{c}ons $X$ par $\lambda X$, avec $\lambda\in F^{\times 2}$ et faisons tendre $\lambda$ vers $0$. D'apr\`es 2.6(1), la fonction 
 $$\lambda\mapsto \hat{j}^{G''_{\flat,T}}({\cal O},\lambda X)D^{G''_{\flat,T}}(\lambda X)^{1/2}$$
 tend vers $0$ si ${\cal O}$ n'est pas r\'eguli\`ere. Supposons d'abord $d$ impair ou $dim(V''_{\flat,T})\leq2$. Choisissons pour $X$ un \'el\'ement  de la forme $X_{qd}$. D'apr\`es le lemme 13.4, la fonction ci-dessus est constante de valeur $\vert W^{G''_{\flat,T}})\vert $ pour ${\cal O}={\cal O}_{reg}$. On obtient
 $$c_{\pi}(t)=\vert W^{G''_{\flat,T}})\vert^{-1}lim_{\lambda\to 0}\theta_{\pi}(texp(\lambda X_{qd}))D^{G''_{\flat,T}}(\lambda X_{qd})^{1/2}.$$
 Si maintenant $d$ est pair et $dim(V''_{\flat,T})\geq4$, on introduit dans l'alg\`ebre $\mathfrak{g}''_{\flat,T}(F)$ des \'el\'ements $X_{F_{1}}^{\pm}$. Le lemme 13.4 conduit alors \`a l'\'egalit\'e
 $$c_{\pi}(t)=lim_{\lambda\to 0}(\vert W^{G''_{\flat,T}})\vert^{-1}\theta_{\pi}(texp(\lambda X_{qd}))D^{G''_{\flat,T}}(\lambda X_{qd})^{1/2}$$
$$+\sum_{F_{1}\in {\cal F}^{V''_{\flat,T}}}\vert W(G''_{\flat,T},T_{F_{1}})\vert ^{-1}\epsilon_{F_{1}}\chi_{F_{1}}(\nu_{0}\eta)(\theta_{\pi}(texp(\lambda X_{F_{1}}^+))-\theta_{\pi}(texp(\lambda X_{F_{1}}^-))) D^{G''_{\flat,T}}(\lambda X_{F_{1}})^{1/2}.$$
 Pour calculer $c_{\Pi_{\flat}}(t)$, on somme ces expressions sur $\pi\in \Pi_{\flat}$. Cela revient \`a remplacer les caract\`eres $\theta_{\pi}$ par $\theta_{\Pi_{\flat}}$. Dans le cas o\`u $d$ est pair et $dim(V''_{\flat,T})\geq4$, on  remarque que, pour $F_{1}\in {\cal F}^{V''_{\flat,T}}$, les points $texp(\lambda X_{F_{1}}^+)$ et $texp(\lambda X_{F_{1}}^-)$ sont stablement conjugu\'es. Or $\theta_{\Pi_{\flat}}$ est une distribution stable. Donc 
 $$\theta_{\Pi_{\flat}}(texp(\lambda X_{F_{1}}^+))-\theta_{\Pi_{\flat}}(texp(\lambda X_{F_{1}}^-))=0.$$ 
 En tout cas, on obtient l'\'egalit\'e
 $$(2) \qquad c_{\Pi_{\flat}}(t)=\vert W^{G''_{\flat,T}})\vert^{-1}lim_{\lambda \to 0}\theta_{\Pi_{\flat}}(texp(\lambda X_{qd}))D^{G''_{\flat,T}}(\lambda X_{qd})^{1/2}.$$
 Soit $T'\in{\cal T}_{\flat}$ stablement conjugu\'e \`a $T$. D'apr\`es le lemme 13.5(i), on peut choisir un \'el\'ement $h\in H$ v\'erifiant la condition 13.5(1) et tel que sa restriction \`a $V''_{\flat,T}$ soit un isomorphisme d\'efini sur $F$ de cet espace sur $V''_{\flat,T'}$. La conjugaison par $h$ identifie $T$ \`a $T'$, $V''_{\flat,T}$ \`a $V''_{\flat,T'}$ et $G''_{\flat,T}$ \`a $G''_{\flat,T'}$. Soit $t'=hth^{-1}$. Posons $X'_{qd}=hX_{qd}h^{-1}$. C'est un \'el\'ement de $\mathfrak{g}''_{\flat,T'}(F)$ qui a les m\^emes propri\'et\'es que $X_{qd}$. On peut calculer $c_{\Pi_{\flat}}(t')$ en rempla\c{c}ant $t$ par $t'$ et $X_{qd}$ par $X_{qd}'$ dans la formule (2). Mais les \'el\'ements $texp(\lambda X_{qd})$ et $t'exp(\lambda X'_{qd})$ sont stablement conjugu\'es. Puisque $\theta_{\Pi_{\flat}}$ est stable, cette fonction prend la m\^eme valeur en ces deux \'el\'ements et on en d\'eduit $c_{\Pi_{\flat}}(t)=c_{\Pi_{\flat}}(t')$. On d\'emontre de m\^eme l'\'egalit\'e $c_{\check{\Sigma}_{\flat}}(t)=c_{\check{\Sigma}_{\flat}}(t)$. On a aussi $D^H(t)\Delta(t)^r=D^H(t')\Delta(t)^r$. Alors l'int\'egrale index\'ee par $T'$ dans la formule (1) a la m\^eme valeur que celle index\'ee par $T$.
 
    Notons ${\cal T}_{\flat,st}$ un sous-ensemble de repr\'esentants dans ${\cal T}_{\flat}$ des classes de conjugaison stable coupant $\underline{\cal T}_{\flat}$. Alors
 $$m(\Sigma_{\flat},\Pi_{\flat})=\sum_{T\in {\cal T}_{\flat,st}} m(T)\nu(T)\int_{T(F)}c_{\check{\Sigma}_{\flat}}(t)c_{\Pi_{\flat}}(t)D^H(t)\Delta(t)^rdt,$$
  o\`u
 $$m (T)=\sum_{T'\in {\cal T}_{\flat}; T'\sim_{st}T}\vert W(H_{\flat},T')\vert ^{-1}.$$
 Le lemme 13.5(iii) d\'efinit une injection ${\cal T}_{a,st}\to {\cal T}_{i,st}$. Soient $T_{a}\in {\cal T}_{a,st}$ et $T_{i}\in {\cal T}_{i,st}$ son image. Choisissons $h\in H_{i}$ v\'erifiant la condition 13.5(2) et tel que la restriction de $h\circ \Phi$ \`a $W''_{a,T_{a}}$ soit d\'efinie sur $F$ (lemme 13.5(i)). Soit $t\in T_{a,\sharp}(F)$, posons $t'=h\phi(t)h^{-1}$. Un argument similaire \`a celui ci-dessus montre que
 $$c_{\Pi_{a}}(t)=(-1)^dc_{\Pi_{i}}(t'),\,\,c_{\check{\Sigma}_{a}}(t)=(-1)^{d_{W}}c_{\check{\Sigma}_{i}}(t').$$
 Les signes proviennent de la relation 13.2(2). Leur produit est $-1$. On a aussi $m(T_{a})=m(T_{i})$ d'apr\`es le lemme 13.5(ii) et, bien s\^ur, $D^{H_{a}}(t)\Delta(t)^r=D^{H_{i}}(t')\Delta(t')^r$. Alors la contribution de $T_{a}$ \`a $m(\Sigma_{a},\Pi_{a})$ est l'oppos\'ee de celle de $T_{i}$ \`a $m(\Sigma_{i},\Pi_{i})$. La somme $m(\Sigma_{a},\Pi_{a})+m(\Sigma_{i},\Pi_{i})$ se r\'eduit donc \`a la contribution de l'unique classe de conjugaison stable qui ne coupe pas ${\cal T}_{a,st}$, c'est-\`a-dire \`a la contribution du tore $\{1\}$ de ${\cal T}_{i}$. On obtient
$$(3) \qquad m(\Sigma_{a},\Pi_{a})+m(\Sigma_{i},\Pi_{i})=c_{\check{\Sigma}_{i}}(1)c_{\Pi_{i}}(1).$$ 
Rappelons un r\'esultat de Rodier ([R] th\'eor\`eme p.161 et remarque 2 p.162). Soient $\pi$ une repr\'esentation admissible irr\'eductible de $G_{i}(F)$ et ${\cal O}$ une orbite nilpotente r\'eguli\`ere de $\mathfrak{g}_{i}(F)$. Alors $c_{\theta_{\pi},{\cal O}}(1)$ vaut $1$ si $\pi$ poss\`ede un mod\`ele de Whittaker relatif \`a ${\cal O}$ et $0$ sinon. On a
$$c_{\Pi_{i}}(1)=\sum_{\pi\in \Pi_{i}}c_{\theta_{\pi},{\cal O}}(1),$$
o\`u ${\cal O}={\cal O}_{reg}$ ou ${\cal O}_{\nu_{0}}$ selon le cas. Le r\'esultat ci-dessus et la propri\'et\'e 13.2(3) entra\^{\i}nent que cette somme ne contient qu'un terme non nul, qui vaut $1$. Donc $c_{\Pi_{i}}(1))=1$ et, de m\^eme, $c_{\check{\Sigma}_{i}}=1$. Le membre de gauche de (3) est la somme des $m(\sigma,\pi)$ pour $(\sigma,\pi)\in (\Sigma_{i}\times \Pi_{i})\cup(\Sigma_{a}\times \Pi_{a})$. Puisqu'elle vaut $1$, il y a un unique couple $(\sigma,\pi)$ pour lequel $m(\sigma,\pi)=1$. $\square$

\bigskip
 {\bf Bibliographie}
 
  [AGRS] A. Aizenbud, D. Gourevitch, S. Rallis, G. Schiffmann: {\it Multiplicity one theorems}, pr\'epublication 2007

[A1] J. Arthur: {\it The trace formula in invariant form}, Annals of Math. 114 (1981), p.1-74

[A2] ...............: {\it The invariant trace formula I. Local theory}, J. AMS 1 (1988), p.323-383

[A3] ...............: {\it  A local trace formula}, Publ. Math. IHES 73 (1991), p.5-96

[A4] ...............: {\it On the transfer of distributions: weighted orbital integrals}, Duke Math. J. 99 (1999), p.209-283

 [A5]  .............: {\it The local behaviour of weighted orbital integrals}, Duke Math. J. 56 (1988), p.223-293

 [A6]  .............: {\it The characters of supercuspidal representations as weighted orbital integrals}, Proc. Indian Acad. Sci . 97 (1987), p.3-19
 
 [GGP] W. T. Gan, B. Gross, D. Prasad: {\it Symplectic local root numbers, central critical $L$-values and restriction problems in the representation theory of classical groups}, pr\'epublication 2008
 
 [GP] B. Gross, D. Prasad: {\it On irreducible representations of $SO_{2n+1}\times SO_{2m}$}, Can. J. Math. 46 (1994), p.930-950
 
[HCDS]  Harish-Chandra: {\it Admissible invariant distributions on reductive $p$-adic groups}, notes par S. DeBacker et P. Sally, University Lecture series 16, AMS (1999)

[HCvD] .........................: {\it Harmonic analysis on reductive $p$-adic groups}, notes par G. van Dijk, Springer Lecture Notes 162 (1970)

 [Konno] T. Konno: {\it Twisted endoscopy implies the generic packet conjecture}, Isra\" el J.  of Math. 129 (2002), p.253-289

[Kot] R. Kottwitz: {\it Transfer factors for Lie algebras}, Representation Th. 3 (1999), p.127-138

 [R] F. Rodier: {\it Mod\`ele de Whittaker et caract\`eres de repr\'esentations}, in Non commutative harmonic analysis, J. Carmona, J. Dixmier, M. Vergne ed. Springer LN 466 (1981), p.151-171

[W1] J.-L. Waldspurger: {\it Une formule des traces locale pour les alg\`ebres de Lie $p$-adiques}, J. f.  reine u.  ang. Math. 465 (1995), p.41-99

[W2] ................................: {\it D\'emonstration d'une conjecture de dualit\'e de Howe dans le cas $p$-adique,$p\not=2$}, in Festschrift in honor of I.I. Piatetski-Shapiro, S. Gelbart, R. Howe, P. Sarnak ed., the Weizmann science press of Isra\" el (1990)

[W3] ................................: {\it Int\'egrales orbitales nilpotentes et endoscopie pour les groupes classiques non ramifi\'es}, Ast\'erisque 269 (2001)
 
  [W4] ............................... : {\it Le lemme fondamental implique le transfert}, Compositio Mathematica 105 (1997), p.153-236
 
 \bigskip
 
 Institut de math\'ematiques de Jussieu- CNRS
 
 175 rue du Chevaleret
 
 75013 Paris
 
 e-mail: waldspur@math.jussieu.fr

\end{document}